\documentclass[a4paper,12pt]{article}%
\usepackage[a4paper]{geometry}
\makeatletter
\def\@makefnmark{\hbox{\@textsuperscript{\normalfont\@thefnmark}}}
\renewcommand\@makefntext[1]%
    {\noindent\makebox[0pt][r]{\textsuperscript{\@thefnmark}\,}#1}
\makeatother
\usepackage{graphicx}
\usepackage{scrextend}
\deffootnote{0em}{1.6em}{\thefootnotemark\ }
\usepackage{setspace}
\usepackage{etoolbox}

\makeatletter
\patchcmd{\@makecaption}
  {\parbox}
  {\advance\@tempdima-\fontdimen2} 
  {}{}
\makeatother 
\usepackage{graphicx}

\usepackage{amssymb,amsmath,amsthm,amsfonts,bbm,bm,mathrsfs,xcolor,array}
\usepackage{mhequ}
\usepackage{mathtools}
\usepackage[colorlinks,citecolor=blue,urlcolor=blue,filecolor=blue,backref=page]{hyperref}
\usepackage{accents,url}
\usepackage[colorinlistoftodos,prependcaption,textsize=tiny]{todonotes}
\usepackage{float}  
\usepackage{subcaption}
\usepackage[utf8]{inputenc} 
\usepackage[T1]{fontenc}    
\usepackage{xurl}
\usepackage{hyperref}       
\usepackage{url}            
\usepackage{booktabs}       
\usepackage{nicefrac}       
\usepackage{microtype}      
\usepackage{lipsum}		
 \usepackage{fixltx2e}
 \usepackage{appendix}
 \usepackage{algorithmicx}
 \usepackage{multirow}
  \usepackage{natbib}
  \usepackage[ruled,vlined]{algorithm2e}

\def\m{\mathcal} \def\mb{\mathbb}   \def\dd{{\rm d}}  \def\wt{\widetilde} \def\wh{\widehat} \def\ov{\overline}   \def \bX {{\bold{X}}}  \def \wX {\widetilde{\bold{X}}}   
\newcommand{\mnorm}[1]{{\vert\kern-0.25ex\vert\kern-0.25ex\vert #1 \vert\kern-0.25ex\vert\kern-0.25ex\vert}}
\newcommand{\bmnorm}[1]{{\big\vert\kern-0.25ex\big\vert\kern-0.25ex\big\vert #1 \big\vert\kern-0.25ex\big\vert\kern-0.25ex\big\vert}}
\newcommand{\Bmnorm}[1]{{\Big\vert\kern-0.25ex\Big\vert\kern-0.25ex\Big\vert #1 \Big\vert\kern-0.25ex\Big\vert\kern-0.25ex\Big\vert}}
\newcommand{\bbmnorm}[1]{{\bigg\vert\kern-0.25ex\bigg\vert\kern-0.25ex\bigg\vert #1 \bigg\vert\kern-0.25ex\bigg\vert\kern-0.25ex\bigg\vert}}
\newcommand{\BBmnorm}[1]{{\Bigg\vert\kern-0.25ex\Bigg\vert\kern-0.25ex\Bigg\vert #1 \Bigg\vert\kern-0.25ex\Bigg\vert\kern-0.25ex\Bigg\vert}}

\newcommand{\be}{\begin{equs}}
\newcommand{\ee}{\end{equs}}

\theoremstyle{plain}
\newtheorem{theorem}{Theorem}

\newtheorem{remark}{Remark}
\newtheorem{corollary}{Corollary}
\newtheorem{lemma}{Lemma}

\newtheorem*{definition*}{Definition}

\title{Subsampling Approximations to Gibbs Posteriors with Non-smooth Losses}

\author{%
  Wenshuo Zhao\thanks{Email: \texttt{wzhaoar@connect.ust.hk}}
  \quad and \quad
  Rong Tang\thanks{Email: \texttt{martang@ust.hk}}
}
\date{Department of Mathematics, The Hong Kong University of Science and Technology, Hong Kong.}

\begin{document}
\maketitle
\begin{abstract}
Gibbs posteriors provide a general framework for loss-based Bayesian inference, but evaluating the full empirical risk can make posterior computation costly for large datasets. We study how accurately subsampling preserves these posteriors when individual losses may not be differentiable. Building on the augmented-state subsampling Metropolis--Hastings construction, we derive finite-sample total variation bounds for a subgradient-based first-order control variate. Under smoothness and identification conditions on the population risk, our analysis quantifies the effects of minibatch size, risk scaling, and reference-point accuracy. With a sufficiently accurate reference point, these bounds yield consistent approximation of the untempered Gibbs posterior using sublinear minibatches. For a lazy version of the subsampling Metropolis--Hastings algorithm, we also establish warm-start convergence bounds under explicit geometric and proposal-comparison conditions. Together, these results connect stationary approximation error with finite-iteration convergence and provide conditions for accurate subsampling with non-smooth losses.

\end{abstract}

\par\medskip

\section{Introduction}

Bayesian inference provides a principled way to combine prior information with data and to quantify uncertainty about unknown parameters through the posterior distribution. In the standard parametric setting, given i.i.d.\ data $\bX^{(n)}=\{X_1,\ldots,X_n\}$ from $\m P^*$, a likelihood $p(x\mid\theta)$, and a prior distribution with density $\pi$ on $\Theta$, the posterior density is
 $\pi_n(\theta)\propto
\pi(\theta)\exp\{\sum_{i=1}^n\log p(X_i\mid\theta)\}$. In many applications, however, specifying an appropriate probabilistic model is difficult, and parametric likelihoods may be misspecified or overly restrictive. This has motivated generalized Bayesian and Gibbs posterior methods~\citep{syring2019gibbs,martin2022direct,1614067}, which retain the exponential form but replace the negative log-likelihood with a user-chosen loss function. Specifically, given a loss $\ell(X,\theta)$ and a prior density $\pi$, the Gibbs posterior has density $\pi_n(\theta)\propto
\pi(\theta)\exp\{-\sum_{i=1}^n\ell(X_i,\theta)\}$. This construction recovers the standard Bayesian posterior when
$\ell(X_i,\theta)=-\log p(X_i\mid\theta)$. In both settings, posterior expectations and credible regions are often unavailable in closed form. Markov chain Monte Carlo (MCMC) methods~\citep{b1878965-730b-33c7-b1ad-dfd3acb6f61b,brooks2011handbook,robert1999monte} are therefore widely used to approximate posterior quantities. However, standard full-data implementations can be computationally expensive for large datasets, since each iteration may require evaluating the loss, and possibly its gradient, over all $n$ observations. Developing scalable algorithms is consequently important for making posterior computation practical in large-data applications.

Subsampling MCMC methods~\citep{bardenet2014towards,bardenet2017markov,10.5555/3020751.3020808,Quiroz03042019,nemeth2021stochastic,welling2011bayesian,zhang2020asymptotically,10.5555/2832747.2832855,pmlr-v32-korattikara14} offer one approach to reducing this computational burden. A common strategy replaces full-data quantities, such as $\sum_{i=1}^n\ell(X_i,\theta)$, with estimates computed from a subset of the observations. Control variates play an important role in reducing the variability of these estimates and can improve sampling efficiency~\citep{Quiroz03042019,bardenet2017markov,nemeth2021stochastic,baker2019control}. Common parameter-expanded control variates use Taylor approximations of $\ell(X,\cdot)$ around a reference point in the parameter space. In particular, second-order constructions use gradients and Hessians, and their associated approximation theory imposes bounds on the third-order derivatives of the individual loss contributions~\citep{Quiroz03042019,bardenet2017markov}.

Many Gibbs posteriors of practical interest, however, are built from non-smooth losses. Examples include absolute and check losses in robust and quantile regression~\citep{yang2016posterior,martin2022direct}, hinge losses in large-margin classification~\citep{10.1214/11-BA601}, and other non-smooth losses used in multivariate quantile inference~\citep{bhattacharya2022gibbs}. For such problems, Gibbs posteriors are particularly appealing because they provide uncertainty summaries directly through the posterior, without constructing a normal approximation to the sampling distribution of an M-estimator. Scalable computation nevertheless presents additional challenges, since the usual smooth Taylor-remainder arguments do not directly apply to non-smooth losses. This motivates studying when subgradient-based control variates provide accurate approximations to the full-data Gibbs posterior.

\subsection{Related Work and Our Contribution}

Gibbs posteriors have a well-developed statistical foundation. \citet{bissiri2016general} provide a coherent framework for loss-based belief updating. Subsequent work studies the statistical properties of these posteriors. For example, \citet{syring2023concentration} establish concentration rates under sub-exponential-type loss conditions, and \citet{syring2019calibrating} develop learning-rate selection for calibrating posterior credible regions. \citet{martin2022direct} provide a broader account of Gibbs posterior construction, concentration, and calibration. This body of work motivates computational methods that accurately approximate the full-data Gibbs posterior.

 There is also a substantial literature on scaling MCMC to large datasets. Stochastic gradient Langevin dynamics and its variants use minibatch gradient estimates to drive discretized Langevin dynamics~\citep{welling2011bayesian,teh2016consistency,dubey2016variance,nemeth2021stochastic}. Subsampling Metropolis--Hastings (MH) methods incorporate minibatches through approximate acceptance tests, augmented-state constructions, or exact factorized acceptance rules~\citep{Quiroz03042019,Wu02012022,bardenet2017markov,zhang2020asymptotically,cornish2019scalable,pmlr-v32-korattikara14,10.5555/2832747.2832855}. Several of these methods use control variates to reduce minibatch variability and improve computational efficiency. In particular, \citet{Quiroz03042019} analyze Taylor-based control variates and quantify how posterior approximation error depends on reference-point accuracy, while \citet{cornish2019scalable} combine control-variate ideas with factorized acceptance rules to preserve the full-data posterior as the invariant distribution. In a related augmented-state approach, \citet{Wu02012022} study a minibatch MH scheme whose invariant distribution approximates a tempered posterior. More recently, \citet{frazier2026exact} analyze the distortion of Gibbs measures caused by estimated losses and develop an exact zig-zag sampler.

In this work, we develop finite-sample approximation theory for an augmented-state subsampling MH scheme for Gibbs posteriors with non-smooth losses. For subgradient-based first-order control variates, our bounds quantify how minibatch size, risk scaling, and reference-point accuracy affect posterior approximation error. Under suitable conditions, sufficiently accurate reference points enable consistent approximation of the original, untempered Gibbs posterior using sublinear minibatches.  We also provide a theoretical basis for constructing accurate reference points through staged refinement. For uniformly bounded losses, suitable risk scaling  introduces a tempered Gibbs target~\citep{10.1214/18-AOS1712,10.1214/19-AOS1855} and yields stationary means that accurately estimate the risk minimizer even without control variates. These estimates provide initial reference points for subsequent refinement using control variates.   A detailed summary of these results is given in Section~\ref{sec:summary}.  Finally, we establish warm-start total variation convergence bounds for a lazy version of the algorithm. Under explicit geometric and proposal-comparison conditions, these bounds retain the same dependence on geometric parameters as the corresponding full-data MH bounds.

\subsection{Notation}
Throughout, we use $n$ to denote the full-data size and $m$ the subsample/minibatch size. In the theoretical analysis we assume $2\leq m\leq n$. We use $\Theta$ to denote the parameter space and $\m X$ to denote the data space. The symbols $\lesssim$ and $\gtrsim$ mean the corresponding inequality up to a  multiplicative constant independent of $(m,n)$. For real numbers $a,b$, we use $a \land b$ and $a \lor b$ to denote the minimum and maximum of them. For two distributions $p$ and $q$, we use $\|p - q\|_{\mathrm{TV}} = \frac12 \int |p(x) - q(x)| \dd x$ to denote the total variation distance. We use $T_{\#}\mu$  to denote the pushforward  measure of $\mu$ under map $T$. We use ${B}(x,r)$ to denote the closed ball centered at $x$ with radius $r$ and $\operatorname{diam}(X)$ to denote diameter of $X$ in $\ell_2$ norm. For a function $f: \mathbb{R}^d \to \mathbb{R}$, we use $\nabla f(x)$ to denote the $d$-dimensional gradient vector of $f$ at $x$ and $\operatorname{Hess}(f(x))$ to denote the Hessian matrix of $f$ at $x$. For a matrix $J$, we use $\|J\|_{\mathrm{F}}$ to denote its Frobenius norm. For two symmetric matrices $A$ and $B$, we use $A \succcurlyeq B$ to mean that $A - B$ is a positive semi-definite matrix. We use the standard $\mathcal{O}$ notation: $f = \mathcal{O}(g)$ means $f \leq C \cdot g$ for some constant $C>0$ (independent of $m,n$) and $\tilde{\mathcal{O}}$ further suppresses logarithmic factors:
$f = \tilde{\mathcal{O}}(g)$ iff $f = \mathcal{O}(g \cdot (\log m)^a (\log n)^b)$ for some $a,b \in \mathbb{R}$.  Similarly,  $f=\Omega(g)$  means  $f \geq C \cdot g$ for some constant $C>0$  and $\tilde{\Omega}(\cdot)$ suppresses logarithmic factors in the same way. 
\subsection{Organization}
The rest of the paper is organized as follows. In Section~\ref{Problem Setup and Overview of the Main Results}, we specify the augmented-state subsampling Metropolis-Hastings (SMH) construction and summarize the approximation results. Section~\ref{Construction} presents our theoretical analysis of the stationary distribution of SMH. In Section~\ref{Mixing time analysis}, we study the mixing time of SMH, including its behavior under a minibatch MALA proposal. Section~\ref{Numerical Study} contains four numerical examples, with further implementation details in Appendix~\ref{detail}. The proofs of the main results are given in Appendix~\ref{Proof of Main Results}, and the proofs of technical lemmas are collected in Appendix~\ref{Proof of technical lemmas}.

\section{Problem Setup and Overview of the Main Results}\label{Problem Setup and Overview of the Main Results}
\subsection{Augmented-State Subsampling Metropolis-Hastings Algorithm}\label{sec:smh}
Consider a Gibbs posterior with density $\pi_n(\theta)\propto \pi(\theta)\exp(-\sum_{i=1}^n \ell(X_i,\theta))$, where $\ell:\m X\times \Theta\to \mb R$ is a loss function.  A standard MH algorithm for sampling from $\pi_n$  can be computationally expensive, as each iteration requires a full pass over the data. An appealing remedy is to employ a subsampling MH algorithm~\citep{Wu02012022,quiroz2018subsampling}.  Given a subsample size $m$, we   construct a mini-batch surrogate $\m L: \m X^{\otimes m}\times \Theta\to \mb R$, so that for subsamples $\wX^{(m)} = \{\wt X_1, \wt X_2, \dots, \wt X_m\}$ drawn from $\bX^{(n)}$, we treat $\m L(\wX^{(m)}, \theta)$ as an estimator of $\frac{1}{n}\sum_{i=1}^n \ell(X_i, \theta)$, up to additive constants independent of $\theta$.  Also, we consider a proposal distribution $q(\cdot \mid \theta', \wX'^{(m)})$ that may depend on the current state $\theta'$ and  subsample $\wX'^{(m)}$. Following~\cite{Wu02012022}, we additionally introduce an inverse-temperature parameter $\tau\in (0,1]$ and the associated risk-scaling factor $\alpha=n\tau$ to mitigate potential variance explosion in $n\cdot\m L(\wX^{(m)}, \theta)$.  Let $\nu_0$ be the distribution that generates the initial state of $\theta$, and let  $\mu_{n}^{(m)}$  be the distribution on $(\bX^{(n)})^{\otimes m}$ that generates subsamples. The SMH algorithm then produces a sequence $\{\theta_k\}_{k\geq 0}$ as follows: for $k = 0, 1, 2, \dots$,\\
\noindent \emph{\underline{1. Initialization.} If $k = 0$, sample $\theta_0 \sim \nu_0$ and $\wX^{(m)}_0 = \{\wt X_{01}, \wt X_{02}, \dots, \wt X_{0m}\}\sim\mu_{n}^{(m)}$}.\\
 \emph{\underline{2. Proposal.} If $k \geq 1$, given the previous state $(\theta_{k-1}, \wX^{(m)}_{k-1})$, generate a candidate $\theta$ from the proposal $q\big(\cdot|\theta_{k-1}, \wX^{(m)}_{k-1}\big)$. Then with probability $1-\zeta$, generate $\wX^{(m)} = \{\wt X_{1}, \dots, \wt X_{m}\}$ from $\mu_{n}^{(m)}$; otherwise set $\wX^{(m)}=\wX^{(m)}_{k-1}$.}\\
  \emph{\underline{3. Metropolis-Hastings Rejection Step.} Define the acceptance probability:} 
  \vspace{-1em}
    \begin{equation*}
      A\big((\theta_{k-1}, \wX^ {(m)}_{k-1}), (\theta, \wX^ {(m)})\big)= 1\wedge  \frac{\pi(\theta)\exp\Big(-\alpha\cdot\m L(\wX^ {(m)}, \theta)\Big) \cdot q\big(\theta_{k-1} \mid \theta, \wX^ {(m)}\big)}{\pi(\theta_{k-1})\exp\Big(-\alpha\cdot\m L(\wX^ {(m)}_{k-1}, \theta_{k-1})\Big) \cdot q\big(\theta \mid \theta_{k-1}, \wX^ {(m)}_{k-1} \big)}.
    \end{equation*}
      \vspace{-1em}
      
    \noindent \emph{Flip a coin and set $( \theta_k,\wX_{k}^{(m)})=(\theta,\wX^{(m)})$ with probability $A\big((\theta_{k-1}, \wX^ {(m)}_{k-1}), (\theta, \wX^ {(m)})\big)$; otherwise, set $( \theta_k,\wX_{k}^{(m)})=(\theta_{k-1},\wX_{k-1}^{(m)})$.}
 

 \vspace{0.5em}
The minibatch-retention parameter $\zeta\in[0,1)$ controls how often a new minibatch is proposed. This retention mechanism is used in the numerical implementations; the fully lazy variant used in the mixing analysis is specified in Section~\ref{Mixing time analysis}. To understand the motivation, suppose the chain reaches a state $(\theta,\wX^{(m)})$ at which $\m L(\wX^{(m)},\theta)$ is unusually small relative to its values on newly sampled minibatches. Proposals that refresh the minibatch may then be rejected frequently, causing the chain to remain at this state for many iterations. Retaining the current minibatch allows parameter updates without simultaneously requiring a minibatch refresh. 
The construction follows the augmented-state principle of pseudo-marginal MCMC~\citep{508a6095-c33f-3d5d-9a9a-f004c0c0a59a}, with  $\wX^{(m)}$  treated as an auxiliary variable. In particular, if $q(\theta|\theta',\wX'{}^{(m)})>0$ for every current state and every $\theta\in \Theta$ with $\pi(\theta)>0$ so that the chain is irreducible, then the joint stationary distribution of $\{(\theta_k, \wX^ {(m)}_k)\}_{k\geq 0}$ is $$\wt\pi_{\m L,\alpha}^{\rm joint}(\theta,\wX^ {(m)})\propto\pi(\theta)\exp\!\big(-\alpha\cdot\m L(\wX^ {(m)},\theta)\big)\mu_{n}^{(m)}(\wX^ {(m)})$$ with marginal in $\theta$:
\begin{equation*}
\wt\pi_{\m L,\alpha}(\theta)
=
\frac{
\pi(\theta)\,
\mb{E}_{\wX^{(m)}\sim\mu_n^{(m)}}
\left[\exp\!\bigl(-\alpha\m L(\wX^{(m)},\theta)\bigr)\right]
}{
\mb{E}_{\wX^{(m)}\sim\mu_n^{(m)}}
\left[
\int_{\Theta}
\pi(\vartheta)\exp\!\bigl(-\alpha\m L(\wX^{(m)},\vartheta)\bigr)
\,\dd\vartheta
\right]
}.
\end{equation*}
In general, $\wt\pi_{\m L,\alpha}$ differs from the tempered Gibbs posterior
\[
\pi_{\alpha}(\theta)
\propto
\pi(\theta)\exp\Big(
-\frac{\alpha}{n}\sum_{i=1}^{n}\ell(X_i,\theta)
\Big).
\]
The surrogate $\m L$ and the risk-scaling factor $\alpha$ must therefore be chosen jointly to control the approximation error. In particular, setting $\alpha=n$ recovers the original target $\pi_{\alpha}=\pi_n$, but the minibatch marginal $\wt\pi_{\m L,n}$ generally remains an approximation.

  \subsection{Summary of Main Results}\label{sec:summary}
Our main approximation result concerns a first-order subgradient control variate. For comparison, we also analyze a zero-order control variate and, under a bounded-loss assumption, the naive surrogate  $\m L_{\rm naive}(\wX^{(m)},\theta)=\frac{1}{m}\sum_{j=1}^m \ell(\wt X_j,\theta)$.   The first-order surrogate is 
$ \m L_{\rm fcv}(\wX^{(m)},\theta)=\max\big(m^{-1}\sum_{j=1}^m \big[\ell(\wt X_j,\theta)-\ell\big(\wt X_j,\theta^{\dagger}\big)-g(\wt X_j,\theta^{\dagger})(\theta-\theta^{\dagger})\big],0\big) +\frac{1}{n}\sum_{i=1}^n g(X_i,\theta^{\dagger})(\theta-\theta^{\dagger})$, where $g(X,\cdot)$ is a subgradient of $\ell(X,\cdot)$ and $\theta^\dagger$ is a reference point  close to the empirical risk minimizer $\wh\theta=\arg\min_{\theta\in \Theta}\sum_{i=1}^n \ell(X_i,\theta)$. Here,
 $\ell(X,\theta)$ is uniformly Lipschitz in $\theta$ but may fail to be everywhere differentiable. The population risk $\m R(\theta)=\mb{E}_{X\sim \m P^*}[\ell(X,\theta)]$ is assumed to be twice continuously differentiable, with the Hessian regularity specified in Assumption A.1. The subgradient $g(X,\theta)$ of $\ell(X,\theta)$ with respect to $\theta$ is assumed to be $\beta_1$-smooth ``in expectation''; in particular, $\mb{E}_{X\sim \m P^*}[\|g(X,\theta)-g(X,\theta')\|_2^2]\leq L\|\theta-\theta'\|_2^{2\beta_1}$ with $\beta_1\in (0,1]$. These averaged smoothness conditions avoid imposing pointwise differentiability assumptions on \(g(X,\cdot)\) and \(\ell(X,\cdot)\) for every \(X\). For simplicity of theoretical analysis,  we fix the subsampling to be with replacement from $\bX^{(n)}$, i.e., $\mu_n^{(m)}=\mu_n^{\otimes m}$ with $\mu_n$ the empirical distribution of $\bX^{(n)}$. Under the full assumptions and scaling conditions of Theorems~\ref{th:1} and~\ref{th:2}, the resulting stationary approximation guarantees include the following.

\begin{enumerate}
    \item If we take $\m L=\m L_{\rm naive}$ and choose $\alpha=\wt{\m O}(m^{2/5})$, then the mean  of  $\wt\pi_{\m L,\alpha}(\theta)$ lies within $\wt {\m O}(m^{-\frac{2}{5}})$  of the empirical risk minimizer $\wh\theta$. Using this estimator as a reference point $\theta^\dagger$ so that $\|\theta^\dagger-\wh\theta\|_2=\wt {\m O}(m^{-\frac{2}{5}})$, and  switching to $\m L=\m L_{\rm fcv}$ with $\alpha=\wt{\m O} (m^{\frac{4}{5}})$,  we obtain, for $m=\wt \Omega(n^{\max(\frac{5}{6+8\beta_1},\frac{5}{8})})$,   a posterior mean under $\wt\pi_{\m L,\alpha}(\theta)$ that is ${\m O}(n^{-\frac{1}{2}})$-close  to $\wh\theta$.
    
    \item If  a reference point $\theta^\dagger$  is available with $\|\theta^\dagger-\wh\theta\|_2={\m O}(n^{-\frac{1}{2}})$, then taking $\m L=\m L_{\rm fcv}$ and $\alpha=n$ allows us to approximate the untempered posterior. Specifically, for $m\geq\sqrt{n}$,
  
\begin{enumerate}
    \item If $g(X,\theta)$ is Lipschitz in $\theta$ for every $X$, then ${\rm TV}(\wt \pi_{\m L,n}(\theta),\pi_n(\theta))=\wt {\m O}(\frac{1}{m})$.
    \item If $g(X,\cdot)$ is discontinuous for some $X$,  then ${\rm TV}(\wt \pi_{\m L,n}(\theta),\pi_n(\theta))=\wt {\m O}(\frac{n^{1-\beta_1}}{m}+\frac{n}{m^2})$. 
\end{enumerate}
 \end{enumerate}
 
These results show how the approximation error rate and the minibatch size sufficient for consistency depend on the smoothness of $\ell(X,\theta)$ and $g(X,\theta)$.  Under the stated assumptions, if $\ell(X,\theta)$  is differentiable in $\theta$ with a uniformly Lipschitz continuous gradient, choosing $m=\sqrt{n}$ ensures that $\wt \pi_{\m L,n}(\theta)$ is within $\wt {\m O}(n^{-\frac{1}{2}})$ total variation distance of the target posterior. On the other hand, if $\ell(X,\theta)$  is not everywhere  continuously differentiable, we need $m\geq n^{\gamma}$ with $\gamma>\max((1-\beta_1),\frac{1}{2})$  to achieve a consistent approximation.

\section{Control Variates and Posterior Approximation}\label{Construction}
We study two complementary constructions of the minibatch surrogate $\m L$. The zero-order construction uses only loss evaluations and provides a simple approach to tempered posterior approximation and point estimation. Under an additional bounded-loss assumption, its bounds with an uninformative reference point also cover the no-control-variate case. The first-order construction incorporates subgradient information and, with an accurate reference point, enables approximation of the untempered Gibbs posterior using sublinear minibatches. We derive total variation and stationary-mean error bounds for both constructions and show how they support staged reference-point refinement. Throughout this section, we take $\mu_n^{(m)}=\mu_n^{\otimes m}$.

\subsection{Zero-Order Control Variates and Point Estimation}\label{sec:zero}

We begin with a construction based entirely on loss evaluations, which is particularly useful when derivative information is unavailable or costly to compute. A natural choice is to replace the full-data empirical risk by its minibatch counterpart, $\m L_{\rm naive}(\wX^{(m)},\theta)
=\frac1m\sum_{j=1}^m\ell(\wt X_j,\theta)$. Conditional on the observed data, its variance is $m^{-1}\operatorname{Var}_{X\sim\mu_n}[\ell(X,\theta)]$.
Since the surrogate enters the Gibbs exponent multiplied by $\alpha$, minibatch variability can lead to substantial posterior approximation error when $\alpha$ is large.  This motivates reducing the variability of the surrogate. Suppose we have a reference point $\theta^\dagger$, we can instead use a zero-order control variate and consider
\begin{equation*}
\underbrace{m^{-1}\sum_{j=1}^m \big[\ell(\wt X_j,\theta)-\ell(\wt X_j,\theta^\dagger)\big]}_{\text{estimate of }n^{-1}\sum_{i=1}^n [\ell(X_i,\theta)-\ell(X_i,\theta^\dagger)]}+\underbrace{\frac{1}{n}\sum_{i=1}^n \ell(X_i,\theta^\dagger)}_{\text{independent of }\theta}.
\end{equation*}
This remains an unbiased estimator of the full-data empirical risk. The first term averages loss differences evaluated on the same observations, which can be much less variable than the individual losses when $\theta$ is close to $\theta^\dagger$. The second term is independent of both $\theta$ and the minibatch and cancels from the MH acceptance ratio, giving the zero-order surrogate
\[
\m L_{\rm zcv}(\wX^{(m)},\theta)
=\frac1m\sum_{j=1}^m
\big[\ell(\wt X_j,\theta)-\ell(\wt X_j,\theta^\dagger)\big].
\]

Under a uniform Lipschitz condition on the loss, the conditional variance of this surrogate is bounded by a constant multiple of $\|\theta-\theta^\dagger\|_2^2/m$. When the Gibbs posterior concentrates around the empirical risk minimizer $\wh\theta$, a reference point close to $\wh\theta$ can therefore reduce minibatch variability throughout the region carrying most of the posterior mass. To quantify this effect, we assume
$\|\theta^\dagger-\wh\theta\|_2\leq C_{m,n}$, where $C_{m,n}$ may depend on $m$ and $n$. An accurate reference point may not be available initially. We therefore also allow an uninformative reference point with $C_{m,n}=\m O(1)$. In this regime, suitable risk scaling becomes particularly useful, and the results below show how a tempered target can still yield a stationary mean close to $\wh\theta$. Such a mean provides a point estimator in its own right and can also serve as a more informative reference point for subsequent refinement. Both the zero-order surrogate and its use with a derivative-free proposal require only loss evaluations. The no-control-variate case fits into the same analysis under the additional assumption
$\sup_{x\in\m X,\theta\in\Theta}|\ell(x,\theta)|\leq L_0$,
where $L_0$ is independent of $(m,n)$. Under this condition, the bounds of Theorem~\ref{th:1} also hold for $\m L_{\rm naive}$, with $C_{m,n}$ replaced by a fixed positive constant and the reference-point condition omitted. 

We now introduce the regularity conditions used in the approximation analysis. A key observation is that averaging over observations can smooth the population risk
$\m R(\theta)=\mb E_{X\sim\m P^*}[\ell(X,\theta)]$,
even when individual losses and the empirical risk are non-smooth. Our analysis exploits this population-level smoothness, together with identification conditions,  prior regularity, and averaged subgradient regularity. These conditions allow us to quantify the approximation error while retaining the loss-only construction of the zero-order surrogate.

 \vspace{0.5em}
\noindent\textbf{Assumption A.1 (Risk function)} \emph{ (1)   The population risk $\mathcal{R}(\theta) = \mb{E}_{X\sim \m P^*}[\ell(X,\theta)]: \Theta \to \mathbb{R}$ has a unique minimizer $\theta^*$ over $\Theta$, and $\theta^*$ is in the interior of $\Theta$. (2) $\m R$ is twice continuously differentiable on $\Theta$, with uniformly bounded mixed partial derivatives up to the second order. Moreover, there exist positive constants $(r, L)$ such that for all $\theta \in  B(\theta^*,r)$, it holds that $ \|\mathcal{H}({\theta}) - \mathcal{H}({\theta^*})\|_{\rm F} \leq L \|\theta - \theta^*\|_2$, where $\mathcal{H}({\theta}) = \operatorname{Hess}\mathcal{R}(\theta)$. (3)  There exist positive constants $\left(a, b\right)$ such that $\mathcal{H}({\theta^*}) \succcurlyeq a \mathbf{I}_d$ and $\mathcal{R}(\theta) - \mathcal{R}(\theta^*) \geq b \left\|\theta - \theta^*\right\|_2^2$, for all $\theta \in \Theta$.
}

Additionally, we impose the following regularity conditions on the smoothness of the prior distribution and the compactness of the parameter space.
\\\textbf{Assumption A.2 (Prior and parameter space)} \emph{(1) The parameter space $\Theta \subset \mb{R}^d$ is compact. (2)  There exist positive constants $\left(l, L\right)$ such that $\pi(\theta^*) \geq l$, and for all $\theta \in \Theta$, it holds that    $\left|\pi(\theta) - \pi(\theta^*)\right| \leq L \left\|\theta - \theta^*\right\|_2$. }

The compactness assumption on the parameter space is imposed primarily for technical convenience in our non-asymptotic analysis and is standard in the Bayesian literature~\citep{10.1214/12-EJS675,yang2012bayesian,MR4758128}. Next, we introduce the following Lipschitz continuity assumption on the loss function $\ell$.

\noindent\textbf{Assumption B.1 (Loss function)}\emph{ There exists a positive constant $L$, such that for any $x \in \mathcal{X}$ and $\theta, \theta' \in \Theta$, it holds that $
    \left|\ell(x, \theta) - \ell(x, \theta')\right| \leq L \left\|\theta - \theta'\right\|_2.$}

The global Lipschitz condition in $\theta$ implies that  $\ell(x,\cdot)$ admits a generalized (Clarke) subgradient. Motivated by this,  we introduce a function $g: \mathcal{X} \times \Theta \to \mathbb{R}^d$ such that $g(x,\theta)$ can be taken as a  subgradient of $\ell(x,\cdot)$ at $\theta$. We then assume that $g$ satisfies the following conditions.

\noindent\textbf{Assumption B.2 (Subgradient of loss function) }\emph{There exist positive constants $(L,L_1,L_2)$ such that (1) For any $\theta \in \Theta$, $\mb E_{X\sim \m P^*}[g(X, \theta)] = \nabla \mathcal{R}(\theta)$ and $\sup_{x \in \mathcal{X}} \|g(x, \theta)\|_2 \leq L$. (2) There exists a constant $\beta_1 \in (0,1]$ such that, for any $\theta, \theta'\in \Theta$, it holds that $\mb{E}_{X\sim \m P^*}[\|g(X, \theta) - g(X, \theta')\|_2^2] \leq L \left\|\theta - \theta'\right\|_2^{2\beta_1}$ and $\mb{E}_{X\sim \m P^*}\big[(\ell(X, \theta) - \ell(X, \theta') - g(X, \theta')\cdot (\theta - \theta'))^2\big] \leq L \left\|\theta - \theta'\right\|_2^{2 + 2\beta_1}$. (3) Define the pseudo-metric $d_n^g$ as
    \(d_n^g(\theta, \theta') = \big(\frac{1}{n} \sum_{i=1}^n \|g(X_i, \theta) - g(X_i, \theta')\|_2^2\big)^{\frac{1}{2}}.
    \) Then for any $n\in \mb N_{+}$ and $\bX^{(n)}\in \m X^{\otimes n}$, the $\varepsilon$-covering number of $\Theta$ with respect to $d_n^g$ is upper bounded by $1\vee (L_1\left(\frac{n}{\varepsilon}\right)^{L_2})$.}

\medskip
Note that the uniform boundedness  of $\|g(x, \theta)\|_2$ and of the Lipschitz constant of $\ell(x,\cdot)$ over $x\in \m X$ is imposed mainly for technical simplicity and can be relaxed. For example, one could instead assume $\|g(x, \theta)\|_2\leq b(x)$  with $\mb E_{X\sim \m P^*}[\exp(|b(X)|^{\gamma})]<\infty$ for some  constant $\gamma>0$. This relaxation would only introduce  additional $(\log n)^{{\rm poly}(\frac{1}{\gamma})}$ factors in the final rates.  Assumption $B.2(2)$ encodes a population-level smoothness condition on the subgradient, with $\beta_1$ quantifying the degree of smoothness. After taking expectations over $X$, the first term controls the Lipschitz behavior of $g(X,\cdot)$, while the second term governs the remainder term in the first-order Taylor expansion of $\ell(X,\cdot)$, where the gradient is replaced by $g(X,\cdot)$. The covering number condition in Assumption $B.2(3)$ then enables uniform control of the random fluctuation of the empirical average $\frac{1}{n}\sum_{i=1}^n g(X_i,\cdot)$ around $\nabla \m R(\cdot)$. These assumptions are trivially satisfied with $\beta_1=1$ when $\ell(X,\cdot)$ is differentiable with a Lipschitz continuous gradient. For non-smooth losses such as the check loss in quantile regression and the hinge loss in SVMs,~\cite{10.1111/rssb.12510, MR4758128} show that analogous conditions can hold with $\beta_1=0.5$. We are now ready to state our main result for  $\wt \pi_{\m L,\alpha}$ under the zero-order control variate.  Define the mapping $\wh\theta:\m X^{\otimes n}\to \Theta$ by $\wh\theta(\bX^{(n)})\in \arg\min_{\theta\in \Theta} \sum_{i=1}^n \ell(X_i,\theta)$ and let   $\theta^\dagger:\m X^{\otimes n}\to \Theta$ be the mapping  that assigns the reference point $\theta^\dagger=\theta^\dagger(\bX^{(n)})$ used in $\m L_{\rm zcv}$.

\begin{theorem}\label{th:1}
Suppose Assumptions $A.1,A.2$ and $B.1, B.2$ hold. Moreover, assume there exist positive constants $(C_1, C_2, C_3)$ such that
 $C_1 \log m \leq \alpha \leq C_2 \frac{m}{(\log m)^2}$ and $C_{m,n} \leq C_3 \sqrt{\frac{m}{\alpha^2 \log m}}$. Then,  there exists a set $\m B\subset \m X^{\otimes n}$ with $\m P^*{}^{\otimes n}(\m B)\geq 1-\frac{1}{n^2}$ such that, for any $\bX^{(n)}\in \m B$ satisfying $\|\wh\theta(\bX^{(n)})-\theta^\dagger(\bX^{(n)})\|_2\leq C_{m,n}$, we have
\[
\left\|\wt{\pi}_{\mathcal{L}_{\rm zcv},\alpha}(\theta) - \pi_\alpha(\theta)\right\|_{\mathrm{TV}} \lesssim \alpha \frac{(\log m)^2}{m} + \alpha^2 \frac{\log m}{m} C_{m,n}^2.
\]
Moreover, if in addition $ \alpha\geq C_4 (\log n)^{\frac{3}{2}}$ for some constant $C_4>0$, then
\[
\Big\|\int_{\Theta} \theta \, \wt{\pi}_{\mathcal{L}_{\rm zcv},\alpha}(\theta) \dd\theta - \wh\theta(\bX^{(n)})\Big\|_2 \lesssim \sqrt{\alpha} \frac{(\log m)^2}{m} + \alpha^{\frac{3}{2}} \frac{\log m}{m} C_{m,n}^2 + \frac{(\log n)^{\frac{3}{2}}}{\alpha} + \frac{(\log n)^{\frac{\beta_1 + 2}{2}}}{\sqrt{n} \alpha^{\frac{\beta_1}{2}}}.
\]
 \end{theorem}
 \begin{remark}\label{re:1}
If we take an uninformative $\theta^\dagger$ so that $C_{m,n}=\m O(1)$, then the point estimation error bound in Theorem~\ref{th:1} simplifies to $\|\int_{\Theta} \theta \, \wt{\pi}_{\mathcal{L}_{\rm zcv},\alpha}(\theta) \dd\theta - \wh\theta(\bX^{(n)})\|_2 \lesssim \alpha^{\frac{3}{2}} \frac{\log m}{m} + \frac{(\log n)^{\frac{3}{2}}}{\alpha}$. Thus, if no control variate is used, or a zero-order control variate is used with an uninformative reference point,  choosing $\alpha = \tilde{\mathcal{O}} (m^{\frac{2}{5}} )$ balances the terms in this upper bound, giving an order $\tilde{\mathcal{O}} (m^{-\frac{2}{5}})$ for the deviation between the corresponding stationary mean $\int_{\Theta} \theta \, \wt{\pi}_{\mathcal{L},\alpha}(\theta) \dd\theta$ and  the empirical risk minimizer $\wh\theta$.
\end{remark}

\begin{remark}
        The main technical challenge in proving Theorem~\ref{th:1} is the ``in-expectation'' analysis. Specifically, for $  \wt\pi_{\mathcal{L}_{\rm zcv},\alpha}(\theta)\propto \mb E_{\wX^{(m)} \sim \mu_n^{\otimes m}}\big[\pi(\theta) \exp  \big(-\alpha\cdot \m{L}_{\rm zcv}(\wX^{(m)}, \theta)  \big)\big]$, we analyze the expectation over the subsample $\wX^{(m)}$  rather than taking a supremum over $\wX^{(m)}$ of the discrepancy between the mini-batch surrogate and the full-batch empirical risk . This ``in-expectation'' analysis lifts the leading error to higher order and requires  a careful treatment  of higher-order terms.
 
\end{remark}

\subsection{First-Order Control Variates and Posterior Approximation}\label{sec:first}
The zero-order bounds describe how risk scaling and reference-point accuracy control approximation to a tempered posterior and estimation of its center. We now incorporate subgradient information to improve posterior approximation and identify regimes that permit direct approximation of the untempered Gibbs posterior. 

Given a reference point $\theta^\dagger$, we construct a first-order control variate using the subgradient $g(X,\theta^\dagger)$, by analogy with a first-order Taylor expansion, and consider 
\begin{equation*}
\underbrace{\frac{1}{m}\sum_{j=1}^m \left[\ell(\widetilde{X}_j, \theta) - \ell\left(\widetilde{X}_j, \theta^\dagger\right) - g(\widetilde{X}_j, \theta^\dagger)\cdot(\theta - \theta^\dagger)\right]}_{\text{estimate of }  \frac{1}{n}\sum_{i=1}^n [\ell(X_i, \theta) - \ell\left(X_i, \theta^\dagger\right) - g(X_i, \theta^\dagger)\cdot(\theta - \theta^\dagger)]}+\frac{1}{n} \sum_{i=1}^n \big[\ell(X_i,\theta^\dagger)+  g(X_i, \theta^\dagger)\cdot (\theta - \theta^\dagger)\big]
\end{equation*}
 as a variance-reduced approximation to the empirical risk. Since the term $n^{-1} \sum_{i=1}^n \ell(X_i,\theta^\dagger)$ is independent of $\theta$,  we  define the first-order control variate surrogate as:
\begin{equation*}
 \begin{aligned}
     &\mathcal{L}_{\rm fcv}(\wX^{(m)}, \theta)\\
     &=\Big(\frac{1}{m}\sum_{j=1}^m \left[\ell(\widetilde{X}_j, \theta) - \ell(\widetilde{X}_j, \theta^\dagger) - g(\widetilde{X}_j, \theta^\dagger)\cdot(\theta - \theta^\dagger)\right]\Big)_{+}  +\frac{1}{n}\sum_{i=1}^n  g(X_i, \theta^\dagger)\cdot (\theta - \theta^\dagger),
 \end{aligned}
 \end{equation*}
 where $a_+=\max(a,0)$. For technical reasons, we take the positive part in the first term. This is motivated by the fact that, if $\ell$ is continuously differentiable in $\theta$ and $\theta^{\dagger}=\wh\theta$, then $\sum_{i=1}^n \ell(X_i,\theta)-\sum_{i=1}^n\ell(X_i,\wh\theta)-\sum_{i=1}^n \nabla\ell(X_i,\wh\theta)\cdot(\theta-\wh\theta)$  is non-negative for all $\theta$.  The term $\frac{1}{n}\sum_{i=1}^n  g(X_i, \theta^\dagger) $ is independent of $\theta$ and can therefore be precomputed and treated as a constant in the algorithm. If $\ell(X,\cdot)$ is twice differentiable,  the variance of $\mathcal{L}_{\rm fcv}(\wX^{(m)}, \theta)$ is of order $\m O(\frac{\|\theta-\theta^\dagger\|_2^4}{m})$, which is smaller than that of $\mathcal{L}_{\rm zcv}(\wX^{(m)}, \theta)$ when $\|\theta-\theta^\dagger\|_2$ is small. We also show that, even for non-smooth loss functions, this first-order construction can substantially improve the posterior approximation. In addition, to obtain a sharper bound in the smooth case, we impose the following assumption.

\noindent\textbf{Assumption B.3 (Pointwise smoothness of subgradients)}
 \emph{There exists a constant $\beta_2 \in [0, \beta_1]$ and a positive constant $L$ such that, for all $x \in \mathcal{X}$ and $\theta, \theta' \in \Theta$, it holds that $       \|g(x, \theta) - g(x, \theta')\|_2 \leq L \|\theta - \theta'\|_2^{\beta_2}$  and
 $        |\ell(x, \theta) - \ell(x, \theta') - g(x, \theta')\cdot(\theta - \theta')| \leq L \| \theta - \theta' \|_2^{1+\beta_2}$.}

\vspace{0.5em}
 Note that when $\beta_2=\beta_1$, Assumption B.3 is stronger than Assumption B.2, which only imposes smoothness at the population level. Here, however, we allow $\beta_2<\beta_1$, and even $\beta_2=0$. It is easy to see that Assumption B.3  with $\beta_2=0$ follows automatically from Assumptions B.1 and B.2. On the other hand,  if $\ell(x,\cdot)$ is differentiable with  a $\beta_2$-H\"{o}lder smooth gradient for some $\beta_2\in (0,1]$,  then Assumption B.3 yields sharper bounds than in the non-smooth case $\beta_2=0$. We are now ready to state our main result for the first-order control variate construction.

\begin{figure}[htbp]
\centering
\includegraphics[width=0.43\linewidth]{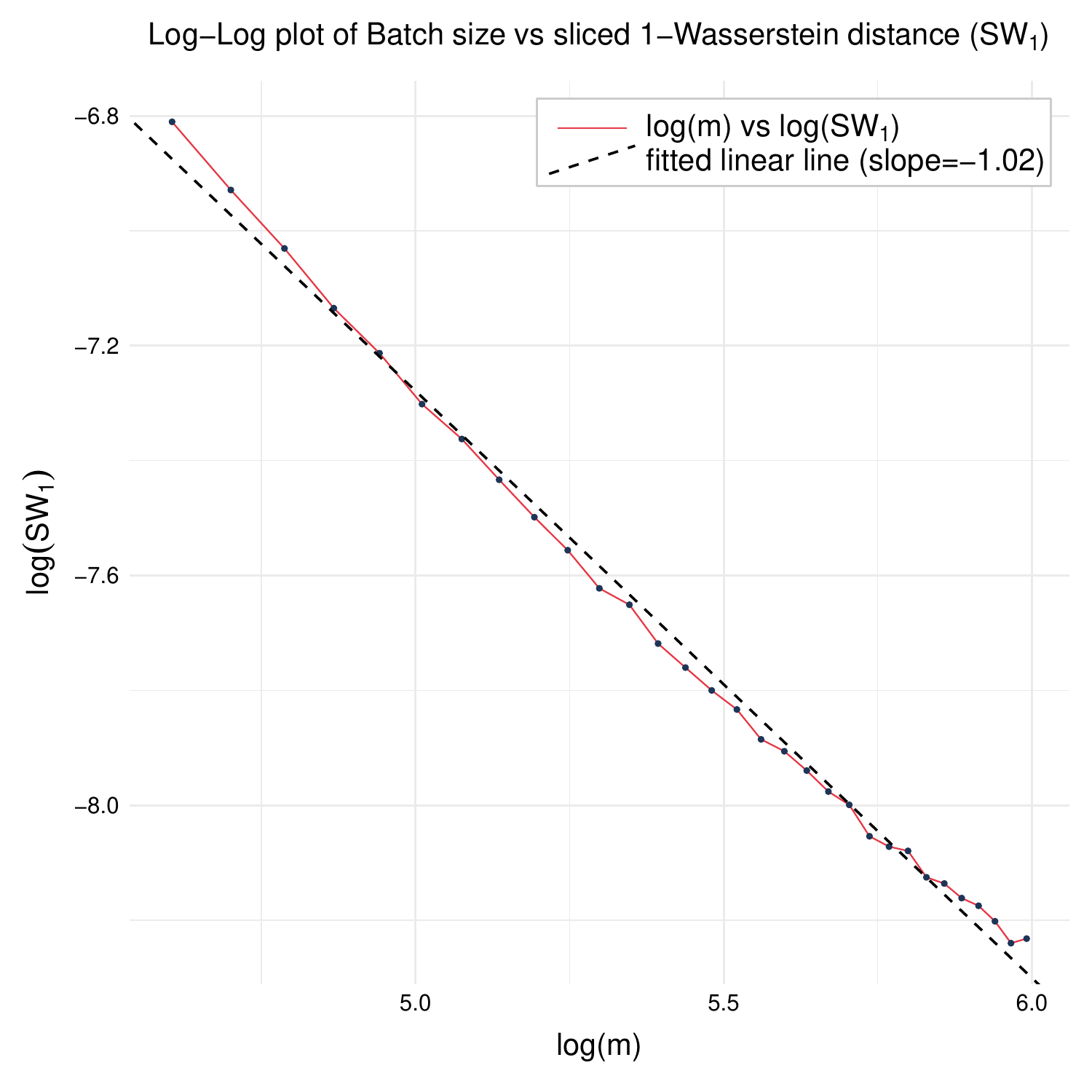}
\caption{To examine whether the empirical approximation error exhibits the $m^{-1}$ decay suggested by Theorem~\ref{th:2}, we use the Huber regression model described in Section~\ref{Bayesian Huber regression}, with parameter dimension $d=20$ and sample size $n=20000$. The log--log plot shows the sliced $1$-Wasserstein distance ($\mathrm{SW}_1$) between the empirical distributions obtained from first-order control-variate SMH and full-batch MH. Distances are estimated using $30000$ MH samples and averaged over $30$ independent runs. We use $\mathrm{SW}_1$ because it is convenient to estimate in high dimensions and satisfies $\mathrm{SW}_1\lesssim\mathrm{TV}$ for distributions supported on a common compact parameter space. The Huber loss has a Lipschitz-continuous gradient, so we compare the results with the theorem’s $\wt{\mathcal O}(m^{-1})$ bound for $\beta_1=\beta_2=1$. The fitted slope of $-1.02$ indicates approximately $m^{-1}$ decay, consistent with this theoretical rate.}
\label{SW Huber}
\end{figure}

\begin{theorem}\label{th:2}
Suppose Assumptions $A.1-A.2$ and  $B.1-B.3$ hold and let $\alpha \leq n$. Assume there exist positive constants $\left(\gamma, C_5, C_{6}, C_{7},C_{8}\right)$ such that:
$m \geq C_5\max\left(n^\gamma,\alpha^{\frac{1}{2}}\right)$, $C_{6} \log m \leq \alpha$, $C_{7} \max\left(\alpha^{1-\beta_1}, \alpha^{\frac{1-\beta_2}{2}}\right) \leq \frac{m}{(\log m)^3}$, and $C_{m,n} \leq \frac{C_{8}}{\sqrt{\alpha}}$.
Then,  there exists a set $\m B\subset \m X^{\otimes n}$ with $\m P^*{}^{\otimes n}(\m B)\geq 1-\frac{1}{n^2}$ such that, for any $\bX^{(n)}\in \m B$ satisfying $\|\wh\theta(\bX^{(n)})-\theta^\dagger(\bX^{(n)})\|_2\leq C_{m,n}$,  
\begin{equation*}
\left\|\wt{\pi}_{\mathcal{L}_{\rm fcv},\alpha}(\theta) - \pi_\alpha(\theta)\right\|_{\mathrm{TV}}  \lesssim \alpha^{1-\beta_1}\frac{(\log m)^{\beta_1+2}}{m} + \alpha^{1-\beta_2}\frac{(\log m)^{\beta_2+3}}{m^2}.
\end{equation*}
Moreover, if in addition $\alpha\geq C_{9} (\log n)^{\frac{3}{2}} $ for some constant $C_{9}>0$ , then 

\begin{equation*}
\begin{aligned}
     &\big\|\int_{\Theta} \theta \, \wt{\pi}_{\mathcal{L}_{\rm fcv},\alpha}(\theta) \dd\theta - \wh\theta(\bX^{(n)})\big\|_2\\
     &\lesssim \alpha^{\frac{1-2\beta_1}{2}}\frac{(\log m)^{\beta_1+2}}{m} + \alpha^{\frac{1-2\beta_2}{2}}\frac{(\log m)^{\beta_2+3}}{m^2} + \frac{(\log n)^{\frac{3}{2}}}{\alpha} + \frac{(\log n)^{\frac{\beta_1+2}{2}}}{\sqrt{n}\alpha^{\frac{\beta_1}{2}}}. 
\end{aligned}
\end{equation*}

\end{theorem}
\begin{remark}\label{re:untempered}
The condition $m\geq n^\gamma$ for some fixed $\gamma>0$ ensures that $n$ is at most polynomial in $m$. If an accurate reference point is available with $C_{m,n}=\m O(n^{-1/2})$, then, under the remaining conditions of Theorem~\ref{th:2}, setting $\alpha=n$ gives an approximation to the original, untempered Gibbs posterior $\pi_n$. The resulting error bound depends on the regularity of the subgradient:
\begin{enumerate}
\item \textbf{Uniformly Lipschitz gradients.}
When $\beta_1=\beta_2=1$, the bound becomes
\[
\|\wt\pi_{\m L_{\rm fcv},n}-\pi_n\|_{\rm TV}
=\wt{\m O}(m^{-1}).
\]
In particular, $m\asymp\sqrt n$ yields a total variation error of order $\wt{\m O}(n^{-1/2})$.   The numerical results in Figure~\ref{SW Huber} show approximately $m^{-1}$ decay in sliced Wasserstein distance, consistent with the minibatch-size dependence of this bound.

\item \textbf{Possibly discontinuous subgradients.}
Taking $\beta_2=0$, while retaining the averaged regularity indexed by $\beta_1>0$, gives
\[
\|\wt\pi_{\m L_{\rm fcv},n}-\pi_n\|_{\rm TV}
=\wt{\m O}\left(\frac{n^{1-\beta_1}}m+\frac{n}{m^2}\right).
\]
Both terms vanish if $m\asymp n^\delta$ with
$\max\{1-\beta_1,1/2\}<\delta<1$.
Thus consistent approximation remains possible with a minibatch size that grows more slowly than $n$, even when the subgradient is discontinuous.
\end{enumerate}
 
\end{remark}

 \subsection{Reference-Point Refinement and Comparison}\label{sec:refinement}
The first-order approximation results rely on an accurate reference point, which can be obtained by directly optimizing the empirical risk. Here, we develop a sampling-based alternative through staged refinement. When posterior sampling is the goal, this approach reuses pilot runs to estimate reference points and initialize subsequent stages.

We begin with the naive surrogate or a zero-order surrogate using an uninformative reference point. Suitable risk scaling allows this initial stage to produce a useful stationary-mean estimate using only loss evaluations. Taking this estimate as the reference point for a first-order stage reduces minibatch variability and permits a larger value of $\alpha$. The stationary-mean bounds in Theorems~\ref{th:1} and~\ref{th:2} quantify how this refinement improves reference-point accuracy and ultimately supports approximation of the untempered Gibbs posterior.

To make this refinement precise, suppose that $\theta^\dagger$ is the stationary mean obtained from the initial stage in Remark~\ref{re:1}, so that
$C_{m,n}=\wt{\m O}(m^{-2/5})$.
Choosing $\alpha$ of order $m^{4/5}$ gives
\[
\Big\|
\int_\Theta\theta\,
\wt\pi_{\m L_{\rm fcv},\alpha}(\theta)\,\dd\theta
-\wh\theta
\Big\|_2
=
\wt{\m O}\left(
m^{-(3+4\beta_1)/5}
+m^{-4/5}
+n^{-1/2}m^{-2\beta_1/5}
\right).
\]
Consequently, taking $m=\wt\Omega(
n^{5/8}+n^{5/(6+8\beta_1)})$ with suitable logarithmic factors, yields a stationary mean within $\m O(n^{-1/2})$ of $\wh\theta$. Using this refined mean as the reference point supplies the accuracy required in Remark~\ref{re:untempered}. A final first-order stage can then target the untempered posterior by setting $\alpha=n$.  

An accurate reference point also supports an alternative route to the untempered posterior using the zero-order construction. Here, the idea is to first approximate a tempered posterior and then adjust its spread by rescaling samples around their stationary mean. The transformation preserves the estimated center while contracting deviations from it by a factor of $\sqrt{\alpha/n}$. The following corollary quantifies the error of this rescaled approximation and allows us to compare it with the direct first-order approach.

\begin{corollary}\label{co:2}
Consider the transformation $G(\theta) = \frac{\sqrt{\alpha} \left(\theta - \theta_P\right)}{\sqrt{n}} + \theta_P$, where $\theta_P = \int_{\Theta} \theta \, \wt{\pi}_{\mathcal{L}_{\rm zcv},\alpha}(\theta) \dd\theta$. Under Assumptions $A.1-A.2$ and $B.1-B.2$, if there exist positive constants $\left(C_{10}, C_{11}, C_{12}\right)$ such that $C_{10} n^{\frac{1}{2}} (\log n)^{\frac{3}{2}} \leq \alpha \leq C_{11} \frac{m^2}{n (\log m)^4}$ and $C_{m,n} \leq C_{12} \frac{m^{\frac{1}{2}}}{(\log m)^{\frac{1}{2}} \alpha^{\frac{3}{4}} n^{\frac{1}{4}}}$. Then,  there exists a set $\m B\subset \m X^{\otimes n}$ with $\m P^*{}^{\otimes n}(\m B)\geq 1-\frac{1}{n^2}$ such that, for any $\bX^{(n)}\in \m B$ satisfying $\|\wh\theta(\bX^{(n)})-\theta^\dagger(\bX^{(n)})\|_2\leq C_{m,n}$,  
 
\[
\left\|G_{\#}\wt{\pi}_{\mathcal{L}_{\rm zcv},\alpha}(\theta) - \pi_n(\theta)\right\|_{\mathrm{TV}} \lesssim \sqrt{n\alpha} \frac{(\log m)^2}{m} + \sqrt{n} \alpha^{\frac{3}{2}} \frac{\log m}{m} C_{m,n}^2 + \sqrt{n} \frac{(\log n)^{\frac{3}{2}}}{\alpha} + \frac{(\log n)^{\frac{\beta_1 + 2}{2}}}{\alpha^{\frac{\beta_1}{2}}}.
\]
 
\end{corollary}
\begin{remark}
For a concrete comparison, consider $\beta_1=1$ and
$C_{m,n}=n^{-1/2}$. The bound in Corollary~\ref{co:2} simplifies to
\[
\|G_\#\wt\pi_{\m L_{\rm zcv},\alpha}-\pi_n\|_{\rm TV}
\lesssim
\sqrt{n\alpha}\frac{(\log m)^2}{m}
+\sqrt n\,\frac{(\log n)^{3/2}}{\alpha}.
\]
 The two terms reflect a tradeoff in choosing $\alpha$. Increasing $\alpha$ reduces the contribution from tempering and rescaling, but increases the contribution from subsampling, including its effect on the estimated center.
 Choosing $\alpha$ of order $m^{2/3}$, up to logarithmic factors, balances them and gives the upper bound
$\wt{\m O}(\sqrt n/m^{2/3})$.
 Thus $m\asymp n^\delta$ with $\delta>3/4$ suffices for consistency in this regime. With the same reference-point accuracy and uniformly Lipschitz gradients, the first-order construction permits $\alpha=n$ directly and gives the sharper bound
$\wt{\m O}(m^{-1})$.  This comparison quantifies how incorporating subgradient information improves posterior approximation and reduces the required minibatch size. The zero-order construction, meanwhile, provides both a loss-only starting point for reference-point refinement and, through rescaling, a route to posterior approximation when only loss evaluations are available.
    
\end{remark}

\section{Mixing time analysis}\label{Mixing time analysis}

Subsampling MH reduces the per-iteration cost, but a key question is whether this saving comes at the expense of slower convergence. Our numerical experiments, presented in Section~\ref{Numerical Study}, suggest that the effective sample size over a fixed number of iterations can remain close to that of full-batch MH when $m$ is sufficiently large, even if $m\ll n$; see Table~\ref{Huber cov CI ess} for example. Intuitively, if the full-batch chain mixes well and the minibatch approximations are sufficiently accurate, the subsampling chain should retain comparable mixing behavior while requiring less computation per iteration. Motivated by this observation, we establish warm-start total variation convergence bounds for a lazy version of subsampling MH under explicit geometric and proposal-comparison conditions.

We begin with the full-batch chain, which provides the benchmark for our analysis. Let
$q^*(\theta\mid\theta')=q(\theta\mid\theta',\bX^{(n)})$
denote the full-batch proposal distribution. The MH algorithm targeting the tempered Gibbs posterior
 $\pi_\alpha(\theta)
\propto
\pi(\theta)\exp\left(
-\frac{\alpha}{n}\sum_{i=1}^n\ell(X_i,\theta)
\right)$ has transition kernel
\begin{equation*}
\begin{aligned}
\m T^*(\theta',\dd\theta)
= 
A^*(\theta',\theta)q^*(\theta\mid\theta')\,\dd\theta +
\Big(
1-\int_\Theta
A^*(\theta',\vartheta)q^*(\vartheta\mid\theta')\,\dd\vartheta
\Big)\delta_{\theta'}(\dd\theta),
\end{aligned}
\end{equation*}
where $\delta_{\theta'}$ is the point mass at $\theta'$ and $A^*(\theta',\theta)
=
1\wedge
\frac{\pi_\alpha(\theta)q^*(\theta'\mid\theta)}
{\pi_\alpha(\theta')q^*(\theta\mid\theta')}$. Assumption C.1 below collects standard conditions for controlling the mixing time of this full-batch chain~\citep{pmlr-v134-chewi21a,10.2307/30243694,JMLR:v20:19-306,JMLR:v21:19-441,MR4758128}. Following common practice~\citep{lovasz1993random,JMLR:v21:19-441,MR4758128,pmlr-v134-chewi21a}, we work with its $\zeta$-lazy version and fix $\zeta=1/2$ for simplicity. At each iteration, this chain remains at its current state with probability $1/2$ and otherwise applies the MH transition above. Under Assumption C.1, a standard $s$-conductance argument~\citep{pmlr-v134-chewi21a,JMLR:v20:19-306,JMLR:v21:19-441} gives
\begin{equation}\label{full mixing}
\|\nu^*_k-\pi_\alpha\|_{\rm TV}
\leq
C M_0\frac{p_1}{\lambda\psi\wedge 1}
+
M_0\exp\left(
-c(\lambda\psi\wedge 1)^2\omega^2k
\right),
\end{equation}
where $C,c>0$ are universal constants and $\nu^*_k$ is the law of the $1/2$-lazy full-batch chain after $k$ iterations, initialized from the $M_0$-warm start specified in Assumption C.1. In large-sample settings admitting a Bernstein--von Mises approximation, conditions of this type have been verified for random-walk Metropolis (RWM)~\citep{10.2307/30243694} and the Metropolis-adjusted Langevin algorithm (MALA)~\citep{MR4758128}.

We now turn to the subsampling chain and explain how its convergence can be related to the full-batch bound. The main complication is that the chain evolves jointly in the parameter $\theta$ and the minibatch $\wX^{(m)}$. This dependence becomes manageable when the joint stationary distribution is close to the product distribution
$\pi_\alpha(\theta)\mu_n^{(m)}(\wX^{(m)})$.
Under this approximate product structure, the parameter and minibatch are approximately independent at stationarity, and the parameter marginal is close to the desired tempered posterior. If, in addition, the minibatch proposal
$q(\cdot\mid\theta,\wX^{(m)})$
is close to its full-batch counterpart
$q^*(\cdot\mid\theta)=q(\cdot\mid\theta,\bX^{(n)})$
on a high-probability region, the subsampling updates should behave similarly to the full-batch updates. The density-ratio, proposal-comparison, and localization conditions in Assumption C.2 make these ideas precise and allow us to transfer conductance bounds between the chains.   
As in the full-batch analysis, we consider a fully lazy variant. With probability $\zeta$, the entire augmented state remains unchanged. Otherwise, a parameter proposal is drawn using the current minibatch, a new minibatch is drawn independently, and the two proposals are accepted or rejected jointly. This variant preserves the same joint stationary distribution and allows us to apply standard conductance arguments to establish convergence bounds. 
Throughout the analysis, we take $\mu_n^{(m)}=\mu_n^{\otimes m}$. The transition kernel of the fully lazy subsampling chain is
\begin{equation}\label{eqn:tlazy}
\begin{aligned}
&\m T\big(
(\theta',\wX'^{(m)}),
(\dd\theta,\wX^{(m)})
\big)\\
&\quad=
(1-\zeta)
A\big(
(\theta',\wX'^{(m)}),
(\theta,\wX^{(m)})
\big)
\mu_n^{\otimes m}(\wX^{(m)})
q(\theta\mid\theta',\wX'^{(m)})\,\dd\theta\\
&\qquad+
\Big(
1-(1-\zeta)
\sum_{\mathsf x\in(\bX^{(n)})^{\otimes m}}
\int_\Theta
A\big(
(\theta',\wX'^{(m)}),
(\vartheta,\mathsf x)
\big)
\mu_n^{\otimes m}(\mathsf x)
q(\vartheta\mid\theta',\wX'^{(m)})\,\dd\vartheta
\Big)\\
&\qquad\qquad\cdot
\delta_{(\theta',\wX'^{(m)})}
(\dd\theta,\wX^{(m)}),
\end{aligned}
\end{equation}
where $A$ is the joint MH acceptance probability defined in Section~\ref{sec:smh}. We retain $\zeta$ in this expression and set $\zeta=1/2$ in the convergence analysis. 

To separate the effect of the approximate target from that of the minibatch proposal, we introduce a comparison chain with the same joint proposal and laziness but with stationary distribution
$\pi_\alpha(\theta)\mu_n^{\otimes m}(\wX^{(m)})$.
Its acceptance probability is
\begin{equation*}
\ov A\big(
(\theta',\wX'^{(m)}),
(\theta,\wX^{(m)})
\big)
=
1\wedge
\frac{
\pi_\alpha(\theta)
q(\theta'\mid\theta,\wX^{(m)})
}{
\pi_\alpha(\theta')
q(\theta\mid\theta',\wX'^{(m)})
}.
\end{equation*}
The minibatch probabilities cancel from this ratio because the independent proposal and the minibatch marginal of the comparison target both use $\mu_n^{\otimes m}$. The comparison chain therefore retains the minibatch proposal while preserving the desired parameter marginal $\pi_\alpha$. To relate the subsampling chain to this comparison chain, we control the difference between their joint stationary distributions. Specifically, we seek to establish
\[
\begin{aligned}
\wt\pi_{\m L,\alpha}^{\rm joint}
(\theta,\wX^{(m)})
&=
\frac{
\pi(\theta)
\exp\big(-\alpha\m L(\wX^{(m)},\theta)\big)
\mu_n^{\otimes m}(\wX^{(m)})
}{
\mb E_{\wX^{(m)}\sim\mu_n^{\otimes m}}
\left[
\int_\Theta
\pi(\vartheta)
\exp\big(-\alpha\m L(\wX^{(m)},\vartheta)\big)
\,\dd\vartheta
\right]
}\\
&\approx
\pi_\alpha(\theta)\mu_n^{\otimes m}(\wX^{(m)}).
\end{aligned}
\]
Assumption C.2 quantifies this approximation on a high-probability region and controls the difference between the minibatch proposal $q(\cdot\mid\theta,\wX^{(m)})$ and the full-batch proposal $q^*(\cdot\mid\theta)$. These conditions address both the change in the joint stationary distribution and the change in the proposal introduced by subsampling. We collect the assumptions below, with C.1 governing the full-batch chain and C.2 specifying the comparisons needed to transfer its convergence bounds to the subsampling chain.

\medskip
\noindent\textbf{Assumption C.}
There exist sets $K,E_1,E_2\subset\Theta$,
$\m A\subseteq(\bX^{(n)})^{\otimes m}$, positive constants
$M_0,\lambda,\psi,\varepsilon_1,p_1,p_2$, and
$\omega\in(0,1]$ such that the following conditions hold.
\begin{enumerate}
\item \textbf{Full-batch chain.}
\begin{enumerate}
\item The initial distribution $\nu_0$ is an $M_0$-warm start of
$\pi_\alpha$, that is, $\nu_0(E)\leq M_0\pi_\alpha(E)$ for every Borel set $E\subset\mathbb R^d$.
The initial joint distribution is
$\nu_0\otimes\mu_n^{\otimes m}$.

\item For every $\theta,\theta'\in E_1$ satisfying
$\|\theta-\theta'\|_2\leq\psi$, $\|\m T^*(\theta,\cdot)-\m T^*(\theta',\cdot)\|_{\rm TV}
\leq 1-\omega.$

\item The restriction of $\pi_\alpha$ to $K$, defined by $\pi_\alpha|_K(\cdot)
=
\frac{\pi_\alpha(\cdot\cap K)}{\pi_\alpha(K)}$, satisfies
\[
\pi_\alpha|_K(S_3)
\geq
\lambda t
\min\{\pi_\alpha|_K(S_1),\pi_\alpha|_K(S_2)\}
\]
for every partition $K=S_1\cup S_2\cup S_3$ with
$\inf_{x\in S_1,z\in S_2}\|x-z\|_2\geq t$.

\item The sets $K$ and $E_1$ satisfy $\pi_\alpha(K)\wedge\pi_\alpha(E_1)\geq 1-p_1.$
 
\end{enumerate}

\item \textbf{Joint-target and proposal comparisons.}
\begin{enumerate}
\item The set $E_2\times\m A$ has high probability under both
joint distributions,
\[
\wt\pi_{\m L,\alpha}^{\rm joint}(E_2\times\m A)
\wedge
(\pi_\alpha\otimes\mu_n^{\otimes m})(E_2\times\m A)
\geq 1-p_2.
\]

\item There exists a set $E_3\supseteq E_2$ such that
\begin{enumerate}
\item For every $\theta\in E_2$ and $\wX^{(m)}\in\m A$,
\[
\mu_n^{\otimes m}(\m A)
\int_{E_3}
q(\theta'\mid\theta,\wX^{(m)})\,\dd\theta'
\geq 1-\frac{\omega}{4}.
\]

\item For every $\theta\in E_3$ and $\wX^{(m)}\in\m A$,
\begin{enumerate}
\item For every $\theta'\in E_3$,
\[
\Big|
\log
\frac{q^*(\theta\mid\theta')}
{q(\theta\mid\theta',\wX^{(m)})}
\Big|
\leq \frac{\omega}{100}.
\]

\item The joint densities satisfy
\[
\big|
\log\wt\pi_{\m L,\alpha}^{\rm joint}(\theta,\wX^{(m)})
-
\log\big(
\pi_\alpha(\theta)\mu_n^{\otimes m}(\wX^{(m)})
\big)
\big|
\leq\varepsilon_1.
\]
\end{enumerate}
\end{enumerate}
\end{enumerate}
\end{enumerate}
Assumption C.2 specifies the region and accuracy of the comparison between the two chains. Part C.2.a ensures that $E_2\times\m A$ has high probability under both joint targets, while Part C.2.b.i ensures that proposals from this region remain in $E_3\times\m A$ with high probability. When the proposal
makes local moves, $E_3$ can be viewed as an enlargement of $E_2$ that accommodates these moves. Within this enlarged region, Parts C.2.b.ii.A and C.2.b.ii.B control the proposal and joint-target density ratios respectively. Lemma~\ref{le:1} in the supplement establishes the joint-target comparisons for the zero-order and first-order constructions under their respective assumptions and scaling conditions.
Together, these conditions allow us to transfer the full-batch conductance bound to the subsampling chain. Let $\nu_k$ denote the parameter marginal after $k$ iterations of the lazy subsampling MH algorithm with transition kernel~\eqref{eqn:tlazy} and $\zeta=1/2$.

\begin{theorem}\label{th:3}
Suppose Assumption C holds with positive constants
$(M_0,\lambda,\psi,\varepsilon_1,p_1,p_2,\omega)$ and $q(\theta\mid\theta',\wX'^{(m)})>0$
 for every $\theta,\theta'\in\Theta$ and
$\wX'^{(m)}\in(\bX^{(n)})^{\otimes m}$.
Then there exists a universal constant $c>0$ such that
\begin{equation}\label{mixing time}
\begin{aligned}
\|\nu_k-\wt\pi_{\m L,\alpha}\|_{\rm TV}
&\leq
M_0\frac{24\exp(2\varepsilon_1)}
{\lambda\psi\wedge 1}
\big(\exp(\varepsilon_1)p_1+p_2\big)\\
&\quad+
M_0\exp(\varepsilon_1)
\exp\left(
-c k\exp(-8\varepsilon_1)
(\lambda\psi\wedge 1)^2\omega^2
\right).
\end{aligned}
\end{equation}
\end{theorem}

\noindent The first term accounts for probability mass outside the comparison regions, while the second decays exponentially with the number of iterations. When $\varepsilon_1=\m O(1)$ and $p_2=\m O(p_1)$, the bound matches the order of the full-batch bound~\eqref{full mixing}, including its dependence on the geometric parameters $(\lambda,\psi,\omega)$. Combining this convergence bound with the stationary approximation bounds in Section~\ref{Construction} controls the distance to $\pi_\alpha$ after finitely many iterations.

We next specialize this result to minibatch versions of a
MALA-type proposal based on the empirical risk. The full-batch
proposal is $q^*(\cdot\mid\theta')
=
N\big(
\theta'-\frac{h}{n}\sum_{i=1}^n g(X_i,\theta'),
\frac{2h}{\alpha}I_d
\big)$, where $h>0$ is bounded above by a constant independent of $(n,m)$.
The factor $1/\alpha$ in the covariance adjusts the proposal
variance to the scale of the tempered posterior. We consider two
minibatch proposals corresponding to the two surrogate
constructions.

\medskip
\noindent\textbf{Case (i) (Zero-order).}\label{case:zero_order}
For $\m L=\m L_{\rm zcv}$, we pair the zero-order surrogate with
the uncorrected minibatch subgradient proposal
\[
q_{\rm zcv}(\cdot\mid\theta',\wX^{(m)})
=
N\big(
\theta'-\frac{h}{m}\sum_{j=1}^m g(\wt X_j,\theta'),
\frac{2h}{\alpha}I_d
\big).
\]
\noindent\textbf{Case (ii) (First-order).}\label{case:first_order}
For $\m L=\m L_{\rm fcv}$, we also use the reference point to reduce
the variability of the proposal drift. Specifically, let
\[
q_{\rm fcv}(\cdot\mid\theta',\wX^{(m)})
=
N\big(
\mu(\wX^{(m)},\theta'),
\frac{2h}{\alpha}I_d
\big),
\]
where $\mu(\wX^{(m)},\theta')
=
\theta'-h\big[
 \frac1m\sum_{j=1}^m
\big(
g(\wt X_j,\theta')-g(\wt X_j,\theta^\dagger)
\big) +\frac1n\sum_{i=1}^n g(X_i,\theta^\dagger)
\big].$
 
 \medskip\noindent
The minibatch term estimates the change in the full-data
subgradient between $\theta^\dagger$ and $\theta'$, while the
full-data term at $\theta^\dagger$ can be precomputed.
Consequently, $\mu(\wX^{(m)},\theta')$ is a variance-reduced
approximation to the full-batch proposal center
$\theta'-h n^{-1}\sum_{i=1}^n g(X_i,\theta')$. For these proposals, the concentration bounds used in the
approximation analysis also control the proposal discrepancies.
Verifying Assumption C.2 and applying Theorem~\ref{th:3} gives the
following convergence result.

 \begin{theorem}\label{th:4}
Consider positive constants $(M_0,\lambda,\psi,p_1,\omega)$, where $p_1\geq m^{-\gamma_1}$ for some $\gamma_1>0$, and let $\mu_n^{(m)}=\mu_n^{\otimes m}$. There exists a set $\m B\subset\m X^{\otimes n}$ with $\m P^*{}^{\otimes n}(\m B)\geq 1-n^{-2}$ such that, for every $\bX^{(n)}\in\m B$ satisfying $\|\wh\theta(\bX^{(n)})-\theta^\dagger(\bX^{(n)})\|_2\leq C_{m,n}$ and for which Assumption C.1 holds with $(M_0,\lambda,\psi,p_1,\omega)$, there exist universal constants $c_0,c_1>0$ such that
\[
\|\nu_k-\wt\pi_{\m L,\alpha}\|_{\rm TV}
\leq
\frac{c_0M_0p_1}{\lambda\psi\wedge 1}
+c_0M_0\exp\left(-c_1(\lambda\psi\wedge 1)^2\omega^2k\right)
\]
in each of the following cases.
\begin{enumerate}
\item \textbf{Zero-order control variate.}
Take $\m L=\m L_{\rm zcv}$ and $q(\cdot\mid\theta',\wX^{(m)})=q_{\rm zcv}(\cdot\mid\theta',\wX^{(m)})$. Suppose Assumptions A.1--A.2 and B.1--B.2 hold, with $m\geq C_0/\sqrt\omega$, $C_{m,n}\leq \sqrt m/(C_0\alpha\sqrt{\log m})$, and $C_1\log m\leq\alpha\leq \omega^2m/[C_1(\log m)^2]$, for sufficiently large constants $C_0,C_1$ independent of $(n,m)$.

\item \textbf{First-order control variate.}
Take $\m L=\m L_{\rm fcv}$ and $q(\cdot\mid\theta',\wX^{(m)})=q_{\rm fcv}(\cdot\mid\theta',\wX^{(m)})$. Suppose Assumptions A.1--A.2 and B.1--B.3 hold, with $m\geq n^\gamma$ for some $\gamma>0$. Moreover, assume $C_2\log m\leq\alpha\leq n$, $m\geq C_2\max\{\sqrt\alpha,1/\sqrt\omega\}$, $C_{m,n}\leq 1/(C_2\sqrt\alpha)$, and
$C_3\max\{\alpha^{1-\beta_1}/\omega^2,\alpha^{(1-\beta_2)/2}/\omega\}\leq m/(\log m)^3$,
for sufficiently large constants $C_2,C_3$ independent of $(n,m)$.
\end{enumerate}
\end{theorem}
Theorem~\ref{th:4} shows that, when the minibatch size, risk scaling, and reference-point accuracy satisfy the stated conditions, the convergence bound for the  lazy subsampling chain retains the full-batch dependence on $(\lambda,\psi,\omega)$. This is achieved at a lower per-iteration cost, reduced from $\m O(n)$ to $\m O(m)$ for fixed dimension after precomputing the reference-point quantities.
In particular, suppose $\omega$ is bounded away from zero, $\alpha=n$, and the reference point satisfies the required $\m O(n^{-1/2})$ accuracy condition. A sufficient minibatch order for the first-order construction is
\[
m\gtrsim(\log n)^3 n^{\max\{1/2,\,1-\beta_1\}}.
\]
Thus, sublinear minibatches can preserve the order of the full-batch convergence bound while reducing the cost per iteration.

\section{Numerical Studies}\label{Numerical Study}

Our numerical studies examine the practical implications of the approximation theory and explore the constructions in broader settings. Huber regression compares zero-order rescaling with the first-order method when the loss has a Lipschitz-continuous gradient. Ranking from noisy pairwise comparisons assesses first-order control variates for a non-smooth hinge loss. Robust classification evaluates posterior accuracy and
sampling efficiency for a discontinuous indicator loss,
using a smooth surrogate to construct control variates. The simulator-based mixture example explores zero-order
control variates for model fitting with a discontinuous loss.

The theory guides the scaling of $\alpha$ with the minibatch size $m$ but leaves the multiplicative constant $k_1$ unspecified. A practical approach is to tune $k_1$ jointly with the scaled proposal step size $h$ through alternating pilot runs. Starting from an initial value of $k_1$, we calibrate $h$ using the acceptance rate,  then hold $h$ fixed and reduce $k_1$ if the estimated ESS falls below a prespecified threshold. We repeat these updates to account for the interaction between risk scaling and proposal size. When $\alpha$ is fixed by the target distribution, only $h$ is tuned.  We illustrate the procedure with Huber regression and provide implementation details in Appendix~\ref{details of Huber Regression}.

\subsection{Bayesian Huber regression}\label{Bayesian Huber regression}
We consider the tempered Gibbs posterior for Huber regression with $\theta\in \mb R^d$:  $\pi_\alpha(\theta ) \propto \pi(\theta)\exp\big(- \frac{\alpha}{n}\sum_{i=1}^n H(Y_i - {X}_i^T \theta) \big)$, where $H(z)=\frac{1}{2}z^2\textbf{1}(|z|\leq 2)+2(|z|-1)\textbf{1}(|z|>2)$. The Huber loss is not twice continuously differentiable, but its gradient is Lipschitz, so
this example illustrates the smooth-gradient regime.  We set $n=20000$, $m\in \{200, 400, 800, 1200, 1600, 2000\}$ and $d=2,5,10$, with a uniform prior on   $[-10,10]^d$. Data-generation details are given in Appendix~\ref{details of Huber Regression}. We compare zero-order control-variate subsampling MH with rescaling, using $\alpha=m^{2/3}$ as suggested by Corollary~\ref{co:2}, and first-order control-variate subsampling MH with $\alpha=n$. Both methods use minibatch MALA proposals and the ordinary least squares (OLS) estimator as the reference point $\theta^\dagger$.  
Table~\ref{Huber cov CI ess} reports the empirical coverage and length of the $95\%$ credible interval for the first coordinate of $\theta$, together with ESS averaged across coordinates and ESS per second, based on $1000$ replications. With $m=200$, the first-order method already yields coverage, interval length, and ESS comparable to those of full-batch MH. The zero-order method with rescaling requires larger minibatches for comparable performance and tends to produce longer credible intervals.
 
\begin{table}[h]
\centering
 
\renewcommand{\arraystretch}{1.0}
\setlength{\tabcolsep}{3pt}
\setlength{\abovecaptionskip}{2pt}
\setlength{\belowcaptionskip}{3pt}

\caption{Coverage probability (\%), 95\% CI length $(\times 10^{-2})$, ESS, and ESS generated per second in Bayesian Huber regression. Subscripts $_0$ and $_1$ denote the zero-order and the first-order control variate methods, respectively. ``Full MH'' refers to the full-batch MH algorithm.}

\begin{tabular}{|lccccccc|}
\hline
& $m=200$ & $m=400$ & $m=800$ & $m=1200$
& $m=1600$ & $m=2000$ & Full MH \\
\hline
\hline
\multicolumn{8}{|c|}{$d=2$}\\
Coverage$_0$ & 93.1 & 94.4 & 95.7 & 95.0 & 95.0 & 94.9
& \multirow{2}{*}{94.3} \\
Coverage$_1$ & 94.9 & 94.9 & 94.8 & 94.7 & 94.9 & 94.6 & \\
\hline
Length$_0$ & 3.20 & 3.10 & 3.04 & 3.01 & 2.99 & 2.97
& \multirow{2}{*}{2.85} \\
Length$_1$ & 2.89 & 2.88 & 2.87 & 2.86 & 2.86 & 2.86 & \\
\hline
ESS$_0$ & 5045 & 6245 & 7417 & 8044 & 8393 & 8771
& \multirow{2}{*}{12196} \\
ESS$_1$ & 11550 & 11875 & 12009 & 12115 & 12081 & 12151 & \\
\hline
ESS\_per\_sec$_0$ & 989 & 1177 & 1298 & 1319 & 1296 & 1282
& \multirow{2}{*}{1026} \\
ESS\_per\_sec$_1$ & 2118 & 2062 & 1879 & 1737 & 1597 & 1494 & \\
\hline
\hline
\multicolumn{8}{|c|}{$d=5$}\\
Coverage$_0$ & 91.7 & 95.1 & 94.2 & 95.7 & 95.0 & 95.1
& \multirow{2}{*}{94.2} \\
Coverage$_1$ & 94.9 & 94.7 & 94.8 & 94.4 & 94.2 & 94.2 & \\
\hline
Length$_0$ & 3.24 & 3.13 & 3.05 & 3.02 & 3.00 & 2.98
& \multirow{2}{*}{2.86} \\
Length$_1$ & 2.92 & 2.89 & 2.87 & 2.87 & 2.86 & 2.86 & \\
\hline
ESS$_0$ & 2644 & 3434 & 4205 & 4650 & 4930 & 5151
& \multirow{2}{*}{7768} \\
ESS$_1$ & 6724 & 7218 & 7511 & 7584 & 7644 & 7664 & \\
\hline
ESS\_per\_sec$_0$ & 510 & 634 & 707 & 725 & 713 & 701
& \multirow{2}{*}{529} \\
ESS\_per\_sec$_1$ & 1211 & 1219 & 1118 & 1018 & 940 & 868 & \\
\hline
\hline
\multicolumn{8}{|c|}{$d=10$}\\
Coverage$_0$ & 88.2 & 92.4 & 94.6 & 95.1 & 95.6 & 95.5
& \multirow{2}{*}{96.1} \\
Coverage$_1$ & 96.3 & 96.4 & 96.0 & 95.8 & 95.9 & 95.9 & \\
\hline
Length$_0$ & 3.31 & 3.16 & 3.07 & 3.03 & 3.01 & 2.99
& \multirow{2}{*}{2.86} \\
Length$_1$ & 2.96 & 2.91 & 2.88 & 2.87 & 2.87 & 2.87 & \\
\hline
ESS$_0$ & 1436 & 2025 & 2601 & 2929 & 3152 & 3320
& \multirow{2}{*}{5544} \\
ESS$_1$ & 3995 & 4691 & 5092 & 5261 & 5319 & 5361 & \\
\hline
ESS\_per\_sec$_0$ & 268 & 357 & 416 & 429 & 428 & 422
& \multirow{2}{*}{289} \\
ESS\_per\_sec$_1$ & 691 & 751 & 715 & 660 & 610 & 567 & \\
\hline
\end{tabular}
\label{Huber cov CI ess}
\end{table}


 
 

 

\subsection{Ranking from Noisy Pairwise Comparisons}\label{ranking}
We consider the preference-graph hinge objective~\citep{pahikkala2009preference}, based on observed pairwise comparisons $C_{ij} \in \{-1, 1\}$ ($i < j$) satisfying $C_{ji} = -C_{ij}$. Here, $C_{ij}=1$ indicates that item $i$ is preferred to item $j$. The associated Gibbs posterior is
\[
\pi_\alpha(\theta)\propto\pi(\theta)\exp(-\alpha R_n(\theta)),
\qquad
R_n(\theta)=\frac{6}{n(n-1)}
\sum_{i<j}[1-C_{ij}\theta^\top(X_i-X_j)]_+,
\]
where $X_i\in\mb R^d$ is the covariate vector for item $i$, $\theta^\top X_i$ is its model score, and $[u]_+=\max\{u,0\}$. We set $n=1000$, $d=3$, and use a uniform prior on $[-10,10]^3$. The data generation details are provided in Appendix \ref{details of Ranking problem}.  The hinge loss is non-differentiable at unit margin, making this a non-smooth example. The empirical risk also has a $U$-statistic-type structure, averaging pairwise losses rather than individual-observation losses. We implement the SMH construction by replacing observation-wise averages in the surrogate with the corresponding pairwise averages within each minibatch. Direct pairwise evaluation involves $\m O(n^2)$ terms for the full dataset and $\m O(m^2)$ terms for a minibatch of $m$ items, providing additional computational motivation for subsampling.

For each $m\in\{100,200,250,300\}$, we follow the staged-refinement strategy in Section~\ref{sec:refinement}. We first obtain an initial reference point from a no-control-variate pilot with $\alpha=10m^{2/5}$, then refine it using a first-order pilot with $\alpha=10m^{4/5}$. The resulting reference point is used in the final first-order chain with $\alpha=n$. Each pilot runs for $2000$ iterations, discarding the first $1000$, and its retained sample mean provides both the reference point and the initial state for the next stage. The final chain and full MH each run for $10000$ iterations, discarding the first $2000$. All methods use the random-walk proposal $N(\theta,1.5^2\mathbf{I}_3/\alpha)$.
Table~\ref{tab:ranking-refinement} illustrates the improvement across stages. We compute the absolute coordinate-wise differences between each stage’s sample mean and the full-MH sample mean (with $\alpha=n$), divide by the corresponding full-MH posterior standard deviations, and average across coordinates. The reported values are averaged over $1000$ independent replications. At $m=200$, the standardized deviation decreases from $3.992$ for the initial reference point to $0.233$ after refinement, and reaches $0.058$ for the final posterior mean. Similar improvements occur across all four minibatch sizes, supporting the use of preliminary sampling stages to construct more accurate reference points. The performance of the final chains is summarized in Table~\ref{tab:ranking-fast}, which reports coverage, credible-interval length, ESS, and ESS per second over the same $1000$ replications. At $m=200$, the subsampling method achieves coverage, interval length, and ESS comparable to those of full-batch MH, while requiring substantially less total computation time.

\begin{table}[htbp]
\centering
\small
\caption{ Ranking with noisy pairwise comparisons. Mean deviations at each refinement and final sampling stage, standardized by the corresponding posterior standard deviations from full MH with $\alpha=n$ and averaged across coordinates and $1000$ independent replications.}
\label{tab:ranking-refinement}
\begin{tabular}{lrrrr}
\toprule
 & $m=100$ & $m=200$ & $m=250$ & $m=300$\\
\midrule
Initial reference (stage 1) & 6.725 & 3.992 & 3.481 & 2.978\\
Refined reference (stage 2) & 0.736 & 0.233 & 0.171 & 0.144\\
Final posterior mean       & 0.304 & 0.058 & 0.056 & 0.055\\
\bottomrule
\end{tabular}
\end{table}

\begin{table}[htbp]
\renewcommand{\arraystretch}{1.0}
\setlength{\tabcolsep}{3pt}
\setlength{\abovecaptionskip}{2pt}
\setlength{\belowcaptionskip}{3pt}
\centering
\caption{ Ranking with noisy pairwise comparisons. Coverage and credible-interval length concern the first coordinate. Coverage is evaluated at the population hinge-risk minimizer, $\theta^*\approx0.86\boldsymbol{1}_3$, obtained by numerically solving the population first-order condition. ESS is averaged over the three coordinates, and total runtime includes all initialization and sampling stages.}
\label{tab:ranking-fast}
 
\setlength{\tabcolsep}{4pt}
\begin{tabular}{lrrrrr}
\toprule
 & $m=100$ & $m=200$ & $m=250$ & $m=300$ & Full MH\\
\midrule
Coverage (\%) &94.9&97.5&97.6&97.7&97.6\\
95\% interval length &0.1579&0.1467&0.1462&0.1460&0.1453\\
ESS &332.1&419.8&420.5&420.1&424.3\\
ESS\_per\_second   & 275.2 & 119.1 & 80.3 & 60.5 & 12.1 \\
\bottomrule
\end{tabular}
 
\end{table}

\subsection{Robust Classification under Gaussian Perturbations}

This example examines latent-variable subsampling with a discontinuous
indicator loss. We consider the Gibbs posterior
\[
\pi_{\alpha}(\theta)
\propto
\pi(\theta)\exp\Big(
-\frac{\alpha}{n}\sum_{i=1}^{n}
\mathbb{E}_{Z}
\left[
\mathbf{1}\{Y_i f_{\theta}(X_i+Z)<0\}
\right]
\Big),
\]
where $\theta=(\theta_1,\theta_2,\theta_3)^\top\in\mathbb{R}^3$,
$f_{\theta}(x)=x_1+\theta_1x_2+\theta_2x_1^2+\theta_3x_2^2$,
and $Z\sim N(0,0.1\mathbf{I}_2)$.
We generate $n=2000$ independent covariates
$X_i\sim N(0,\mathbf{I}_2)$ and responses satisfying
$\mathbb{P}(Y_i=1\mid X_i)=\sigma(f_{\theta^*}(X_i))$,
where $\theta^*=(1,1,-1)^\top$,
$\sigma(t)=1/(1+\exp(-t))$, and $Y_i\in\{-1,1\}$.
We use the Gaussian prior $\theta\sim N(0,\mathbf{I}_3)$
and set $\alpha=n$ for the final sampling stage.

Computing this posterior involves an expectation over Gaussian
perturbations and a discontinuous indicator integrand.
To approximate the expectation, for a perturbation batch
$\mathcal Z_s=(Z_1,\ldots,Z_s)$, define
\[
\widehat R_s(\theta;\mathcal Z_s)
=
\frac{1}{ns}\sum_{i=1}^{n}\sum_{j=1}^{s}
\mathbf{1}\{Y_i f_{\theta}(X_i+Z_j)<0\}.
\]
Within each replication, we generate a fixed bank $\mathcal Z_N$
of $N=20000$ independent perturbations. Full-batch MH
(denoted as full MH) targets the resulting finite-bank posterior
with density proportional to
$\pi(\theta)\exp(-n\widehat R_N(\theta;\mathcal Z_N))$.
The subsampling method reduces the number of perturbations used
per iteration while retaining all $n$ observed data points.

To construct a first-order control variate for this discontinuous
loss, we obtain gradient information from the smooth surrogate
\[
\widehat R_{\tau,s}(\theta;\mathcal Z_s)
=
\frac{1}{ns}\sum_{i=1}^{n}\sum_{j=1}^{s}
\sigma\left(-\frac{Y_i f_{\theta}(X_i+Z_j)}{\tau}\right),
\qquad \tau=0.3.
\]
This surrogate is used only to construct the linear correction.
The risk evaluations retain the original indicator loss. For each minibatch size $m\in\{200,400,600,800,1000\}$,
we first obtain a reference point through a loss-only pilot.
Starting at $(5,5,5)^\top$, we run subsampling MH without
control variates for $2000$ iterations with
$\alpha=50m^{2/5}$ and $\zeta=1/2$, discard the first $1000$ iterations,
and take the retained sample mean as $\theta^\dagger$.
We then use the minibatch surrogate
\[
\begin{aligned}
\mathcal L_m(\theta;\mathcal Z_m)
={}&
\widehat R_m(\theta;\mathcal Z_m)
-\widehat R_m(\theta^\dagger;\mathcal Z_m)\\
&-
\left[
\nabla_\theta\widehat R_{\tau,m}
    (\theta^\dagger;\mathcal Z_m)
-
\nabla_\theta\widehat R_{\tau,N}
    (\theta^\dagger;\mathcal Z_N)
\right]^\top
(\theta-\theta^\dagger).
\end{aligned}
\]
The resulting first-order control-variate method, denoted as SMH,
starts at $\theta^\dagger$. Full MH starts at $(5,5,5)^\top$.
Both main chains run for $10000$ iterations with $\alpha=n$,
discarding the first $1000$ iterations as burn-in.
All runs use the random-walk proposal
$q(\theta'\mid\theta)=N(\theta,20\mathbf{I}_3/n)$. Table~\ref{tab:robust-classification-results} reports the length
of the $95\%$ credible interval for the first coordinate,
ESS, ESS per second, and standardized mean deviation.
Within each replication, we define the mean deviation as $D_{\mathrm{mean}}
=
\frac{1}{3}\sum_{j=1}^{3}
\frac{
\left|
\bar\theta_{\mathrm{SMH},j}
-
\bar\theta_{\mathrm{full},j}
\right|
}{
s_{\mathrm{full},j}
},$ where $\bar\theta_{\mathrm{SMH},j}$ and
$\bar\theta_{\mathrm{full},j}$ are the retained sample means
for coordinate $j$, and $s_{\mathrm{full},j}$ is the
corresponding posterior standard deviation estimated from
the full-MH samples at $\alpha=n$.
ESS and ESS per second are also averaged across coordinates.
All reported metrics are averaged over $100$ independent
replications. The SMH runtime includes both reference-point
construction and the final sampling stage.
Across the minibatch sizes considered, SMH produces interval
lengths and ESS similar to those of full MH while achieving
higher ESS per second.

To further examine the role of the smooth surrogate, we compare
this approach with directly smoothing the loss in the posterior.
Specifically, we run full-batch MH targeting the density
proportional to
$\pi(\theta)\exp(-n\widehat R_{\tau,N}(\theta;\mathcal Z_N))$,
and denote this method by Smooth MH. Figure~\ref{fig:posterior-kde} in the supplement compares
the marginal density estimates from Smooth MH, full MH,
and SMH. At $\tau=0.3$, directly smoothing the loss noticeably
shifts the marginal distributions, whereas SMH, which uses
the smooth surrogate only to construct the control variate,
more closely matches the full-MH benchmark.
These results illustrate how smoothing can aid computation
while retaining the original loss.

\begin{table}[!htbp]
\centering
\caption{\textit{95\% CI length for the first coordinate, ESS, ESS per second, and standardized mean deviation for robust classification}}
\label{tab:robust-classification-results}
\renewcommand{\arraystretch}{1.25}
\setlength{\tabcolsep}{4pt}
\small
\begin{tabular}{lcccccc}
\toprule
 & $m=200$ & $m=400$ & $m=600$ & $m=800$ & $m=1000$ & Full MH \\
\midrule
CI length
& 0.4565 & 0.4558 & 0.4545 & 0.4542 & 0.4564 & 0.4524 \\
 
ESS & 300.3 & 325.8 & 336.7 & 339.7 & 346.3 & 348.6 \\
ESS\_per\_second & 16.0 & 17.6 & 14.4 & 12.1 & 10.8 & 3.6 \\
 
Mean deviation
& 0.0841 & 0.0677 & 0.0630 & 0.0757 & 0.0739 & / \\
\bottomrule
\end{tabular}
\end{table}

\subsection{G-and-K Mixture Model Fitting}
\label{sec:gk}

This example illustrates the use of zero-order control variates
for fitting a simulator-based model with a discontinuous loss.
We compare full-batch MH (Full MH), naive minibatching (Naive),
and the zero-order control-variate method (Sub$_0$) using $n = 1234$ gamma-ray burst $T_{90}$ durations from the
BATSE 4B catalog \citep{paciesas1999fourth}. Let $t_i$ denote the observed duration of burst $i$ in seconds.
We log-transform these measurements as $x_i=\log_{10}(t_i)$
and standardize them to $y_i=(x_i-\bar x)/s_x$,
where $\bar x$ and $s_x$ are the sample mean and standard
deviation of the log-durations. All MCMC updates and Wasserstein-loss evaluations
use these standardized observations. We model the standardized observations using a two-component
g-and-k mixture, represented by
\[
Y=
\{1-\mathbf 1(U<\pi)\}Q_1(Z_1)
+\mathbf 1(U<\pi)Q_2(Z_2),
\]
where $U\sim\operatorname{Unif}(0,1)$ and
$Z_1,Z_2\sim N(0,1)$ are mutually independent.
For component $j\in\{1,2\}$, the quantile function is
\[
Q_j(z)
=
A_j+B_j
\left\{1+0.8\tanh\left(\frac{g_jz}{2}\right)\right\}
(1+z^2)^{k_j}z.
\]
The parameter vector is $\theta=(\pi,\eta)$, with
$\eta=(A_1,B_1,g_1,k_1,A_2,B_2,g_2,k_2)$.
We assign independent uniform priors to these parameters,
with ranges given in Appendix~\ref{details of gk}.

For a latent bank
$\mathcal Z_K=\{(U_i,Z_{1i},Z_{2i})\}_{i=1}^{K}$,
we measure model fit using the squared 2-Wasserstein distance
\[
\ell_K(\theta;\mathcal Z_K)
=
W_2^2(\widehat P_{\mathrm{obs}},
      \widehat P_{\theta,\mathcal Z_K}),
\]
where $\widehat P_{\mathrm{obs}}$ is the empirical distribution
of the standardized observations and
$\widehat P_{\theta,\mathcal Z_K}$ is that of the simulated
observations.  For fixed latent draws, the mixture allocation changes
discontinuously with $\pi$. We use the zero-order surrogate
to handle this loss through function evaluations alone.

Full MH uses a fixed latent bank of size $N=20000$ and
targets the Gibbs posterior proportional to
$p(\theta)\exp(-N\ell_N(\theta;\mathcal Z_N))$,
where $p(\theta)$ denotes the prior density.
Naive and Sub$_0$ instead augment the state with a latent
minibatch $\mathcal Z_m$ of size $m=1000$.
Their respective surrogates are
\[
\begin{aligned}
&\mathcal L_{\rm naive}(\theta;\mathcal Z_m)
 =\ell_m(\theta;\mathcal Z_m),\\
&\mathcal L_{\rm zcv}(\theta;\mathcal Z_m)
 =\ell_m(\theta;\mathcal Z_m)
-\ell_m(\theta_{\rm ref};\mathcal Z_m).
\end{aligned}
\]
Thus, Sub$_0$ centers the minibatch loss at a fixed reference
point $\theta_{\rm ref}$, evaluating both terms with the
same latent draws.  
Both minibatch methods use a minibatch-retention probability
of $\zeta=1/5$, with risk-scaling factors
$\alpha_{\mathrm{Naive}}=50$ and
$\alpha_{\mathrm{Sub0}}=125$.

 All methods use a hybrid proposal that combines a random-walk
update for the mixture weight $\pi$ with a MALA update for
the eight component parameters $\eta$. Full MH and Naive SMH start from the same state
$\theta_{\mathrm{init}}$.
For Sub$0$, we use the mean of a $40000$-iteration Naive chain, after discarding the first $10000$ iterations, as both the reference point $\theta_{\mathrm{ref}}$ and the initial state.
Each chain runs for $200000$ iterations.
The initialization and proposal settings are specified in
Appendix~\ref{details of gk}.

\begin{figure}[htbp]
\centering
\begin{minipage}{0.32\textwidth}
  \centering
  \includegraphics[width=\linewidth]
  {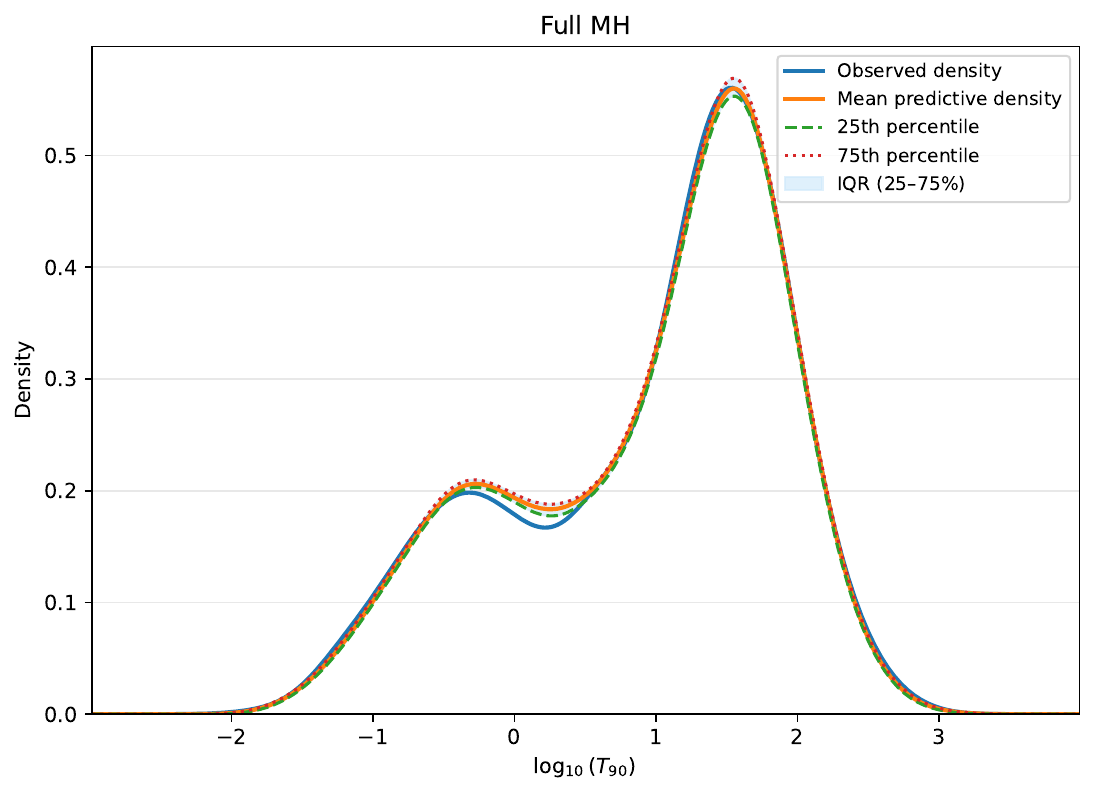}
\end{minipage}
\hfill
\begin{minipage}{0.32\textwidth}
  \centering
  \includegraphics[width=\linewidth]
  {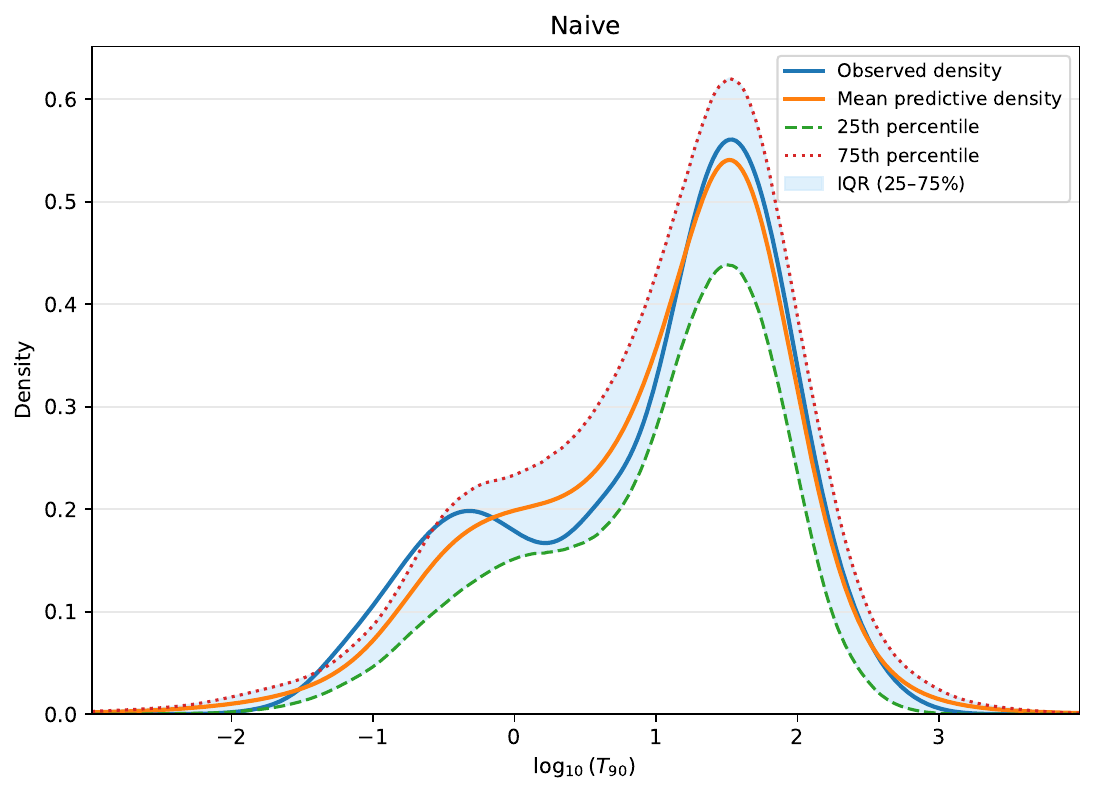}
\end{minipage}
\hfill
\begin{minipage}{0.32\textwidth}
  \centering
  \includegraphics[width=\linewidth]
  {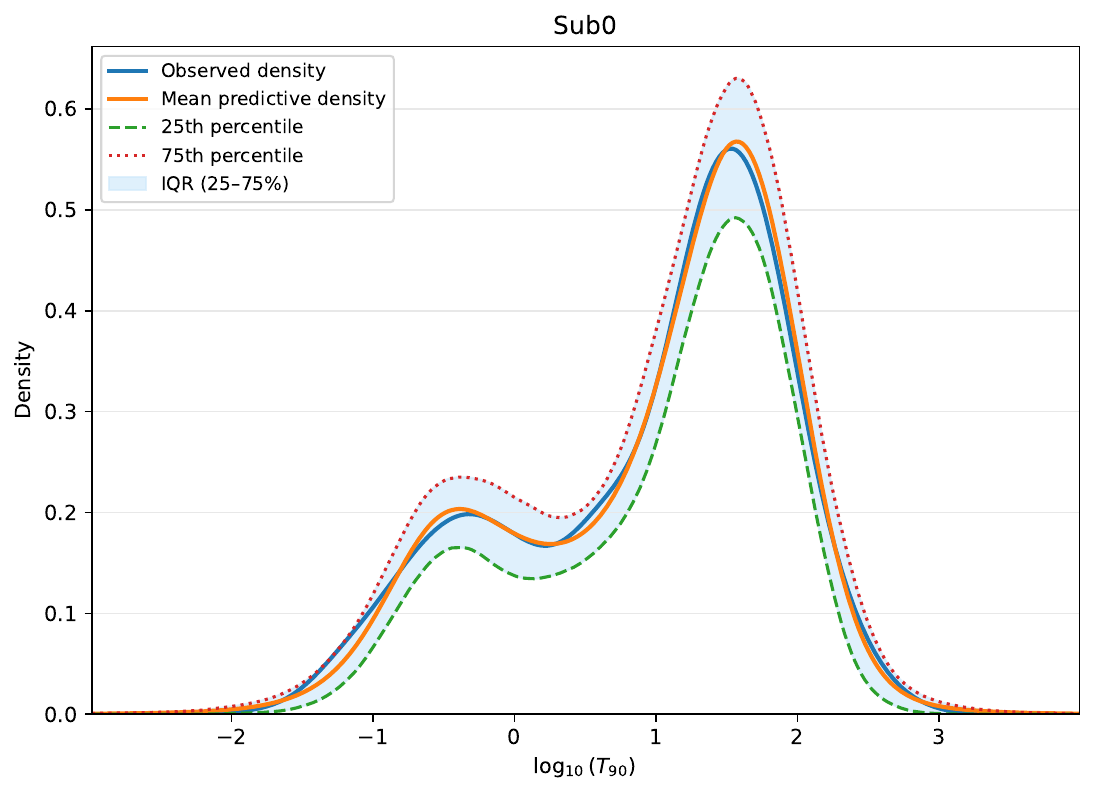}
\end{minipage}
\caption{Predictive-density comparisons for the BATSE
g-and-k mixture model, using Full MH, Naive, and Sub$_0$
from left to right.
The blue curve is the observed-data KDE, the orange curve
is the mean predictive KDE, and the shaded region spans
the pointwise $25$th to $75$th percentiles of the
state-specific predictive KDEs.
All densities are shown on the log-duration scale.}
\label{fig:batse_gk_density}
\end{figure}

To assess the fitted distributions, we select every tenth
state from each chain, yielding $S=20000$ parameter draws.
For each draw, we generate $M_{\mathrm{pred}}=2000$ observations,
transform them back to the log-duration scale, and compute
a normalized predictive KDE $\widehat f_s(x)$.
We summarize these densities by their mean
$\bar f(x)=S^{-1}\sum_{s=1}^{S}\widehat f_s(x)$ and assess
agreement with the observed-data KDE $f_{\mathrm{obs}}$ using
\[
L_1=\int|\bar f(x)-f_{\mathrm{obs}}(x)|\,\mathrm dx,
\qquad
L_2=
\left[
\int|\bar f(x)-f_{\mathrm{obs}}(x)|^2\,\mathrm dx
\right]^{1/2}.
\]
Figure~\ref{fig:batse_gk_density} compares the observed-data
KDE with the mean predictive KDE and pointwise interquartile
band from each method.
Naive SMH underestimates the dominant mode, whereas zero-order SMH
captures the two-mode structure more closely.
Table~\ref{tab:batse_gk_density} reports the density discrepancies
and computation times, averaged over $30$ independent
replications. For Sub$0$, the reported computation time includes the time required to obtain the reference point $\theta_{\mathrm{ref}}$. Sub$_0$ reduces both discrepancies relative
to Naive. Full MH achieves the smallest discrepancies,
while Sub$_0$ provides a closer fit than Naive at substantially
less reported computation time than Full MH.

\begin{table}[htbp]
\centering
\caption{G-and-k mixture model fitting.
Distances between the mean predictive KDE and the observed-data
KDE, together with reported MCMC computation times,
averaged over $30$ independent replications.}
\label{tab:batse_gk_density}
\begin{tabular}{lccc}
\toprule
& Full MH & Naive & Sub$_0$ \\
\midrule
$L_1$ & 0.0419 & 0.1272 & 0.0758 \\
$L_2$ & 0.0240 & 0.0647 & 0.0410 \\
Time (s) & 326.35 & 46.42  & 77.05 \\
\bottomrule
\end{tabular}
\end{table}

 \section{Conclusion}
In this paper, we develop finite-sample approximation theory
for subsampling Gibbs posteriors with non-smooth losses.
Our bounds for subgradient-based first-order control variates
quantify how minibatch size, risk scaling, and reference-point
accuracy affect posterior approximation, and identify conditions
for consistent approximation of the untempered posterior using
sublinear minibatches. Complementary zero-order results support
loss-only estimation and staged reference-point refinement.
Under explicit geometric and proposal-comparison conditions,
our warm-start convergence bounds connect stationary
approximation guarantees with finite-iteration accuracy.
While the approximation theory assumes Lipschitz continuity
of the individual losses, our numerical results also suggest
that the constructions can be useful for discontinuous losses.
Extending the theory to this setting is a direction for future work.

 \section*{Statement on AI use}
The authors used OpenAI's ChatGPT to assist with language editing
and code refinement. They reviewed all outputs and take full
responsibility for the content of the manuscript.

 \bibliography{references}{}
\bibliographystyle{abbrvnat} 
\newpage

\appendix
\begin{center}
{\bf\Large Appendix}
\end{center}
\textbf{Notations:} We use the notations defined in the main text, and introduce the following additional notations used in our proofs. We use $N(\mu, \Sigma)$ to denote the $d$-dimensional multivariate Gaussian distribution with mean vector $\mu \in \mathbb{R}^d$ and covariance matrix $\Sigma \in \mathbb{R}^{d \times d}$.
We use $\m N(\mathcal{F}, d_n, \varepsilon)$ to denote the $\varepsilon$-covering number of the function class $\mathcal{F}$ with respect to the pseudo-metric $d_n$. 

\section{Details for Numerical Studies}\label{detail}
To support the large number of replications required to estimate coverage probabilities, the experiments in Sections~\ref{details of Huber Regression} and~\ref{details of Ranking problem} were conducted on a CPU cluster with Intel Xeon Gold 6348 processors. The experiment in Section~\ref{details:robust-classification}, which involves tensor computations, was conducted on NVIDIA Quadro RTX 8000 GPUs. The experiment in Section~\ref{details of gk} was conducted on an Apple M4 chip. The code for the numerical experiments is available at \url{https://github.com/DMASTIMA/Subsampling-Approximations-to-Gibbs-Posteriors-with-Non-smooth-Losses}.
\subsection{Bayesian Huber Regression}\label{details of Huber Regression}
We consider the tempered Gibbs posterior for Huber regression 
\[
\pi_\alpha(\theta) \propto \pi(\theta)\exp\Big(- \frac{\alpha}{n}\sum_{i=1}^n H_{\delta}\big(Y_i - {X}_i^T \theta\big) \Big) , \quad \theta \in \mathbb{R}^d,
\]where 
\[
H_{\delta}(z) =
\begin{cases}
\frac{1}{2}z^2, & |z| \leq \delta \\
\delta\left(|z| - \frac{1}{2}\delta\right), & |z| > \delta
\end{cases},
\]
and compare three algorithms:  the zero-order control variate SMH  algorithm  with rescaling (denoted as sub0), the first-order control variate SMH algorithm (denoted as sub1) and the full-batch  MH  algorithm (denoted as full MH),  all using (mini-batch) MALA proposals. Here instead of sampling $\wX^{(m)}$ with replacement as in the theoretical analysis, we sample $\wX^{(m)}$ without replacement from $\bX^{(n)}$. Sampling with replacement is adopted in the theory for technical convenience, and we expect both schemes to behave similarly when $m\ll n$. However, when $m$ is comparable to $n$, sampling without replacement may perform better as it avoids duplicating data points. Therefore, in all simulations we fix the design to sample $\wX^{(m)}$ without replacement. Additionally, we introduce \(1/2\)-laziness in the subsample generation: at each SMH step, with probability \(1/2\) we draw a new subsample, and with probability \(1/2\) we retain the current subsample.

We set the sample size $n = 20000$, batch size $m\in\{200,400,800,1200,1600,2000\}$ and choose the parameter dimension $d$ from  $\{2, 5, 10\}$, with hyperparameter $\delta=2$ and the prior is given by the uniform distribution over the hypercube $[-10,10]^d \subset \mathbb{R}^d$. The covariates $X$ are generated from the standard multivariate Gaussian distribution with zero mean and identity covariance matrix. The response variable $Y$ is defined via the linear model $Y = {X}^\top \boldsymbol{1}_d + \varepsilon$, where the noise $\varepsilon$ is formed by concatenating $\lfloor 4n/5 \rfloor$ standard normal samples $N(0,1)$ and $n-\lfloor 4n/5 \rfloor$ Student's $t$-samples with $4$ degrees of freedom $t(4)$, with the two noise components being mutually independent. All numerical results are based on a total of $1000$ independent experimental replications, and each replication runs $10000$ iterations of the (subsampling) MH algorithm.  
All methods use the (mini-batch) Metropolis-adjusted Langevin algorithm (MALA) proposals described in Section~\ref{Mixing time analysis}, with $h = \frac{1.35}{d^{1/3}}$.

We use the ordinary least squares (OLS) estimator as the
reference point for both control-variate methods and as the
initial state for all three methods. Under the present
data-generating model and for fixed dimension, both OLS and
the empirical Huber-risk minimizer $\widehat\theta$ are
root-$n$ consistent for the regression parameter. Consequently,
$\|\theta_{\mathrm{OLS}}-\widehat\theta\|_2
=\mathcal O_{\mathbb P}(n^{-1/2})$.

We set $\alpha = n$ for full MH and the first-order method so that they target and approximate, respectively, the untempered Gibbs posterior $\pi_n$. For the zero-order method with rescaling, balancing the
leading terms in the bound of Corollary~\ref{co:2}
motivates the scaling $\alpha\asymp m^{2/3}$, up to
logarithmic factors. We therefore parameterize its
risk-scaling factor as $\alpha=k_1m^{2/3}$ and use
the proposal step size $h=k_2d^{-1/3}$.

To assess the multiplicative constants, we conducted an
exploratory alternating tuning procedure with $m=200$
and $d=2$, starting from $(k_1,k_2)=(5,0.3)$.
Holding $k_1$ fixed, we updated $k_2$ every $500$ iterations
using the acceptance rate $A_t$ in the current block,
\[
k_2^{(t+1)}
=
\begin{cases}
k_2^{(t)}, & A_t\in[0.55,0.65],\\[3pt]
k_2^{(t)}
\exp\!\left(
\max\!\left\{-\log 1.5,\,
\min\!\left\{\log 1.5,\,1.5(A_t-0.6)\right\}\right\}
\right), & \text{otherwise}.
\end{cases}
\]
This update moves the acceptance rate toward the target
window $[0.55,0.65]$, while limiting each multiplicative
change in the step size to a factor of $1.5$.

When two consecutive blocks generated under the same
$(k_1,k_2)$ both had acceptance rates within this window,
we pooled their $1000$ samples and computed the mean
relative ESS,
$E_t=(1000d)^{-1}\sum_{j=1}^{d}
\widehat{\mathrm{ESS}}_{t,j}$,
where $\widehat{\mathrm{ESS}}_{t,j}$ is the ESS for
coordinate $j$. We stopped tuning when $E_t\geq0.2$.
Otherwise, we reduced the risk-scaling constant according
to $k_1^{(t+1)}=\max\{0.05,\,0.85k_1^{(t)}\}$ and
resumed step-size tuning from the current value of $k_2$.
This alternation seeks adequate sampling efficiency
while accounting for the interaction between risk scaling
and proposal size.

\begin{figure}[H]
\centering
\includegraphics[width=0.8\textwidth]{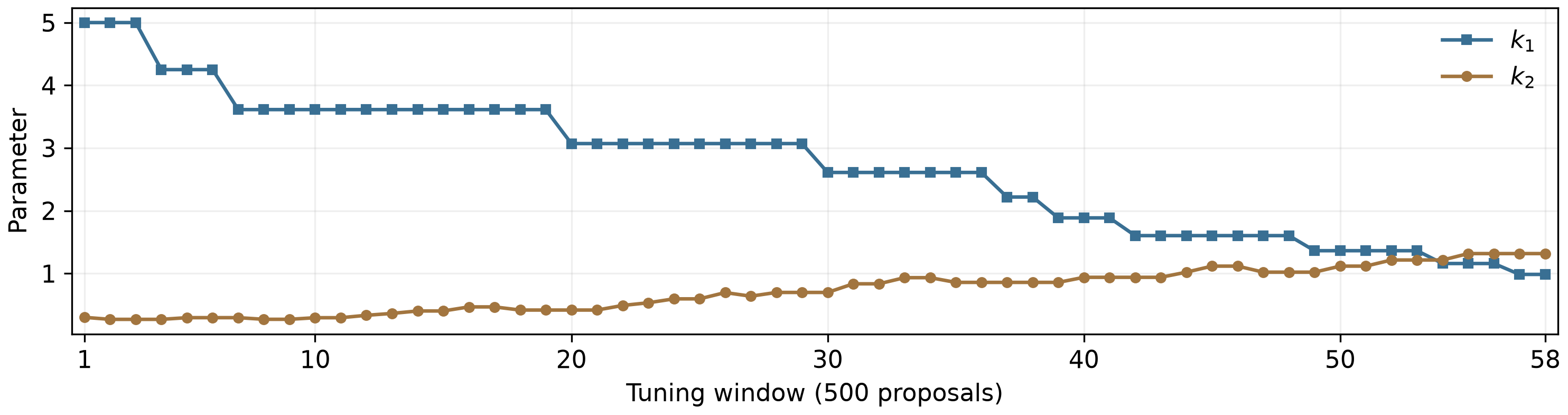}
\caption{Tuning trajectories of $k_1$ and $k_2$ with
$m=200$ and $d=2$}
\label{fig:sub0-tuning}
\end{figure}

 Figure~\ref{fig:sub0-tuning} shows the tuning trajectories.
The pilot stopped at $(k_1,k_2)=(0.9844,1.3166)$.
We adopted the nearby values $(1,1.35)$, using $k_1=1$
for the zero-order method and $k_2=1.35$ for all three
methods throughout the main Huber regression experiments.
The comparisons with full MH below assess posterior
approximation accuracy under these fixed settings.

\begin{figure}[!htbp]
\centering
\includegraphics[
    width=0.8\textwidth,
    height=0.36\textheight,
    keepaspectratio
]{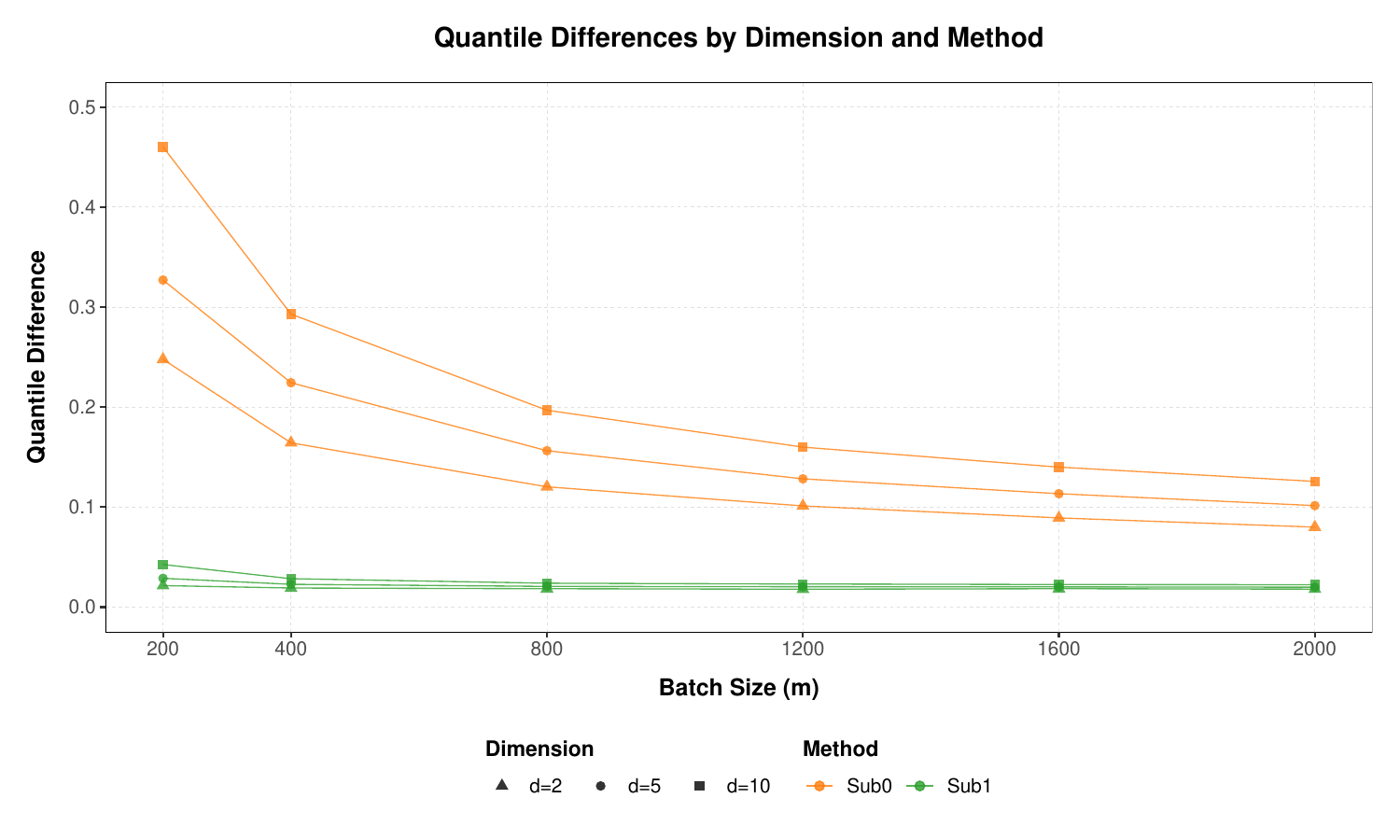}
\caption{Quantile difference for Bayesian Huber regression}
\label{Huber diff}
\end{figure}

\begin{figure}[H]
\centering
\includegraphics[width=0.8\textwidth]{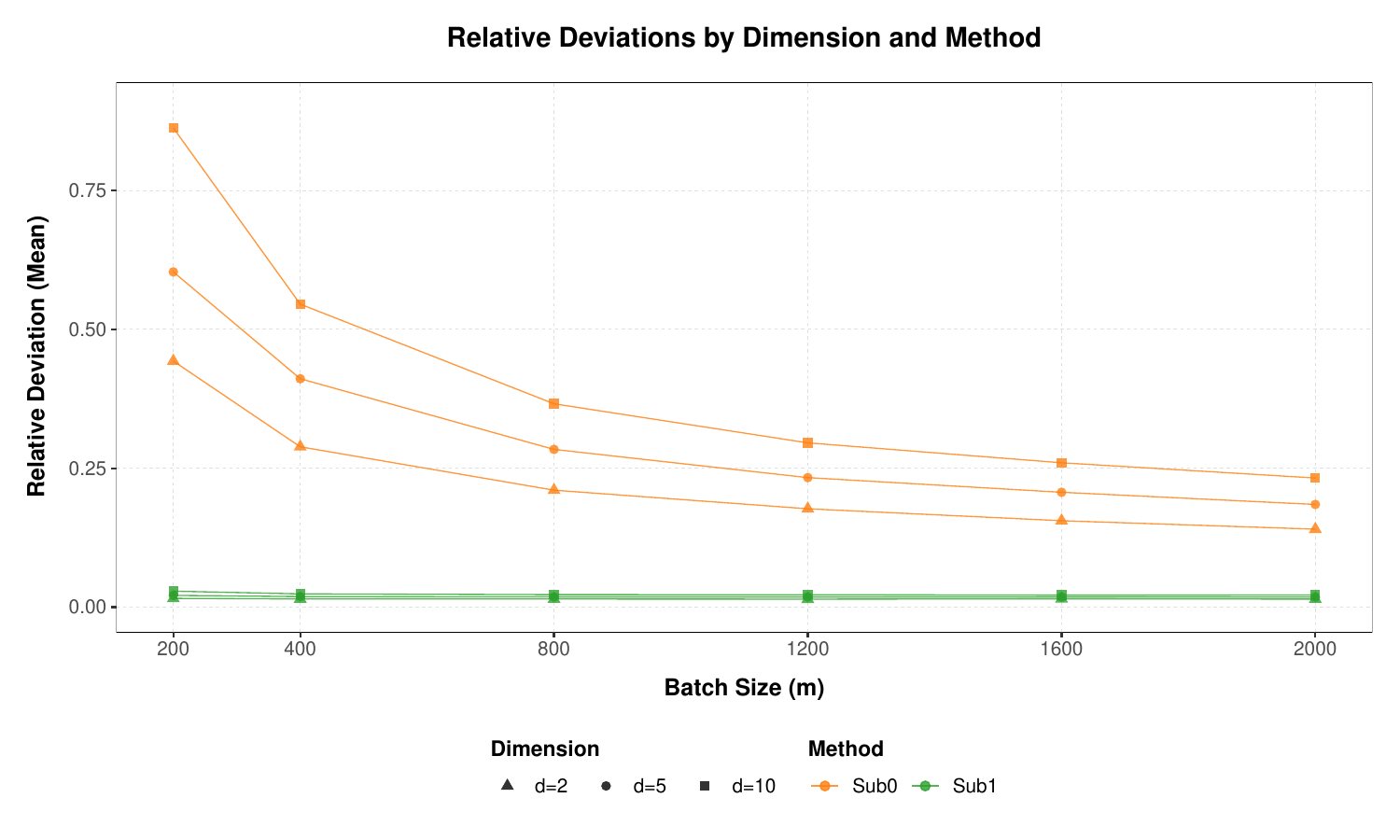}
\caption{Relative bias for Bayesian Huber regression}
\label{Huber Rdev}
\end{figure}

Table~\ref{Huber cov CI ess} reports the  empirical coverage and length of the 95\% credible interval for the first coordinate of $\theta$, together with the effective sample size (ESS)  and the ESS generated per second (both averaged over coordinates). To assess the approximation quality of sub0 and sub1 relative to full MH, we consider the following relative quantile error:
\begin{equation}\label{def:quantile}
\text{quantile\_diff}_i = \frac{1}{d} \sum_{j=1}^{d} \frac{ \left| q^{(i)}_{j,97.5} - q^{\text{full}}_{j,97.5} \right| + \left|q^{\text{full}}_{j,2.5} - q^{(i)}_{j,2.5} \right| }{ q^{\text{full}}_{j,97.5} - q^{\text{full}}_{j,2.5} }, i=0,1,
\end{equation}
where $q^{\text{full}}_{j,\alpha}$ is the $\alpha\%$-quantile of the $j$th coordinate under the full-batch MH algorithm, and $q^{(i)}_{j,\alpha}$ is the corresponding quantile under SMH with an $i$th-order control variate. The results for this metric are shown in Figure~\ref{Huber diff}. In addition, we examine the relative deviation (bias) of the posterior means produced by sub0 and sub1, defined as
\begin{equation}\label{def:mean}
\text{mean\_diff}_i = \frac{1}{d} \sum_{j=1}^{d} \frac{ \left| u_{i,j} - u_{\text{full},j} \right|}{ sd_{\text{full},j}}, i=0,1,
\end{equation}
where $u_{\text{full},j}$  and $sd_{\text{full},j}$  
  denote the mean and standard deviation of the  $j$th marginal of the  untempered Gibbs posterior $\pi_n$, estimated from the full-batch MH output, while $ u_{i,j}$ is the sample mean of the $j$th coordinate under SMH with the $i$th-order control variate. The  results of this metric are presented in Figure \ref{Huber Rdev}.

\subsection{Ranking from Noisy Pairwise Comparisons}\label{details of Ranking problem}

We consider the tempered Gibbs posterior
\[
\pi_\alpha(\theta)
\propto
\pi(\theta)\exp\left(
-\frac{6\alpha}{n(n-1)}
\sum_{i<j}
\bigl[1-C_{ij}\theta^\top(X_i-X_j)\bigr]_+
\right),
\]
where $[u]_+=\max\{u,0\}$.
To generate the data, we draw covariates
$X_i\sim N(0,\mathbf I_d)$ and assign each item a latent score
$U_i=X_i^\top\mathbf 1_d+\epsilon_i$,
where $\epsilon_i\sim N(0,1)$.
For each pair $i<j$, we generate the comparison
$C_{ij}=\operatorname{sign}(U_i-U_j+\eta_{ij})\in\{-1,1\}$,
where $\eta_{ij}\sim N(0,0.5^2)$.
All covariates and noise variables are mutually independent.
We set $n=1000$, $d=3$, and
$m\in\{100,200,250,300\}$, and use a uniform prior on
$[-10,10]^d$.

The empirical risk is a pairwise average, so we construct
the minibatch surrogates using all comparisons among
$m$ sampled items. Specifically, we select $m$ distinct
indices uniformly without replacement and write them as
$I_1<\cdots<I_m$.
Let $\wt X_j=X_{I_j}$ and
$\wt C_{j_1j_2}=C_{I_{j_1}I_{j_2}}$, and let
$\wX^{(m)}$ denote the selected covariates.
Without control variates, we use
\[
\m L_{\rm naive}(\wX^{(m)},\theta)
=
\frac{1}{m(m-1)}
\sum_{1\leq j_1<j_2\leq m}
\wt\ell_{j_1,j_2}(\theta),
\]
where
$\wt\ell_{j_1,j_2}(\theta)
=6[1-\wt C_{j_1j_2}
\theta^\top(\wt X_{j_1}-\wt X_{j_2})]_+$.
For a reference point $\theta^\dagger$, the first-order
surrogate is
\[
\begin{aligned}
\m L_{\rm fcv}(\wX^{(m)},\theta)
={}&
\Big[
\frac{1}{m(m-1)}
\sum_{1\leq j_1<j_2\leq m}
\left\{
\wt\ell_{j_1,j_2}(\theta)
-\wt\ell_{j_1,j_2}(\theta^\dagger)
-\wt g_{j_1,j_2}(\theta^\dagger)^\top
(\theta-\theta^\dagger)
\right\}
\Big]_+\\
&+
\frac{1}{n(n-1)}
\sum_{1\leq i_1<i_2\leq n}
g_{i_1,i_2}(\theta^\dagger)^\top
(\theta-\theta^\dagger).
\end{aligned}
\]
Here, the pairwise subgradients are chosen as
\[
\begin{aligned}
g_{i_1,i_2}(\theta)
&=
-6C_{i_1i_2}(X_{i_1}-X_{i_2})
\mathbf 1\!\left\{
C_{i_1i_2}\theta^\top(X_{i_1}-X_{i_2})<1
\right\},\\
\wt g_{j_1,j_2}(\theta)
&=
-6\wt C_{j_1j_2}(\wt X_{j_1}-\wt X_{j_2})
\mathbf 1\!\left\{
\wt C_{j_1j_2}
\theta^\top(\wt X_{j_1}-\wt X_{j_2})<1
\right\}.
\end{aligned}
\]
The full-data reference gradient is computed once at
the beginning of each stage.

Given the small parameter dimension, all methods use
the random-walk proposal
$q(\theta\mid\theta')
=N(\theta',\sigma^2\mathbf I_d/\alpha)$
with $\sigma=1.5$.
Full MH uses $\alpha=n$ and starts from the pairwise
least-squares estimate $\theta_{\rm init}
=
\operatorname*{arg\,min}_{\theta}
\sum_{i<j}
\left[C_{ij}-(X_i-X_j)^\top\theta\right]^2.$ We run this chain for $K=10000$ iterations and discard
the first $2000$ as burn-in.  For each minibatch size, the subsampling method proceeds
through three stages.
First, we run $K_1=2000$ iterations without control
variates at $\alpha=10m^{2/5}$, starting from
$\theta_{\rm init}$.
After discarding the first $1000$ iterations, we take
the retained sample mean $\theta_1$ as both the reference
point and initial state for the next stage.
We then run the first-order method for another
$K_1=2000$ iterations at $\alpha=10m^{4/5}$.
Discarding the first $1000$ iterations gives the refined
sample mean $\theta_2$.
Finally, we use $\theta_2$ as both the reference point
and initial state for a first-order chain with
$\alpha=n$, run for $K=10000$ iterations with the first
$2000$ discarded as burn-in. The choices of $\alpha$ are motivated by the initial
estimation and refinement strategies in
Theorems~\ref{th:1},~\ref{th:2}, and Remark~\ref{re:1}.
This experiment evaluates their performance for the
pairwise objective.
At each subsampling iteration, the minibatch proposal
retains the current minibatch with probability $1/2$
and otherwise draws a new subset of $m$ distinct items.
We evaluate empirical coverage relative to the population
hinge-risk minimizer
$\theta^*=\operatorname*{arg\,min}_{\theta}
\mathbb E\bigl\{[1-C_{12}\theta^\top(X_1-X_2)]_+\bigr\}$.
Under the stated data-generating model,
$\theta^*=c\mathbf 1_3$ with $c\approx0.8610$.
All numerical results are based on $1000$ independent
experimental replications.

\subsection{Robust Classification under Gaussian Perturbations}\label{details:robust-classification}

We consider the tempered Gibbs posterior
$$
\pi_{\alpha}(\theta) \propto \pi(\theta)\exp\big(-\alpha R_n(\theta)\big),
$$
where the empirical risk under perturbation is defined as
$$
R_n(\theta) = \frac{1}{n}\sum_{i=1}^{n} \mathbb{E}_{Z} \left[ \mathbf{1} \left\{ Y_i f_{\theta}(X_i+Z) < 0 \right\} \right].
$$
Here, the parameter vector is $\theta=(\theta_1,\theta_2,\theta_3)^\top\in\mathbb{R}^3$, and the decision function is given by $f_{\theta}(x) = x_1+\theta_1x_2+\theta_2x_1^2+\theta_3x_2^2$.
The simulated variable follows $Z\sim N(0,0.1\mathbf{I}_2)$. Approximating the expectation over $Z$ using many Monte Carlo
draws can be computationally expensive.  To reduce this cost, we consider applying the first-order control variate method to the latent variable. However, our subsampling strategy relies on gradient information, while the original loss is nonsmooth. We therefore introduce a smooth approximation to facilitate gradient-based subsampling. We also investigate how smoothing itself affects the sampling mechanism.

We generate $n=2000$ independent covariates
$X_i\sim N(0,\mathbf I_2)$ and binary responses satisfying
\[
\mathbb P(Y_i=1\mid X_i)
=
\sigma\bigl(f_{\theta^*}(X_i)\bigr),
\qquad
\theta^*=(1,1,-1)^\top,
\]
where $\sigma(t)=1/(1+\exp(-t))$ and
$Y_i\in\{-1,1\}$.
We use the Gaussian prior $\theta\sim N(0,\mathbf I_3)$
and set $\alpha=n$ for the main sampling runs. For a latent-variable batch $\mathcal Z_m=(Z_1,\ldots,Z_m)$,
define the Monte Carlo risk and its smooth counterpart by
\[
\begin{aligned}
\widehat R_m(\theta;\mathcal Z_m)
&=
\frac{1}{nm}\sum_{i=1}^{n}\sum_{j=1}^{m}
\mathbf 1\{Y_i f_\theta(X_i+Z_j)<0\},\\
\widehat R_{\tau,m}(\theta;\mathcal Z_m)
&=
\frac{1}{nm}\sum_{i=1}^{n}\sum_{j=1}^{m}
\sigma\left(-\frac{Y_i f_\theta(X_i+Z_j)}{\tau}\right),
\end{aligned}
\]
with $\tau=0.3$.
Within each replication, we generate a fixed bank
$\mathcal Z_N$ of $N=20000$ independent perturbations.
Full MH uses $\widehat R_N(\theta;\mathcal Z_N)$,
providing the benchmark based on the original indicator loss.
To examine the effect of directly smoothing the loss,
we also run Smooth Full MH using
$\widehat R_{\tau,N}(\theta;\mathcal Z_N)$.
Both chains start at $(5,5,5)^\top$ and run for $10000$
iterations, discarding the first $1000$ as burn-in.

For each minibatch size $m\in\{200,400,600,800,1000\}$,
we obtain the reference point for SMH through a loss-only
pilot. Starting at $(5,5,5)^\top$, we run subsampling MH
without control variates for $2000$ iterations at
$\alpha=50m^{2/5}$.
We discard the first $1000$ iterations and take the
retained sample mean as $\theta^\dagger$.
The subsequent first-order method uses the surrogate
\[
\begin{aligned}
\mathcal L_m(\theta;\mathcal Z_m)
={}&
\widehat R_m(\theta;\mathcal Z_m)
-\widehat R_m(\theta^\dagger;\mathcal Z_m)\\
&-
\left[
\nabla_\theta\widehat R_{\tau,m}
    (\theta^\dagger;\mathcal Z_m)
-
\nabla_\theta\widehat R_{\tau,N}
    (\theta^\dagger;\mathcal Z_N)
\right]^\top
(\theta-\theta^\dagger).
\end{aligned}
\]
Thus, smoothing supplies the gradient information for
the linear correction, while the loss differences are
computed using the original indicator function.
The SMH chain starts at $\theta^\dagger$ and runs for
$10000$ iterations at $\alpha=n$, discarding the first
$1000$ as burn-in. We set the minibatch-retention
parameter to $\zeta=1/2$.
All methods use the random-walk proposal
$q(\theta'\mid\theta)=N(\theta,20\mathbf I_3/n)$.
The numerical results are based on $100$ independent
experimental replications.

\begin{figure}[t]
    \centering
    \begin{subfigure}[t]{0.32\textwidth}
        \centering
        \includegraphics[width=\linewidth]{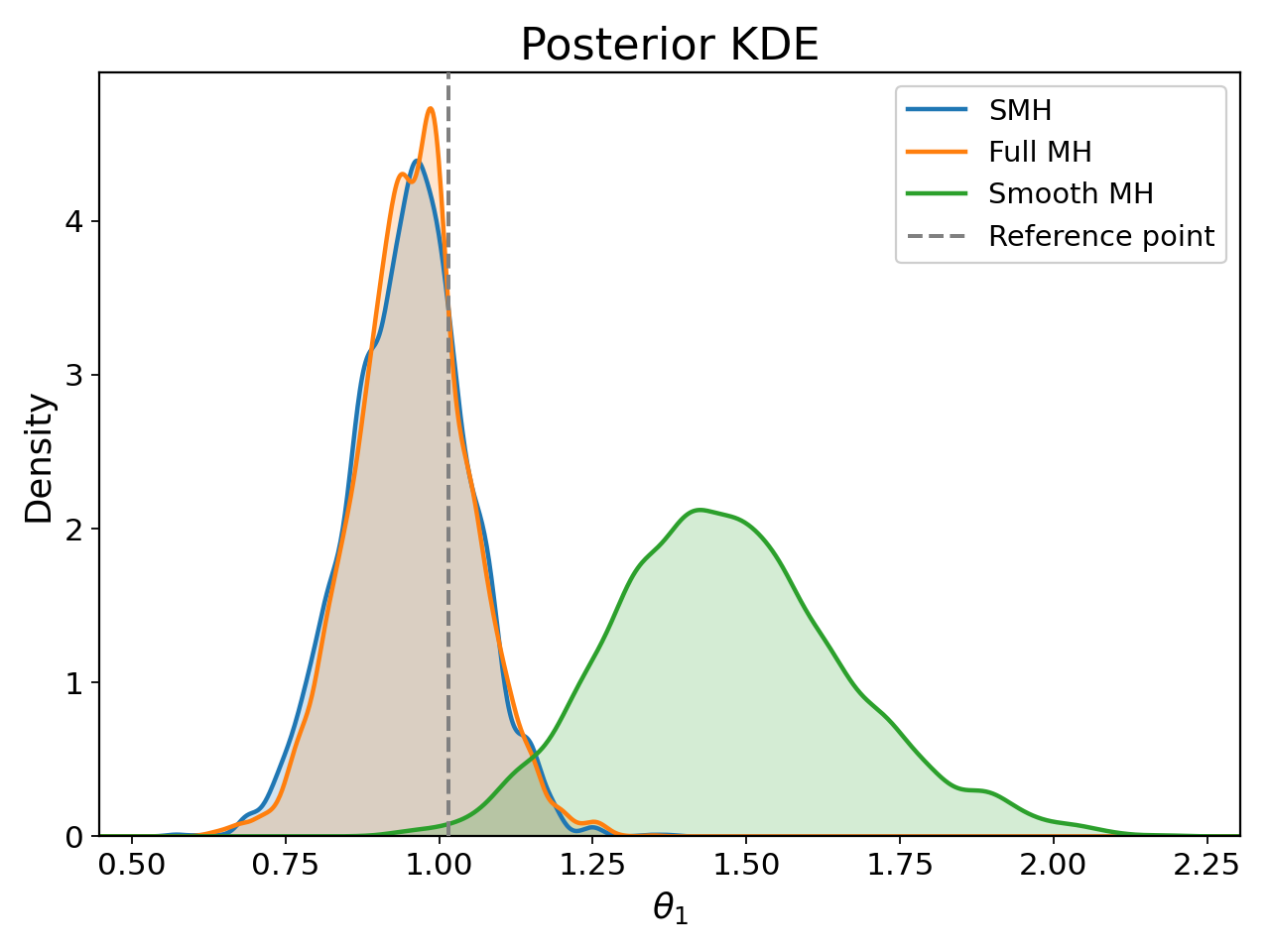}
        \caption{$\theta_1$}
    \end{subfigure}\hfill
    \begin{subfigure}[t]{0.32\textwidth}
        \centering
        \includegraphics[width=\linewidth]{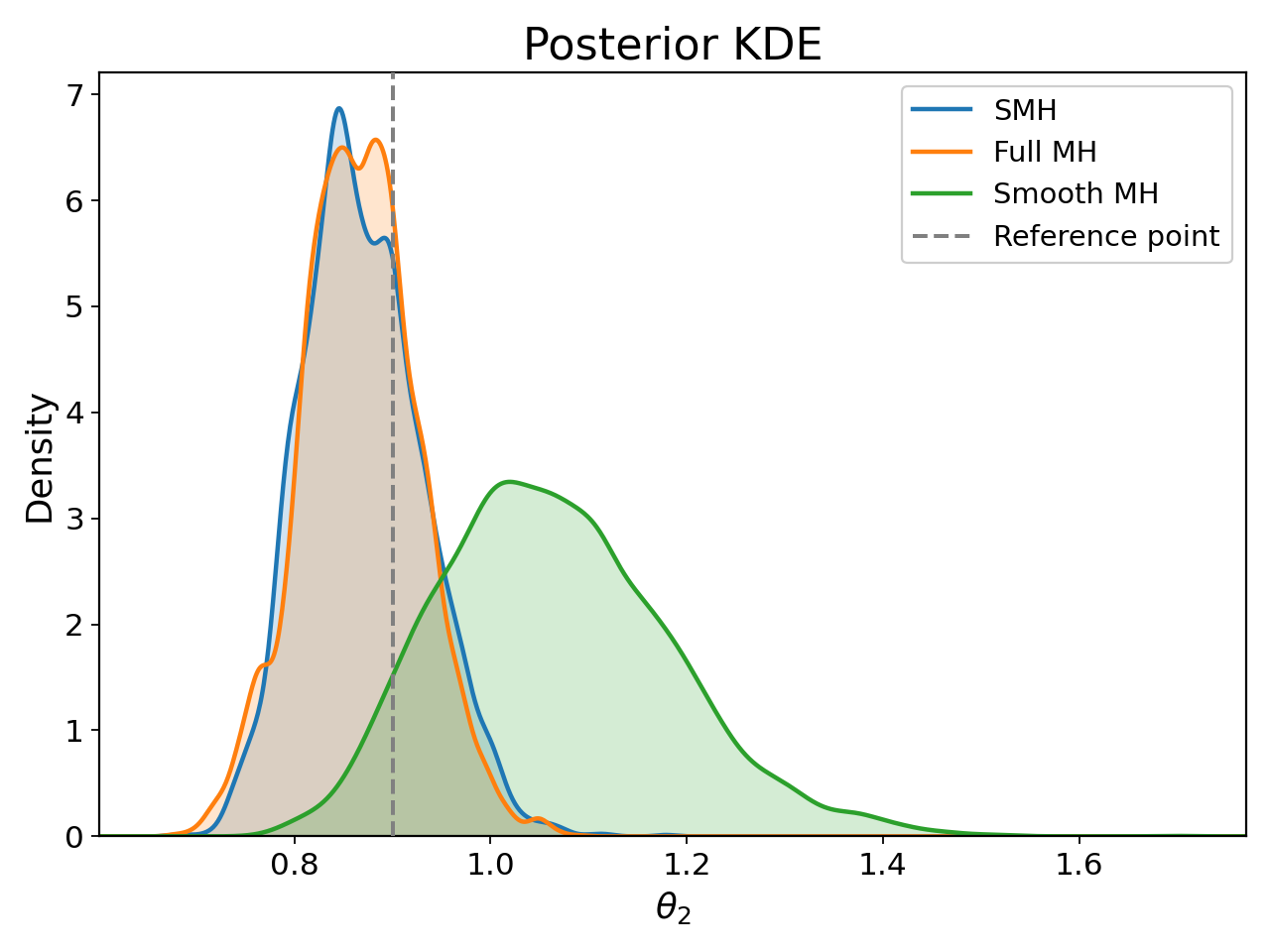}
        \caption{$\theta_2$}
    \end{subfigure}\hfill
    \begin{subfigure}[t]{0.32\textwidth}
        \centering
        \includegraphics[width=\linewidth]{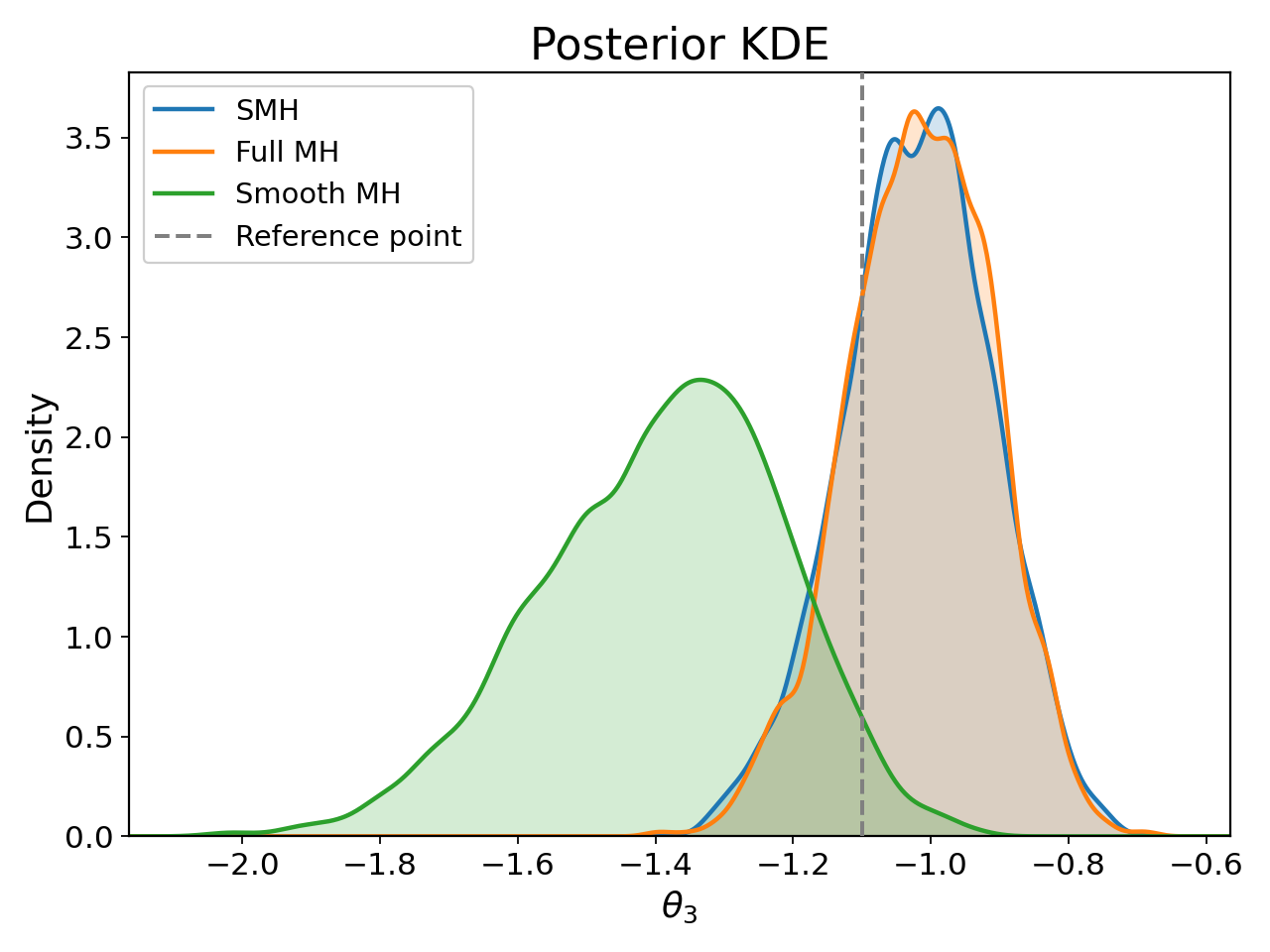}
        \caption{$\theta_3$}
    \end{subfigure}
    \caption{Effect of smoothing in robust classification.
Marginal KDEs from Full MH, Smooth Full MH with $\tau=0.3$,
and SMH with $m=1000$, all at $\alpha=n$.
Smooth Full MH replaces the indicator loss by its smooth
approximation, whereas SMH uses smoothing only to construct
the control variate. The dashed line in each panel marks
the corresponding coordinate of the reference point
$\theta^\dagger$.}
    \label{fig:posterior-kde}
\end{figure}

 Figure~\ref{fig:posterior-kde} compares the marginal density
estimates from Full MH, Smooth Full MH, and SMH with $m=1000$.
At $\tau=0.3$, directly smoothing the loss noticeably
shifts the marginal distributions, whereas SMH follows
the Full-MH benchmark more closely.
Table~\ref{tab:robust-classification-results} further
shows that SMH produces credible-interval lengths and
posterior mean estimates close to those of Full MH,
while achieving higher ESS per second.

\subsection{G-and-K Mixture Model Fitting}
\label{details of gk}

This example examines zero-order control variates for fitting
a simulator-based mixture model with a discontinuous loss.
We compare full-batch MH (Full MH), naive minibatching (Naive),
and the zero-order control-variate method (Sub$_0$), assessing
both predictive fit and computational cost. We use $n=1234$ finite, positive gamma-ray burst duration
measurements from the BATSE 4B catalog~\citep{paciesas1999fourth}. Let $t_i$ denote the observed duration of burst $i$ in seconds.
We set $x_i=\log_{10}(t_i)$ and standardize the log-durations
as $y_i=(x_i-\bar x)/s_x$, where $\bar x$ and $s_x$ are their
sample mean and standard deviation.
All MCMC updates and Wasserstein-loss evaluations are
performed on this standardized scale.

We model the standardized observations using a two-component
g-and-k mixture,
\[
Y=
\{1-\mathbf 1(U<\pi)\}Q_1(Z_1)
+\mathbf 1(U<\pi)Q_2(Z_2),
\]
where $U\sim\operatorname{Unif}(0,1)$ and
$Z_1,Z_2\sim N(0,1)$ are mutually independent.
The component weights are $1-\pi$ and $\pi$, respectively,
and the quantile function for component $j\in\{1,2\}$ is
\[
Q_j(z)
=
A_j+B_j
\left\{1+0.8\tanh\left(\frac{g_jz}{2}\right)\right\}
(1+z^2)^{k_j}z.
\]
The parameter vector is $\theta=(\pi,\eta)$, where
$\eta=(A_1,B_1,g_1,k_1,A_2,B_2,g_2,k_2)$.
We assign independent priors
$\pi\sim\operatorname{Unif}(0,0.5)$ and, for $j=1,2$,
$A_j\sim\operatorname{Unif}(-2,2)$,
$B_j\sim\operatorname{Unif}(0,3)$,
$g_j\sim\operatorname{Unif}(-1,1)$, and
$k_j\sim\operatorname{Unif}(-0.04,2)$.

Given $\theta$ and a latent bank
$\mathcal Z_K=\{(U_i,Z_{1i},Z_{2i})\}_{i=1}^{K}$,
we generate $K$ observations using the mixture representation
above. We measure model fit using
\[
\ell_K(\theta;\mathcal Z_K)
=
W_2^2\bigl(
\widehat P_{\mathrm{obs}},
\widehat P_{\theta,\mathcal Z_K}
\bigr),
\]
where $W_2$ denotes the 2-Wasserstein distance,
$\widehat P_{\mathrm{obs}}$ is the empirical distribution
of the standardized observations, and
$\widehat P_{\theta,\mathcal Z_K}$ is that of the simulated
observations.
For fixed latent draws, changes in $\pi$ can switch
component allocations and cause jumps in this loss.

At the beginning of each replication, Full MH generates
a latent bank $\mathcal Z_N$ of size $N=20000$ and holds
it fixed throughout the chain. Its target is
\[
\pi_{\mathrm{Full}}(\theta\mid\mathcal Z_N)
\propto
p(\theta)\exp\bigl(-N\ell_N(\theta;\mathcal Z_N)\bigr).
\]
Naive and Sub$_0$ instead augment the Markov state with
a latent minibatch, giving the joint state
$(\theta,\mathcal Z_m)$ with $m=1000$.
Both methods use a minibatch-retention probability of $1/5$.
Their risk-scaling factors are
$\alpha_{\mathrm{Naive}}=50$ and
$\alpha_{\mathrm{Sub0}}=120$, respectively.

All methods use a hybrid proposal combining a random-walk
update for the mixture weight $\pi$ with a MALA update for
the eight component parameters $\eta$,
\[
\begin{aligned}
\pi'
&=
\pi+\sqrt{\frac{2h_\pi}{\alpha}}\,\varepsilon_\pi,\\
\eta'
&=
\eta-\boldsymbol h_\eta\odot
\nabla_\eta\ell_K(\theta;\mathcal Z_K)
+\sqrt{\frac{2\boldsymbol h_\eta}{\alpha}}
\odot\varepsilon_\eta,
\end{aligned}
\]
where $\varepsilon_\pi\sim N(0,1)$ and
$\varepsilon_\eta\sim N(0,\mathbf I_8)$ are independent.
Here, $K=N$ for Full MH and $K=m$ for Naive and Sub$_0$.
The zero-order reference term in Sub$_0$ is independent
of $\theta$ and therefore does not affect this gradient.
The symbol $\odot$ denotes elementwise multiplication,
and the square root is taken elementwise.
We use $\alpha=N$ for Full MH and the corresponding
risk-scaling factor for each minibatch method.
The step sizes are $h_\pi=0.020$ and
$\boldsymbol h_\eta
=(0.030,0.020,0.020,0.015,0.030,0.020,0.020,0.015)$.

Full MH and Naive SMH start from the same initial state
$\theta_{\mathrm{init}}
= (0.25,0.25,0.60,0,0.25,\\-0.25,0.60,0,0.25)$.
For Sub$0$, we use the mean of a $40000$-iteration Naive chain, after discarding the first $10000$ iterations, as both the reference point $\theta_{\mathrm{ref}}$ and the initial state.
Each chain runs for $200000$ iterations.

For predictive evaluation, we select the parameter states
at iterations $10,20,\ldots,200000$, yielding $S=20000$
draws per method.
For each selected state $\theta^{(s)}$, we generate
$M_{\mathrm{pred}}=2000$ observations on the standardized
scale. We use a common predictive latent bank across
all selected states and all three methods to reduce
Monte Carlo variability in the comparisons.
The simulated observations are then transformed back
to the log-duration scale through
$X_r^{(s)}=\bar x+s_xY_r^{(s)}$ and used to construct
a normalized kernel density estimate (KDE)
$\widehat f_s(x)$. We summarize the predictive densities by their mean
$\overline f(x)=S^{-1}\sum_{s=1}^{S}\widehat f_s(x)$.
Agreement with the observed-data KDE $f_{\mathrm{obs}}$
is measured by the $L_1$ and $L_2$ distances,
\[
L_1
=
\int|\overline f(x)-f_{\mathrm{obs}}(x)|\,\mathrm dx,
\qquad
L_2
=
\left[
\int|\overline f(x)-f_{\mathrm{obs}}(x)|^2\,\mathrm dx
\right]^{1/2}.
\]

\section{Proof of Main Results}\label{Proof of Main Results}
\subsection{Auxiliary Lemmas: Uniform Concentration for the Full Sample}
To guarantee the stability of the full-sample objective while accommodating potential non-smoothness in the moment function, we restrict attention to a high-probability event on which the full-sample empirical moments remain uniformly close to their population analogues. The required bounds are
summarized below.
\begin{lemma}
\label{lemmaprobB1}
Under Assumptions~$A.1-A.2$ and $B.1-B.2$, for any positive constant $c$, there exists a positive constant $c_0$ such that for $\bX^{(n)}\sim \m P^*{}^{\otimes n}$, it holds with probability at least $1-\frac{1}{n^c}$ that 
\begin{enumerate}
\item 
\begin{equation*}
    \begin{aligned}
        \underset{\theta\in \Theta} {\sup}\, \Big\|\frac{1}{n}\sum_{i=1}^n g(X_i,\theta)-\mb E_{X\sim \m P^*} [g(X,\theta)]\Big\|_2 \leq c_0\sqrt{\frac{\log n}{n}}.    \end{aligned}
\end{equation*}
\item 
\begin{equation*}
    \begin{aligned}    
    \underset{\theta, \theta^{\prime}       \in \Theta} {\sup} \frac{\left| \frac{1}{n} \sum_{i=1}^n \ell(X_i, \theta) - \frac{1}{n} \sum_{i=1}^n \ell(X_i, \theta') - \mb E_{X\sim \m P^*}[\ell(X, \theta)] + \mb E_{X\sim \m P^*}[\ell(X, \theta')] \right| }{
\sqrt{\frac{\log n}{n}}  \left\| \theta - \theta' \right\|_2 + \frac{\log n}{n} 
}\leq c_0.
\end{aligned}
      \end{equation*}
 \item 
\begin{equation*}
    \begin{aligned}    
    \underset{\theta, \theta^{\prime}       \in \Theta} {\sup} \frac{\left\| \frac{1}{n} \sum_{i=1}^n g(X_i, \theta) - \frac{1}{n} \sum_{i=1}^n g(X_i, \theta')- \mb E_{X\sim \m P^*}[g(X, \theta)] + \mb E_{X\sim \m P^*}[g(X, \theta')]    \right\|_2 }{
\sqrt{\frac{\log n}{n}}  \| \theta - \theta' \|_2^{\beta_1} + \frac{\log n}{n}  
}\leq c_0.
\end{aligned}
\end{equation*}
\item 
\begin{equation*}
{\footnotesize
    \begin{aligned}    
    \sup_{\theta, \theta' \in \Theta} 
\frac{
\Bigl| 
\frac{1}{n} \sum_{i=1}^n \bigl[ \ell(X_i,\theta) - \ell(X_i,\theta') - g(X_i,\theta')\cdot(\theta-\theta') \bigr]
- \mb E_{X\sim \m P^*}\bigl[ \ell(X,\theta) - \ell(X,\theta') - g(X,\theta')\cdot(\theta-\theta') \bigr]
\Bigr|
}{
\sqrt{\frac{\log n}{n}} \|\theta-\theta'\|_2^{\beta_1+1}
+ \frac{\log n}{n} \|\theta-\theta'\|_2
+ \left( \frac{\log n}{n} \right)^2
}
\le c_0 .
\end{aligned}}
      \end{equation*}\end{enumerate}
\end{lemma}
To study the stationary distribution of the subsampling Markov chain with control variates, we need stronger concentration bounds.
\begin{lemma}\label{lemmaprobB2}
Under the Assumptions~$A.1-A.2$ and $B.1-B.3$, for any positive constant $c$, there exists a positive constant $c_0$ such that for $\bX^{(n)}\sim \m P^*{}^{\otimes n}$, it holds with probability at least $1-\frac{1}{n^c}$ that 
\begin{enumerate}
\item \begin{equation*}
    \begin{aligned}
\underset{\theta, \theta^{\prime}       \in \Theta} {\sup}\frac{\, \frac{1}{n}\sum_{i=1}^n      \left\|g(X_i,\theta)-g(X_i,\theta')\right\|^2_2}{\|\theta-\theta'\|^{2\beta_1}_2+\sqrt{\frac{\log n}{n}} \,\|\theta-\theta'\|^{\beta_1+\beta_2}_2+\frac{\log n}{n}\,\|\theta-\theta'\|^{2\beta_2}_2+(\frac{\log n}{n})^2
} \leq c_0.
      \end{aligned}
      \end{equation*}           \item 
\begin{equation*}
    \begin{aligned}
\underset{\theta, \theta^{\prime}       \in \Theta} {\sup}\frac{\, \frac{1}{n}\sum_{i=1}^n      \left|\ell(X_i,\theta)- \ell(X_i,\theta')-g(X_i,\theta')\cdot(\theta-\theta')\right|^2}{\|\theta-\theta'\|^{2\beta_1+2}_2+\sqrt{\frac{\log n}{n}} \,\|\theta-\theta'\|^{\beta_1+\beta_2+2}_2+\frac{\log n}{n}\,\|\theta-\theta'\|^{2\beta_2+2}_2+(\frac{\log n}{n})^4
} \leq c_0.
      \end{aligned}
      \end{equation*} 
      \item 
\begin{equation*}
    \begin{aligned}    
    \underset{\theta, \theta^{\prime}       \in \Theta} {\sup} \frac{\Big\| \frac{1}{n} \sum_{i=1}^n g(X_i, \theta) - \frac{1}{n} \sum_{i=1}^n g(X_i, \theta')- \mb E_{X\sim \m P^*}[g(X, \theta)] + \mb E_{X\sim \m P^*}[g(X, \theta')]    \Big\|_2 }{
\sqrt{\frac{\log n}{n}}  \| \theta - \theta' \|_2^{\beta_1} + \frac{\log n}{n}   \| \theta - \theta' \|_2^{\beta_2}+( \frac{\log n}{n})^2}\leq c_0.
\end{aligned}
      \end{equation*} \item 
{\footnotesize
\begin{equation*}
\sup_{\theta,\theta'\in\Theta}
\frac{
    \Bigl|
    \frac{1}{n}\sum_{i=1}^n
    \bigl[
        \ell(X_i,\theta)-\ell(X_i,\theta')
        -g(X_i,\theta')(\theta-\theta')
    \bigr]
    -
    \mb E_{X\sim\m P^*}
    \bigl[
        \ell(X,\theta)-\ell(X,\theta')
        -g(X,\theta')(\theta-\theta')
    \bigr]
    \Bigr|
}{
    \sqrt{\frac{\log n}{n}}\,
    \|\theta-\theta'\|_2^{\beta_1+1}
    +
    \frac{\log n}{n}\,
    \|\theta-\theta'\|_2^{\beta_2+1}
    +
    \Big(\frac{\log n}{n}\Big)^3
}
\leq c_0.
\end{equation*}
}
      \end{enumerate}
\end{lemma}

\subsection{Proof of Theorem \ref{th:1}}
Let $\mathcal{B}\subset \m X^{\otimes n}$ denote the event on which all conclusions in Lemma \ref{lemmaprobB1} hold. For sufficiently large $c_0$, we have  $\m P^*{}^{\otimes n}(\m B^c)\leq\frac{1}{n^2}$. Unless otherwise specified, all subsequent analysis is performed under $\bX^{(n)}\in \m B$.
\subsubsection{
Distance between
\texorpdfstring{$\theta^{\diamond}(\bX^{(n)})$}{theta-diamond(X)}
and
\texorpdfstring{$\widehat{\theta}(\boldsymbol{X}^{(n)})$}{theta-hat(X)}
}
Denote
$ \theta^{\diamond}(\bX^{(n)})= \theta^*-\m H(\theta^*)^{-1}\frac{1}{n}\sum_{i=1}^ng( X_i,\theta^*)$. For notational simplicity, we abbreviate $\theta^\diamond(\boldsymbol{X}^{(n)})$ to $\theta^\diamond$, $\widehat{\theta}(\boldsymbol{X}^{(n)})$ to $\widehat{\theta}$, and $\theta^\dagger(\boldsymbol{X}^{(n)})$ to $\theta^\dagger$.
By definition, $\wh\theta$ minimizes the empirical risk
$\frac{1}{n}\sum_{i=1}^n \ell(X_i,\theta)$. Hence,
\[
\frac{1}{n}\sum_{i=1}^n \ell(X_i,\wh\theta)
\le
\frac{1}{n}\sum_{i=1}^n \ell(X_i,\theta^*).
\]
Moreover, by Assumption~$A.1$,
\[
\mb E_{X\sim \m P^*}[\ell(X,\wh\theta)]
-
\mb E_{X\sim \m P^*}[\ell(X,\theta^*)]
\gtrsim
\|\wh\theta-\theta^*\|_2^2.
\]

On the other hand, the second statement of Lemma \ref{lemmaprobB1} yields
\begin{equation*}
\Big|
\frac{1}{n}\sum_{i=1}^n \ell(X_i,\wh\theta)
-
\frac{1}{n}\sum_{i=1}^n \ell(X_i,\theta^*)
-
\mb E_{X\sim \m P^*}[\ell(X,\wh\theta)]
+
\mb E_{X\sim \m P^*}[\ell(X,\theta^*)]
\Big|
\lesssim
\sqrt{\frac{\log n}{n}}\,
\|\wh\theta-\theta^*\|_2
+
\frac{\log n}{n},
\end{equation*}
we obtain
\[
\|\wh\theta-\theta^*\|_2^2
\lesssim
\mb E_{X\sim \m P^*}[\ell(X,\wh\theta)]
-
\mb E_{X\sim \m P^*}[\ell(X,\theta^*)]
\lesssim
\sqrt{\frac{\log n}{n}}\,
\|\wh\theta-\theta^*\|_2
+
\frac{\log n}{n}.
\]
Solving the resulting quadratic inequality gives
\begin{equation*}
\|\wh\theta-\theta^*\|_2
\lesssim
\sqrt{\frac{\log n}{n}}.
\end{equation*}
 The first statement of Lemma \ref{lemmaprobB1} 
 gives $\|\theta^* - \theta^{\diamond}\|_2=\Big\|\m H(\theta^*)^{-1}\frac{1}{n}\sum_{i=1}^ng( X_i,\theta^*)\Big\|_2\lesssim\sqrt{\frac{\log n}{n}}.$ Therefore, we get
\begin{equation}\label{estimate 0.1}
\|\wh\theta - \theta^{\diamond}\|_2 
\leq\|\wh\theta - \theta^*\|_2 + \|\theta^* - \theta^{\diamond}\|_2\lesssim \sqrt{\frac{\log n}{n}}.
\end{equation}
By Assumption A.1, for $u$ and $v$ near $\theta^*$, we have
\[
\|\mathcal{H}(u)-\mathcal{H}(v)\|_{\mathrm{F}}
\leq L\left(\|u-v\|_2+2\|v-\theta^*\|_2\right).
\]
hence
\begin{equation}\label{Hessian Taylor}
\begin{aligned}
&\Big|\mathcal{R}(\wh\theta) - \mathcal{R}(\theta^{\diamond}) 
- \nabla\mathcal{R}(\theta^{\diamond}) \cdot (\wh\theta - \theta^{\diamond})
- \frac{1}{2}(\wh\theta - \theta^{\diamond})^T 
\m H(\theta^{\diamond}) (\wh\theta - \theta^{\diamond})\Big| \\
&= \Big|\frac{1}{2}(\wh\theta - \theta^{\diamond})^T 
\m H(\theta^{\diamond}) (\wh\theta - \theta^{\diamond}) - (\wh\theta - \theta^{\diamond})^T 
\int_{0}^{1}(1-t)\, \m H\!\big(\theta^{\diamond} + t(\wh\theta - \theta^{\diamond})\big) dt \,
(\wh\theta - \theta^{\diamond})\Big| \\
&\lesssim \|\wh\theta - \theta^{\diamond}\|_2^3+\|\theta^{\diamond}-\theta^*\|_2
\|\wh\theta - \theta^{\diamond}\|_2^2,
\end{aligned}
\end{equation}
and by the last statement in Lemma \ref{lemmaprobB1}, it holds that
\begin{equation*}
\begin{aligned}
&\Big| \big[ \frac{1}{n} \sum_{i=1}^n \ell(X_i,\wh\theta) 
- \frac{1}{n} \sum_{i=1}^n \ell(X_i,\theta^{\diamond}) 
- \frac{1}{n} \sum_{i=1}^n g(X_i,\theta^{\diamond}) \cdot 
(\wh\theta-\theta^{\diamond}) \big] \\
&\quad - \big[ \mathcal{R}(\wh\theta) - \mathcal{R}(\theta^{\diamond}) 
- \nabla\mathcal{R}(\theta^{\diamond}) \cdot (\wh\theta - \theta^{\diamond}) \big] \Big|
\\
&\lesssim \sqrt{\frac{\log n}{n}} \|\wh\theta-\theta^{\diamond}\|_2^{\beta_1+1}
+ \frac{\log n}{n} \|\wh\theta-\theta^{\diamond}\|_2
+ \Big( \frac{\log n}{n}\Big)^2.
\end{aligned}
\end{equation*}
As we have 
\[
(\wh\theta - \theta^{\diamond})^T \m{H}(\theta^{\diamond}) 
(\wh\theta - \theta^{\diamond}) \gtrsim \|\wh\theta-\theta^{\diamond}\|^2_2,
\]
the following holds 
\begin{equation}\label{estimate 0.2}
\begin{aligned}
&\Big|\frac{1}{n}\sum_{i=1}^n g(X_i,\theta^{\diamond})\cdot(\wh\theta-\theta^{\diamond})\Big|
\\
&\geq -\frac{1}{n}\sum_{i=1}^n g(X_i,\theta^{\diamond})\cdot(\wh\theta-\theta^{\diamond})
\\
&\geq \frac{1}{2}(\wh\theta - \theta^{\diamond})^T \m H(\theta^{\diamond})(\wh\theta - \theta^{\diamond})
\\
&\quad - \Big|\mathcal{R}(\wh\theta) - \mathcal{R}(\theta^{\diamond}) 
- \nabla\mathcal{R}(\theta^{\diamond})\cdot (\wh\theta - \theta^{\diamond})
- \frac{1}{2}(\wh\theta - \theta^{\diamond})^T \m H(\theta^{\diamond})(\wh\theta - \theta^{\diamond})\Big|
\\
&\quad - \Big|\Big[\frac{1}{n}\sum_{i=1}^n \ell(X_i,\wh\theta)-\frac{1}{n}\sum_{i=1}^n \ell(X_i,\theta^{\diamond})
-\frac{1}{n}\sum_{i=1}^n g(X_i,\theta^{\diamond})\cdot(\wh\theta-\theta^{\diamond})\Big]
\\
&\quad -\Big[\mathcal{R}(\wh\theta) - \mathcal{R}(\theta^{\diamond}) 
- \nabla\mathcal{R}(\theta^{\diamond})\cdot(\wh\theta - \theta^{\diamond})\Big]\Big|
 - \Big(\frac{1}{n}\sum_{i=1}^n \ell(X_i,\wh\theta)-\frac{1}{n}\sum_{i=1}^n \ell(X_i,\theta^{\diamond})\Big)
\\
&\gtrsim \|\wh\theta-\theta^{\diamond}\|^2_2 
- \|\theta^{\diamond}-\theta^*\|_2
\|\wh\theta - \theta^{\diamond}\|_2^2-\|\wh\theta-\theta^{\diamond}\|^3_2 - \Big(\sqrt{\frac{\log n}{n}}\|\wh\theta-\theta^{\diamond}\|_2^{\beta_1+1}
+ \frac{\log n}{n}\|\wh\theta-\theta^{\diamond}\|_2
+ \Big(\frac{\log n}{n}\Big)^2\Big)
\\
&\gtrsim \|\wh\theta-\theta^{\diamond}\|^2_2 
- \Big(\sqrt{\frac{\log n}{n}}\|\wh\theta-\theta^{\diamond}\|_2^{\beta_1+1}
+ \frac{\log n}{n}\|\wh\theta-\theta^{\diamond}\|_2
+ \Big(\frac{\log n}{n}\Big)^2\Big).
\end{aligned}
\end{equation}
The last inequality follows by absorbing both remainder terms into the leading quadratic term, using the bound $\|\theta^\diamond-\theta^*\|_2 + \|\widehat{\theta}-\theta^\diamond\|_2 \lesssim \sqrt{\frac{\log n}{n}}$. As $\|\theta^*-\theta^{\diamond}\|\lesssim\sqrt{\frac{\log n}{n}}$, the third statement of Lemma \ref{lemmaprobB1} gives 
\begin{equation*}
\begin{aligned}
&\Big\| \frac{1}{n} \sum_{i=1}^n g(X_i, \theta^*) 
- \frac{1}{n} \sum_{i=1}^n g(X_i, \theta^{\diamond})
- \nabla\mathcal{R}(\theta^*) 
+ \nabla\mathcal{R}(\theta^{\diamond}) \Big\|_2
\\
&\lesssim \sqrt{\frac{\log n}{n}} \| \theta^* - \theta^{\diamond} \|_2^{\beta_1} 
+ \frac{\log n}{n}
\\
&\lesssim \Big( \frac{\log n}{n} \Big)^{\frac{\beta_1+1}{2}}.
\end{aligned}
\end{equation*}
Moreover
\begin{equation*}
    \begin{aligned}
    \nabla\m R(\theta^{\diamond}) 
    &= \nabla\m R(\theta^{\diamond}) - \nabla\m R(\theta^*) \\
    &= \m H(\theta^*) \bigl(\theta^{\diamond}-\theta^*\bigr) 
       + \m O\bigl(\|\theta^{\diamond}-\theta^*\|^2_2\bigr) \\
    &= -\frac{1}{n}\sum_{i=1}^n g(X_i,\theta^*) 
       + \m O\left(\frac{\log n}{n}\right),
    \end{aligned}
\end{equation*}
we get
\begin{equation}\label{estimate 0.3}\begin{aligned}
\Big\|\frac{1}{n}\sum_{i=1}^n g\big(X_i,\theta^{\diamond}\big)\Big\|_2\lesssim\Big(\frac{\log n}{n}\Big)^{\frac{\beta_1+1}{2}}.
\end{aligned}
\end{equation}
Combining \eqref{estimate 0.2} and \eqref{estimate 0.3}, we get
\begin{equation*}
\begin{aligned}
& \|\wh\theta-\theta^{\diamond}\|^2_2
- \Big(\sqrt{\frac{\log n}{n}}\|\wh\theta-\theta^{\diamond}\|_2^{\beta_1+1}
+ \frac{\log n}{n}\|\wh\theta-\theta^{\diamond}\|_2
+ \Big(\frac{\log n}{n}\Big)^2\Big)
\\
&\lesssim \Big|\frac{1}{n}\sum_{i=1}^n g(X_i,\theta^{\diamond})
\cdot(\wh\theta-\theta^{\diamond})\Big|
\\
&\lesssim \Big(\frac{\log n}{n}\Big)^{\frac{\beta_1+1}{2}} \|\wh\theta-\theta^{\diamond}\|_2,
\end{aligned}
\end{equation*}
hence \begin{equation}\label{estimate 0.4}
    \|\wh\theta-\theta^{\diamond}\|_2 \lesssim \Big(\frac{\log n}{n}\Big)^{\frac{\beta_1+1}{2}}.    
\end{equation}
\subsubsection{Proof of the Total Variation Bound}
By definition
\begin{equation*}
\begin{aligned}
&\wt{\pi}_{\mathcal{L}_{\rm zcv},\alpha}(\theta)=
\frac{\mb{E}_{\wX^{(m)} \sim  \mu_n^{\otimes m} }\Big[\pi(\theta) \exp\Big(-\alpha\Big(\frac{1}{m}\sum_{j=1}^m\ell(\wt X_j,\theta)- \frac{1}{m}\sum_{j=1}^m\ell (\wt X_j,\theta^{\dagger})\Big)\Big)\Big]}{\mb{E}_{\wX^{(m)} \sim  \mu_n^{\otimes m} }\Big[\int\pi(\theta) \exp\Big(-\alpha\Big(\frac{1}{m}\sum_{j=1}^m\ell(\wt X_j,\theta)- \frac{1}{m}\sum_{j=1}^m\ell(\wt X_j,\theta^{\dagger})\Big)\Big)\,\dd\theta\Big]}\\&=\frac{\mb{E}_{\wX^{(m)} \sim  \mu_n^{\otimes m} }\Big[\pi(\theta) \exp\Big(-\alpha\Big(\frac{1}{m}\sum_{j=1}^m[\ell(\wt X_j,\theta)- \ell (\wt X_j,\theta^{\dagger})]-\frac{1}{n}\sum_{i=1}^n [\ell(X_i,\theta^{\diamond})-\ell(X_i,\theta^{\dagger})]
\Big)\Big)\Big]}{\mb{E}_{\wX^{(m)} \sim  \mu_n^{\otimes m} }\Big[\int\pi(\theta) \exp\Big(-\alpha\Big(\frac{1}{m}\sum_{j=1}^m[\ell(\wt X_j,\theta)- \ell (\wt X_j,\theta^{\dagger})]-\frac{1}{n}\sum_{i=1}^n [\ell(X_i,\theta^{\diamond})-\ell(X_i,\theta^{\dagger})]\Big)\Big)\,\dd\theta\Big]},
\end{aligned}
\end{equation*}and
\begin{equation*}
 \begin{aligned}
     \pi_\alpha(\theta)&=
\frac{\pi(\theta)\exp\Big(-\frac{\alpha
}{n}\sum_{i=1}^n \ell(X_i,\theta)\Big)}{\int\pi(\theta)\exp\Big(-\alpha\Big( \frac{1}{n}\sum_{i=1}^n \ell(X_i,\theta)\Big)\Big)\,\dd\theta}\\&=
\frac{\pi(\theta)\exp\Big(- \alpha\Big(\frac{1}{n}\sum_{i=1}^n \ell(X_i,\theta)-\frac{1}{n}\sum_{i=1}^n \ell(X_i,\theta^{\diamond})\Big)\Big)}{\int\pi(\theta)\exp\Big(-\alpha\Big( \frac{1}{n}\sum_{i=1}^n \ell(X_i,\theta)-\frac{1}{n}\sum_{i=1}^n \ell(X_i,\theta^{\diamond})\Big)\Big)\,\dd\theta} ,\end{aligned}   
\end{equation*} 
The key step is to bound the deviation between the full-sample and
minibatch loss. To this end, we next derive
a concentration inequality for $\widetilde{\bX}^{(m)}\sim\mu_n^{\otimes m}$.
\begin{lemma}\label{lemmaprobA1}
    Given $\bX^{(n)}\in \m X^{\otimes n} $, under Assumptions~$A.1-A.2$ and $B.1-B.2$, for any positive constant $c$, there exists a constant $c_0$ such that, for
$\wX^{(m)}=
(\wt X_1,\wt X_2,\ldots,\wt X_m)
\sim \mu_n^{\otimes m}$,
the following holds with probability at least
$1-m^{-c}$.
    \begin{enumerate}
\item 
\begin{equation*}
     \begin{aligned}
        \underset{\theta\in \Theta}{\sup}\Big\|\frac{1}{m}\sum_{j=1}^mg (\wt X_j,\theta)-\frac{1}{n}\sum_{i=1}^ng(X_i,\theta)\Big\|_2\leq c_0\sqrt{\frac{\log m}{m}}.  
        \end{aligned}
 \end{equation*}
 
 \item 
\begin{equation*}
         \underset{\theta, \theta^{\prime}\in \Theta}{\sup}\frac{\left|\frac{1}{m}\sum_{j=1}^m\left[\ell (\wt X_j,\theta)-\ell (\wt X_j,\theta')\right]-\frac{1}{n}\sum_{i=1}^n \left[\ell(X_i,\theta)-\ell(X_i,\theta')\right]\right|}{
         \sqrt{\frac{\log m}{m}}\|\theta-\theta'\|_2+\frac{\log m}{m}         }\leq c_0.
\end{equation*}
\end{enumerate}
\end{lemma}
In the remainder of the analysis, we fix an $\bX^{(n)} \in \m B$. Let $\m A$ denote the event of $\widetilde{\bX}^{(m)}$ on which all
conclusions of Lemma \ref{lemmaprobA1} hold. By choosing $c_0$ to be
sufficiently large, we have
 \begin{equation}\label{eqn:prob1}   \mb P_{\wX^{(m)}\sim \mu_n^{\otimes m}}(\m A)\geq 1-\frac{1}{m^2\alpha^{d+1}}.
\end{equation}
Our proof proceeds by bounding
\begin{equation}\label{goal}
    \begin{aligned}
       & \alpha^{\frac{d}{2}}\int_{\Theta}\Big|\mb{E}_{\wX^{(m)}\sim \mu_n^{\otimes m}}\Big[\pi(\theta)\exp\Big(-\alpha\Big(\frac{1}{m}\sum_{j=1}^m\ell^* (\wt X_j,\theta)- \frac{1}{n}\sum_{i=1}^n\ell^* (X_i,\theta^{\diamond})\Big)\Big)\cdot\bold{1}(\wX^{(m)}\in \m A)\Big] \\
         &\qquad-\mb{E}_{\wX^{(m)}\sim \mu_n^{\otimes m}}\Big[\pi(\theta)\exp\Big(-\alpha\Big(  \frac{1}{n}\sum_{i=1}^n\ell^* (X_i,\theta)-\frac{1}{n}\sum_{i=1}^n\ell^* (X_i,\theta^{\diamond})\Big)\Big)\cdot\bold{1}(\wX^{(m)}\in \m A)\Big]\Big|\,\dd\theta
         \end{aligned}
         \end{equation}
Next we  show that the contribution to the expectation from its complement \(\mathcal A^c\) is negligible.
 We bound the integrals \eqref{goal} by splitting the domain of integration into three regions. Specifically, define
\begin{equation}\label{setstheta1}
\begin{aligned}
&\Theta_1=\Big\{\theta\in \Theta : \|\theta-\theta^{\diamond}\|_2\leq \delta_1\sqrt{\frac{\log m}{\alpha}}\Big\},\\
&\Theta_2=\Big\{\theta\in \Theta : \delta_1\sqrt{\frac{\log m}{\alpha}}<\|\theta-\theta^{\diamond}\|_2\leq \delta_2\Big\},\\
&\Theta_3=\Big\{\theta\in \Theta : \|\theta-\theta^{\diamond}\|_2> \delta_2\Big\}.
\end{aligned}
\end{equation}
Here, $\delta_1,\delta_2>0$ are suitably chosen constants. By the last statement in Lemma \ref{lemmaprobB1}, we have \begin{equation}\label{estimate 1.1}
\begin{aligned}
&\Big| \frac{1}{n} \sum_{i=1}^n \ell(X_i, \theta) 
      - \frac{1}{n} \sum_{i=1}^n \ell(X_i, \theta^{\diamond}) 
      - \frac{1}{n} \sum_{i=1}^n g(X_i, \theta^{\diamond}) 
        \cdot (\theta - \theta^{\diamond}) \\
&\quad - \m{R}(\theta) + \m{R}(\theta^{\diamond}) 
      + \nabla\m{R}(\theta^{\diamond}) 
        \cdot (\theta - \theta^{\diamond}) \Big| \\
&\lesssim \sqrt{\frac{\log n}{n}} \, \| \theta - \theta^{\diamond} \|_2^{\beta_1+1} 
        + \frac{\log n}{n} \| \theta - \theta^{\diamond} \|_2 
        + \Big( \frac{\log n}{n} \Big)^2.
\end{aligned}
\end{equation}
Combining \eqref{estimate 0.3} with \eqref{estimate 1.1}, we get
\begin{equation*}
\begin{aligned}
&\Big| \frac{1}{n} \sum_{i=1}^n \ell(X_i, \theta) 
     - \frac{1}{n} \sum_{i=1}^n \ell(X_i, \theta^{\diamond}) 
     - \m R(\theta) + \m R(\theta^{\diamond}) + \nabla\m R(\theta^{\diamond})\cdot(\theta - \theta^{\diamond})\Big| \\
&\lesssim \sqrt{\frac{\log n}{n}} \, \| \theta - \theta^{\diamond} \|_2^{\beta_1+1} 
     + \Big( \frac{\log n}{n} \Big)^{\frac{\beta_1+1}{2}}\|\theta - \theta^{\diamond} \|_2
     + \Big( \frac{\log n}{n} \Big)^2,
\end{aligned}
\end{equation*}
note that\begin{equation*}
    \begin{aligned}
    &\qquad\m R(\theta) - \m R(\theta^{\diamond}) 
     - \nabla\m R(\theta^{\diamond}) \cdot (\theta - \theta^{\diamond}) \\
    &\quad = \int_0^1 (1 - t) \, \left(\theta - \theta^{\diamond}\right)^T 
       \m H\big(\theta^{\diamond} + t(\theta - \theta^{\diamond})\big)
       (\theta - \theta^{\diamond}) \, dt,
    \end{aligned}
\end{equation*}
and by the Cauchy--Schwarz inequality we have
\begin{align*}
    \sqrt{\frac{\log n}{n}} \, \| \theta - \theta^{\diamond} \|_2^{\beta_1+1}
    &\leq 4 \frac{\log n}{n} \| \theta - \theta^{\diamond} \|_2^{2\beta_1} 
    + \frac{ \| \theta - \theta^{\diamond} \|_2^2 }{8}, \\[6pt]
    \Big( \frac{\log n}{n} \Big)^{\frac{\beta_1+1}{2}} \|\theta - \theta^{\diamond} \|_2
    &\leq 4 \Big( \frac{\log n}{n} \Big)^{\beta_1+1}   
    + \frac{ \| \theta - \theta^{\diamond} \|_2^2 }{8}.
\end{align*}
Then, for $\theta\in\Theta_1\cup\Theta_2$, we have
\begin{equation}\label{estimate 1.2}
M_1\Big(\|\theta - \theta^{\diamond} \|^2_2 
    - \frac{\log n}{n}\| \theta -\theta^{\diamond} \|_2^{2\beta_1}
    - \Big( \frac{\log n}{n} \Big)^{\beta_1+1}\Big)
    \leq \frac{1}{n} \sum_{i=1}^n \ell(X_i, \theta) 
    - \frac{1}{n} \sum_{i=1}^n \ell(X_i, \theta^{\diamond}),
\end{equation}
and
\begin{equation}\label{estimate 1.3}
\frac{1}{n} \sum_{i=1}^n \ell(X_i, \theta) 
    - \frac{1}{n} \sum_{i=1}^n \ell(X_i, \theta^{\diamond})
    \leq M_2\Bigl( \|\theta -\theta^{\diamond} \|^2_2 
    + \frac{\log n}{n}\| \theta - \theta^{\diamond} \|_2^{2\beta_1}
    + \Big( \frac{\log n}{n} \Big)^{\beta_1+1}\Bigr).
\end{equation}
Moreover, it holds that 
\begin{equation*}
    \begin{aligned}
       \m R(\theta)-\m R(\theta^*) \gtrsim  \|\theta-\theta^*\|^2_2.   \end{aligned}
\end{equation*}
With the help of $\|\theta^*-\theta^{\diamond}\|_2\lesssim \sqrt{\frac{\log n}{n}}$ and $\|\nabla\m R(\theta^{\diamond})\|_2\lesssim\sqrt{\frac{\log n}{n}}$, we obtain a positive constant $M_3$ such that for $\theta\in\Theta_3$,
\begin{equation}\label{global hessian}
    M_3\|\theta - \theta^{\diamond} \|^2_2\leq\frac{1}{n} \sum_{i=1}^n \ell(X_i, \theta) - \frac{1}{n} \sum_{i=1}^n \ell(X_i, \theta^{\diamond}).\end{equation}
\paragraph{1. \underline{On $\Theta_1$.}}
For $\wX^{(m)}\in \m A$ and $\theta\in \Theta_1$, combining the second statement of Lemma \ref{lemmaprobA1} with \eqref{estimate 0.4}, we obtain
\begin{equation}\label{estimate 1.4}
    \begin{aligned}
        &\Big|\frac{1}{m}\sum_{j=1}^m\Big(\ell (\wt X_j,\theta)- \ell(\wt X_j,\theta^\dagger)\Big)
          -\frac{1}{n}\sum_{i=1}^n\Big(\ell (X_i,\theta)- \ell(X_i,\theta^{\dagger})\Big)\Big|\\
        &\lesssim \sqrt{\frac{\log m}{m}}\|\theta-\theta^{\dagger}\|_2+\frac{\log m}{m}\\
        &\lesssim\sqrt{\frac{\log m}{m}}\bigg(\|\theta-\theta^{\diamond}\|_2
          +\|\wh\theta-\theta^{\diamond}\|_2+C_{m,n}\bigg)+\frac{\log m}{m}\\
        &\lesssim\sqrt{\frac{\log m}{m}}\bigg(\|\theta-\theta^{\diamond}\|_2 +C_{m,n}\bigg)+\frac{\log m}{m}.
    \end{aligned}
\end{equation}
To simplify the notation, we denote 
\begin{equation*}
\ell^* (x,\theta)=\ell (x,\theta)- \ell(x,\theta^\dagger),
\end{equation*}thus we have 
\begin{equation*}
\begin{aligned}
&\Big|\frac{1}{m}\sum_{j=1}^m\ell^*(\wt X_j,\theta)-\frac{1}{n}\sum_{i=1}^n\ell^*(X_i,\theta^{\diamond})-\frac{1}{n}\sum_{i=1}^n\ell^*(X_i,\theta)+\frac{1}{n}\sum_{i=1}^n\ell^*(X_i,\theta^{\diamond})\Big|\\
&=\Big|\frac{1}{m}\sum_{j=1}^m\ell^*(\wt X_j,\theta)-\frac{1}{n}\sum_{i=1}^n\ell^*(X_i,\theta)\Big|\\
&\lesssim\sqrt{\frac{\log m}{m}}\Big(\|\theta-\theta^{\diamond}\|_2+C_{m,n}\Big)+\frac{\log m}{m}.
\end{aligned}
\end{equation*}
    Hence
\begin{equation*}
\begin{aligned}
&\exp\Big(-\alpha\Big(\frac{1}{m}\sum_{j=1}^m\ell^* (\wt X_j,\theta)- \frac{1}{n}\sum_{i=1}^n\ell^* (X_i,\theta^{\diamond})\Big)\Big)-\exp\Big(-\alpha\Big( \frac{1}{n}\sum_{i=1}^n\ell (X_i,\theta)-\frac{1}{n}\sum_{i=1}^n\ell(X_i,\theta^{\diamond})\Big)\Big)\\
&=\exp\Big(-\alpha\Big(\frac{1}{m}\sum_{j=1}^m\ell^* (\wt X_j,\theta)- \frac{1}{n}\sum_{i=1}^n\ell^* (X_i,\theta^{\diamond})\Big)\Big)-\exp\Big(-\alpha\Big( \frac{1}{n}\sum_{i=1}^n\ell^* (X_i,\theta)-\frac{1}{n}\sum_{i=1}^n\ell^* (X_i,\theta^{\diamond})\Big)\Big)\\
&=\exp\Big(-\alpha\Big( \frac{1}{n}\sum_{i=1}^n\ell^* (X_i,\theta)-\frac{1}{n}\sum_{i=1}^n\ell^* (X_i,\theta^{\diamond})\Big)\Big)\\
&\quad\cdot\Big[\m O\Big(\alpha\frac{(\log m)^2}{m}+\alpha^2\frac{\log m}{m}C^2_{m,n}+\alpha^2\Big(\frac{\log m}{m}\Big)^{2}\Big)-\alpha\Big(\frac{1}{m}\sum_{j=1}^m \ell^*(\wt X_j,\theta)-\frac{1}{n}\sum_{i=1}^n \ell^*(X_i,\theta)\Big) \Big]\\
&=\exp\Big(-\alpha\Big( \frac{1}{n}\sum_{i=1}^n\ell (X_i,\theta)-\frac{1}{n}\sum_{i=1}^n\ell(X_i,\theta^{\diamond})\Big)\Big)\\
&\quad\cdot\Big[\m O\Big(\alpha\frac{(\log m)^2}{m}+\alpha^2\frac{\log m}{m}C^2_{m,n}\Big)-\alpha\Big(\frac{1}{m}\sum_{j=1}^m \ell^*(\wt X_j,\theta)-\frac{1}{n}\sum_{i=1}^n \ell^*(X_i,\theta)\Big) \Big],
\end{aligned}
\end{equation*}
where the equality is due to the fact that \[
\alpha^2\Big(\frac{\log m}{m}\Big)^{2}\lesssim\alpha\frac{(\log m)^2}  
    {m}.\]Therefore, we can obtain
\begin{equation}\label{estimate 1.5}
    \begin{aligned}
       & \alpha^{\frac{d}{2}}\int_{\Theta_1}\Big|\mb{E}_{\wX^{(m)}\sim \mu_n^{\otimes m}}\Big[\pi(\theta)\exp\Big(-\alpha\Big(\frac{1}{m}\sum_{j=1}^m\ell^* (\wt X_j,\theta)- \frac{1}{n}\sum_{i=1}^n\ell^* (X_i,\theta^{\diamond})\Big)\Big)\cdot\bold{1}(\wX^{(m)}\in \m A)\Big] \\
         &\qquad-\mb{E}_{\wX^{(m)}\sim \mu_n^{\otimes m}}\Big[\pi(\theta)\exp\Big(-\alpha\Big(  \frac{1}{n}\sum_{i=1}^n\ell^* (X_i,\theta)-\frac{1}{n}\sum_{i=1}^n\ell^* (X_i,\theta^{\diamond})\Big)\Big)\cdot\bold{1}(\wX^{(m)}\in \m A)\Big]\Big|\,\dd\theta \\
         &\lesssim \alpha^{\frac{d}{2}}\int_{\Theta_1}\pi(\theta)\exp\Big(-\alpha\Big( \frac{1}{n}\sum_{i=1}^n \ell(X_i,\theta)-\frac{1}{n}\sum_{i=1}^n \ell(X_i,\theta^{\diamond})\Big)\Big) \\
         &\qquad\cdot \Big(\alpha\frac{(\log m)^2}{m}+\alpha^2\frac{\log m}{m}C^2_{m,n}+\alpha\cdot\Big|\mb{E}_{\wX^{(m)}\sim \mu_n^{\otimes m}}\Big[\bold{1}(\wX^{(m)}\in \m A) \\
         &\qquad\qquad\cdot \Big(\frac{1}{m}\sum_{j=1}^m \ell(\wt X_j,\theta)-\frac{1}{m}\sum_{j=1}^m \ell(\wt X_j,\theta^{\dagger})-\frac{1}{n}\sum_{i=1}^n \ell(X_i,\theta)+\frac{1}{n}\sum_{i=1}^n \ell(X_i,\theta^\dagger)\Big)\Big]\Big|\Big)\,\dd\theta.
    \end{aligned}
\end{equation}
By assumption $\ell(x,\cdot)$ is Lipschitz and  $\Theta\subset \mb{R}^d$ is compact, we have 
\begin{equation}\label{estimate 1.6}
    \begin{aligned}
        &\Big|\mb{E}_{\wX^{(m)}\sim \mu_n^{\otimes m}}\Big[\bold{1}(\wX^{(m)}\in \m A)\\&\quad \cdot \Big(\frac{1}{m}\sum_{j=1}^m \ell(\wt X_j,\theta)-\frac{1}{m}\sum_{j=1}^m \ell(\wt X_j,\theta^\dagger)-\frac{1}{n}\sum_{i=1}^n \ell(X_i,\theta)+\frac{1}{n}\sum_{i=1}^n \ell(X_i,\theta^\dagger)\Big)\Big]\Big|\\
        &\leq \Big|\mb{E}_{\wX^{(m)}\sim \mu_n^{\otimes m}}\Big[ \Big(\frac{1}{m}\sum_{j=1}^m \ell(\wt X_j,\theta)-\frac{1}{m}\sum_{j=1}^m \ell(\wt X_j,\theta^\dagger)-\frac{1}{n}\sum_{i=1}^n \ell(X_i,\theta)+\frac{1}{n}\sum_{i=1}^n \ell(X_i,\theta^\dagger)\Big)\Big]\Big|\\
        &+\quad \Big|\mb{E}_{\wX^{(m)}\sim \mu_n^{\otimes m}}\Big[\bold{1}(\wX^{(m)}\in \m A^c)\\&\qquad\cdot \Big(\frac{1}{m}\sum_{j=1}^m \ell(\wt X_j,\theta)-\frac{1}{m}\sum_{j=1}^m \ell(\wt X_j,\theta^\dagger)-\frac{1}{n}\sum_{i=1}^n \ell(X_i,\theta)+\frac{1}{n}\sum_{i=1}^n \ell(X_i,\theta^\dagger)\Big)\Big]\Big|\\
        &\lesssim\mb P_{\wX^{(m)}\sim \mu_n^{\otimes m}}(\m A^c)\cdot 
        \underset{x\in \m X,\theta\in\Theta}{\sup}|\ell(x,\theta)-\ell(x,\theta^{\dagger})|\\&\lesssim\mb P_{\wX^{(m)}\sim \mu_n^{\otimes m}}(\m A^c)\cdot\operatorname{diam}(\Theta)        \\&\lesssim \frac{1}{m^2\alpha^{d+1}}.
    \end{aligned}
\end{equation}
Hence, by combining \eqref{estimate 1.2}, \eqref{estimate 1.5}, and \eqref{estimate 1.6} and using $\alpha \lesssim m$, we obtain
\begin{equation}\label{step 1.1}
    \begin{aligned}
       & \alpha^{\frac{d}{2}}\int_{\Theta_1}\Big|\mb{E}_{\wX^{(m)}\sim \mu_n^{\otimes m}}\Big[\pi(\theta)\exp\Big(-\alpha\Big(\frac{1}{m}\sum_{j=1}^m\ell^* (\wt X_j,\theta)- \frac{1}{n}\sum_{i=1}^n\ell^* (X_i,\theta^{\diamond})\Big)\Big)\cdot\bold{1}(\wX^{(m)}\in \m A)\Big] \\
         &\qquad-\mb{E}_{\wX^{(m)}\sim \mu_n^{\otimes m}}\Big[\pi(\theta)\exp\Big(-\alpha\Big(  \frac{1}{n}\sum_{i=1}^n\ell^* (X_i,\theta)-\frac{1}{n}\sum_{i=1}^n\ell^* (X_i,\theta^{\diamond})\Big)\Big)\cdot\bold{1}(\wX^{(m)}\in \m A)\Big]\Big|\,\dd\theta \\
         &\lesssim \alpha^{\frac{d}{2}}\int_{\Theta_1}\pi(\theta)\exp\Big(-\alpha M_1\|\theta - \theta^{\diamond} \|^2_2 
           +\alpha M_1\frac{\log n}{n}\| \theta - \theta^{\diamond} \|_2^{2\beta_1}
           +\alpha M_1\Big( \frac{\log n}{n} \Big)^{\beta_1+1}\Big) \\
         &\qquad\cdot \Big(\alpha\frac{(\log m)^2}{m}+\alpha^2\frac{\log m}{m}C^2_{m,n}+ \frac{1}{m^2\alpha^{d}}\Big)\,\dd\theta \\
         &\lesssim \alpha^{\frac{d}{2}}\int_{\Theta_1}\pi(\theta)\exp\Big(-\alpha M_1\|\theta - \theta^{\diamond} \|^2_2\Big)
           \cdot \Big(\alpha\frac{(\log m)^2}{m}+\alpha^2\frac{\log m}{m}C^2_{m,n}\Big)\,\dd\theta \\
         &\lesssim \alpha\frac{(\log m)^2}{m}+ \alpha^2\frac{\log m}{m}C^2_{m,n},
    \end{aligned}
\end{equation}
Moreover, since $\pi(\theta^*) \ge l$ and $\pi$ is locally Lipschitz in a neighborhood of $\theta^*$, we have $\pi(\theta) \ge \frac{l}{2}$ for all $\theta \in \Theta_1$; together with $\alpha \lesssim m$ and \eqref{estimate 1.3}, this implies that there exists a constant $k_1 > 0$ such that
\begin{equation}\label{lowerbound 1}
\begin{aligned}
&\alpha^{\frac{d}{2}}\int_{\Theta_1}\pi(\theta)\exp\Big(-\alpha\Big( \frac{1}{n}\sum_{i=1}^n \ell^*(X_i,\theta)-\frac{1}{n}\sum_{i=1}^n \ell^*(X_i,\theta^{\diamond})\Big)\Big)\,\dd\theta\\
&= \alpha^{\frac{d}{2}}\int_{\Theta_1}\pi(\theta)\exp\Big(-\alpha\Big( \frac{1}{n}\sum_{i=1}^n \ell(X_i,\theta)-\frac{1}{n}\sum_{i=1}^n \ell(X_i,\theta^{\diamond})\Big)\Big)\,\dd\theta\\
&\geq \alpha^{\frac{d}{2}}\int_{\Theta_1}\pi(\theta)\exp\Big(-\alpha M_2\|\theta - \theta^{\diamond} \|^2_2 
     -\alpha M_2\frac{\log n}{n}\| \theta - \theta^{\diamond} \|_2^{2\beta_1}
     -\alpha M_2\Big( \frac{\log n}{n} \Big)^{\beta_1+1}\Big)\,\dd\theta\\
&\gtrsim \alpha^{\frac{d}{2}}\int_{\Theta_1}\pi(\theta)\exp\Big(-\alpha M_2\|\theta - \theta^{\diamond} \|^2_2\Big)\,\dd\theta \\&\geq k_1.
\end{aligned}
\end{equation}
\paragraph{2. \underline{On $\Theta_2$.}}
For $\theta\in \Theta_2$ and $\wX^{(m)}\in \m A$, by the second statement of Lemma \ref{lemmaprobA1}, we have 
\begin{equation*}
    \begin{aligned}
        &-\Big(\frac{1}{m}\sum_{j=1}^m\ell^* (\wt X_j,\theta)- \frac{1}{n}\sum_{i=1}^n\ell^* (X_i,\theta^{\diamond})\Big)\\
        &=-\Big(\frac{1}{m}\sum_{j=1}^m \ell(\wt X_j,\theta)-\frac{1}{m}\sum_{j=1}^m \ell(\wt X_j,\theta^{\dagger})-\frac{1}{n}\sum_{i=1}^n \ell(X_i,\theta^{\diamond})+\frac{1}{n}\sum_{i=1}^n \ell(X_i,\theta^{\dagger})\Big)\\
        &\leq-\Big(\frac{1}{n}\sum_{i=1}^n \ell(X_i,\theta)-\frac{1}{n}\sum_{i=1}^n \ell(X_i,\theta^{\diamond})\Big)\\
        &\quad+\Big|\frac{1}{m}\sum_{j=1}^m \ell(\wt X_j,\theta)-\frac{1}{m}\sum_{j=1}^m \ell(\wt X_j,\theta^{\dagger})-\frac{1}{n}\sum_{i=1}^n \ell(X_i,\theta)+\frac{1}{n}\sum_{i=1}^n \ell(X_i,\theta^\dagger)\Big|\\
        &\leq-\Big(\frac{1}{n}\sum_{i=1}^n \ell(X_i,\theta)-\frac{1}{n}\sum_{i=1}^n \ell(X_i,\theta^{\diamond})\Big)
        +c_0\sqrt{\frac{\log m}{m}}\|\theta-\theta^\dagger\|_2+c_0\frac{\log m}{m}.
    \end{aligned}
\end{equation*}
Then, by \eqref{estimate 1.2} and noting that $
    \sqrt{\frac{\log m}{\alpha}}\gtrsim\sqrt{\frac{\log m }{m}}$ and $
       C_{m,n}\lesssim \sqrt{\frac{m}{\alpha^2\log m}}$, for sufficiently large $\delta_1$, we have
\begin{equation}\label{step 1.2.1}
\begin{aligned}
&\alpha^{\frac{d}{2}}\int_{\Theta_2} \pi(\theta)\exp\Big(-\alpha\Big(\frac{1}{m}\sum_{j=1}^m\ell^* (\wt X_j,\theta)- \frac{1}{n}\sum_{i=1}^n\ell^* (X_i,\theta^{\diamond})\Big)\Big)\,\dd\theta\\
&\lesssim \alpha^{\frac{d}{2}}\int_{\Theta_2} \pi(\theta)\exp\Big(-\alpha M_1\|\theta -\theta^{\diamond} \|^2_2 
     +\alpha M_1\frac{\log n}{n}\| \theta - \theta^{\diamond} \|_2^{2\beta_1}
     +\alpha M_1\Big( \frac{\log n}{n} \Big)^{\beta_1+1}\\
&\quad+\alpha c_0\sqrt{\frac{\log m}{m}}\|\theta-\theta^{\dagger}\|_2
     +\alpha c_0\frac{\log m}{m}\Big)\,\dd\theta\\
&\lesssim \alpha^{\frac{d}{2}}\int_{\Theta_2}\pi(\theta)\exp\Big(-\frac{\alpha M_1}{2}\|\theta-\theta^{\diamond}\|_2^2\Big)\,\dd\theta\\     
&\lesssim \frac{1}{m^2},
\end{aligned}
\end{equation}
and
\begin{equation}\label{step 1.2.2} 
\begin{aligned}
&\alpha^{\frac{d}{2}}\int_{\Theta_2} \pi(\theta)\exp\Big(-\alpha\Big( \frac{1}{n}\sum_{i=1}^n \ell^*(X_i,\theta)-\frac{1}{n}\sum_{i=1}^n \ell^*(X_i,\theta^{\diamond})\Big)\Big)\,\dd\theta\\
&= \alpha^{\frac{d}{2}}\int_{\Theta_2} \pi(\theta)\exp\Big(-\alpha\Big( \frac{1}{n}\sum_{i=1}^n \ell(X_i,\theta)-\frac{1}{n}\sum_{i=1}^n \ell(X_i,\theta^{\diamond})\Big)\Big)\,\dd\theta\\
&\leq \alpha^{\frac{d}{2}}\int_{\Theta_2}\pi(\theta)\exp\Big(-\alpha M_1\|\theta - \theta^{\diamond}\|^2_2 
     +\alpha M_1\frac{\log n}{n}\| \theta - \theta^{\diamond} \|_2^{2\beta_1}
     +\alpha M_1\Big( \frac{\log n}{n} \Big)^{\beta_1+1}\Big)\,\dd\theta\\&\lesssim \frac{1}{m^2}.
\end{aligned}
\end{equation}
\paragraph{3. \underline{On $\Theta_3$.}}
 For $\theta\in \Theta_3$, by \eqref{global hessian}, there exists a positive constant $K_1$ such that
 \begin{equation}\label{estimate 1.7}
\begin{aligned}
&\frac{1}{n}\sum_{i=1}^n \ell(X_i,\theta)-\frac{1}{n}\sum_{i=1}^n \ell(X_i,\theta^{\diamond})\geq K_1.
\end{aligned} 
\end{equation}
Then, the second statement of Lemma \ref{lemmaprobA1} gives
 \begin{equation}\label{estimate 1.8}
\begin{aligned}
&\frac{1}{m}\sum_{j=1}^m\ell^* (\wt X_j,\theta)- \frac{1}{n}\sum_{i=1}^n\ell^* (X_i,\theta^{\diamond})\\
&=\frac{1}{m}\sum_{j=1}^m \ell(\wt X_j,\theta)-\frac{1}{m}\sum_{j=1}^m \ell(\wt X_j,\theta^{\dagger})-\frac{1}{n}\sum_{i=1}^n \ell(X_i,\theta^{\diamond})+\frac{1}{n}\sum_{i=1}^n \ell(X_i,\theta^\dagger)\\
&\geq \frac{1}{n}\sum_{i=1}^n \ell( X_i,\theta)-\frac{1}{n}\sum_{i=1}^n \ell(X_i,\theta^{\diamond})
\\&\quad-\Big|\frac{1}{m}\sum_{j=1}^m \ell(\wt X_j,\theta)-\frac{1}{m}\sum_{j=1}^m \ell(\wt X_j,\theta^{\dagger})-\frac{1}{n}\sum_{i=1}^n \ell(X_i,\theta)+\frac{1}{n}\sum_{i=1}^n \ell(X_i,\theta^\dagger)\Big|\\
&\geq \frac{K_1}{2}.
\end{aligned}
\end{equation}
Therefore, by the compactness of $\Theta$, we can obtain
 \begin{equation}\label{step 1.3.1}
\begin{aligned}
& \alpha^{\frac{d}{2}}\int_{\Theta_3} \pi(\theta)\exp\Big(-\alpha\Big(\frac{1}{m}\sum_{j=1}^m\ell^* (\wt X_j,\theta)- \frac{1}{n}\sum_{i=1}^n\ell^* (X_i,\theta^{\diamond})\Big)\Big)\,\dd\theta \lesssim\alpha^{\frac{d}{2}}\exp\Big(-\frac{\alpha K_1}{2}\Big)\lesssim\frac{1}{m^2}
,\end{aligned}
\end{equation}
and
\begin{equation}\label{step 1.3.2}
\begin{aligned}
&\alpha^{\frac{d}{2}}\int_{\Theta_3} \pi(\theta)\exp\Big(-\alpha\Big( \frac{1}{n}\sum_{i=1}^n \ell^*(X_i,\theta)-\frac{1}{n}\sum_{i=1}^n \ell^*(X_i,\theta^{\diamond})\Big)\Big)\,\dd\theta\lesssim\alpha^{\frac{d}{2}}\exp(-\alpha K_1)\lesssim\frac{1}{m^2}.
\end{aligned}
\end{equation}
Combining the estimates \eqref{step 1.1}, \eqref{step 1.2.1}, \eqref{step 1.2.2}, \eqref{step 1.3.1}, and \eqref{step 1.3.2}, and using the lower bound in \eqref{lowerbound 1} for the denominator, we obtain
 \begin{equation*} 
\begin{aligned}
& \int\Big|\frac{\mb{E}_{\wX^{(m)} \sim  \mu_n^{\otimes m} }\Big[\pi(\theta) \exp\Big(-\alpha\Big(\frac{1}{m}\sum_{j=1}^m\ell^* (\wt X_j,\theta)- \frac{1}{n}\sum_{i=1}^n\ell^* (X_i,\theta^{\diamond})\Big)\Big)\cdot\bold{1}(\wX^{(m)}\in \m A)\Big]}{\mb{E}_{\wX^{(m)} \sim  \mu_n^{\otimes m} }\Big[\int\pi(\theta) \exp\Big(-\alpha\Big(\frac{1}{m}\sum_{j=1}^m\ell^* (\wt X_j,\theta)- \frac{1}{n}\sum_{i=1}^n\ell^* (X_i,\theta^{\diamond})\Big)\Big)\,\dd\theta\cdot\bold{1}(\wX^{(m)}\in \m A)\Big]}\\
&-\frac{\mb{E}_{\wX^{(m)} \sim  \mu_n^{\otimes m} }\Big[\pi(\theta)\exp\Big(-\alpha\Big( \frac{1}{n}\sum_{i=1}^n \ell(X_i,\theta)-\frac{1}{n}\sum_{i=1}^n \ell(X_i,\theta^{\diamond})\Big)\Big)\cdot\bold{1}(\wX^{(m)}\in \m A)\Big]}{\mb{E}_{\wX^{(m)} \sim  \mu_n^{\otimes m} }\Big[\int\pi(\theta)\exp\Big(-\alpha\Big( \frac{1}{n}\sum_{i=1}^n \ell(X_i,\theta)-\frac{1}{n}\sum_{i=1}^n \ell(X_i,\theta^{\diamond})\Big)\Big)\,\dd\theta\cdot\bold{1}(\wX^{(m)}\in \m A)\Big]}\Big|\,\dd\theta\\
&\lesssim \alpha^{\frac{d}{2}}\cdot \Big|\mb{E}_{\wX^{(m)} \sim  \mu_n^{\otimes m} }\Big[\int\pi(\theta) \exp\Big(-\alpha\Big(\frac{1}{m}\sum_{j=1}^m\ell^* (\wt X_j,\theta)- \frac{1}{n}\sum_{i=1}^n\ell^* (X_i,\theta^{\diamond})\Big)\Big)\\
&\quad- \pi(\theta)\exp\Big(-\alpha\Big( \frac{1}{n}\sum_{i=1}^n \ell(X_i,\theta)-\frac{1}{n}\sum_{i=1}^n \ell(X_i,\theta^{\diamond})\Big)\Big)\,\dd\theta \cdot\bold{1}(\wX^{(m)}\in \m A)\Big]\Big|\\
&\lesssim \alpha\frac{(\log m)^2}{m}+\alpha^2\frac{\log m}{m}C^2_{m,n}.
\end{aligned}
\end{equation*}
Thus,
\begin{equation}\label{total 1}
\begin{aligned}
& \int \left|\wt{\pi}_{\mathcal{L}_{\rm zcv},\alpha}(\theta)-\pi_{\alpha}(\theta)\right|\,\dd\theta \\
&= \int\Big|\frac{\mb{E}_{\wX^{(m)} \sim  \mu_n^{\otimes m} }\Big[\pi(\theta)\exp\Big(-\alpha\Big(\frac{1}{m}\sum_{j=1}^m\ell^* (\wt X_j,\theta)- \frac{1}{n}\sum_{i=1}^n\ell^* (X_i,\theta^{\diamond})\Big)\Big) \Big]}{\mb{E}_{\wX^{(m)} \sim  \mu_n^{\otimes m} }\Big[\int\pi(\theta)\exp\Big(-\alpha\Big(\frac{1}{m}\sum_{j=1}^m\ell^* (\wt X_j,\theta)- \frac{1}{n}\sum_{i=1}^n\ell^* (X_i,\theta^{\diamond})\Big)\Big)\,\dd\theta\Big]}\\
&-\frac{\mb{E}_{\wX^{(m)} \sim  \mu_n^{\otimes m} }\Big[\pi(\theta)\exp\Big(-\alpha\Big(\frac{1}{n}\sum_{i=1}^n \ell(X_i,\theta)-\frac{1}{n}\sum_{i=1}^n \ell( X_i,\theta^{\diamond})\Big)\Big)\Big]}{\mb{E}_{\wX^{(m)} \sim  \mu_n^{\otimes m} }\Big[\int \pi(\theta)\exp\Big(-\alpha\Big(\frac{1}{n}\sum_{i=1}^n \ell(X_i,\theta)-\frac{1}{n}\sum_{i=1}^n \ell( X_i,\theta^{\diamond})\Big)\Big)\,\dd\theta \Big]}\Big|\,\dd\theta\\
&\lesssim \alpha\frac{(\log m)^2}{m}+\alpha^2\frac{\log m}{m}C^2_{m,n}\\
&\quad+\alpha^{\frac{d}{2}}\cdot\underset{\theta\in \Theta}{\sup}\,\mb{E}_{\wX^{(m)} \sim  \mu_n^{\otimes m} }\Big[\exp\Big(-\alpha\Big(\frac{1}{n}\sum_{i=1}^n \ell(X_i,\theta)-\frac{1}{n}\sum_{i=1}^n \ell( X_i,\theta^{\diamond})\Big)\Big)\cdot\bold{1}(\wX^{(m)}\in \m A^c)\Big]\\
&\quad+\alpha^{\frac{d}{2}}\cdot\underset{\theta\in \Theta}{\sup}\,\mb{E}_{\wX^{(m)}\sim \mu_n^{\otimes m}}\Big[\exp\Big(-\alpha\Big(\frac{1}{m}\sum_{j=1}^m\ell^* (\wt X_j,\theta)- \frac{1}{n}\sum_{i=1}^n\ell^* (X_i,\theta^{\diamond})\Big)\Big)\cdot \bold{1}(\wX^{(m)}\in \m A^c)\Big].
\end{aligned}
\end{equation}
\paragraph{4. \underline{Control of the Expectation over $\mathcal{A}^c$.}}
For $\theta\in \Theta_1\cup\Theta_2$, by \eqref{estimate 1.2}, we have
\begin{equation*}
\begin{aligned}
&\exp\Big(-\alpha\Big(
\frac{1}{n}\sum_{i=1}^n \ell(X_i,\theta)
-\frac{1}{n}\sum_{i=1}^n \ell(X_i,\theta^{\diamond})
\Big)\Big)\\
&\leq
\exp\Big(
-\alpha M_1\|\theta-\theta^{\diamond}\|_2^2
+\alpha M_1\frac{\log n}{n}
\|\theta-\theta^{\diamond}\|_2^{2\beta_1}
+\alpha M_1\Big(\frac{\log n}{n}\Big)^{\beta_1+1}
\Big)\\
&\lesssim
\exp\Big(-\alpha M_1\|\theta-\theta^{\diamond}\|_2^2\Big).
\end{aligned}
\end{equation*}
For $\theta\in\Theta_3$, by \eqref{estimate 1.7}, we have
\begin{equation*}
\exp\Big(-\alpha\Big(
\frac{1}{n}\sum_{i=1}^n \ell(X_i,\theta)
-\frac{1}{n}\sum_{i=1}^n \ell(X_i,\theta^{\diamond})
\Big)\Big)
\leq \exp(-\alpha K_1).
\end{equation*}
Combining the two cases, for any $\theta\in\Theta$, we obtain
\begin{equation}\label{A complement 1.1}
\begin{aligned}
&\mb{E}_{\wX^{(m)}\sim\mu_n^{\otimes m}}
\Big[
\exp\Big(-\alpha\Big(
\frac{1}{n}\sum_{i=1}^n \ell(X_i,\theta)
-\frac{1}{n}\sum_{i=1}^n \ell(X_i,\theta^{\diamond})
\Big)\Big)
\cdot\bold{1}\big(\wX^{(m)}\in\m A^c\big)
\Big]\\
&\qquad\lesssim
\mb{P}_{\wX^{(m)}\sim\mu_n^{\otimes m}}
\big(\m A^c\big)
\leq \frac{1}{m^2\alpha^{d+1}}.
\end{aligned}
\end{equation}
 For the last term of \eqref{total 1},
\begin{equation*}
\begin{aligned}
& \mb{E}_{\wX^{(m)}\sim \mu_n^{\otimes m}}\Big[\exp\Big(-\alpha\cdot \Big(\frac{1}{m}\sum_{j=1}^m\ell^* (\wt X_j,\theta)- \frac{1}{n}\sum_{i=1}^n\ell^* (X_i,\theta^{\diamond})\Big)\Big)\cdot \bold{1}(\wX^{(m)}\in \m A^c)\Big]\\
&=\int_{0}^{+\infty}\mb{P}_{\wX^{(m)}\sim \mu_n^{\otimes m}}\Big(\exp\Big(-\alpha\cdot \Big(\frac{1}{m}\sum_{j=1}^m\ell^* (\wt X_j,\theta)- \frac{1}{n}\sum_{i=1}^n\ell^* (X_i,\theta^{\diamond})\Big)\Big)\cdot \bold{1}(\wX^{(m)}\in \m A^c)>t\Big)\,\dd t\\
&=\int_{0}^{m\alpha^{\frac{d+1}{2}}}\mb{P}_{\wX^{(m)}\sim \mu_n^{\otimes m}}\Big(\exp\Big(-\alpha\cdot \Big(\frac{1}{m}\sum_{j=1}^m\ell^* (\wt X_j,\theta)- \frac{1}{n}\sum_{i=1}^n\ell^* (X_i,\theta^{\diamond})\Big)\Big)\cdot \bold{1}(\wX^{(m)}\in \m A^c)>t\Big)\,\dd t\\
&\quad+\int_{m\alpha^{\frac{d+1}{2}}}^{+\infty}\mb{P}_{\wX^{(m)}\sim \mu_n^{\otimes m}}\Big(\exp\Big(-\alpha\cdot \Big(\frac{1}{m}\sum_{j=1}^m\ell^* (\wt X_j,\theta)- \frac{1}{n}\sum_{i=1}^n\ell^* (X_i,\theta^{\diamond})\Big)\Big)\cdot \bold{1}(\wX^{(m)}\in \m A^c)>t\Big)\,\dd t\\
&\leq\int_{0}^{m\alpha^{\frac{d+1}{2}}}\mb{P}_{\wX^{(m)}\sim \mu_n^{\otimes m}}\big(\bold{1}(\wX^{(m)}\in \m A^c)>0\big)\,\dd t\\
&\quad+\int_{m\alpha^{\frac{d+1}{2}}}^{+\infty}\mb{P}_{\wX^{(m)}\sim \mu_n^{\otimes m}}\Big(\exp\Big(-\alpha\cdot \Big(\frac{1}{m}\sum_{j=1}^m \ell^*(\wt X_j,\theta)-\frac{1}{n}\sum_{i=1}^n\ell^*(X_i,\theta^{\diamond})\Big)\Big)\cdot \bold{1}(\wX^{(m)}\in \m A^c)>t\Big)\,\dd t\\
&\leq \frac{1}{m\alpha^{\frac{d+1}{2}}}+\int_{m\alpha^{\frac{d+1}{2}}}^{+\infty}\mb{P}_{\wX^{(m)}\sim \mu_n^{\otimes m}}\Big(\exp\Big(-\alpha\cdot \Big(\frac{1}{m}\sum_{j=1}^m \ell^*(\wt X_j,\theta)-\frac{1}{n}\sum_{i=1}^n\ell^*(X_i,\theta^{\diamond})\Big)\Big)>t\Big)\,\dd t,\end{aligned}\end{equation*}
as $\ell(x,\cdot)$ is Lipschitz, there exists a positive constant $M_4$ such that 
   \begin{equation*}
    \begin{aligned}
       \Big|\frac{1}{m}\sum_{j=1}^m \ell^*(\wt X_j,\theta)-\frac{1}{n}\sum_{i=1}^n\ell^*(X_i,\theta^{\diamond})\Big|\leq M_4,   
       \end{aligned}
\end{equation*}
and with the help of
\begin{equation*}
\min\{M_1, M_3\}\Big(\|\theta - \theta^{\diamond} \|^2_2 
- \frac{\log n}{n}\| \theta - \theta^{\diamond} \|_2^{2\beta_1}
- \Big( \frac{\log n}{n} \Big)^{\beta_1+1}\Big)
\leq \frac{1}{n} \sum_{i=1}^n \ell(X_i, \theta) 
- \frac{1}{n} \sum_{i=1}^n \ell(X_i, \theta^{\diamond}),
\end{equation*}
we have
\begin{equation*}
\begin{aligned}
&\int_{m\alpha^{\frac{d+1}{2}}}^{+\infty}\mb{P}_{\wX^{(m)}\sim \mu_n^{\otimes m}}\Big(\exp\Big(-\alpha\Big(\frac{1}{m}\sum_{j=1}^m \ell^*(\wt X_j,\theta)-\frac{1}{n}\sum_{i=1}^n\ell^*(X_i,\theta^{\diamond})\Big)\Big)>t\Big)\,\dd t\\
&= \int_{m\alpha^{\frac{d+1}{2}}}^{\exp(\alpha M_4)}\mb{P}_{\wX^{(m)}\sim \mu_n^{\otimes m}}\Big(\exp\Big(-\alpha\Big(\frac{1}{m}\sum_{j=1}^m \ell^*(\wt X_j,\theta)-\frac{1}{n}\sum_{i=1}^n\ell^*(X_i,\theta^{\diamond})\Big)\Big)>t\Big)\,\dd t\\        
&=\int_{m\alpha^{\frac{d+1}{2}}}^{\exp(\alpha M_4)}\mb{P}_{\wX^{(m)}\sim \mu_n^{\otimes m}}\Big(\exp\Big(-\alpha\Big(\frac{1}{m}\sum_{j=1}^m \ell^*(\wt X_j,\theta)-\frac{1}{n}\sum_{i=1}^n\ell^*(X_i,\theta)\Big)\Big)\\
&\quad\cdot\exp\Big(-\alpha\Big(\frac{1}{n}\sum_{i=1}^n \ell^*( X_i,\theta)-\frac{1}{n}\sum_{i=1}^n\ell^*(X_i,\theta^{\diamond})\Big)\Big) >t\Big)\,\dd t \\
&\leq\int_{m\alpha^{\frac{d+1}{2}}}^{\exp(\alpha M_4)}\mb{P}_{\wX^{(m)}\sim \mu_n^{\otimes m}}\Big( \Big|\frac{1}{m}\sum_{j=1}^m \ell^*(\wt X_j,\theta)-\frac{1}{n}\sum_{i=1}^n\ell^*(X_i,\theta)\Big|\\
&\quad>\frac{\log t}{\alpha}+\min\{M_1, M_3\}\Big(\|\theta-\theta^{\diamond}\|^2_2-\frac{\log n}{n}\| \theta - \theta^{\diamond} \|_2^{2\beta_1}-\Big( \frac{\log n}{n} \Big)^{\beta_1+1}\Big)\Big)\,\dd t\\
&\leq\int_{m\alpha^{\frac{d+1}{2}}}^{\exp(\alpha M_4)}\mb{P}_{\wX^{(m)}\sim \mu_n^{\otimes m}}\Big( \Big|\frac{1}{m}\sum_{j=1}^m \ell^*(\wt X_j,\theta)-\frac{1}{n}\sum_{i=1}^n\ell^*(X_i,\theta)\Big|>\frac{\log t}{2\alpha}+\min\{M_1, M_3\}\|\theta-\theta^{\diamond}\|_2^2\Big)\,\dd t,
\end{aligned}
\end{equation*}
where the last inequality is due to the fact that when $t\geq m\alpha^{\frac{d+1}{2}}$, it holds that for $\theta\in\Theta $, 
\begin{equation*}
\frac{\log t}{\alpha}\gtrsim\frac{\log n}{n}\| \theta - \theta^{\diamond} \|_2^{2\beta_1}+\Big( \frac{\log n}{n} \Big)^{\beta_1+1}.
\end{equation*}
In the following $C$ is a  positive constant which may change from line to line.\\ 
Since 
\begin{equation*}
    \begin{aligned}
        \Big| \ell^*(\wt X_j,\theta)-\frac{1}{n}\sum_{i=1}^n\ell^*(X_i,\theta)\Big|\leq C\|\theta-\theta^\dagger\|_2,\end{aligned}
\end{equation*}
by Bernstein's inequality (for example Theorem 8.4 in \cite{vershynin2018high}), we have 
\begin{equation*}
\begin{aligned}
&\mb{P}_{\wX^{(m)}\sim \mu_n^{\otimes m}}\Big( \Big|\frac{1}{m}\sum_{j=1}^m \ell^*(\wt X_j,\theta)-\frac{1}{n}\sum_{i=1}^n\ell^*(X_i,\theta)\Big|>\frac{\log t}{2\alpha}+\min\{M_1, M_3\}\|\theta-\theta^{\diamond}\|_2^2\Big)\\
&\lesssim \max\Big(\exp\Big(-mC\frac{\big(\frac{\log t}{\alpha}+\|\theta-\theta^{\diamond}\|_2^2\big)^2}{\|\theta-\theta^\dagger\|^2_2}\Big),\exp\Big(-mC\frac{\frac{\log t}{\alpha}+\|\theta-\theta^{\diamond}\|^2_2}{\|\theta-\theta^\dagger\|_2}\Big)\Big).
\end{aligned}
\end{equation*}
\textbf{Case (i):}
If $\frac{\log t}{\alpha}\leq\|\theta-\theta^{\diamond}\|_2^2$, then
\begin{equation*}
\begin{aligned}
&\exp\Big(-mC\frac{\big(\frac{\log t}{\alpha}+\|\theta-\theta^{\diamond}\|_2^2\big)^2}{\|\theta-\theta^\dagger\|_2^2}\Big)\\
&\leq\exp\Big(-mC\frac{\big(\frac{\log t}{\alpha}+\|\theta-\theta^{\diamond}\|_2^2\big)^2}{\|\theta-\theta^{\diamond}\|_2^2+\|\theta^{\diamond}-\wh\theta\|^2_2+\|\wh\theta-\theta^{\dagger}\|^2_2}\Big)\\
&\leq\exp\Big(-mC\frac{\|\theta-\theta^{\diamond}\|_2^4}{\|\theta-\theta^{\diamond}\|_2^2+C_{m,n}^2+\|\theta^{\diamond}-\wh\theta\|^2_2}\Big)\\
&\leq\max\Big(\exp\big(-mC\|\theta-\theta^{\diamond}\|_2^2\big),\exp\Big(-mC\frac{\|\theta-\theta^{\diamond}\|_2^4}{C_{m,n}^2+\|\theta^{\diamond}-\wh\theta\|^2_2}\Big)\Big)\\
&\leq\max\Big(\exp\big(-mC\frac{\log t}{\alpha}\big),\exp\Big(-mC\frac{\big(\frac{\log t}{\alpha}\big)^2}{C_{m,n}^2+\|\theta^{\diamond}-\wh\theta\|^2_2}\Big)\Big).
\end{aligned}
\end{equation*}
Since $C_{m,n}\lesssim \sqrt{\frac{m}{\alpha^2\log m}}$, $\|\theta^{\diamond}-\wh\theta\|^2_2\lesssim\Big(\frac{\log n}{n}\Big)^{1+\beta_1}$, $\alpha\lesssim \frac{m}{(\log m)^2}$ and $t\geq m\alpha^{\frac{d+1}{2}}$, we have
\begin{equation*}
\max\Big(\exp\big(-mC\frac{\log t}{\alpha}\big),\exp\Big(-mC\frac{\big(\frac{\log t}{\alpha}\big)^2}{C_{m,n}^2+\|\theta^{\diamond}-\wh\theta\|^2_2}\Big)\Big)\lesssim\exp(-2\log t)=\frac{1}{t^2}.
\end{equation*}
Similarly,
\begin{equation*}
\begin{aligned}
&\exp\Big(-mC\frac{\frac{\log t}{\alpha}+\|\theta-\theta^{\diamond}\|_2^2}{\|\theta-\theta^\dagger\|_2}\Big)\\
&\leq \exp\Big(-mC\frac{\|\theta-\theta^{\diamond}\|_2^2}{\|\theta-\theta^{\diamond}\|_2+C_{m,n}+\|\theta^{\diamond}-\wh\theta\|_2}\Big)\\
&\leq\max\Big(\exp\big(-mC\|\theta-\theta^{\diamond}\|_2\big),\exp\Big(-mC\frac{\|\theta-\theta^\diamond\|_2^2}{C_{m,n}+\|\theta^{\diamond}-\wh\theta\|_2}\Big)\Big)\\
&\leq\max\Big(\exp\big(-mC\sqrt{\frac{\log t}{\alpha}}\big),\exp\Big(-mC\frac{\frac{\log t}{\alpha}}{C_{m,n}+\|\theta^{\diamond}-\wh\theta\|_2}\Big)\Big)\\
&\lesssim\max\Big(\exp\big(-2\sqrt{m\log m \log t}\big),\frac{1}{t^2}\Big).
\end{aligned}
\end{equation*}
\\\textbf{Case (ii):}
If $\frac{\log t}{\alpha}\geq\|\theta-\theta^{\diamond}\|_2^2$, then
\begin{equation*}
\begin{aligned}
&\exp\Big(-mC\frac{\big(\frac{\log t}{\alpha}+\|\theta-\theta^{\diamond}\|_2^2\big)^2}{\|\theta-\theta^\dagger\|_2^2}\Big)\\
&\leq\exp\Big(-mC\frac{\big(\frac{\log t}{\alpha}\big)^2}{\frac{\log t}{\alpha}+C_{m,n}^2+\|\theta^{\diamond}-\wh\theta\|^2_2}\Big)\\
&\leq\max\Big(\exp\big(-mC\frac{\log t}{\alpha}\big),\exp\Big(-mC\frac{\big(\frac{\log t}{\alpha}\big)^2}{C_{m,n}^2+\|\theta^{\diamond}-\wh\theta\|^2_2}\Big)\Big)\\
&\lesssim\frac{1}{t^2}.
\end{aligned}
\end{equation*}
Similarly,
\begin{equation*}
\begin{aligned}
&\exp\Big(-mC\frac{\frac{\log t}{\alpha}+\|\theta-\theta^{\diamond}\|_2^2}{\|\theta-\theta^\dagger\|_2}\Big)\\
&\leq \exp\Big(-mC\frac{\frac{\log t}{\alpha}}{\sqrt{\frac{\log t}{\alpha}}+C_{m,n}+\|\theta^{\diamond}-\wh\theta\|_2}\Big)\\
&\leq\max\Big(\exp\big(-mC\sqrt{\frac{\log t}{\alpha}}\big),\exp\Big(-mC\frac{\frac{\log t}{\alpha}}{C_{m,n}+\|\theta^{\diamond}-\wh\theta\|_2}\Big)\Big)\\
&\lesssim\max\Big(\exp\big(-2\sqrt{m\log m\log t}\big),\frac{1}{t^2}\Big).
\end{aligned}
\end{equation*}
Since
$\sqrt{m\log m\log t}\geq \log t$
for all \(t\in(1,\exp(\alpha M_4))\) when \(m\) is sufficiently large, we have
\begin{equation}\label{estimate 1.9}
\begin{aligned}
&\int_{m\alpha^{\frac{d+1}{2}}}^{+\infty}\mb{P}_{\wX^{(m)}\sim \mu_n^{\otimes m}}\Big(\exp\Big(-\alpha\Big(\frac{1}{m}\sum_{j=1}^m \ell^*(\wt X_j,\theta)-\frac{1}{n}\sum_{i=1}^n\ell^*(X_i,\theta^{\diamond})\Big)\Big)>t\Big)\,\dd t\\
&=\int_{m\alpha^{\frac{d+1}{2}}}^{\exp(\alpha M_4)}\mb{P}_{\wX^{(m)}\sim \mu_n^{\otimes m}}\Big(\exp\Big(-\alpha\Big(\frac{1}{m}\sum_{j=1}^m \ell^*(\wt X_j,\theta)-\frac{1}{n}\sum_{i=1}^n\ell^*(X_i,\theta^{\diamond})\Big)\Big)>t\Big)\,\dd t\\
&\lesssim\int_{m\alpha^{\frac{d+1}{2}}}^{\exp(\alpha M_4)}\max\Big(\exp\big(-2\sqrt{m\log m \log t}\big),\frac{1}{t^2}\Big)\,\dd t\\
&\lesssim\int_{m\alpha^{\frac{d+1}{2}}}^{\exp(\alpha M_4)}\frac{1}{t^2}\,\dd t \\
&\leq \frac{1}{m\alpha^{\frac{d+1}{2}}}.
\end{aligned}
\end{equation}
Hence,
\begin{equation}\label{A complement 1.2}
\alpha^{\frac{d}{2}} \mb{E}_{\wX^{(m)}\sim \mu_n^{\otimes m}}\Big[\exp\Big(-\alpha\Big(\frac{1}{m}\sum_{j=1}^m\ell^* (\wt X_j,\theta)- \frac{1}{n}\sum_{i=1}^n\ell^* (X_i,\theta^{\diamond})\Big)\Big)\cdot \bold{1}(\wX^{(m)}\in \m A^c)\Big]\lesssim \frac{1}{m\sqrt{\alpha}}.
\end{equation}
Combining all pieces, we can obtain 
\begin{equation}\label{TV bound 1.1}
    \int \big| \wt{\pi}_{\mathcal{L}_{\rm zcv},\alpha}(\theta)-\pi_\alpha(\theta)\big|\,\dd\theta \lesssim \alpha\frac{(\log m)^2}{m}+\alpha^2\frac{\log m}{m}C^2_{m,n}.
\end{equation}

\subsubsection{Proof of the Mean Deviation}
\paragraph{1. \underline{Difference Between the Means of 
$\widetilde{\pi}_{\mathcal{L}_{\mathrm{zcv}},\alpha}$ and $\pi_\alpha$}}

Proceeding as in \eqref{step 1.1}, we obtain
\begin{equation}\label{mean step 1.1}
    \begin{aligned}
       & \alpha^{\frac{d}{2}}\int_{\Theta_1}\Big\|\mb{E}_{\wX^{(m)}\sim \mu_n^{\otimes m}}\Big[\sqrt{\alpha}\cdot(\theta-\theta^{\diamond})\pi(\theta)\exp\Big(-\alpha\cdot\Big(\frac{1}{m}\sum_{j=1}^m \ell^*(\wt X_j,\theta)-\frac{1}{n}\sum_{i=1}^n \ell^*(X_i,\theta^{\diamond})\Big)\Big)\\&\quad\cdot\bold{1}(\wX^{(m)}\in \m A)\Big] \\
         &\quad-\mb{E}_{\wX^{(m)}\sim \mu_n^{\otimes m}}\Big[\sqrt{\alpha}\cdot(\theta-\theta^{\diamond})\pi(\theta)\exp\Big(-\alpha\cdot\Big( \frac{1}{n}\sum_{i=1}^n \ell^*(X_i,\theta)-\frac{1}{n}\sum_{i=1}^n \ell^*(X_i,\theta^{\diamond})\Big)\Big)\\&\quad\cdot\bold{1}(\wX^{(m)}\in \m A)\Big]\Big\|_2\,\dd\theta\\&\lesssim \alpha^{\frac{d+1}{2}}\int_{\Theta_1}
         \pi(\theta)\|\theta - \theta^{\diamond} \|_2         \exp\Big(-\alpha M_1\|\theta - \theta^{\diamond} \|^2_2\Big)
           \cdot \Big(\alpha\frac{(\log m)^2}{m}+\alpha^2\frac{\log m}{m}C^2_{m,n}\Big)\,\dd\theta \\         &\lesssim \alpha\frac{(\log m)^2}{m}+ \alpha^2\frac{\log m}{m}C^2_{m,n}.
    \end{aligned}
\end{equation}
By arguments analogous to those used in \eqref{step 1.2.1} and \eqref{step 1.2.2}, for every $\wX^{(m)}\in \m A$, we have
\begin{equation}\label{mean step 1.2.1}
\begin{aligned}
&\alpha^{\frac{d}{2}}\Big\|\int_{\Theta_2}\sqrt{\alpha}\cdot\big(\theta-\theta^{\diamond}\big) \pi(\theta)\exp\Big(-\alpha\cdot\Big(\frac{1}{m}\sum_{j=1}^m \ell^*(\wt X_j,\theta)-\frac{1}{n}\sum_{i=1}^n \ell^*(X_i,\theta^{\diamond})\Big)\Big)\,\dd\theta\Big\|_2\lesssim\frac{1}{m^2},
\end{aligned}
\end{equation}
and
\begin{equation}\label{mean step 1.2.2} 
\begin{aligned}       
&\alpha^{\frac{d}{2}} \Big\|\int_{\Theta_2}\sqrt{\alpha}\cdot\big(\theta-\theta^{\diamond}\big) \pi(\theta)\exp\Big(-\alpha\cdot\Big( \frac{1}{n}\sum_{i=1}^n \ell^*(X_i,\theta)-\frac{1}{n}\sum_{i=1}^n \ell^*(X_i,\theta^{\diamond})\Big)\Big)\,\dd\theta\Big\|_2\lesssim\frac{1}{m^2}.
\end{aligned}
\end{equation}
Similarly, following the arguments in \eqref{step 1.3.1} and \eqref{step 1.3.2}, for every $\wX^{(m)}\in \m A$, we have
\begin{equation}\label{mean step 1.3.1} 
\begin{aligned}
&\alpha^{\frac{d}{2}}\Big\|\int_{\Theta_3}\sqrt{\alpha}\cdot\big(\theta-\theta^{\diamond}\big) \pi(\theta)\exp\Big(-\alpha\cdot\Big(\frac{1}{m}\sum_{j=1}^m \ell^*(\wt X_j,\theta)-\frac{1}{n}\sum_{i=1}^n \ell^*(X_i,\theta^{\diamond})\Big)\Big)\,\dd\theta\Big\|_2\lesssim \frac{1}{m^2},
\end{aligned}
\end{equation}
and
\begin{equation}\label{mean step 1.3.2} 
\begin{aligned}       
&\alpha^{\frac{d}{2}} \Big\|\int_{\Theta_3}\sqrt{\alpha}\cdot\big(\theta-\theta^{\diamond}\big) \pi(\theta)\exp\Big(-\alpha\cdot\Big( \frac{1}{n}\sum_{i=1}^n \ell^*(X_i,\theta)-\frac{1}{n}\sum_{i=1}^n \ell^*(X_i,\theta^{\diamond})\Big)\Big)\,\dd\theta\Big\|_2\lesssim\frac{1}{m^2}.
\end{aligned}
\end{equation}
Combining \eqref{mean step 1.1}, \eqref{mean step 1.2.1}, \eqref{mean step 1.2.2}, \eqref{mean step 1.3.1}, and \eqref{mean step 1.3.2} with \eqref{A complement 1.1} and \eqref{A complement 1.2}, we obtain
\begin{equation}\label{omit 1.1}
    \begin{aligned}
&\Big\|\alpha^{\frac{d}{2}}\mb{E}_{\wX^{(m)} \sim  \mu_n^{\otimes m} }\Big[\int(\theta-\theta^{\diamond})\pi(\theta) \exp\Big(-\alpha\Big(\frac{1}{m}\sum_{j=1}^m \ell^*(\wt X_j,\theta)-\frac{1}{n}\sum_{i=1}^n \ell^*(X_i,\theta^{\diamond})\Big)\Big)\,\dd\theta \Big]\\&-\alpha^{\frac{d}{2}}\int(\theta-\theta^{\diamond})\pi(\theta)\exp\Big(-\alpha\Big(\frac{1}{n}\sum_{i=1}^n \ell( X_i,\theta)-\frac{1}{n}\sum_{i=1}^n \ell( X_i,\theta^{\diamond})\Big)\Big)\,\dd\theta\Big\|_2
\\&\lesssim \alpha^{\frac{d}{2}}\int_{\Theta}\Big\|\mb{E}_{\wX^{(m)}\sim \mu_n^{\otimes m}}\Big[(\theta-\theta^{\diamond})\pi(\theta)\exp\Big(-\alpha\Big(\frac{1}{m}\sum_{j=1}^m \ell^*(\wt X_j,\theta)-\frac{1}{n}\sum_{i=1}^n \ell^*(X_i,\theta^{\diamond})\Big)\Big)\cdot\bold{1}(\wX^{(m)}\in \m A)\Big]\\
&\qquad-\mb{E}_{\wX^{(m)}\sim \mu_n^{\otimes m}}\Big[(\theta-\theta^{\diamond})\pi(\theta)\exp\Big(-\alpha\Big( \frac{1}{n}\sum_{i=1}^n \ell(X_i,\theta)-\frac{1}{n}\sum_{i=1}^n \ell(X_i,\theta^{\diamond})\Big)\Big)\cdot\bold{1}(\wX^{(m)}\in \m A)\Big]\Big\|_2\,\dd\theta\\
&\quad+\alpha^{\frac{d}{2}}\operatorname{diam}(\Theta)\underset{\theta\in \Theta}{\sup}\,\mb{E}_{\wX^{(m)} \sim  \mu_n^{\otimes m} }\Big[\exp\Big(-\alpha\Big(\frac{1}{n}\sum_{i=1}^n \ell(X_i,\theta)-\frac{1}{n}\sum_{i=1}^n \ell( X_i,\theta^{\diamond})\Big)\Big)\cdot\bold{1}(\wX^{(m)}\in \m A^c)\Big]\\
&\quad+\alpha^{\frac{d}{2}}\operatorname{diam}(\Theta)\underset{\theta\in \Theta}{\sup}\,\mb{E}_{\wX^{(m)}\sim \mu_n^{\otimes m}}\Big[\exp\Big(-\alpha\Big(\frac{1}{m}\sum_{j=1}^m\ell^* (\wt X_j,\theta)- \frac{1}{n}\sum_{i=1}^n\ell^* (X_i,\theta^{\diamond})\Big)\Big)\cdot \bold{1}(\wX^{(m)}\in \m A^c)\Big]\\         
&\lesssim\sqrt{\alpha}\frac{(\log m)^2}{m}+ \alpha^{\frac{3}{2}}\frac{\log m}{m}C^2_{m,n}. 
\end{aligned}
\end{equation}
Using the denominator lower bound established in \eqref{lowerbound 1}, we further obtain
 \begin{equation}\label{omit 1.2}
 {\small
\begin{aligned}
&\Big\|\int \theta\,\wt{\pi}_{\mathcal{L}_{\rm zcv},\alpha}(\theta)\,\dd\theta-\int \theta\,\pi_{\alpha}(\theta)\,\dd\theta\Big\|_2\\
&= \Big\|\int (\theta-\theta^{\diamond}) \wt{\pi}_{\mathcal{L}_{\rm zcv},\alpha}(\theta)\,\dd\theta-\int (\theta-\theta^{\diamond})\pi_{\alpha}(\theta)\,\dd\theta\Big\|_2\\
&=\Big\|\int \frac{\alpha^{\frac{d}{2}}\mb{E}_{\wX^{(m)} \sim  \mu_n^{\otimes m} }\Big[(\theta-\theta^{\diamond})\pi(\theta) \exp\Big(-\alpha\Big(\frac{1}{m}\sum_{j=1}^m \ell^*(\wt X_j,\theta)-\frac{1}{n}\sum_{i=1}^n \ell^*(X_i,\theta^{\diamond})\Big)\Big) \Big]}{\alpha^{\frac{d}{2}}\mb{E}_{\wX^{(m)} \sim  \mu_n^{\otimes m} }\Big[\int\pi(\theta) \exp\Big(-\alpha\Big(\frac{1}{m}\sum_{j=1}^m \ell^*(\wt X_j,\theta)-\frac{1}{n}\sum_{i=1}^n \ell^*(X_i,\theta^{\diamond})\Big)\Big)\,\dd\theta\Big]}\,\dd\theta\\
&\quad-\int \frac{\alpha^{\frac{d}{2}}(\theta-\theta^{\diamond})\pi(\theta)\exp\Big(-\alpha\Big(\frac{1}{n}\sum_{i=1}^n \ell( X_i,\theta)-\frac{1}{n}\sum_{i=1}^n \ell( X_i,\theta^{\diamond})\Big)\Big)}{\alpha^{\frac{d}{2}}\int\pi(\theta)\exp\Big(-\alpha\Big(\frac{1}{n}\sum_{i=1}^n \ell(X_i,\theta)-\frac{1}{n}\sum_{i=1}^n \ell( X_i,\theta^{\diamond})\Big)\Big)\,\dd\theta}\,\dd\theta\Big\|_2\\&\lesssim\Big\|\alpha^{\frac{d}{2}}\mb{E}_{\wX^{(m)} \sim  \mu_n^{\otimes m} }\Big[\int(\theta-\theta^{\diamond})\pi(\theta) \exp\Big(-\alpha\Big(\frac{1}{m}\sum_{j=1}^m \ell^*(\wt X_j,\theta)-\frac{1}{n}\sum_{i=1}^n \ell^*(X_i,\theta^{\diamond})\Big)\Big)\,\dd\theta \Big]\\&-\alpha^{\frac{d}{2}}\int(\theta-\theta^{\diamond})\pi(\theta)\exp\Big(-\alpha\Big(\frac{1}{n}\sum_{i=1}^n \ell( X_i,\theta)-\frac{1}{n}\sum_{i=1}^n \ell( X_i,\theta^{\diamond})\Big)\Big)\,\dd\theta\Big\|_2\\&+\Big\|\alpha^{\frac{d}{2}}\int(\theta-\theta^{\diamond})\pi(\theta)\exp\Big(-\alpha\Big(\frac{1}{n}\sum_{i=1}^n \ell( X_i,\theta)-\frac{1}{n}\sum_{i=1}^n \ell( X_i,\theta^{\diamond})\Big)\Big)\,\dd\theta
\Big\|_2\\&\cdot\Big|\alpha^{\frac{d}{2}}\mb{E}_{\wX^{(m)} \sim  \mu_n^{\otimes m} }\Big[\int\pi(\theta) \exp\Big(-\alpha\Big(\frac{1}{m}\sum_{j=1}^m \ell^*(\wt X_j,\theta)-\frac{1}{n}\sum_{i=1}^n \ell^*(X_i,\theta^{\diamond})\Big)\Big)\,\dd\theta \Big]\\&-\alpha^{\frac{d}{2}}\int\pi(\theta)\exp\Big(-\alpha\Big(\frac{1}{n}\sum_{i=1}^n \ell( X_i,\theta)-\frac{1}{n}\sum_{i=1}^n \ell( X_i,\theta^{\diamond})\Big)\Big)\,\dd\theta\Big|.
\end{aligned}}
\end{equation}
Combining \eqref{step 1.1}, \eqref{step 1.2.1}, \eqref{step 1.3.1}, and \eqref{step 1.3.2} with \eqref{A complement 1.1} and \eqref{A complement 1.2}, we obtain
\begin{equation}\label{omit 1.3}
    \begin{aligned}
        &\Big|\alpha^{\frac{d}{2}}\mb{E}_{\wX^{(m)} \sim  \mu_n^{\otimes m} }\Big[\int\pi(\theta) \exp\Big(-\alpha\Big(\frac{1}{m}\sum_{j=1}^m \ell^*(\wt X_j,\theta)-\frac{1}{n}\sum_{i=1}^n \ell^*(X_i,\theta^{\diamond})\Big)\Big)\,\dd\theta \Big]\\&-\alpha^{\frac{d}{2}}\int\pi(\theta)\exp\Big(-\alpha\Big(\frac{1}{n}\sum_{i=1}^n \ell( X_i,\theta)-\frac{1}{n}\sum_{i=1}^n \ell( X_i,\theta^{\diamond})\Big)\Big)\,\dd\theta\Big|\lesssim\alpha\frac{(\log m)^2}{m}+\alpha^2\frac{\log m}{m}C^2_{m,n},
    \end{aligned}
\end{equation}
We will establish in \eqref{omit 1.4} below that
\[\Big\|\alpha^{\frac{d}{2}}\int(\theta-\theta^{\diamond})\pi(\theta)\exp\Big(-\alpha\Big(\frac{1}{n}\sum_{i=1}^n \ell( X_i,\theta)-\frac{1}{n}\sum_{i=1}^n \ell( X_i,\theta^{\diamond})\Big)\Big)\,\dd\theta
\Big\|_2\lesssim\frac{(\log n)^\frac{3}{2}}{\alpha}+\frac{(\log n)^{\frac{\beta_1+2}{2}}}{\sqrt{n}\alpha^{\frac{\beta_1}{2}}}.\]
It then follows from this bound, together with \eqref{omit 1.1} and \eqref{omit 1.3}, that
\begin{equation}\label{mean diff 1.1}
    \begin{aligned}
       \Big\|\int \theta\,\wt{\pi}_{\mathcal{L}_{\rm zcv},\alpha}(\theta)\,\dd\theta-\int \theta\,\pi_{\alpha}(\theta)\,\dd\theta\Big\|_2 \lesssim\sqrt{\alpha}\frac{(\log m)^2}{m}+ \alpha^{\frac{3}{2}}\frac{\log m}{m}C^2_{m,n}
    \end{aligned}
\end{equation}
\paragraph{2. \underline{Difference Between the Mean Under $\pi_\alpha$ and $\theta^{\diamond}(\bX^{(n)})$}}\mbox{}\\
We introduce the following two regions:
\begin{equation*}
\begin{aligned}
&\Theta'_1=\left\{\theta\in\mb R^d : \|\theta-\theta^{\diamond}\|_2\leq \delta'_1\sqrt{\frac{\log n}{\alpha}}\right\},\\
&\Theta'_2=\left\{\theta\in \mb R^d : \delta'_1\sqrt{\frac{\log n}{\alpha}}<\|\theta-\theta^{\diamond}\|_2\right\}.
\end{aligned}
\end{equation*}
Here, $\delta'_1 > 0$ is a suitably chosen constant. Since $\theta^*$ lies in the interior of $\Theta$ and $\|\theta^\diamond - \theta^*\|_2 \lesssim \sqrt{\frac{\log n}{n}}$ on $\mathcal{B}$, the assumed lower bound on $\alpha$ ensures that $\Theta_1' \subset \Theta$ for sufficiently large $n$. In the following, $\pi(\theta)$ is extended by zero outside $\Theta$.

We first consider the contribution from $\Theta'_1$. Analogously to \eqref{Hessian Taylor}, by Taylor's theorem with integral remainder and the local Lipschitz continuity of $\m{H}$ around $\theta^*$, we have
\begin{equation*}
\begin{aligned}
&\Big|\m R(\theta) - \m R(\theta^{\diamond}) - \nabla \m R(\theta^{\diamond})\cdot (\theta - \theta^{\diamond}) -\frac{1}{2}(\theta-\theta^{\diamond})^T\m H(\theta^{\diamond})(\theta-\theta^{\diamond})\Big|\\
&=\Big|\int_0^1 (1-t) \, (\theta-\theta^{\diamond})^T \m H\big(\theta^{\diamond} + t(\theta-\theta^{\diamond})\big) (\theta-\theta^{\diamond}) \, \dd t
-\frac{1}{2}(\theta-\theta^{\diamond})^T\m H(\theta^{\diamond})(\theta-\theta^{\diamond})\Big|\\
&\lesssim
\left(
\|\theta-\theta^\diamond\|_2
+\|\theta^\diamond-\theta^*\|_2
\right)
\|\theta-\theta^\diamond\|_2^2\\
&\lesssim
\sqrt{\frac{\log n}{\alpha}}
\|\theta-\theta^\diamond\|_2^2,
\end{aligned}
\end{equation*}
where the last inequality follows from the bounds $\|\theta^\diamond - \theta^*\|_2 \lesssim \sqrt{\frac{\log n}{n}}$ and $\alpha \le n$.
Moreover, the final assertion of Lemma \ref{lemmaprobB1}, together with \eqref{estimate 0.3}, yields
\begin{equation*}
\begin{aligned}
&\Big| \frac{1}{n} \sum_{i=1}^n \ell(X_i, \theta) - \frac{1}{n} \sum_{i=1}^n \ell(X_i, \theta^{\diamond}) - \m R(\theta) + \m R(\theta^{\diamond})
+ \nabla\m R(\theta^{\diamond})\cdot(\theta - \theta^{\diamond})\Big|\\
&\lesssim \sqrt{\frac{\log n}{n}} \| \theta - \theta^{\diamond} \|_2^{\beta_1+1} 
+ \Big( \frac{\log n}{n} \Big)^{\frac{\beta_1+1}{2}}\|\theta - \theta^{\diamond} \|_2
+ \Big( \frac{\log n}{n} \Big)^2.
\end{aligned}
\end{equation*}
Since $\alpha\lesssim n$, the preceding estimates imply
\begin{equation*}
\begin{aligned}
&\Big| \frac{1}{n} \sum_{i=1}^n \ell(X_i, \theta) - \frac{1}{n} \sum_{i=1}^n \ell(X_i, \theta^{\diamond})
- \frac{1}{2}(\theta-\theta^{\diamond})^T\m H(\theta^{\diamond})(\theta-\theta^{\diamond})\Big|\\
&\lesssim \Big(\frac{\log n}{\alpha}\Big)^{\frac{3}{2}}
+ \sqrt{\frac{\log n}{n}}\Big(\frac{\log n}{\alpha}\Big)^{\frac{\beta_1+1}{2}}
+ \Big( \frac{\log n}{n} \Big)^{\frac{\beta_1+1}{2}}\sqrt{\frac{\log n}{\alpha}}
+\Big( \frac{\log n}{n} \Big)^2\\
&\lesssim \Big(\frac{\log n}{\alpha}\Big)^{\frac{3}{2}}
+ \sqrt{\frac{\log n}{n}}\Big(\frac{\log n}{\alpha}\Big)^{\frac{\beta_1+1}{2}}.
\end{aligned}
\end{equation*}
By the inequality $\left|e^{-a}-e^{-b}\right|
\leq |a-b|\left(e^{-a}+e^{-b}\right)$,
we have
\begin{equation}\label{estimate 1.10}
\begin{aligned}
&\Big|\exp\Big(-\alpha\Big(\frac{1}{n}\sum_{i=1}^n \ell(X_i,\theta)-\frac{1}{n}\sum_{i=1}^n \ell( X_i,\theta^{\diamond})\Big)\Big)
-\exp\Big(-\frac{\alpha}{2}(\theta-\theta^{\diamond})^T\m H(\theta^{\diamond})(\theta-\theta^{\diamond})\Big)\Big|\\
&\lesssim\Bigg[\exp\Big(-\alpha\Big(\frac{1}{n}\sum_{i=1}^n \ell(X_i,\theta)-\frac{1}{n}\sum_{i=1}^n \ell( X_i,\theta^{\diamond})\Big)\Big)+\exp\Big(-\frac{\alpha}{2}(\theta-\theta^{\diamond})^T\m H(\theta^{\diamond})(\theta-\theta^{\diamond})\Big)\Bigg]
\\&\quad\cdot\Big(\frac{(\log n)^\frac{3}{2}}{\sqrt{\alpha}}
+\frac{(\log n)^{\frac{\beta_1+2}{2}}}{\sqrt{n}\alpha^{\frac{\beta_1-1}{2}}}\Big).
\end{aligned}
\end{equation}
By Assumption A.2, the local Lipschitz continuity of $\pi$ around $\theta^*$ implies 
\begin{equation}\label{pi lip}
\left|\pi(\theta)-\pi(\theta^\diamond)\right|
\le L\|\theta-\theta^\diamond\|_2
+2L\|\theta^\diamond-\theta^*\|_2
\lesssim \sqrt{\frac{\log n}{\alpha}}
\qquad \text{for all }\theta\in\Theta_1'.
\end{equation}
As $\alpha\lesssim m$, combining \eqref{estimate 1.2} with \eqref{estimate 1.10}, we obtain
\begin{equation*}
\begin{aligned}
& \alpha^{\frac{d}{2}}\int_{\Theta'_1}\Big\|\sqrt{\alpha}\cdot(\theta-\theta^{\diamond})\pi(\theta)
\exp\Big(-\alpha\Big(\frac{1}{n}\sum_{i=1}^n \ell(X_i,\theta)-\frac{1}{n}\sum_{i=1}^n \ell( X_i,\theta^{\diamond})\Big)\Big)\\
&\quad-\sqrt{\alpha}\cdot(\theta-\theta^{\diamond})\pi(\theta^{\diamond})
\exp\Big(-\frac{\alpha}{2}(\theta-\theta^{\diamond})^T\m H(\theta^{\diamond})(\theta-\theta^{\diamond})\Big)\Big\|_2\,\dd\theta\\
&\leq\alpha^{\frac{d}{2}}\int_{\Theta'_1}\Big\|\sqrt{\alpha}\cdot(\theta-\theta^{\diamond})\pi(\theta)
\exp\Big(-\alpha\Big(\frac{1}{n}\sum_{i=1}^n \ell(X_i,\theta)-\frac{1}{n}\sum_{i=1}^n \ell( X_i,\theta^{\diamond})\Big)\Big)\\
&\quad-\sqrt{\alpha}\cdot(\theta-\theta^{\diamond})\pi(\theta)
\exp\Big(-\frac{\alpha}{2}(\theta-\theta^{\diamond})^T\m H(\theta^{\diamond})(\theta-\theta^{\diamond})\Big)\Big\|_2\,\dd\theta\\ 
&\quad+\alpha^{\frac{d}{2}}\int_{\Theta'_1}\sqrt{\alpha}\|\theta-\theta^{\diamond}\|_2
\cdot|\pi(\theta)-\pi(\theta^{\diamond})|
\exp\Big(-\frac{\alpha}{2}(\theta-\theta^{\diamond})^T\m H(\theta^{\diamond})(\theta-\theta^{\diamond})\Big)\,\dd\theta\\ 
&\lesssim \alpha^{\frac{d}{2}}\int_{\Theta'_1}\sqrt{\alpha}\|\theta-\theta^{\diamond}\|_2 \pi(\theta)
\exp\Big(-\alpha\Big( \frac{1}{n}\sum_{i=1}^n \ell(X_i,\theta)-\frac{1}{n}\sum_{i=1}^n \ell(X_i,\theta^{\diamond})\Big)\Big) \cdot\Big(\frac{(\log n)^\frac{3}{2}}{\sqrt{\alpha}}
+\frac{(\log n)^{\frac{\beta_1+2}{2}}}{\sqrt{n}\alpha^{\frac{\beta_1-1}{2}}}\Big)\,\dd\theta
\\&\quad+\alpha^{\frac{d}{2}}\int_{\Theta'_1}\sqrt{\alpha}\|\theta-\theta^{\diamond}\|_2 \pi(\theta)
\exp\Big(-\frac{\alpha}{2}(\theta-\theta^{\diamond})^T\m H(\theta^{\diamond})(\theta-\theta^{\diamond})\Big) \cdot\Big(\frac{(\log n)^\frac{3}{2}}{\sqrt{\alpha}}
+\frac{(\log n)^{\frac{\beta_1+2}{2}}}{\sqrt{n}\alpha^{\frac{\beta_1-1}{2}}}\Big)\,\dd\theta
\\&\quad+\frac{\log n}{\sqrt{\alpha}}\\
&\lesssim\frac{(\log n)^\frac{3}{2}}{\sqrt{\alpha}}
+\frac{(\log n)^{\frac{\beta_1+2}{2}}}{\sqrt{n}\alpha^{\frac{\beta_1-1}{2}}}.
\end{aligned}
\end{equation*}

We next consider the contribution from $\Theta'_2$. Proceeding as in \eqref{mean step 1.2.2} and \eqref{mean step 1.3.2}, and choosing $\delta'_1$ sufficiently large, we have
\begin{equation*}
\begin{aligned}
\alpha^{\frac{d}{2}} \int_{\Theta'_2}\sqrt{\alpha}\|\theta-\theta^{\diamond}\|_2\pi(\theta)
\exp\Big(-\alpha\Big( \frac{1}{n}\sum_{i=1}^n \ell(X_i,\theta)
-\frac{1}{n}\sum_{i=1}^n \ell(X_i,\theta^{\diamond})\Big)\Big)\,\dd\theta
\lesssim \frac{1}{n^2}.
\end{aligned}
\end{equation*}
Likewise, the corresponding Gaussian tail satisfies
\begin{equation*}
\begin{aligned}
\alpha^{\frac{d}{2}}\int_{\Theta'_2}\sqrt{\alpha}\|\theta-\theta^{\diamond}\|_2\pi(\theta^{\diamond})
\exp\Big(-\frac{\alpha}{2}(\theta-\theta^{\diamond})^T\m H(\theta^{\diamond})(\theta-\theta^{\diamond})\Big)\,\dd\theta
\lesssim \frac{1}{n^2}.
\end{aligned}
\end{equation*}
Taken together, the estimates over $\Theta'_1$ and $\Theta'_2$ yield
\begin{equation}
\begin{aligned}\label{omit 1.5}
&\alpha^{\frac{d}{2}}\Big\|\int_{\mathbb{R}^d}(\theta-\theta^{\diamond})\pi(\theta)
\exp\Big(-\alpha\Big(\frac{1}{n}\sum_{i=1}^n \ell(X_i,\theta)
-\frac{1}{n}\sum_{i=1}^n \ell( X_i,\theta^{\diamond})\Big)\Big)\\
&\quad-(\theta-\theta^{\diamond})\pi(\theta^{\diamond})
\exp\Big(-\frac{\alpha}{2}(\theta-\theta^{\diamond})^T
\m H(\theta^{\diamond})(\theta-\theta^{\diamond})\Big)\,\dd\theta\Big\|_2\\
&\lesssim \frac{(\log n)^\frac{3}{2}}{\alpha}
+\frac{(\log n)^{\frac{\beta_1+2}{2}}}{\sqrt{n}\alpha^{\frac{\beta_1}{2}}}.
\end{aligned}
\end{equation}
By symmetry,
\[
\int_{\mathbb{R}^d}(\theta-\theta^{\diamond})\pi(\theta^{\diamond})
\exp\Big(-\frac{\alpha}{2}(\theta-\theta^{\diamond})^T
\m H(\theta^{\diamond})(\theta-\theta^{\diamond})\Big)\,\dd\theta=0.
\]
It therefore follows that
\begin{equation}\label{omit 1.4}
\begin{aligned}
&\Big\|\alpha^{\frac{d}{2}}\int(\theta-\theta^{\diamond})\pi(\theta)
\exp\Big(-\alpha\Big(\frac{1}{n}\sum_{i=1}^n \ell( X_i,\theta)
-\frac{1}{n}\sum_{i=1}^n \ell( X_i,\theta^{\diamond})\Big)\Big)\,\dd\theta
\Big\|_2\\
&\lesssim\frac{(\log n)^\frac{3}{2}}{\alpha}
+\frac{(\log n)^{\frac{\beta_1+2}{2}}}{\sqrt{n}\alpha^{\frac{\beta_1}{2}}}.
\end{aligned}
\end{equation}
Applying the denominator lower bound established in \eqref{lowerbound 1} to \eqref{omit 1.5}, we obtain
\begin{equation}\label{mean diff 1.2}
\begin{aligned}
&\Big\|\int \theta\,\pi_{\alpha}(\theta)\,\dd\theta-\theta^{\diamond}\Big\|_2\\
&=\Big\|\int \theta\,\pi_{\alpha}(\theta)\,\dd\theta
-\int \theta\, N\Big(\theta^{\diamond},\frac{1}{\alpha}\m H^{-1}(\theta^{\diamond})\Big)\,\dd\theta\Big\|_2\\
&= \Big\|\int (\theta-\theta^{\diamond})\, \pi_{\alpha}(\theta)\,\dd\theta-\int (\theta-\theta^{\diamond})
N\Big(\theta^{\diamond},\frac{1}{\alpha}\m H^{-1}(\theta^{\diamond})\Big)\,\dd\theta\Big\|_2\\
&= \Big\|\int \frac{(\theta-\theta^{\diamond})\pi(\theta)
\exp\Big(-\alpha\Big(\frac{1}{n}\sum_{i=1}^n \ell(X_i,\theta)
-\frac{1}{n}\sum_{i=1}^n \ell( X_i,\theta^{\diamond})\Big)\Big)}
{\int\pi(\theta)\exp\Big(-\alpha\Big(\frac{1}{n}\sum_{i=1}^n \ell(X_i,\theta)
-\frac{1}{n}\sum_{i=1}^n \ell( X_i,\theta^{\diamond})\Big)\Big)\,\dd\theta}\,\dd\theta \\
&\quad-\int \frac{(\theta-\theta^{\diamond})\pi(\theta^{\diamond})
\exp\Big(-\frac{\alpha}{2}(\theta-\theta^{\diamond})^T
\m H(\theta^{\diamond})(\theta-\theta^{\diamond})\Big)}
{\int\pi(\theta^{\diamond})
\exp\Big(-\frac{\alpha}{2}(\theta-\theta^{\diamond})^T
\m H(\theta^{\diamond})(\theta-\theta^{\diamond})\Big)\,\dd\theta}\,\dd\theta\Big\|_2\\    
&\lesssim\frac{(\log n)^\frac{3}{2}}{\alpha}
+\frac{(\log n)^{\frac{\beta_1+2}{2}}}{\sqrt{n}\alpha^{\frac{\beta_1}{2}}}.
\end{aligned}
\end{equation}
Finally, the triangle inequality, \eqref{mean diff 1.1}, \eqref{mean diff 1.2}, and the previously established bound on $\|\theta^{\diamond}-\wh\theta\|_2$ yield
\begin{equation*}
\begin{aligned}
&\Big\|\int \theta \wt{\pi}_{\mathcal{L}_{\rm zcv},\alpha}(\theta)\,\dd\theta-\wh\theta\Big\|_2\\
&\lesssim\Big\|\int \theta \wt{\pi}_{\mathcal{L}_{\rm zcv},\alpha}(\theta)\,\dd\theta
-\theta^{\diamond}\Big\|_2+\Big\|\theta^{\diamond}-\wh\theta\Big\|_2\\
&\lesssim \sqrt{\alpha}\frac{(\log m)^2}{m}
+ \alpha^{\frac{3}{2}}\frac{\log m}{m}C^2_{m,n}
+ \frac{(\log n)^\frac{3}{2}}{\alpha}
+\frac{(\log n)^{\frac{\beta_1+2}{2}}}{\sqrt{n}\alpha^{\frac{\beta_1}{2}}}
+\Big(\frac{\log n}{n}\Big)^{\frac{\beta_1+1}{2}}\\
&\lesssim \sqrt{\alpha}\frac{(\log m)^2}{m}
+ \alpha^{\frac{3}{2}}\frac{\log m}{m}C^2_{m,n}
+ \frac{(\log n)^\frac{3}{2}}{\alpha}
+\frac{(\log n)^{\frac{\beta_1+2}{2}}}{\sqrt{n}\alpha^{\frac{\beta_1}{2}}}.
\end{aligned}
\end{equation*}
The last inequality follows from $\alpha\lesssim n$.
\subsection{Proof of Corollary \ref{co:2}}
\begin{proof}
Let $\m B$ be the set of $\bX^{(n)}$ defined in the proof of Theorem~\ref{th:1}, and fix an arbitrary realization $\bX^{(n)}\in\m B$. Recall that  $ \theta_P=\int_{\Theta}\theta\,\wt{\pi}_{\mathcal{L}_{\rm zcv},\alpha}(\theta)\,\dd\theta$. As in the proof of Theorem~\ref{th:1}, we decompose the parameter space into the following two regions:
\begin{equation*}
\begin{aligned}
&\Theta'_1
=
\Big\{
\theta\in\mb R^d:
\|\theta-\theta^{\diamond}\|_2
\leq
\delta'_1\sqrt{\frac{\log n}{\alpha}}
\Big\},\\
&\Theta'_2
=
\Big\{
\theta\in\mb R^d:
\delta'_1\sqrt{\frac{\log n}{\alpha}}
<
\|\theta-\theta^{\diamond}\|_2
\Big\},
\end{aligned}
\end{equation*}
where $\delta'_1>0$ is a suitably chosen constant. We begin with the local region $\Theta'_1$. Applying \eqref{estimate 1.10} and using the local Lipschitz continuity of $\pi$ as in \eqref{pi lip}, we obtain
\begin{equation*}
\begin{aligned}
&\alpha^{\frac{d}{2}}
\int_{\Theta'_1}
\Big|
\pi(\theta)
\exp\Big(
-\alpha\Big(
\frac{1}{n}\sum_{i=1}^n\ell(X_i,\theta)
-
\frac{1}{n}\sum_{i=1}^n\ell(X_i,\theta^{\diamond})
\Big)
\Big)\\
&\qquad\qquad
-
\pi(\theta^{\diamond})
\exp\Big(
-\frac{\alpha}{2}
(\theta-\theta^{\diamond})^T
\m H(\theta^{\diamond})
(\theta-\theta^{\diamond})
\Big)
\Big|
\,\dd\theta\\
&\leq
\alpha^{\frac{d}{2}}
\int_{\Theta'_1}
\Big|
\pi(\theta)
\exp\Big(
-\alpha\Big(
\frac{1}{n}\sum_{i=1}^n\ell(X_i,\theta)
-
\frac{1}{n}\sum_{i=1}^n\ell(X_i,\theta^{\diamond})
\Big)
\Big)\\
&\qquad\qquad
-
\pi(\theta)
\exp\Big(
-\frac{\alpha}{2}
(\theta-\theta^{\diamond})^T
\m H(\theta^{\diamond})
(\theta-\theta^{\diamond})
\Big)
\Big|
\,\dd\theta\\
&\quad+
\alpha^{\frac{d}{2}}
\int_{\Theta'_1}
|\pi(\theta)-\pi(\theta^{\diamond})|
\exp\Big(
-\frac{\alpha}{2}
(\theta-\theta^{\diamond})^T
\m H(\theta^{\diamond})
(\theta-\theta^{\diamond})
\Big)
\,\dd\theta\\
&\lesssim
\alpha^{\frac{d}{2}}
\int_{\Theta'_1}
\pi(\theta)
\exp\Big(
-\alpha\Big(
\frac{1}{n}\sum_{i=1}^n\ell(X_i,\theta)
-
\frac{1}{n}\sum_{i=1}^n\ell(X_i,\theta^{\diamond})
\Big)
\Big)\cdot
\Big(
\frac{(\log n)^{\frac{3}{2}}}{\sqrt{\alpha}}
+
\frac{(\log n)^{\frac{\beta_1+2}{2}}}
{\sqrt{n}\alpha^{\frac{\beta_1-1}{2}}}
\Big)
\,\dd\theta
\\&\quad+\alpha^{\frac{d}{2}}\int_{\Theta'_1}\pi(\theta)
\exp\Big(-\frac{\alpha}{2}(\theta-\theta^{\diamond})^T\m H(\theta^{\diamond})(\theta-\theta^{\diamond})\Big) \cdot\Big(\frac{(\log n)^\frac{3}{2}}{\sqrt{\alpha}}
+\frac{(\log n)^{\frac{\beta_1+2}{2}}}{\sqrt{n}\alpha^{\frac{\beta_1-1}{2}}}\Big)\,\dd\theta
+\sqrt{\frac{\log n}{\alpha}}\\
&\lesssim
\frac{(\log n)^{\frac{3}{2}}}{\sqrt{\alpha}}
+
\frac{(\log n)^{\frac{\beta_1+2}{2}}}
{\sqrt{n}\alpha^{\frac{\beta_1-1}{2}}}.
\end{aligned}
\end{equation*}
We next turn to the  region $\Theta'_2$. By the arguments used in \eqref{step 1.2.2} and \eqref{step 1.3.2}, with $\delta'_1$ chosen sufficiently large, we have
\begin{equation*}
\begin{aligned}
&\alpha^{\frac{d}{2}}
\int_{\Theta'_2}
\pi(\theta)
\exp\Big(
-\alpha\Big(
\frac{1}{n}\sum_{i=1}^n\ell(X_i,\theta)
-
\frac{1}{n}\sum_{i=1}^n\ell(X_i,\theta^{\diamond})
\Big)
\Big)
\,\dd\theta
\lesssim
\frac{1}{n^2}.
\end{aligned}
\end{equation*}
An analogous estimate holds for the Gaussian term:
\begin{equation*}
\begin{aligned}
&\alpha^{\frac{d}{2}}
\int_{\Theta'_2}
\pi(\theta^{\diamond})
\exp\Big(
-\frac{\alpha}{2}
(\theta-\theta^{\diamond})^T
\m H(\theta^{\diamond})
(\theta-\theta^{\diamond})
\Big)
\,\dd\theta
\lesssim
\frac{1}{n^2}.
\end{aligned}
\end{equation*}
Taken together, the estimates over $\Theta'_1$ and $\Theta'_2$ imply
\begin{equation*}
\begin{aligned}
&\alpha^{\frac{d}{2}}
\int_{\mathbb{R}^d}
\Big|
\pi(\theta)
\exp\Big(
-\alpha\Big(
\frac{1}{n}\sum_{i=1}^n\ell(X_i,\theta)
-
\frac{1}{n}\sum_{i=1}^n\ell(X_i,\theta^{\diamond})
\Big)
\Big)\\
&\qquad\qquad
-
\pi(\theta^{\diamond})
\exp\Big(
-\frac{\alpha}{2}
(\theta-\theta^{\diamond})^T
\m H(\theta^{\diamond})
(\theta-\theta^{\diamond})
\Big)
\Big|
\,\dd\theta\\
&\lesssim
\frac{(\log n)^{\frac{3}{2}}}{\sqrt{\alpha}}
+
\frac{(\log n)^{\frac{\beta_1+2}{2}}}
{\sqrt{n}\alpha^{\frac{\beta_1-1}{2}}}.
\end{aligned}
\end{equation*}
Using this estimate together with the denominator lower bound established in \eqref{lowerbound 1}, we obtain
\begin{equation}\label{estimate 1.11}
\Big\|
\pi_{\alpha}
-
N\Big(
\theta^{\diamond},
\frac{1}{\alpha}\m H^{-1}(\theta^{\diamond})
\Big)
\Big\|_{\mathrm{TV}}
\lesssim
\frac{(\log n)^{\frac{3}{2}}}{\sqrt{\alpha}}
+
\frac{(\log n)^{\frac{\beta_1+2}{2}}}
{\sqrt{n}\alpha^{\frac{\beta_1-1}{2}}}.
\end{equation}
Setting $\alpha=n$ in the same argument gives
\begin{equation}\label{estimate 1.12}
\begin{aligned}
\Big\|
\pi_n
-
N\Big(
\theta^{\diamond},
\frac{1}{n}\m H^{-1}(\theta^{\diamond})
\Big)
\Big\|_{\mathrm{TV}}
&\lesssim
\frac{(\log n)^{\frac{3}{2}}}{\sqrt{n}}
+
\frac{(\log n)^{\frac{\beta_1+2}{2}}}{n^{\frac{\beta_1}{2}}}\lesssim
\frac{(\log n)^{\frac{\beta_1+2}{2}}}{n^{\frac{\beta_1}{2}}}.
\end{aligned}
\end{equation}
Since $G$ is an invertible affine transformation, the total variation distance is invariant under the pushforward by $G$. Therefore, by \eqref{TV bound 1.1},
\begin{equation}\label{estimate 1.13}
\begin{aligned}
\Big\|
G_{\#}\wt{\pi}_{\mathcal{L}_{\rm zcv},\alpha}
-
G_{\#}\pi_{\alpha}
\Big\|_{\mathrm{TV}}
&=
\Big\|
\wt{\pi}_{\mathcal{L}_{\rm zcv},\alpha}
-
\pi_{\alpha}
\Big\|_{\mathrm{TV}}\lesssim
\alpha\frac{(\log m)^2}{m}
+
\alpha^2\frac{\log m}{m}C^2_{m,n}.
\end{aligned}
\end{equation}
Moreover, the pushforward of the Gaussian distribution $N\Big(
\theta^{\diamond},
\frac{1}{\alpha}\m H^{-1}(\theta^{\diamond})
\Big)
$ under $G$ is given by
\[
G_{\#}
N\Big(
\theta^{\diamond},
\frac{1}{\alpha}\m H^{-1}(\theta^{\diamond})
\Big)
=
N\Big(
\theta_P
+
\sqrt{\frac{\alpha}{n}}
(\theta^{\diamond}-\theta_P),
\frac{1}{n}\m H^{-1}(\theta^{\diamond})
\Big).
\]
Consequently, \eqref{estimate 1.11} and the invariance of the total variation distance yield
\begin{equation}\label{estimate 1.14}
\begin{aligned}
\Big\|
G_{\#}\pi_{\alpha}
-
N\Big(
\theta_P
+
\sqrt{\frac{\alpha}{n}}
(\theta^{\diamond}-\theta_P),
\frac{1}{n}\m H^{-1}(\theta^{\diamond})
\Big)
\Big\|_{\mathrm{TV}}&=
\Big\|
G_{\#}\pi_{\alpha}
-
G_{\#}
N\Big(
\theta^{\diamond},
\frac{1}{\alpha}\m H^{-1}(\theta^{\diamond})
\Big)
\Big\|_{\mathrm{TV}}\\
&=
\Big\|
\pi_{\alpha}
-
N\Big(
\theta^{\diamond},
\frac{1}{\alpha}\m H^{-1}(\theta^{\diamond})
\Big)
\Big\|_{\mathrm{TV}}\\
&\lesssim
\frac{(\log n)^{\frac{3}{2}}}{\sqrt{\alpha}}
+
\frac{(\log n)^{\frac{\beta_1+2}{2}}}
{\sqrt{n}\alpha^{\frac{\beta_1-1}{2}}}.
\end{aligned}
\end{equation}
Since the two Gaussian distributions in the following  have the same covariance matrix, we have
\begin{equation}\label{estimate 1.15}
\begin{aligned}
&\Big\|
G_{\#}
N\Big(
\theta^{\diamond},
\frac{1}{\alpha}\m H^{-1}(\theta^{\diamond})
\Big)
-
N\Big(
\theta^{\diamond},
\frac{1}{n}\m H^{-1}(\theta^{\diamond})
\Big)
\Big\|_{\mathrm{TV}}\\
&=
\Big\|
N\Big(
\theta_P
+
\sqrt{\frac{\alpha}{n}}
(\theta^{\diamond}-\theta_P),
\frac{1}{n}\m H^{-1}(\theta^{\diamond})
\Big)
-
N\Big(
\theta^{\diamond},
\frac{1}{n}\m H^{-1}(\theta^{\diamond})
\Big)
\Big\|_{\mathrm{TV}}\\
&\lesssim
\sqrt{n}
\Big|
1-\sqrt{\frac{\alpha}{n}}
\Big|
\|\theta_P-\theta^{\diamond}\|_2\\
&\lesssim
\sqrt{n}\,
\|\theta_P-\theta^{\diamond}\|_2.
\end{aligned}
\end{equation}
By the definition of $\theta_P$, together with \eqref{mean diff 1.1} and \eqref{mean diff 1.2}, we further have
\begin{equation*}
\begin{aligned}
\|\theta_P-\theta^{\diamond}\|_2
&\leq
\Big\|
\int_{\Theta}
\theta\,
\wt{\pi}_{\mathcal{L}_{\rm zcv},\alpha}(\theta)
\,\dd\theta
-
\int_{\Theta}
\theta\,\pi_{\alpha}(\theta)
\,\dd\theta
\Big\|_2+
\Big\|
\int_{\Theta}
\theta\,\pi_{\alpha}(\theta)
\,\dd\theta
-
\theta^{\diamond}
\Big\|_2\\
&\lesssim
\sqrt{\alpha}\frac{(\log m)^2}{m}
+
\alpha^{\frac{3}{2}}
\frac{\log m}{m}C^2_{m,n}
+
\frac{(\log n)^{\frac{3}{2}}}{\alpha}
+
\frac{(\log n)^{\frac{\beta_1+2}{2}}}
{\sqrt{n}\alpha^{\frac{\beta_1}{2}}}.
\end{aligned}
\end{equation*}
Finally, \eqref{estimate 1.12}--\eqref{estimate 1.15}, together with $\alpha\lesssim n$, yield
\begin{equation*}
\begin{aligned}
\Big\|
G_{\#}\wt{\pi}_{\mathcal{L}_{\rm zcv},\alpha}
-
\pi_n
\Big\|_{\mathrm{TV}}
&\leq
\Big\|
G_{\#}\wt{\pi}_{\mathcal{L}_{\rm zcv},\alpha}
-G_{\#}\pi_{\alpha}\Big\|_{\mathrm{TV}}+\Big\|
G_{\#}\pi_{\alpha}
-
G_{\#}
N\Big(
\theta^{\diamond},
\frac{1}{\alpha}\m H^{-1}(\theta^{\diamond})
\Big)
\Big\|_{\mathrm{TV}}\\
&\quad+
\Big\|
G_{\#}
N\Big(
\theta^{\diamond},
\frac{1}{\alpha}\m H^{-1}(\theta^{\diamond})
\Big)
-
N\Big(
\theta^{\diamond},
\frac{1}{n}\m H^{-1}(\theta^{\diamond})
\Big)
\Big\|_{\mathrm{TV}}\\
&\quad+
\Big\|
N\Big(
\theta^{\diamond},
\frac{1}{n}\m H^{-1}(\theta^{\diamond})
\Big)
-
\pi_n
\Big\|_{\mathrm{TV}}\\
&\lesssim
\sqrt{n\alpha}\frac{(\log m)^2}{m}
+
\sqrt{n}\alpha^{\frac{3}{2}}
\frac{\log m}{m}C^2_{m,n}
+\sqrt{n}\frac{(\log n)^{\frac{3}{2}}}{\alpha}+
\frac{(\log n)^{\frac{\beta_1+2}{2}}}
{\alpha^{\frac{\beta_1}{2}}}.
\end{aligned}
\end{equation*}
\end{proof}
\subsection{Proof of Theorem \ref{th:2}}
Denote
$ \theta^{\diamond}(\bX^{(n)})= \theta^*-\m H(\theta^*)^{-1}\frac{1}{n}\sum_{i=1}^ng( X_i,\theta^*)$. As in the proof of Theorem \ref{th:1}, we use the same simplified notation. Let $\mathcal{B}\subset \m X^{\otimes n}$ denote the event on which all conclusions of Lemma \ref{lemmaprobB1} and Lemma \ref{lemmaprobB2} hold. By choosing $c_0$ large enough, we can ensure that
\[
\m P^*{}^{\otimes n}(\m B^c)\leq\frac{1}{n^2}.
\]
Unless otherwise stated, all subsequent analysis is performed under $\mathcal{B}$.
\subsubsection{Proof of the Total Variation Bound}
By definition
\begin{equation*}
 \begin{aligned}
     \wt{\pi}_{\mathcal{L}_{\rm fcv},\alpha}(\theta)&=
\frac{\mb{E}_{\wX^{(m)} \sim  \mu_n^{\otimes m} }\Big[\pi(\theta) \exp\Big(-\alpha\mathcal{L}_{\rm fcv}(\wX^{(m)}, \theta)\Big)\Big]}{\mb{E}_{\wX^{(m)} \sim  \mu_n^{\otimes m} }\Big[\int\pi(\theta) \exp\Big(-\alpha\mathcal{L}_{\rm fcv}(\wX^{(m)}, \theta)\Big)\,\dd\theta\Big]}\\&=\frac{\mb{E}_{\wX^{(m)} \sim  \mu_n^{\otimes m} }\Big[\pi(\theta) \exp\Big(-\alpha\Big(\mathcal{L}_{\rm fcv}(\wX^{(m)}, \theta)+\frac{1}{n}\sum^n_{i=1}[\ell(X_i,\theta^{\dagger})-\ell(X_i,\theta^{\diamond})]\Big)\Big)\Big]}{\mb{E}_{\wX^{(m)} \sim  \mu_n^{\otimes m} }\Big[\int\pi(\theta) \exp\Big(-\alpha\Big(\mathcal{L}_{\rm fcv}(\wX^{(m)}, \theta)+\frac{1}{n}\sum^n_{i=1}[\ell(X_i,\theta^{\dagger})-\ell(X_i,\theta^{\diamond})]\Big)\Big)\,\dd\theta\Big]},\end{aligned}
\end{equation*}
and
\[\pi_\alpha(\theta)=
\frac{\pi(\theta)\exp\Big(- \alpha\Big(\frac{1}{n}\sum_{i=1}^n \ell(X_i,\theta)-\frac{1}{n}\sum_{i=1}^n \ell(X_i,\theta^{\diamond})\Big)\Big)}{\int\pi(\theta)\exp\Big(-\alpha\Big( \frac{1}{n}\sum_{i=1}^n \ell(X_i,\theta)-\frac{1}{n}\sum_{i=1}^n \ell(X_i,\theta^{\diamond})\Big)\Big)\,\dd\theta}, \]
To bound the deviation between the full-sample and minibatch losses, we will derive the following concentration inequality for $\wX^{(m)}\sim\mu_n^{\otimes m}$ conditional on $\bX^{(n)}\in\mathcal{B}$.
\begin{lemma}\label{lemmaprobA2}
Given $\bX^{(n)} \in\mathcal{B}$ defined above, under Assumptions~$A.1-A.2$ and $B.1-B.3$, for any positive constant $c$, there exists a constant $c_0$ such that, for
$\wX^{(m)}=
(\wt X_1,\wt X_2,\ldots,\wt X_m)
\sim \mu_n^{\otimes m}$,
the following holds with probability at least
$1-m^{-c}$
    \begin{enumerate}
    \item 
    \[
    \displaystyle \sup_{\theta, \theta^{\prime}}\frac{\Big\|\frac{1}{m}\sum_{j=1}^m \big[g(\wt X_j,\theta)- g(\wt X_j,\theta')\big] -\frac{1}{n}\sum_{i=1}^n \big[g(X_i,\theta)-g(X_i,\theta')\big]\Big\|_2}{\sqrt{\frac{\log m}{m}}\|\theta-\theta'\|^{\beta_1}_2+\sqrt{\frac{\log m}{m}}(\frac{\log n}{n})^{\frac{1}{4}}\|\theta-\theta'\|^{\frac{\beta_1+\beta_2}{2}}_2+\frac{\log m}{m}\|\theta-\theta'\|^{\beta_2}_2+\sqrt{\frac{\log m}{m}}\frac{\log n}{n}}\leq c_0.
    \]
    \item
{\footnotesize
\[
\displaystyle
\sup_{\theta,\theta'}
\frac{
\Big|
\frac{1}{m}\sum_{j=1}^m
\big[\ell(\wt X_j,\theta)-\ell(\wt X_j,\theta')
-g(\wt X_j,\theta')\cdot(\theta-\theta')\big]
-
\frac{1}{n}\sum_{i=1}^n
\big[\ell(X_i,\theta)-\ell(X_i,\theta')
-g(X_i,\theta')\cdot(\theta-\theta')\big]
\Big|
}{
\sqrt{\frac{\log m}{m}}\|\theta-\theta'\|_2^{\beta_1+1}
+
\sqrt{\frac{\log m}{m}}
(\frac{\log n}{n})^{\frac14}
\|\theta-\theta'\|_2^{\frac{\beta_1+\beta_2}{2}+1}
+
\frac{\log m}{m}\|\theta-\theta'\|_2^{\beta_2+1}
+
\sqrt{\frac{\log m}{m}}(\frac{\log n}{n})^2
}
\leq c_0.
\]
}
\end{enumerate}
\end{lemma}
Throughout the remainder of the proof, we fix an arbitrary $\bX^{(n)}\in\m B$ and let $\m A$ be the event on which all the conclusions of Lemma~\ref{lemmaprobA1} and Lemma~\ref{lemmaprobA2}  hold for $\widetilde{\bX}^{(m)}$. Then, by Lemmas~\ref{lemmaprobA1} and \ref{lemmaprobA2}, with the assumption that $m\gtrsim n^{\gamma}$, if $c_0$ is sufficiently large, we have 
 \begin{equation}\label{eqn:prob2}
\mb P_{\wX^{(m)}\sim\mu_n^{\otimes m}}(\m A)\geq 1-\frac{1}{m^4\alpha^{d+1}}
\quad\text{and}\quad
\mb P_{\wX^{(m)}\sim\mu_n^{\otimes m}}(\m A^c)\exp\Big((1+\delta'\sqrt{\log n})M_5\Big)\leq\frac{1}{m^2\alpha^{d+1}},
\end{equation}
where the constants $M_5$ and $\delta'$ are  given in \eqref{omit 2.7} and \eqref{omit 2.6} respectively. For notational convenience, define
\begin{equation*}
    \begin{aligned}
        \ell^{\triangle}(\wt X_j,\theta) =&\ell(\wt X_j,\theta)-\ell(\wt X_j,\theta^{\dagger})-g(\wt X_j,\theta^{\dagger})\cdot(\theta-\theta^{\dagger}).
    \end{aligned}
\end{equation*}
In the following, we will bound
\begin{equation*}
    \begin{aligned}
       & \alpha^{\frac{d}{2}}\int_{\Theta}\Big|\mb{E}_{\wX^{(m)}\sim \mu_n^{\otimes m}}\Big[\pi(\theta)\exp\Big(-\alpha\Big(\mathcal{L}_{\rm fcv}(\wX^{(m)}, \theta)+\frac{1}{n}\sum^n_{i=1}[\ell(X_i,\theta^{\dagger})-\ell(X_i,\theta^{\diamond})]\Big)\Big)\Big] \\
         &\qquad-\mb{E}_{\wX^{(m)}\sim \mu_n^{\otimes m}}\Big[\pi(\theta)\exp\Big(-\alpha\Big(  \frac{1}{n}\sum_{i=1}^n\ell(X_i,\theta)-\frac{1}{n}\sum_{i=1}^n\ell (X_i,\theta^{\diamond})\Big)\Big)\Big]\Big|\,\dd\theta
         \end{aligned}
         \end{equation*}
 Define the following partition of $\Theta$:
\begin{equation}\label{setstheta2}
\begin{aligned}
&\Theta_1=\Big\{\theta\in \Theta : \|\theta-\theta^{\diamond}\|_2\leq \delta_1\sqrt{\frac{\log m}{\alpha}}\Big\},\\
&\Theta_2=\Big\{\theta\in \Theta : \delta_1\sqrt{\frac{\log m}{\alpha}}<\|\theta-\theta^{\diamond}\|_2\leq \delta_2\Big\},\\
&\Theta_3=\Big\{\theta\in \Theta : \|\theta-\theta^{\diamond}\|_2> \delta_2\Big\}.
\end{aligned}
\end{equation}
Here, $\delta_1,\delta_2>0$ are suitably chosen constants.
\paragraph{1. \underline{On $\Theta_1$.}}
For $\wX^{(m)}\in \m A$ and $\theta\in \Theta_1$,\\
\textbf{Case (i):}
If $\frac{1}{m}\sum_{j=1}^m \ell^{\triangle}(\wt X_j,\theta)\geq 0$, we have
\begin{equation*}
\Big|\Big(\frac{1}{m}\sum_{j=1}^m \ell^{\triangle}(\wt X_j,\theta)\Big)_{+}-\frac{1}{n}\sum_{i=1}^n \ell^{\triangle}(X_i,\theta)\Big|=\Big|\frac{1}{m}\sum_{j=1}^m \ell^{\triangle}(\wt X_j,\theta)-\frac{1}{n}\sum_{i=1}^n \ell^{\triangle}(X_i,\theta)\Big|.
\end{equation*}\textbf{Case (ii):}
If $\frac{1}{m}\sum_{j=1}^m \ell^{\triangle}(\wt X_j,\theta)\leq 0$, we have
\begin{equation*}
\begin{aligned}
\Big|\Big(\frac{1}{m}\sum_{j=1}^m \ell^{\triangle}(\wt X_j,\theta)\Big)_{+}-\frac{1}{n}\sum_{i=1}^n \ell^{\triangle}(X_i,\theta)\Big|
&=\Big|\frac{1}{n}\sum_{i=1}^n \ell^{\triangle}(X_i,\theta)\Big|\\
&\leq \Big|\frac{1}{m}\sum_{j=1}^m \ell^{\triangle}(\wt X_j,\theta)-\frac{1}{n}\sum_{i=1}^n \ell^{\triangle}(X_i,\theta)\Big|\\
&\quad+\Big|\frac{1}{n}\sum_{i=1}^n \ell^{\triangle}(X_i,\theta)- \m R(\theta) +\m R(\theta^{\dagger})+ \nabla\m R(\theta^{\dagger})\cdot(\theta - \theta^{\dagger})\Big|,
\end{aligned}
\end{equation*}
where the last inequality is due to the fact that
\begin{equation}\label{estimate 10.7}
\begin{aligned}
&\m R(\theta) -\m R(\theta^{\dagger})- \nabla\m R(\theta^{\dagger})\cdot(\theta - \theta^{\dagger})\\
&=(\theta - \theta^\dagger)^T \int_{0}^{1}(1-t)\m H\big(\theta^\dagger + t(\theta - \theta^\dagger)\big)dt\,(\theta - \theta^\dagger)\gtrsim\|\theta - \theta^\dagger\|^2_2
\geq0.
\end{aligned}
\end{equation}
 Hence, fix an arbitrary $\wX^{(m)}\in \m A$, when $\theta\in \Theta_1$, by the last statement of Lemma \ref{lemmaprobB2} and the second statement of Lemma \ref{lemmaprobA2}, we have 
\begin{equation}\label{omit 2.1}
{\small
\begin{aligned}
&\Big|\Big(\frac{1}{m}\sum_{j=1}^m \ell^{\triangle}(\wt X_j,\theta)\Big)_{+}-\frac{1}{n}\sum_{i=1}^n \ell^{\triangle}(X_i,\theta)\Big|\\
&\leq\Big|\Big[\frac{1}{m}\sum_{j=1}^m \ell(\wt X_j,\theta)-\frac{1}{m}\sum_{j=1}^m \ell(\wt X_j,\theta^{\dagger})-\frac{1}{m}\sum_{j=1}^mg(\wt X_j,\theta^{\dagger})\cdot\big(\theta-\theta^{\dagger}\big)\Big]\\
&\quad-\Big[\frac{1}{n}\sum_{i=1}^n \ell(X_i,\theta)-\frac{1}{n}\sum_{i=1}^n \ell(X_i,\theta^{\dagger})-\frac{1}{n}\sum_{i=1}^ng(X_i,\theta^{\dagger})\cdot\big(\theta-\theta^{\dagger}\big)\Big]\Big|\\
&\quad+\Big|\Big[\frac{1}{n}\sum_{i=1}^n \ell(X_i,\theta)-\frac{1}{n}\sum_{i=1}^n \ell(X_i,\theta^{\dagger})-\frac{1}{n}\sum_{i=1}^ng(X_i,\theta^{\dagger})\cdot\big(\theta-\theta^{\dagger}\big)\Big]\\
&\quad-\Big[\m R(\theta) -\m R(\theta^{\dagger})- \nabla\m R(\theta^{\dagger})\cdot(\theta - \theta^{\dagger})\Big]\Big|\\
&\lesssim\sqrt{\frac{\log m}{m}}\|\theta-\theta^{\dagger}\|^{\beta_1+1}_2+\sqrt{\frac{\log m}{m}}\Big(\frac{\log n}{n}\Big)^{\frac{1}{4}}\|\theta-\theta^{\dagger}\|^{\frac{\beta_1+\beta_2}{2}+1}_2+\frac{\log m}{m}\|\theta-\theta^{\dagger}\|^{\beta_2+1}_2\\
&\quad+\sqrt{\frac{\log m}{m}}\Big(\frac{\log n}{n}\Big)^{2}+\sqrt{\frac{\log n}{n}}\|\theta-\theta^{\dagger}\|^{\beta_1+1}_2+\frac{\log n}{n}\|\theta-\theta^{\dagger}\|^{\beta_2+1}_2+\Big(\frac{\log n}{n}\Big)^{3}\\&\lesssim\sqrt{\frac{\log m}{m}}\|\theta-\theta^{\dagger}\|^{\beta_1+1}_2+\frac{\log m}{m}\|\theta-\theta^{\dagger}\|^{\beta_2+1}_2+\sqrt{\frac{\log m}{m}}\Big(\frac{\log n}{n}\Big)^{2},
\end{aligned}}
\end{equation}
where the last inequality follows from the Cauchy-Schwarz inequality and
$\frac{\log n}{n}\lesssim\frac{\log m}{m}$. 
By the previously established bound \eqref{estimate 0.4} on
$\|\wh\theta-\theta^{\diamond}\|_2$, together with
$m\gtrsim n^{\gamma}$ and $\alpha\leq n$, we have
\[
\|\wh\theta-\theta^{\diamond}\|_2
\lesssim
\Big(\frac{\log n}{n}\Big)^{\frac{1+\beta_1}{2}}
\lesssim
\sqrt{\frac{\log m}{\alpha}}.
\]
Substituting this bound into the preceding estimate yields
\begin{equation}\label{estimate 2.1}
\begin{aligned}
&\Big|
\Big(\frac{1}{m}\sum_{j=1}^m
\ell^{\triangle}(\wt X_j,\theta)\Big)_{+}
-\frac{1}{n}\sum_{i=1}^n
\ell^{\triangle}(X_i,\theta)
\Big|\\
&\lesssim
\frac{(\log m)^{\frac{\beta_1+2}{2}}}
{m^{\frac{1}{2}}\alpha^{\frac{\beta_1+1}{2}}}
+\frac{(\log m)^{\frac{1}{2}}}{m^{\frac{1}{2}}}
C^{\beta_1+1}_{m,n}
+\frac{(\log m)^{\frac{\beta_2+3}{2}}}
{m\alpha^{\frac{\beta_2+1}{2}}}\\
&\quad
+\frac{\log m}{m}C^{\beta_2+1}_{m,n}
+\sqrt{\frac{\log m}{m}}
\Big(\frac{\log n}{n}\Big)^2\\
&\lesssim
\frac{(\log m)^{\frac{\beta_1+2}{2}}}
{m^{\frac{1}{2}}\alpha^{\frac{\beta_1+1}{2}}}
+\frac{(\log m)^{\frac{\beta_2+3}{2}}}
{m\alpha^{\frac{\beta_2+1}{2}}}.
\end{aligned}
\end{equation}
Here, the last inequality follows from
$C_{m,n}\lesssim\alpha^{-1/2}$
and $\alpha\leq n$.
\\Therefore, for any $\theta\in \Theta_1$ and $\wX^{(m)}\in \m A$, we have 
\begin{equation*}
\begin{aligned}
&\exp\Big(-\alpha\Big(\Big(\frac{1}{m}\sum^m_{j=1}\ell^{\triangle}(\wt X_j,\theta)\Big)_{+}+\frac{1}{n}\sum^n_{i=1}\ell(X_i,\theta^{\dagger})-\frac{1}{n}\sum_{i=1}^n \ell(X_i,\theta^{\diamond})+\frac{1}{n}\sum_{i=1}^ng(X_i,\theta^{\dagger})\cdot(\theta-\theta^{\dagger})\Big)\Big)\\
&-\exp\Big(-\alpha\Big( \frac{1}{n}\sum_{i=1}^n \ell(X_i,\theta)-\frac{1}{n}\sum_{i=1}^n \ell(X_i,\theta^{\diamond})\Big)\Big)\\
&=\exp\Big(-\alpha\Big( \frac{1}{n}\sum_{i=1}^n \ell(X_i,\theta)-\frac{1}{n}\sum_{i=1}^n \ell(X_i,\theta^{\diamond})\Big)\Big)\\
&\quad\cdot\Big[\m O\Big(\alpha^{1-\beta_1}\frac{(\log m)^{\beta_1+2}}{m}+\alpha^{1-\beta_2}\frac{(\log m)^{\beta_2+3}}{m^2}\Big)-\alpha\cdot\Big(\Big(\frac{1}{m}\sum^m_{j=1}\ell^{\triangle}(\wt X_j,\theta)\Big)_{+}
-\frac{1}{n}\sum^n_{i=1}\ell^{\triangle}(X_i,\theta)\Big)\Big].\end{aligned}\end{equation*}
Hence, we can obtain
\begin{equation}\label{estimate 2.2}
\begin{aligned}
& \alpha^{\frac{d}{2}}\int_{\Theta_1}\Big|\mb{E}_{\wX^{(m)}\sim \mu_n^{\otimes m}}\Big[\pi(\theta)\exp\Big(-\alpha\Big(\mathcal{L}_{\rm fcv}(\wX^{(m)}, \theta)+\frac{1}{n}\sum^n_{i=1}[\ell(X_i,\theta^{\dagger})-\ell(X_i,\theta^{\diamond})]\Big)\Big)\cdot\bold{1}(\wX^{(m)}\in \m A)\Big]\\
&-\mb{E}_{\wX^{(m)}\sim \mu_n^{\otimes m}}\Big[\pi(\theta)\exp\Big(-\alpha\Big( \frac{1}{n}\sum_{i=1}^n \ell(X_i,\theta)-\frac{1}{n}\sum_{i=1}^n \ell(X_i,\theta^{\diamond})\Big)\Big)\cdot\bold{1}(\wX^{(m)}\in \m A)\Big]\Big|\,\dd\theta\\
&\lesssim \alpha^{\frac{d}{2}}\int_{\Theta_1}\pi(\theta)\exp\Big(-\alpha\Big( \frac{1}{n}\sum_{i=1}^n \ell(X_i,\theta)-\frac{1}{n}\sum_{i=1}^n \ell(X_i,\theta^{\diamond})\Big)\Big)\\
&\quad\cdot\Big(\alpha^{1-\beta_1}\frac{(\log m)^{\beta_1+2}}{m}+\alpha^{1-\beta_2}\frac{(\log m)^{\beta_2+3}}{m^2}\\
&\quad+\alpha\cdot\Big|\mb{E}_{\wX^{(m)}\sim \mu_n^{\otimes m}}\Big[\bold{1}(\wX^{(m)}\in \m A)\cdot\Big(\Big(\frac{1}{m}\sum^m_{j=1}\ell^{\triangle}(\wt X_j,\theta)\Big)_{+}
-\frac{1}{n}\sum^n_{i=1}\ell^{\triangle}(X_i,\theta)\Big)\Big]\Big|\Big)\,\dd\theta.
\end{aligned}
\end{equation}
As $\ell(x,\cdot)$ is Lipschitz and $g (x,\theta)$ is bounded, we have 
\begin{equation}\label{omit 2.3}
\begin{aligned}
&\Big|\mb{E}_{\wX^{(m)}\sim \mu_n^{\otimes m}}\Big[\bold{1}(\wX^{(m)}\in \m A) \cdot \Big(\frac{1}{m}\sum^m_{j=1}\ell^{\triangle}(\wt X_j,\theta)-\frac{1}{n}\sum^n_{i=1}\ell^{\triangle}(X_i,\theta)\Big)\Big]\Big|\\
&\leq \Big|\mb{E}_{\wX^{(m)}\sim \mu_n^{\otimes m}}\Big[\frac{1}{m}\sum^m_{j=1}\ell^{\triangle}(\wt X_j,\theta)-\frac{1}{n}\sum^n_{i=1}\ell^{\triangle}(X_i,\theta)\Big]\Big|\\
&+\Big|\mb{E}_{\wX^{(m)}\sim \mu_n^{\otimes m}}\Big[\bold{1}(\wX^{(m)}\in \m A^c) \cdot \Big(\frac{1}{m}\sum^m_{j=1}\ell^{\triangle}(\wt X_j,\theta)-\frac{1}{n}\sum^n_{i=1}\ell^{\triangle}(X_i,\theta)\Big)\Big]\Big|\\
&\lesssim\mb P_{\wX^{(m)}\sim \mu_n^{\otimes m}}(\m A^c)\cdot\Big( \underset{x\in \m X,\theta\in\Theta}{\sup}|\ell(x,\theta)-\ell(x,\theta^{\dagger})|
+\underset{x\in \m X,\theta\in\Theta}{\sup}\|g(x,\theta)\|_2\cdot\operatorname{diam}(\Theta)\Big)\\ 
&\lesssim \frac{1}{m^4\alpha^{d+1}}.
\end{aligned}
\end{equation}
Applying the final assertion of Lemma~\ref{lemmaprobB2} and the second assertion of Lemma~\ref{lemmaprobA2}, and arguing as in the derivation of \eqref{omit 2.1}, we obtain
\begin{equation*}
    \begin{aligned}
&\Big|\frac{1}{m}\sum^m_{j=1}\ell^{\triangle}(\wt X_j,\theta)-\m R(\theta) +\m R(\theta^{\dagger})+ \nabla\m R(\theta^{\dagger})\cdot(\theta - \theta^{\dagger})\Big|\\
&\lesssim\sqrt{\frac{\log m}{m}}\,\|\theta-\theta^{\dagger}\|^{\beta_1+1}_2+\frac{\log m}{m}\|\theta-\theta^{\dagger}\|^{\beta_2+1}_2+\sqrt{\frac{\log m}{m}}\Big(\frac{\log n}{n}\Big)^{2}.
\end{aligned}
\end{equation*}
It follows from this estimate and \eqref{estimate 10.7} that
\begin{equation}\label{cutoff}
    \begin{aligned}
    \frac{1}{m}\sum^m_{j=1}\ell^{\triangle}(\wt X_j,\theta)&\gtrsim\|\theta-\theta^{\dagger}\|^{2}_2-\sqrt{\frac{\log m}{m}} \,\|\theta-\theta^{\dagger}\|^{\beta_1+1}_2-\frac{\log m}{m}\|\theta-\theta^{\dagger}\|^{\beta_2+1}_2\\&-\sqrt{\frac{\log m}{m}}\Big(\frac{\log n}{n}\Big)^{2}.
    \end{aligned}
    \end{equation}
 For $\beta_1 = 1$, the term $\sqrt{\frac{\log m}{m}} \|\theta - \theta^{\dagger}\|_2^{\beta_1+1}$ is absorbed into the positive quadratic term, and analogously for $\beta_2 = 1$. Repeating the argument below after these absorptions establishes \eqref{estimate 2.3} without the summand $\alpha\big(\frac{\log m}{m}\big)^{\frac{1}{1-\beta_1}}$ when $\beta_1 = 1$, and without $\alpha\big(\frac{\log m}{m}\big)^{\frac{2}{1-\beta_2}}$ when $\beta_2 = 1$. We may therefore assume $\beta_1 < 1$ and $\beta_2 < 1$. Consequently, there exists a sufficiently large constant $\delta'' > 0$ such that, whenever
\[
\|\theta-\theta^{\dagger}\|_2\geq \delta''\max\Big( \Big(\frac{\log m}{m}\Big)^{\frac{1}{2(1-\beta_1)}},\Big(\frac{\log m}{m}\Big)^{\frac{1}{1-\beta_2}},\Big(\frac{\log m}{m}\Big)^{\frac{1}{4}}\Big(\frac{\log n}{n}\Big)\Big),
\]
the lower bound in \eqref{cutoff} is positive. Therefore,
\begin{equation*}
\Big(\frac{1}{m}\sum^m_{j=1}\ell^{\triangle}(\wt X_j,\theta)\Big)_{+}=\frac{1}{m}\sum^m_{j=1}\ell^{\triangle}(\wt X_j,\theta).
\end{equation*}
On the other hand, suppose that
\begin{equation*}
\|\theta-\theta^{\dagger}\|_2\leq\delta''\max\Big( \Big(\frac{\log m}{m}\Big)^{\frac{1}{2(1-\beta_1)}},\Big(\frac{\log m}{m}\Big)^{\frac{1}{1-\beta_2}},\Big(\frac{\log m}{m}\Big)^{\frac{1}{4}}\Big(\frac{\log n}{n}\Big)\Big),
\end{equation*}
and that
\[
\frac{1}{m}\sum^m_{j=1}\ell^{\triangle}(\wt X_j,\theta)\leq 0,
\]
we have
\begin{equation*}
\begin{aligned}
&\Big|\Big(\frac{1}{m}\sum^m_{j=1}\ell^{\triangle}(\wt X_j,\theta)\Big)_{+}-\frac{1}{m}\sum^m_{j=1}\ell^{\triangle}(\wt X_j,\theta)\Big|\\
&= \Big|\frac{1}{m}\sum^m_{j=1}\ell^{\triangle}(\wt X_j,\theta)\Big|\\&\lesssim\|\theta-\theta^{\dagger}\|^{2}_2+\sqrt{\frac{\log m}{m}} \,\|\theta-\theta^{\dagger}\|^{\beta_1+1}_2+\frac{\log m}{m}\|\theta-\theta^{\dagger}\|^{\beta_2+1}_2+\sqrt{\frac{\log m}{m}}\Big(\frac{\log n}{n}\Big)^{2}\\&\lesssim \Big(\frac{\log m}{m}\Big)^{\frac{1}{1-\beta_1}}+\Big(\frac{\log m}{m}\Big)^{\frac{2}{1-\beta_2}}+\Big(\frac{\log m}{m}\Big)^{\frac{1}{2}}\Big(\frac{\log n}{n}\Big)^{2}.
\end{aligned}
\end{equation*}
Consequently, for every $\theta \in \Theta_1$, we have
\begin{equation}\label{omit 2.2}
    \begin{aligned}
        &\Big|\Big(\frac{1}{m}\sum^m_{j=1}\ell^{\triangle}(\wt X_j,\theta)\Big)_{+}-\frac{1}{m}\sum^m_{j=1}\ell^{\triangle}(\wt X_j,\theta)\Big|\lesssim \Big(\frac{\log m}{m}\Big)^{\frac{1}{1-\beta_1}}+\Big(\frac{\log m}{m}\Big)^{\frac{2}{1-\beta_2}}+\Big(\frac{\log m}{m}\Big)^{\frac{1}{2}}\Big(\frac{\log n}{n}\Big)^{2},
    \end{aligned}
\end{equation}
Together with \eqref{omit 2.3}, the preceding estimate yields
\begin{equation}\label{estimate 2.3}
\begin{aligned}
& \alpha\cdot\Big|\mb{E}_{\wX^{(m)}\sim \mu_n^{\otimes m}}\Big[\bold{1}(\wX^{(m)}\in \m A) \cdot\Big(\Big(\frac{1}{m}\sum^m_{j=1}\ell^{\triangle}(\wt X_j,\theta)\Big)_{+}
-\frac{1}{n}\sum^n_{i=1}\ell^{\triangle}(X_i,\theta)\Big)\Big]\Big|
\\
&\leq  \alpha\cdot\Big|\mb{E}_{\wX^{(m)}\sim \mu_n^{\otimes m}}\Big[\bold{1}(\wX^{(m)}\in \m A) \cdot\Big(\frac{1}{m}\sum^m_{j=1}\ell^{\triangle}(\wt X_j,\theta)
-\frac{1}{n}\sum^n_{i=1}\ell^{\triangle}(X_i,\theta)\Big)\Big]\Big|
\\
&\qquad+\alpha\cdot\Big|\mb{E}_{\wX^{(m)}\sim \mu_n^{\otimes m}}\Big[\bold{1}(\wX^{(m)}\in \m A) \cdot\Big(\Big(\frac{1}{m}\sum^m_{j=1}\ell^{\triangle}(\wt X_j,\theta)\Big)_{+}
-\frac{1}{m}\sum^m_{j=1}\ell^{\triangle}(\wt X_j,\theta)\Big)\Big]\Big|
\\
&\lesssim\alpha\Big(\frac{\log m}{m}\Big)^{\frac{1}{1-\beta_1}}+\alpha\Big(\frac{\log m}{m}\Big)^{\frac{2}{1-\beta_2}}+\alpha\Big(\frac{\log m}{m}\Big)^{\frac{1}{2}}\Big(\frac{\log n}{n}\Big)^{2}+\frac{1}{m^4\alpha^{d}}.
\end{aligned}
\end{equation}
By the choice of $m, n$, and $\alpha$, the following inequalities hold for $\beta_i < 1$ (while the corresponding left-hand side vanishes when $\beta_i = 1$):
\begin{equation}
\begin{aligned}\label{omit 2.4}
\alpha\Big(\frac{\log m}{m}\Big)^{\frac{1}{1-\beta_1}}
&\lesssim
\alpha^{1-\beta_1}\frac{(\log m)^{\beta_1+2}}{m},\\
\alpha\Big(\frac{\log m}{m}\Big)^{\frac{2}{1-\beta_2}}
&\lesssim
\alpha^{1-\beta_2}\frac{(\log m)^{\beta_2+3}}{m^2}.
\end{aligned}
\end{equation}
Moreover, for every $\theta\in\Theta_1$, we have
\begin{equation*}
    \alpha\frac{\log n}{n}\| \theta - \theta^{\diamond} \|_2^{2\beta_1}\lesssim
    \frac{(\log m)^{\beta_1+1}}{n^{\beta_1}}.
\end{equation*}
It then follows from \eqref{estimate 2.2} and \eqref{estimate 2.3} that
\begin{equation}\label{step 2.1}
\begin{aligned}
& \alpha^{\frac{d}{2}}\int_{\Theta_1}\Big|\mb{E}_{\wX^{(m)}\sim \mu_n^{\otimes m}}\Big[\pi(\theta)\exp\Big(-\alpha\Big(\mathcal{L}_{\rm fcv}(\wX^{(m)}, \theta)+\frac{1}{n}\sum^n_{i=1}[\ell(X_i,\theta^{\dagger})-\ell(X_i,\theta^{\diamond})]\Big)\Big)\bold{1}(\wX^{(m)}\in \m A)\Big]\\
&-\mb{E}_{\wX^{(m)}\sim \mu_n^{\otimes m}}\Big[\pi(\theta)\exp\Big(-\alpha\Big(  \frac{1}{n}\sum_{i=1}^n\ell(X_i,\theta)-\frac{1}{n}\sum_{i=1}^n\ell(X_i,\theta^{\diamond})\Big)\Big)\cdot\bold{1}(\wX^{(m)}\in \m A)\Big]\Big|\,\dd\theta
\\
&\lesssim \alpha^{\frac{d}{2}}\int_{\Theta_1}\pi(\theta)\exp\Big(-\alpha M_1\|\theta - \theta^{\diamond}\|^2_2 
+\alpha M_1\frac{\log n}{n}\| \theta - \theta^{\diamond} \|_2^{2\beta_1}
+\alpha M_1\Big( \frac{\log n}{n} \Big)^{\beta_1+1}\Big)\\&\cdot\Big(\alpha^{1-\beta_1}\frac{(\log m)^{\beta_1+2}}{m}
+\alpha^{1-\beta_2}\frac{(\log m)^{\beta_2+3}}{m^2}+\alpha\Big(\frac{\log m}{m}\Big)^{\frac{1}{2}}\Big(\frac{\log n}{n}\Big)^{2}+\frac{1}{m^4\alpha^{d}}
\\&\quad+\alpha\cdot\Big|\mb{E}_{\wX^{(m)}\sim \mu_n^{\otimes m}}\Big[\bold{1}(\wX^{(m)}\in \m A)\Big(\Big(\frac{1}{m}\sum^m_{j=1}\ell^{\triangle}(\wt X_j,\theta)\Big)_{+}
-\frac{1}{n}\sum^n_{i=1}\ell^{\triangle}(X_i,\theta)\Big)\Big]\Big|\Big)\,\dd\theta
\\&\lesssim \alpha^{\frac{d}{2}}\int_{\Theta_1}\pi(\theta)\exp\Big(-\alpha M_1\|\theta - \theta^{\diamond}\|^2_2\Big)\cdot\Big(\alpha^{1-\beta_1}\frac{(\log m)^{\beta_1+2}}{m}+\alpha^{1-\beta_2}\frac{(\log m)^{\beta_2+3}}{m^2}\Big)\,\dd\theta
\\
&\lesssim \alpha^{1-\beta_1}\frac{(\log m)^{\beta_1+2}}{m}+\alpha^{1-\beta_2}\frac{(\log m)^{\beta_2+3}}{m^2},
\end{aligned}
\end{equation}
By \eqref{estimate 1.3}, there exists a constant $k_2>0$ such that
\begin{equation}\label{lowerbound 2}
\begin{aligned}
&\alpha^{\frac{d}{2}}\int_{\Theta_1}\pi(\theta)\exp\Big(-\alpha\Big( \frac{1}{n}\sum_{i=1}^n \ell(X_i,\theta)-\frac{1}{n}\sum_{i=1}^n \ell(X_i,\theta^{\diamond})\Big)\Big)\,\dd\theta\\
&\geq \alpha^{\frac{d}{2}}\int_{\Theta_1}\pi(\theta)\exp\Big(-\alpha M_2\|\theta - \theta^{\diamond}\|^2_2 
-\alpha M_2\frac{\log n}{n}\| \theta - \theta^{\diamond} \|_2^{2\beta_1}
-\alpha M_2\Big( \frac{\log n}{n} \Big)^{\beta_1+1}\Big)\,\dd\theta\\
&\gtrsim \alpha^{\frac{d}{2}}\int_{\Theta_1}\pi(\theta)\exp\Big(-\alpha M_2\|\theta - \theta^{\diamond}\|^2_2\Big)\,\dd\theta \geq k_2 .
\end{aligned}
\end{equation}
\paragraph{2. \underline{On $\Theta_2$.}}
We next turn to the region $\Theta_2$. For any $\theta\in\Theta_2$ and $\wX^{(m)}\in\m A$, the second assertion of Lemma~\ref{lemmaprobA2} yields
\begin{equation*}
\begin{aligned}
&\Big(\frac{1}{m}\sum^m_{j=1}\ell^{\triangle}(\wt X_j,\theta)\Big)_{+}+\frac{1}{n}\sum^n_{i=1}\ell(X_i,\theta^{\dagger})-\frac{1}{n}\sum_{i=1}^n \ell(X_i,\theta^{\diamond})+\frac{1}{n}\sum_{i=1}^ng(X_i,\theta^{\dagger})\cdot(\theta-\theta^{\dagger})\\
&\geq\frac{1}{m}\sum_{j=1}^m \big[\ell(\wt X_j,\theta)-\ell(\wt X_j,\theta^{\dagger})-g(\wt X_j,\theta^{\dagger})\cdot(\theta-\theta^{\dagger})\big]\\
&+\frac{1}{n}\sum_{i=1}^n \big[\ell(X_i,\theta^{\dagger})-\ell(X_i,\theta^{\diamond})+g(X_i,\theta^{\dagger})\cdot(\theta-\theta^{\dagger})\big]\\
&=\frac{1}{n}\sum_{i=1}^n \big[\ell(X_i,\theta^{\dagger})-\ell(X_i,\theta^{\diamond})+g(X_i,\theta^{\dagger})\cdot(\theta-\theta^{\dagger})\big]\\
&+\frac{1}{n}\sum_{i=1}^n \big[\ell(X_i,\theta)-\ell(X_i,\theta^{\dagger})-g(X_i,\theta^{\dagger})\cdot(\theta-\theta^{\dagger})\big]\\
&+
\frac{1}{m}\sum_{j=1}^m \big[\ell(\wt X_j,\theta)-\ell(\wt X_j,\theta^{\dagger})-g(\wt X_j,\theta^{\dagger})\cdot(\theta-\theta^{\dagger})\big]\\
&-\frac{1}{n}\sum_{i=1}^n \big[\ell(X_i,\theta)-\ell(X_i,\theta^{\dagger})-g(X_i,\theta^{\dagger})\cdot(\theta-\theta^{\dagger})\big]\\
&\geq\Big(\frac{1}{n}\sum_{i=1}^n \ell(X_i,\theta)-\frac{1}{n}\sum_{i=1}^n \ell(X_i,\theta^{\diamond})\Big)-c_0\Big(\sqrt{\frac{\log m}{m}}\,\|\theta-\theta^{\dagger}\|^{\beta_1+1}_2\\
&\quad+\sqrt{\frac{\log m}{m}}\Big(\frac{\log n}{n}\Big)^{\frac{1}{4}}\|\theta-\theta^{\dagger}\|^{\frac{\beta_1+\beta_2}{2}+1}_2+\frac{\log m}{m}\|\theta-\theta^{\dagger}\|^{\beta_2+1}_2+\sqrt{\frac{\log m}{m}}\Big(\frac{\log n}{n}\Big)^{2}\Big)\\&\gtrsim\Big(\frac{1}{n}\sum_{i=1}^n \ell(X_i,\theta)-\frac{1}{n}\sum_{i=1}^n \ell(X_i,\theta^{\diamond})\Big)\\
&\quad-\Big(\sqrt{\frac{\log m}{m}}\,\|\theta-\theta^{\dagger}\|^{\beta_1+1}_2+\frac{\log m}{m}\|\theta-\theta^{\dagger}\|^{\beta_2+1}_2+\sqrt{\frac{\log m}{m}}\Big(\frac{\log n}{n}\Big)^{2}\Big).
\end{aligned}
\end{equation*}
Using \eqref{estimate 1.2} and the bounds
\[
\|\wh\theta-\theta^{\diamond}\|_2
\lesssim
\Big(\frac{\log n}{n}\Big)^{\frac{1+\beta_1}{2}}
\lesssim
\sqrt{\frac{\log m}{\alpha}}
\qquad\text{and}\qquad
C_{m,n}\lesssim\alpha^{-1/2},
\]
we obtain 
\begin{equation*}
\begin{aligned}
&\Big(\frac{1}{m}\sum^m_{j=1}\ell^{\triangle}(\wt X_j,\theta)\Big)_{+}+\frac{1}{n}\sum^n_{i=1}\ell(X_i,\theta^{\dagger})-\frac{1}{n}\sum_{i=1}^n \ell(X_i,\theta^{\diamond})+\frac{1}{n}\sum_{i=1}^ng(X_i,\theta^{\dagger})\cdot(\theta-\theta^{\dagger})\\
&\geq \delta'''\Big[\Big(\|\theta - \theta^{\diamond}\|^2_2 -\frac{\log n}{n}\| \theta - \theta^{\diamond} \|_2^{2\beta_1}-\Big( \frac{\log n}{n} \Big)^{\beta_1+1}\Big)\\
&\quad-\Big(\sqrt{\frac{\log m}{m}}\,\|\theta-\theta^{\diamond}\|^{\beta_1+1}_2+\frac{\log m}{m}\|\theta-\theta^{\diamond}\|^{\beta_2+1}_2+\sqrt{\frac{\log m}{m}}\Big(\frac{\log n}{n}\Big)^{2}\Big)\Big],
\end{aligned}
\end{equation*}
where $\delta'''$ is a positive  constant. When $\beta_1 < 1$ and $\beta_2 < 1$, Young's inequality gives
\begin{equation}\label{estimate 2.4}
    \begin{aligned}
    \sqrt{\frac{\log m}{m}} \, \|\theta-\theta^{\diamond}\|^{\beta_1+1}_2
      &\le \frac{1-\beta_1}{2} \Big( \frac{\log m}{m} \Big)^{\frac{1}{1-\beta_1}} 
         \delta^{-\frac{\beta_1+1}{1-\beta_1}}
         + \frac{\beta_1+1}{2}\delta \|\theta-\theta^{\diamond}\|^2_2, \\[6pt]
    \frac{\log m}{m} \|\theta-\theta^{\diamond}\|^{\beta_2+1}_2
      &\le \frac{1-\beta_2}{2} \Big( \frac{\log m}{m} \Big)^{\frac{2}{1-\beta_2}} 
         \delta^{-\frac{\beta_2+1}{1-\beta_2}}
         + \frac{\beta_2+1}{2}\delta \|\theta-\theta^{\diamond}\|^2_2.
    \end{aligned}
\end{equation}
For $\beta_i = 1$ ($i = 1, 2$), given any $\delta > 0$ and sufficiently large $m$, we have
\[
\sqrt{\frac{\log m}{m}}\,\|\theta-\theta^\diamond\|_2^2
\le \delta\|\theta-\theta^\diamond\|_2^2
\qquad \text{and} \qquad
\frac{\log m}{m}\,\|\theta-\theta^\diamond\|_2^2
\le \delta\|\theta-\theta^\diamond\|_2^2.
\]
For every $\theta\in\Theta_2$, we have
\[
\|\theta-\theta^{\diamond}\|_2^2
\gtrsim
\frac{\log n}{n}
\|\theta-\theta^{\diamond}\|_2^{2\beta_1}.
\]
It follows from this observation, together with \eqref{omit 2.4} and the assumed scaling conditions on $m$ and $\alpha$, that, for sufficiently large $\delta_1$ and sufficiently small $\delta$,
\begin{equation}\label{step 2.2.1}
\begin{aligned}
&\alpha^{\frac{d}{2}}\int_{\Theta_2} \pi(\theta)\exp\Big(-\alpha\Big(\mathcal{L}_{\rm fcv}(\wX^{(m)}, \theta)+\frac{1}{n}\sum^n_{i=1}[\ell(X_i,\theta^{\dagger})-\ell(X_i,\theta^{\diamond})]\Big)\Big)\,\dd\theta\\
&= \alpha^{\frac{d}{2}}\int_{\Theta_2} \pi(\theta)\exp\Big(-\alpha\Big(\Big(\frac{1}{m}\sum^m_{j=1}\ell^{\triangle}(\wt X_j,\theta)\Big)_{+}+\frac{1}{n}\sum^n_{i=1}\ell(X_i,\theta^{\dagger})\\
&\quad-\frac{1}{n}\sum_{i=1}^n \ell(X_i,\theta^{\diamond})+\frac{1}{n}\sum_{i=1}^ng(X_i,\theta^{\dagger})\cdot(\theta-\theta^{\dagger})\Big)\Big)\,\dd\theta\\
&\lesssim \alpha^{\frac{d}{2}}\int_{\Theta_2} \pi(\theta)\exp\Big(- \frac{\alpha \delta'''}{2}\|\theta - \theta^{\diamond}\|^2_2 
+\alpha \delta'''\frac{\log n}{n}\| \theta -\theta^{\diamond} \|_2^{2\beta_1}
+\alpha\delta'''\Big( \frac{\log n}{n} \Big)^{\beta_1+1}\Big)\,\dd\theta\\
&\lesssim \alpha^{\frac{d}{2}}\int_{\Theta_2}\pi(\theta)\exp\Big(-\frac{\alpha \delta'''}{4}\|\theta-\theta^{\diamond}\|_2^2\Big)\,\dd\theta\\     
&\lesssim \frac{1}{m^2}.
\end{aligned}
\end{equation}
In addition, \eqref{estimate 1.2} gives
\begin{equation}\label{step 2.2.2} 
\begin{aligned}
&\alpha^{\frac{d}{2}}\int_{\Theta_2} \pi(\theta)\exp\Big(-\alpha\Big( \frac{1}{n}\sum_{i=1}^n \ell(X_i,\theta)-\frac{1}{n}\sum_{i=1}^n \ell(X_i,\theta^{\diamond})\Big)\Big)\,\dd\theta
\\&\leq \alpha^{\frac{d}{2}}\int_{\Theta_2}\pi(\theta)\exp\Big(-\alpha M_1\|\theta - \theta^{\diamond}\|^2_2 
     +\alpha M_1\frac{\log n}{n}\| \theta - \theta^{\diamond} \|_2^{2\beta_1}
     +\alpha M_1\Big( \frac{\log n}{n} \Big)^{\beta_1+1}\Big)\,\dd\theta\\&\lesssim\frac{1}{m^2}
\end{aligned}
\end{equation}
\paragraph{3. \underline{On $\Theta_3$.}}
 We next consider $\theta\in\Theta_3$ and $\wX^{(m)}\in\m A$. By \eqref{estimate 1.7} and the second assertion of Lemma~\ref{lemmaprobA2}, there exists a constant $K_2>0$ such that
 \begin{equation}\label{estimate 2.5}
\begin{aligned}
&\Big(\frac{1}{m}\sum^m_{j=1}\ell^{\triangle}(\wt X_j,\theta)\Big)_{+}+\frac{1}{n}\sum^n_{i=1}\ell(X_i,\theta^{\dagger})-\frac{1}{n}\sum_{i=1}^n \ell(X_i,\theta^{\diamond})+\frac{1}{n}\sum_{i=1}^ng(X_i,\theta^{\dagger})\cdot(\theta-\theta^{\dagger})\\
&\geq\frac{1}{m}\sum_{j=1}^m \big[\ell(\wt X_j,\theta)-\ell(\wt X_j,\theta^{\dagger})-g(\wt X_j,\theta^{\dagger})\cdot(\theta-\theta^{\dagger})\big]\\
&+\frac{1}{n}\sum_{i=1}^n \big[\ell(X_i,\theta^{\dagger})-\ell(X_i,\theta^{\diamond})+g(X_i,\theta^{\dagger})\cdot(\theta-\theta^{\dagger})\big]\\
&\geq\frac{1}{n}\sum_{i=1}^n \ell(X_i,\theta)-\frac{1}{n}\sum_{i=1}^n \ell(X_i,\theta^{\diamond})\\
&\quad-\Big|\frac{1}{m}\sum_{j=1}^m \ell(\wt X_j,\theta)-\frac{1}{m}\sum_{j=1}^m \ell(\wt X_j,\theta^{\dagger})-\frac{1}{m}\sum^m_{j=1}g(\wt X_j,\theta^{\dagger})\cdot(\theta-\theta^{\dagger})\\&\qquad-\frac{1}{n}\sum_{i=1}^n \ell(X_i,\theta)+\frac{1}{n}\sum_{i=1}^n \ell(X_i,\theta^{\dagger})+
\frac{1}{n}\sum^n_{i=1}g(X_i,\theta^{\dagger})\cdot(\theta-\theta^{\dagger})\Big|\\
&\geq K_2,
\end{aligned}
\end{equation}
hence 
\begin{equation}\label{step 2.3.1}
\begin{aligned}
&\alpha^{\frac{d}{2}}\int_{\Theta_3} \pi(\theta)\exp\Big(-\alpha\Big(\mathcal{L}_{\rm fcv}(\wX^{(m)}, \theta)+\frac{1}{n}\sum^n_{i=1}[\ell(X_i,\theta^{\dagger})-\ell(X_i,\theta^{\diamond})]\Big)\Big)\,\dd\theta\lesssim \frac{1}{m^2}.
\end{aligned}
\end{equation}
Using \eqref{estimate 1.8}, we also obtain
\begin{equation} \label{step 2.3.2}
\begin{aligned}
\alpha^{\frac{d}{2}}\int_{\Theta_3} \pi(\theta)\exp\Big(-\alpha\Big( \frac{1}{n}\sum_{i=1}^n \ell(X_i,\theta)-\frac{1}{n}\sum_{i=1}^n \ell(X_i,\theta^{\diamond})\Big)\Big)\,\dd\theta\leq \frac{1}{m^2}.
\end{aligned}
\end{equation}
By combining the estimates \eqref{step 2.1}, \eqref{step 2.2.1}, \eqref{step 2.2.2}, \eqref{step 2.3.1}, \eqref{step 2.3.2} and in view of the lower bound for the denominator given in \eqref{lowerbound 2}, we deduce that
 \begin{equation*}
\begin{aligned}
& \int\Big|\frac{\mb{E}_{\wX^{(m)} \sim  \mu_n^{\otimes m} }\Big[\pi(\theta) \exp\Big(-\alpha\Big(\mathcal{L}_{\rm fcv}(\wX^{(m)}, \theta)+\frac{1}{n}\sum^n_{i=1}[\ell(X_i,\theta^{\dagger})-\ell(X_i,\theta^{\diamond})]\Big)\Big)\cdot\bold{1}(\wX^{(m)}\in \m A)\Big]}{\mb{E}_{\wX^{(m)} \sim  \mu_n^{\otimes m} }\Big[\int\pi(\theta) \exp\Big(-\alpha\Big(\mathcal{L}_{\rm fcv}(\wX^{(m)}, \theta)+\frac{1}{n}\sum^n_{i=1}[\ell(X_i,\theta^{\dagger})-\ell(X_i,\theta^{\diamond})]\Big)\Big)\cdot\bold{1}(\wX^{(m)}\in \m A)\,\dd\theta\Big]}\\
&-\frac{\mb{E}_{\wX^{(m)} \sim  \mu_n^{\otimes m} }\Big[\pi(\theta)\exp\Big(-\alpha\Big( \frac{1}{n}\sum_{i=1}^n \ell(X_i,\theta)-\frac{1}{n}\sum_{i=1}^n \ell(X_i,\theta^{\diamond})\Big)\Big)\cdot\bold{1}(\wX^{(m)}\in \m A)\Big]}{\mb{E}_{\wX^{(m)} \sim  \mu_n^{\otimes m} }\Big[\int\pi(\theta)\exp\Big(-\alpha\Big( \frac{1}{n}\sum_{i=1}^n \ell(X_i,\theta)-\frac{1}{n}\sum_{i=1}^n \ell(X_i,\theta^{\diamond})\Big)\Big)\cdot\bold{1}(\wX^{(m)}\in \m A)\,\dd\theta\Big]}\Big|\,\dd\theta\\
&\lesssim \alpha^{\frac{d}{2}}\cdot \int\Big|\mb{E}_{\wX^{(m)} \sim  \mu_n^{\otimes m} }\Big[\pi(\theta) \exp\Big(-\alpha\Big(\mathcal{L}_{\rm fcv}(\wX^{(m)}, \theta)+\frac{1}{n}\sum^n_{i=1}[\ell(X_i,\theta^{\dagger})-\ell(X_i,\theta^{\diamond})]\Big)\Big)\cdot\bold{1}(\wX^{(m)}\in \m A)\\
&\quad- \pi(\theta)\exp\Big(-\alpha\Big( \frac{1}{n}\sum_{i=1}^n \ell(X_i,\theta)-\frac{1}{n}\sum_{i=1}^n \ell(X_i,\theta^{\diamond})\Big)\Big)\cdot\bold{1}(\wX^{(m)}\in \m A)\Big]\Big|\,\dd\theta\\
&\lesssim \alpha^{1-\beta_1}\frac{(\log m)^{\beta_1+2}}{m}+\alpha^{1-\beta_2}\frac{(\log m)^{\beta_2+3}}{m^2}.
\end{aligned}
\end{equation*}
 Thus,
\begin{equation}\label{total 2}
\begin{aligned}
& \int \big|\wt{\pi}_{\mathcal{L}_{\rm fcv},\alpha}(\theta)-\pi_{\alpha}(\theta)\big|\,\dd\theta \\
&=  \int\Big|\frac{\mb{E}_{\wX^{(m)} \sim  \mu_n^{\otimes m} }\Big[\pi(\theta)\exp\Big(-\alpha\Big(\mathcal{L}_{\rm fcv}(\wX^{(m)}, \theta)+\frac{1}{n}\sum^n_{i=1}[\ell(X_i,\theta^{\dagger})-\ell(X_i,\theta^{\diamond})]\Big)\Big) \Big]}{\mb{E}_{\wX^{(m)} \sim  \mu_n^{\otimes m} }\Big[\int\pi(\theta)\exp\Big(-\alpha\Big(\mathcal{L}_{\rm fcv}(\wX^{(m)}, \theta)+\frac{1}{n}\sum^n_{i=1}[\ell(X_i,\theta^{\dagger})-\ell(X_i,\theta^{\diamond})]\Big)\Big)\,\dd\theta\Big]}\\
&-\frac{\mb{E}_{\wX^{(m)} \sim  \mu_n^{\otimes m} }\Big[\pi(\theta)\exp\Big(-\alpha\Big(\frac{1}{n}\sum_{i=1}^n \ell(X_i,\theta)-\frac{1}{n}\sum_{i=1}^n \ell( X_i,\theta^{\diamond})\Big)\Big)\Big]}{\mb{E}_{\wX^{(m)} \sim  \mu_n^{\otimes m} }\Big[\int \pi(\theta)\exp\Big(-\alpha\Big(\frac{1}{n}\sum_{i=1}^n \ell(X_i,\theta)-\frac{1}{n}\sum_{i=1}^n \ell( X_i,\theta^{\diamond})\Big)\Big)\,\dd\theta \Big]}\Big|\,\dd\theta\\
&\lesssim \alpha^{1-\beta_1}\frac{(\log m)^{\beta_1+2}}{m}+\alpha^{1-\beta_2}\frac{(\log m)^{\beta_2+3}}{m^2}\\
&+\alpha^{\frac{d}{2}}\cdot\underset{\theta\in \Theta}{\sup}\,\mb{E}_{\wX^{(m)}\sim \mu_n^{\otimes m}}\Big[\exp\Big(-\alpha\Big(\frac{1}{n}\sum_{i=1}^n \ell(X_i,\theta)-\frac{1}{n}\sum_{i=1}^n \ell( X_i,\theta^{\diamond})\Big)\Big)\bold{1}(\wX^{(m)}\in \m A^c)\Big]\\
&+\alpha^{\frac{d}{2}}\cdot\underset{\theta\in \Theta}{\sup}\,\mb{E}_{\wX^{(m)} \sim  \mu_n^{\otimes m} }\Big[\exp\Big(-\alpha\Big(\mathcal{L}_{\rm fcv}(\wX^{(m)}, \theta)+\frac{1}{n}\sum^n_{i=1}[\ell(X_i,\theta^{\dagger})-\ell(X_i,\theta^{\diamond})]\Big)\Big)\cdot\bold{1}(\wX^{(m)}\in \m A^c)\Big].
\end{aligned}
\end{equation}  
\paragraph{4. \underline{Control of the Expectation over $\mathcal{A}^c$.}}
Similar to \eqref{A complement 1.1}, we have for any $\theta\in \Theta$,
\begin{equation}\label{A complement 2.1} 
\displaystyle
\mb{E}_{\wX^{(m)} \sim  \mu_n^{\otimes m} }\Big[\exp\Big(-\alpha\Big(\frac{1}{n}\sum_{i=1}^n \ell(X_i,\theta)-\frac{1}{n}\sum_{i=1}^n \ell( X_i,\theta^{\diamond})\Big)\Big)\cdot\bold{1}(\wX^{(m)}\in \m A^c)\Big]\lesssim \mb{P}_{\wX^{(m)} \sim  \mu_n^{\otimes m} }(\m A^c)\leq  \frac{1}{m^4\alpha^{d+1}}.
\end{equation}
For the last term in \eqref{total 2}, the third statement of Lemma~\ref{lemmaprobB1} gives
\begin{equation*}\begin{aligned}
&\Big\|\frac{1}{n} \sum_{i=1}^n g(X_i, \theta^{\diamond}) - \frac{1}{n} \sum_{i=1}^n g(X_i, \theta^{\dagger})- \mb E_{X\sim \m P^*}[g(X, \theta^{\diamond})] + \mb E_{X\sim \m P^*}[g(X, \theta^{\dagger})]\Big\|_2\\&\lesssim \sqrt{\frac{\log n}{n}} \| \theta^{\diamond} - \theta^{\dagger} \|_2^{\beta_1} + \frac{\log n}{n},
\end{aligned}
\end{equation*}
Also, \eqref{estimate 0.4} and $C_{m,n}\lesssim \frac{1}{\sqrt{\alpha}}$ imply
\[\|\theta^{\dagger} -\theta^{\diamond} \|_2\leq C_{m,n}+\|\wh\theta-\theta^{\diamond}\|_2\lesssim\frac{1}{\sqrt{\alpha}}.
\]
Recalling \eqref{estimate 0.3}
\begin{equation*}
\Big\|\frac{1}{n} \sum_{i=1}^n g(X_i, \theta^{\diamond})\Big\|_2\lesssim\Big( \frac{\log n}{n} \Big)^\frac{\beta_1+1}{2}.
\end{equation*}
by the Lipschitz continuity of $\nabla\mathcal{R}(\theta)$ and
the third statement of Lemma~\ref{lemmaprobB1}, we obtain
\begin{equation}\label{estimate 2.6}
\Big\|\frac{1}{n} \sum_{i=1}^n g(X_i, \theta^{\dagger})\Big\|_2\lesssim \|\theta^{\dagger} -\theta^{\diamond} \|_2+\sqrt{\frac{\log n}{n}} \|\theta^{\dagger} -\theta^{\diamond} \|^{\beta_1}_2+\Big( \frac{\log n}{n} \Big)^\frac{\beta_1+1}{2}\lesssim \frac{1}{\sqrt{\alpha}},
\end{equation}
where the last inequality uses $\alpha \le n$ and $\beta_1 > 0$.
Moreover,
by \eqref{estimate 1.2} and \eqref{estimate 1.3}, we have
\begin{equation*}
\Big|\frac{1}{n} \sum_{i=1}^n \ell(X_i, \theta^{\diamond}) - \frac{1}{n} \sum_{i=1}^n \ell(X_i, \theta^{\dagger})\Big|\lesssim \|\theta^{\dagger} -\theta^{\diamond} \|^2_2 +\frac{\log n}{n}\| \theta^{\dagger} - \theta^{\diamond} \|_2^{2\beta_1}+\Big( \frac{\log n}{n} \Big)^{\beta_1+1}\lesssim\frac{1}{\alpha}.
\end{equation*}
Hence, there exists a positive constant $M_5$, such that 
\begin{equation}\label{omit 2.7}
\max\Big( \alpha\Big|\frac{1}{n} \sum_{i=1}^n \ell(X_i, \theta^{\diamond}) - \frac{1}{n} \sum_{i=1}^n \ell(X_i, \theta^{\dagger})\Big|,\sqrt{\alpha}\Big\|\frac{1}{n} \sum_{i=1}^n g(X_i, \theta^{\dagger})\Big\|_2\Big)\leq M_5.
\end{equation}
\textbf{Case (i):} If $\|\theta-\theta^{\dagger}\|_2\leq \delta'\sqrt{\frac{\log n}{\alpha}}$, then
\begin{equation*}
\begin{aligned}
&\exp\Big(-\alpha\Big(\Big(\frac{1}{m}\sum^m_{j=1}\ell^{\triangle}(\wt X_j,\theta)\Big)_{+}
+\frac{1}{n}\sum^n_{i=1}\ell(X_i,\theta^{\dagger})-\frac{1}{n}\sum_{i=1}^n \ell(X_i,\theta^{\diamond})
+\frac{1}{n}\sum_{i=1}^ng(X_i,\theta^{\dagger})\cdot(\theta-\theta^{\dagger})\Big)\Big)\\
&\leq \exp\Big((1+\delta'\sqrt{\log n}) M_5\Big),
\end{aligned}
\end{equation*}
hence, it follows from \eqref{eqn:prob2} that
\begin{equation}\label{estimate 2.7.1}
\begin{aligned}
&\alpha^{\frac{d}{2}}\cdot\mb{E}_{\wX^{(m)} \sim  \mu_n^{\otimes m} }\Big[\exp\Big(-\alpha\Big(\mathcal{L}_{\rm fcv}(\wX^{(m)}, \theta)+\frac{1}{n}\sum^n_{i=1}[\ell(X_i,\theta^{\dagger})-\ell(X_i,\theta^{\diamond})]\Big)\Big)\cdot\bold{1}(\wX^{(m)}\in \m A^c)\Big]\\&\lesssim \frac{1}{m^2\alpha^{\frac{d+2}{2}}}.
\end{aligned}
\end{equation}
\par\noindent\textbf{Case (ii):} If $\|\theta-\theta^{\dagger}\|_2\geq\delta'\sqrt{\frac{\log n}{\alpha}}$, by the Cauchy--Schwarz inequality, we have
\begin{equation*}
\begin{aligned}
&\Big(\mb{E}_{\wX^{(m)} \sim  \mu_n^{\otimes m} }\Big[\exp\Big(-\alpha\Big(\mathcal{L}_{\rm fcv}(\wX^{(m)}, \theta)+\frac{1}{n}\sum^n_{i=1}[\ell(X_i,\theta^{\dagger})-\ell(X_i,\theta^{\diamond})]\Big)\Big)\cdot\bold{1}(\wX^{(m)}\in \m A^c)\Big]\Big)^2\\
&\leq \mb P_{\wX^{(m)}\sim \mu_n^{\otimes m}}(\m A^c)\cdot\mb{E}_{\wX^{(m)} \sim  \mu_n^{\otimes m} }\Big[\exp\Big(-2\alpha\Big(\mathcal{L}_{\rm fcv}(\wX^{(m)}, \theta)+\frac{1}{n}\sum^n_{i=1}[\ell(X_i,\theta^{\dagger})-\ell(X_i,\theta^{\diamond})]\Big)\Big)\Big]\\
&\leq \mb P_{\wX^{(m)}\sim \mu_n^{\otimes m}}(\m A^c)\cdot\mb{E}_{\wX^{(m)} \sim  \mu_n^{\otimes m} }\Big[\exp\Big(-2\alpha\Big(\frac{1}{m}\sum^m_{j=1}\ell^{\triangle}(\wt X_j,\theta)\Big)_{+}+2M_5+2M_5\sqrt{\alpha}\|\theta-\theta^{\dagger}\|_2\Big)\Big]\\
&\lesssim\mb P_{\wX^{(m)}\sim \mu_n^{\otimes m}}(\m A^c)\cdot\mb{E}_{\wX^{(m)} \sim  \mu_n^{\otimes m} }\Big[\exp\Big(-2\alpha\Big(\frac{1}{m}\sum^m_{j=1}\ell^{\triangle}(\wt X_j,\theta)\Big)_{+}+2M_5\sqrt{\alpha}\|\theta-\theta^{\dagger}\|_2\Big)\Big]\\
&=\mb P_{\wX^{(m)}\sim \mu_n^{\otimes m}}(\m A^c)\cdot\int_{0}^{+\infty}\mb{P}_{\wX^{(m)}\sim \mu_n^{\otimes m}}\Big(-2\alpha\Big(\frac{1}{m}\sum^m_{j=1}\ell^{\triangle}(\wt X_j,\theta)\Big)_{+}>\log t-2M_5\sqrt{\alpha}\|\theta-\theta^{\dagger}\|_2\Big)\,\dd t\\
&=\mb P_{\wX^{(m)}\sim \mu_n^{\otimes m}}(\m A^c)\cdot\int_{0}^{\exp(2M_5\sqrt{\alpha}\|\theta-\theta^{\dagger}\|_2)}\mb{P}_{\wX^{(m)}\sim \mu_n^{\otimes m}}\Big(-2\alpha\Big(\frac{1}{m}\sum^m_{j=1}\ell^{\triangle}(\wt X_j,\theta)\Big)_{+}\\
&\qquad>\log t-2M_5\sqrt{\alpha}\|\theta-\theta^{\dagger}\|_2\Big)\,\dd t.
\\&\leq\mb P_{\wX^{(m)}\sim \mu_n^{\otimes m}}(\m A^c)\cdot\int_{1}^{\exp(2M_5\sqrt{\alpha}\|\theta-\theta^{\dagger}\|_2)}\mb{P}_{\wX^{(m)}\sim \mu_n^{\otimes m}}\Big(-2\alpha\Big(\frac{1}{m}\sum^m_{j=1}\ell^{\triangle}(\wt X_j,\theta)\Big)_{+}
\\&\qquad>\log t-2M_5\sqrt{\alpha}\|\theta-\theta^{\dagger}\|_2\Big)\,\dd t
+\mb P_{\wX^{(m)}\sim \mu_n^{\otimes m}}(\m A^c)\cdot\int_{0}^{1}\,\dd t.\\&\end{aligned}
\end{equation*}
And for $1\leq t\leq \exp(2M_5\sqrt{\alpha}\|\theta-\theta^{\dagger}\|_2)$, it holds that 
\begin{equation*}
\begin{aligned}
&\mb{P}_{\wX^{(m)}\sim \mu_n^{\otimes m}}\Big(-2\alpha\cdot\Big(\frac{1}{m}\sum^m_{j=1}\ell^{\triangle}(\wt X_j,\theta)\Big)_{+}>\log t-2M_5\sqrt{\alpha}\|\theta-\theta^{\dagger}\|_2\Big)\\
&\leq\mb{P}_{\wX^{(m)}\sim \mu_n^{\otimes m}}\Big(-2\alpha\cdot\Big(\frac{1}{m}\sum^m_{j=1}\ell^{\triangle}(\wt X_j,\theta)\Big)>\log t-2M_5\sqrt{\alpha}\|\theta-\theta^{\dagger}\|_2\Big)\\&=\mb{P}_{\wX^{(m)}\sim \mu_n^{\otimes m}}\Big(-2\alpha\cdot\Big[\frac{1}{m}\sum^m_{j=1}\ell^{\triangle}(\wt X_j,\theta)-\,\mb{E}_{\wX^{(m)} \sim  \mu_n^{\otimes m} }\Big[\frac{1}{m}\sum^m_{j=1}\ell^{\triangle}(\wt X_j,\theta)\Big]\Big]\\
&\qquad>\log t-2M_5\sqrt{\alpha}\|\theta-\theta^{\dagger}\|_2+\frac{2\alpha}{n}\sum^n_{i=1}\ell^{\triangle}(X_i,\theta)\Big).
\end{aligned}
\end{equation*}
To choose $\delta'$, we first take $\delta_2$ sufficiently small so that,
for every $\theta \in \Theta_1 \cup \Theta_2$,
\begin{equation*}
\begin{aligned}
&\m R(\theta) - \m R(\theta^{\dagger})
- \nabla\m R(\theta^{\dagger})\cdot(\theta-\theta^{\dagger})\\
&\quad =
\int_0^1 (1-t)\,
(\theta-\theta^{\dagger})^T
\m H\bigl(\theta^{\dagger}+t(\theta-\theta^{\dagger})\bigr)
(\theta-\theta^{\dagger})\,dt\\
&\quad \gtrsim \|\theta-\theta^{\dagger}\|_2^2.
\end{aligned}
\end{equation*}
For $\theta \in \Theta_3$, the assumed quadratic growth condition gives
\[
\m R(\theta)-\m R(\theta^*)
\geq b\|\theta-\theta^*\|_2^2.
\]
Combining these bounds with
$\|\theta^{\dagger}-\theta^*\|_2 \lesssim \alpha^{-1/2}+\sqrt{\frac{\log n}{n}}$,
we obtain a constant $M_6>0$ such that, for all $\theta\in\Theta$,
\begin{equation}\label{monit 2.5}
\m R(\theta)-\m R(\theta^{\dagger})
-\nabla\m R(\theta^{\dagger})\cdot(\theta-\theta^{\dagger})
\geq M_6\|\theta-\theta^{\dagger}\|_2^2.
\end{equation}
Moreover, the last assertion of Lemma~\ref{lemmaprobB2} yields
\begin{equation*}
\begin{aligned}
&\bigg|
\Big[
\frac{1}{n}\sum_{i=1}^n \ell(X_i,\theta)
-\frac{1}{n}\sum_{i=1}^n \ell(X_i,\theta^{\dagger})
-\frac{1}{n}\sum_{i=1}^n
g(X_i,\theta^{\dagger})\cdot(\theta-\theta^{\dagger})
\Big]\\
&\qquad -
\Big[
\m R(\theta)-\m R(\theta^{\dagger})
-\nabla\m R(\theta^{\dagger})\cdot(\theta-\theta^{\dagger})
\Big]
\bigg|\\
&\quad \lesssim
\sqrt{\frac{\log n}{n}}\,
\|\theta-\theta^{\dagger}\|_2^{\beta_1+1}
+\frac{\log n}{n}\,
\|\theta-\theta^{\dagger}\|_2^{\beta_2+1}
+\Big(\frac{\log n}{n}\Big)^3.
\end{aligned}
\end{equation*}
Together with \eqref{monit 2.5} and $\alpha\lesssim n$, this bound allows us to choose
a sufficiently large constant $\delta'>0$ such that, whenever
$\theta\in\Theta$ satisfies
$\|\theta-\theta^{\dagger}\|_2
\geq \delta'\sqrt{\frac{\log n}{\alpha}}$, we have
\begin{equation}\label{omit 2.6}
-2M_5\sqrt{\alpha}\,\|\theta-\theta^{\dagger}\|_2
+2\alpha\left(
\frac{1}{n}\sum_{i=1}^n \ell^{\triangle}(X_i,\theta)
\right)
\gtrsim \alpha\|\theta-\theta^{\dagger}\|_2^2.
\end{equation}
In the following $C$ is a positive constant which may change from line to line. By the second assertion of Lemma~\ref{lemmaprobB2}, we have 
\begin{equation*}
\begin{aligned}
&\frac{1}{n}\sum_{i=1}^n      \Big|\ell(X_i,\theta)- \ell(X_i,\theta^{\dagger})-g(X_i,\theta^{\dagger})\cdot(\theta-\theta^{\dagger})\Big|^2\\&\lesssim\|\theta-\theta^{\dagger}\|^{2\beta_1+2}_2+\sqrt{\frac{\log n}{n}} \,\|\theta-\theta^{\dagger}\|^{\beta_1+\beta_2+2}_2+\frac{\log n}{n}\,\|\theta-\theta^{\dagger}\|^{2\beta_2+2}_2+\Big(\frac{\log n}{n}\Big)^4\\&\lesssim \|\theta-\theta^{\dagger}\|^{2\beta_1+2}_2,  
     \end{aligned}
     \end{equation*}
and
\begin{equation*}
\begin{aligned}
\Big|\ell(X_i,\theta)- \ell(X_i,\theta^{\dagger})-g(X_i,\theta^{\dagger})(\theta-\theta^{\dagger})\Big|\lesssim\|\theta-\theta^{\dagger}\|^{1+\beta_2}_2,  
     \end{aligned}
     \end{equation*}
then by Bernstein's inequality, for $t\geq 1$,  we have
\begin{equation*}
\begin{aligned} 
&\mb{P}_{\wX^{(m)}\sim \mu_n^{\otimes m}}\Bigg(-2\alpha\cdot\Big[\frac{1}{m}\sum^m_{j=1}\ell^{\triangle}(\wt X_j,\theta)-\mb{E}_{\wX^{(m)}\sim\mu_n^{\otimes m}}\Big[\frac{1}{m}\sum^m_{j=1}\ell^{\triangle}(\wt X_j,\theta)\Big]\Big]\\
&\qquad>\log t-2M_5\sqrt{\alpha}\|\theta-\theta^{\dagger}\|_2+2\alpha\cdot\Big(\frac{1}{n}\sum^n_{i=1}\ell^{\triangle}(X_i,\theta)\Big)\Bigg)\\
&\leq\mb{P}_{\wX^{(m)}\sim\mu_n^{\otimes m}}\Bigg(\Big|\frac{1}{m}\sum^m_{j=1}\ell^{\triangle}(\wt X_j,\theta)-\mb{E}_{\wX^{(m)}\sim\mu_n^{\otimes m}}\Big[\frac{1}{m}\sum^m_{j=1}\ell^{\triangle}(\wt X_j,\theta)\Big]\Big|>C\|\theta-\theta^{\dagger}\|^2_2\Bigg)\\
&=\mb{P}_{\wX^{(m)}\sim\mu_n^{\otimes m}}\Bigg(\Big|\frac{1}{m}\sum^m_{j=1}\ell^{\triangle}(\wt X_j,\theta)-\frac{1}{n}\sum^n_{i=1}\ell^{\triangle}(X_i,\theta)\Big|>C\|\theta-\theta^{\dagger}\|^2_2\Bigg)\\
&\lesssim\max\Bigg(\exp\Big(-mC\frac{\|\theta-\theta^\dagger\|_2^4}{\|\theta-\theta^\dagger\|^{2\beta_1+2}_2}\Big),\exp\Big(-mC\frac{\|\theta-\theta^\dagger\|^2_2}{\|\theta-\theta^\dagger\|^{\beta_2+1}_2}\Big)\Bigg)\\
&=\max\Big(\exp(-mC\|\theta-\theta^\dagger\|^{2-2\beta_1}_2),\exp(-mC\|\theta-\theta^\dagger\|^{1-\beta_2}_2)\Big),
\end{aligned}
\end{equation*}
hence by the choice of $\alpha$,
\begin{equation*}
\begin{aligned}
&\int_{1}^{\exp(2M_5\sqrt{\alpha}\|\theta-\theta^{\dagger}\|_2)}
\mb{P}_{\wX^{(m)}\sim\mu_n^{\otimes m}}\Big(-2\alpha\cdot\Big(\frac{1}{m}\sum^m_{j=1}\ell^{\triangle}(\wt X_j,\theta)\Big)_{+}>\log t-2M_5\sqrt{\alpha}\|\theta-\theta^{\dagger}\|_2\Big)\,\dd t\\
&\lesssim\max\Big(\exp\big(2M_5\sqrt{\alpha}\|\theta-\theta^{\dagger}\|_2-mC\|\theta-\theta^\dagger\|^{2-2\beta_1}_2\big),\\
&\qquad\exp\big(2M_5\sqrt{\alpha}\|\theta-\theta^{\dagger}\|_2-mC\|\theta-\theta^\dagger\|^{1-\beta_2}_2\big)\Big)\\&\lesssim1,
\end{aligned}
\end{equation*}
which implies, when $\|\theta-\theta^{\dagger}\|_2\geq\delta'\sqrt{\frac{\log n}{\alpha}}$, it holds that 
\begin{equation}\label{estimate 2.7.2}
\displaystyle
\alpha^{\frac{d}{2}}\cdot\mb{E}_{\wX^{(m)}\sim\mu_n^{\otimes m}}\Big[\exp\Big(-\alpha\Big(\mathcal{L}_{\rm fcv}(\wX^{(m)}, \theta)+\frac{1}{n}\sum^n_{i=1}[\ell(X_i,\theta^{\dagger})-\ell(X_i,\theta^{\diamond})]\Big)\Big)\cdot\bold{1}(\wX^{(m)}\in\m A^c)\Big]\lesssim\frac{1}{m^2}.
\end{equation}
Hence,
\begin{equation}\label{A complement 2.2}
\displaystyle
\alpha^{\frac{d}{2}}\cdot\sup_{\theta\in\Theta}\,\mb{E}_{\wX^{(m)}\sim\mu_n^{\otimes m}}\Big[\exp\Big(-\alpha\Big(\mathcal{L}_{\rm fcv}(\wX^{(m)}, \theta)+\frac{1}{n}\sum^n_{i=1}[\ell(X_i,\theta^{\dagger})-\ell(X_i,\theta^{\diamond})]\Big)\Big)\cdot\bold{1}(\wX^{(m)}\in\m A^c)\Big]\lesssim\frac{1}{m^2}.
\end{equation}
Combining all pieces, we can obtain 
\begin{equation*}
\begin{aligned}
    \int \big|\wt{\pi}_{\mathcal{L}_{\rm fcv},\alpha}(\theta)-\pi_\alpha(\theta)\big|\,\dd\theta&
         \lesssim \alpha^{1-\beta_1}\frac{(\log m)^{\beta_1+2}}{m}+\alpha^{1-\beta_2}\frac{(\log m)^{\beta_2+3}}{m^2}.
    \end{aligned}
\end{equation*}
\subsubsection{Proof of the Mean Deviation}
Arguing as in \eqref{step 2.1}, we obtain
\begin{equation}\label{mean step 2.1}
{\small
\begin{aligned}
& \alpha^{\frac{d}{2}}\int_{\Theta_1}\Big\|\mb{E}_{\wX^{(m)}\sim \mu_n^{\otimes m}}\Big[\sqrt{\alpha}\cdot(\theta-\theta^{\diamond})\pi(\theta)\exp\Big(-\alpha\Big(\mathcal{L}_{\rm fcv}(\wX^{(m)}, \theta)+\frac{1}{n}\sum^n_{i=1}[\ell(X_i,\theta^{\dagger})-\ell(X_i,\theta^{\diamond})]\Big)\Big)\bold{1}(\wX^{(m)}\in \m A)\Big]\\
&\qquad-\mb{E}_{\wX^{(m)}\sim \mu_n^{\otimes m}}\Big[\sqrt{\alpha}\cdot(\theta-\theta^{\diamond})\pi(\theta)\exp\Big(-\alpha\Big( \frac{1}{n}\sum_{i=1}^n \ell(X_i,\theta)-\frac{1}{n}\sum_{i=1}^n \ell(X_i,\theta^{\diamond})\Big)\Big)\bold{1}(\wX^{(m)}\in \m A)\Big]\Big\|_2\,\dd\theta\\
&\lesssim \alpha^{1-\beta_1}\frac{(\log m)^{\beta_1+2}}{m}+\alpha^{1-\beta_2}\frac{(\log m)^{\beta_2+3}}{m^2}.
\end{aligned}}
\end{equation}

Following the arguments leading to \eqref{step 2.2.1} and
\eqref{step 2.2.2}, we obtain the following bounds for every
$\wX^{(m)}\in \m A$:
\begin{equation}\label{mean step 2.2.1}
\begin{aligned}
&\alpha^{\frac{d}{2}}\int_{\Theta_2}\Big\|\sqrt{\alpha}\cdot\big(\theta-\theta^{\diamond}\big) \pi(\theta)\exp\Big(-\alpha\cdot\Big( \mathcal{L}_{\rm fcv}(\wX^{(m)}, \theta)+\frac{1}{n}\sum^n_{i=1}[\ell(X_i,\theta^{\dagger})-\ell(X_i,\theta^{\diamond})]\Big)\Big)\Big\|_2\,\dd\theta\lesssim\frac{1}{m^2},
\end{aligned}
\end{equation}
and
\begin{equation}\label{mean step 2.2.2} 
\begin{aligned}       
&\alpha^{\frac{d}{2}} \int_{\Theta_2}\Big\|\sqrt{\alpha}\cdot\big(\theta-\theta^{\diamond}\big) \pi(\theta)\exp\Big(-\alpha\cdot\Big( \frac{1}{n}\sum_{i=1}^n \ell(X_i,\theta)-\frac{1}{n}\sum_{i=1}^n \ell(X_i,\theta^{\diamond})\Big)\Big)\Big\|_2\,\dd\theta\lesssim\frac{1}{m^2}.
\end{aligned}
\end{equation}
Similarly, the arguments used to establish \eqref{step 2.3.1}
and \eqref{step 2.3.2} yield the following bounds for every
$\wX^{(m)}\in \m A$:
\begin{equation}\label{mean step 2.3.1} 
\begin{aligned}
&\alpha^{\frac{d}{2}}\int_{\Theta_3}\Big\|\sqrt{\alpha}\cdot\big(\theta-\theta^{\diamond}\big) \pi(\theta)\exp\Big(-\alpha\Big(\mathcal{L}_{\rm fcv}(\wX^{(m)}, \theta)+\frac{1}{n}\sum^n_{i=1}[\ell(X_i,\theta^{\dagger})-\ell(X_i,\theta^{\diamond})]\Big)\Big)\Big\|_2\,\dd\theta\lesssim \frac{1}{m^2},
\end{aligned}
\end{equation}
and
\begin{equation}\label{mean step 2.3.2}
\begin{aligned}       
&\alpha^{\frac{d}{2}} \int_{\Theta_3}\Big\|\sqrt{\alpha}\cdot\big(\theta-\theta^{\diamond}\big) \pi(\theta)\exp\Big(-\alpha\Big( \frac{1}{n}\sum_{i=1}^n \ell(X_i,\theta)-\frac{1}{n}\sum_{i=1}^n \ell(X_i,\theta^{\diamond})\Big)\Big)\Big\|_2\,\dd\theta\lesssim\frac{1}{m^2}.
\end{aligned}
\end{equation}
Proceeding as in \eqref{omit 1.1}, we combine the estimates in
\eqref{mean step 2.1}, \eqref{mean step 2.2.1},
\eqref{mean step 2.2.2}, \eqref{mean step 2.3.1},
and \eqref{mean step 2.3.2} with those in
\eqref{A complement 2.1} and \eqref{A complement 2.2} to obtain
\begin{equation*}
\begin{aligned}
&\alpha^{\frac{d}{2}}\int_{\Theta}\Big\|\mb{E}_{\wX^{(m)}\sim\mu_n^{\otimes m}}\Big[(\theta-\theta^{\diamond})\pi(\theta)\exp\Big(-\alpha\big(\mathcal{L}_{\rm fcv}(\wX^{(m)}, \theta)+\frac{1}{n}\sum^n_{i=1}[\ell(X_i,\theta^{\dagger})-\ell(X_i,\theta^{\diamond})]\big)\Big)\Big]\\
&\qquad-\mb{E}_{\wX^{(m)}\sim\mu_n^{\otimes m}}\Big[(\theta-\theta^{\diamond})\pi(\theta)\exp\Big(-\alpha\big(\frac{1}{n}\sum_{i=1}^n \ell(X_i,\theta)-\frac{1}{n}\sum_{i=1}^n \ell( X_i,\theta^{\diamond})\big)\Big)\Big]\Big\|_2\,\dd\theta\\
&\lesssim\alpha^{\frac{1-2\beta_1}{2}}\frac{(\log m)^{\beta_1+2}}{m}+\alpha^{\frac{1-2\beta_2}{2}}\frac{(\log m)^{\beta_2+3}}{m^2}.
\end{aligned}
\end{equation*}
Next, arguing as in \eqref{omit 1.2} and using the denominator
lower bound in \eqref{lowerbound 2} together with
\eqref{omit 1.4}, we obtain
\begin{equation}
\begin{aligned}
&\Big\|\int \theta\, \wt{\pi}_{\mathcal{L}_{\rm fcv},\alpha}(\theta)\,\dd\theta-\int \theta\, \pi_{\alpha}(\theta)\,\dd\theta\Big\|_2\\
&=\Big\|\int (\theta-\theta^{\diamond})\wt{\pi}_{\mathcal{L}_{\rm fcv},\alpha}(\theta)\,\dd\theta-\int (\theta-\theta^{\diamond})\pi_{\alpha}(\theta)\,\dd\theta\Big\|_2\\
&=\Bigg\|\frac{\mb{E}_{\wX^{(m)}\sim\mu_n^{\otimes m}}\Big[\int(\theta-\theta^{\diamond})\pi(\theta)\exp\Big(-\alpha\cdot\big(\mathcal{L}_{\rm fcv}(\wX^{(m)}, \theta)+\frac{1}{n}\sum^n_{i=1}[\ell(X_i,\theta^{\dagger})-\ell(X_i,\theta^{\diamond})]\big)\Big)\,\dd\theta\Big]}{\mb{E}_{\wX^{(m)}\sim\mu_n^{\otimes m}}\Big[\int\pi(\theta)\exp\Big(-\alpha\cdot\big(\mathcal{L}_{\rm fcv}(\wX^{(m)}, \theta)+\frac{1}{n}\sum^n_{i=1}[\ell(X_i,\theta^{\dagger})-\ell(X_i,\theta^{\diamond})]\big)\Big)\,\dd\theta\Big]}\\
&\quad-\frac{\int(\theta-\theta^{\diamond})\pi(\theta)\exp\Big(-\alpha\cdot\big(\frac{1}{n}\sum_{i=1}^n\ell(X_i,\theta)-\frac{1}{n}\sum_{i=1}^n\ell(X_i,\theta^{\diamond})\big)\Big)\,\dd\theta}{\int\pi(\theta)\exp\Big(-\alpha\cdot\big(\frac{1}{n}\sum_{i=1}^n\ell(X_i,\theta)-\frac{1}{n}\sum_{i=1}^n\ell(X_i,\theta^{\diamond})\big)\Big)\,\dd\theta}\Bigg\|_2\\&\lesssim\Big\|\alpha^{\frac{d}{2}}\mb{E}_{\wX^{(m)} \sim  \mu_n^{\otimes m} }\Big[\int(\theta-\theta^{\diamond})\pi(\theta) \exp\Big(-\alpha\Big(\mathcal{L}_{\rm fcv}(\wX^{(m)}, \theta)+\frac{1}{n}\sum^n_{i=1}[\ell(X_i,\theta^{\dagger})-\ell(X_i,\theta^{\diamond})]\Big)\Big)\,\dd\theta \Big]\\&-\alpha^{\frac{d}{2}}\int(\theta-\theta^{\diamond})\pi(\theta)\exp\Big(-\alpha\Big(\frac{1}{n}\sum_{i=1}^n \ell( X_i,\theta)-\frac{1}{n}\sum_{i=1}^n \ell( X_i,\theta^{\diamond})\Big)\Big)\,\dd\theta\Big\|_2\\&+\Big\|\alpha^{\frac{d}{2}}\int(\theta-\theta^{\diamond})\pi(\theta)\exp\Big(-\alpha\Big(\frac{1}{n}\sum_{i=1}^n \ell( X_i,\theta)-\frac{1}{n}\sum_{i=1}^n \ell( X_i,\theta^{\diamond})\Big)\Big)\,\dd\theta
\Big\|_2\\&\cdot\Big|\alpha^{\frac{d}{2}}\mb{E}_{\wX^{(m)} \sim  \mu_n^{\otimes m} }\Big[\int\pi(\theta) \exp\Big(-\alpha\Big(\mathcal{L}_{\rm fcv}(\wX^{(m)}, \theta)+\frac{1}{n}\sum^n_{i=1}[\ell(X_i,\theta^{\dagger})-\ell(X_i,\theta^{\diamond})]\Big)\Big)\,\dd\theta \Big]\\&-\alpha^{\frac{d}{2}}\int\pi(\theta)\exp\Big(-\alpha\Big(\frac{1}{n}\sum_{i=1}^n \ell( X_i,\theta)-\frac{1}{n}\sum_{i=1}^n \ell( X_i,\theta^{\diamond})\Big)\Big)\,\dd\theta\Big|\\&\lesssim\alpha^{\frac{1-2\beta_1}{2}}\frac{(\log m)^{\beta_1+2}}{m}+\alpha^{\frac{1-2\beta_2}{2}}\frac{(\log m)^{\beta_2+3}}{m^2}.
\end{aligned}
\end{equation}
Finally, combining this estimate with \eqref{mean diff 1.2},
we obtain
\begin{equation*}
\begin{aligned}
&\Big\|\int \theta\,\wt{\pi}_{\mathcal{L}_{\rm fcv},\alpha}(\theta)\,\dd\theta-\wh\theta\Big\|_2\\
&\leq\Big\|\int \theta\,\wt{\pi}_{\mathcal{L}_{\rm fcv},\alpha}(\theta)\,\dd\theta-\int \theta\,\pi_{\alpha}(\theta)\,\dd\theta\Big\|_2\\
&\quad+\Big\|\int \theta\,\pi_{\alpha}(\theta)\,\dd\theta-\theta^{\diamond}\Big\|_2+\Big\|\theta^{\diamond}-\wh\theta\Big\|_2\\
&\lesssim\alpha^{\frac{1-2\beta_1}{2}}\frac{(\log m)^{\beta_1+2}}{m}+\alpha^{\frac{1-2\beta_2}{2}}\frac{(\log m)^{\beta_2+3}}{m^2}+\frac{(\log n)^\frac{3}{2}}{\alpha}+\frac{(\log n)^{\frac{\beta_1+2}{2}}}{\sqrt{n}\alpha^{\frac{\beta_1}{2}}}+\Big(\frac{\log n}{n}\Big)^{\frac{\beta_1+1}{2}}\\
&\lesssim\alpha^{\frac{1-2\beta_1}{2}}\frac{(\log m)^{\beta_1+2}}{m}+\alpha^{\frac{1-2\beta_2}{2}}\frac{(\log m)^{\beta_2+3}}{m^2}+\frac{(\log n)^\frac{3}{2}}{\alpha}+\frac{(\log n)^{\frac{\beta_1+2}{2}}}{\sqrt{n}\alpha^{\frac{\beta_1}{2}}}.
\end{aligned}
\end{equation*}

\subsection{Proof of Theorem \ref{th:3}}
\begin{proof}
We use an $s$-conductance argument to bound the mixing time of the
Markov chain. We first establish a bound on the total variation
distance between the one-step transition distributions from states
in a high-probability set whose parameter components are
sufficiently close.

\begin{lemma}\label{lemmaTVbound}
Under Assumption~$C$, for any
$(\theta,\wX^{(m)}),(\theta',\wX'^{(m)})
\in (E_1\cap E_2)\times \m A$
satisfying $\|\theta-\theta'\|_2\leq\psi$, we have
\[
\|\m T((\theta,\wX^{(m)}),\cdot)-\m T((\theta',\wX'^{(m)}),\cdot)\|_{\rm TV}\leq 1-(1-\zeta)\exp(-2\varepsilon_1)\frac{\omega}{8}.
\]
\end{lemma}
For the transition kernel $\m T$, define the $s$-conductance by
\[
\Phi_s:=\inf\bigg\{ \frac{\phi(S)}{\wt\pi_{\m L,\alpha}^{\rm joint}(S)-s}:\, s<\wt\pi_{\m L,\alpha}^{\rm joint}(S)\leq \frac{1}{2}\bigg\},\quad\mbox{for }s\in(0,\frac{1}{2}),
\]
where $\phi(S)$ denotes the ergodic flow out of a measurable set $S$:
\[
\phi(S)=\int_{S} \bigg\{\int_{S^c} \m T(\xi,\,\dd y)\bigg\}\,\wt\pi_{\m L,\alpha}^{\rm joint}(\dd\xi).
\]
The following lemma converts the  bound in
Lemma~\ref{lemmaTVbound} into a lower bound on $\Phi_s$.
\begin{lemma}\label{lemmaconduc}
Suppose there exists a constant $\omega_1>0$ such that, for any
$(\theta,\wX^{(m)}),(\theta',\wX'^{(m)})
\in (E_1\cap E_2)\times \m A$
satisfying $\|\theta-\theta'\|_2\leq\psi$, we have
\[
\|\m T((\theta,\wX^{(m)}),\cdot)-\m T((\theta',\wX'^{(m)}),\cdot)\|_{\rm TV}\leq 1-\omega_1.
\]
Then, under Assumption~$C$, for any $s\in(0,\frac{1}{2})$ satisfying
$s\geq\frac{8\exp(\varepsilon_1)}{\lambda\psi\wedge 1}
(\exp(\varepsilon_1)p_1+p_2)$,
the $s$-conductance of $\m T$ satisfies
\begin{equation*}
\Phi_s\geq\min\{\frac{\omega_1}{4},\frac{\omega_1\lambda\psi}{32}\exp(-2\varepsilon_1)\}.
\end{equation*}
\end{lemma}
By Assumption~C.1(a), $\nu_0^{\rm joint} = \nu_0 \otimes \mu_n^{\otimes m} \leq M_0(\pi_\alpha \otimes \mu_n^{\otimes m})$. Combining this with Assumption~C.2(a) yields
\[
\nu_0^{\rm joint}\big((E_2 \times \mathcal{A})^c\big)
\leq
M_0(\pi_\alpha \otimes \mu_n^{\otimes m})\bigl((E_2 \times \mathcal{A})^c\bigr)
\leq M_0 p_2.
\]
If $\nu_0^{\rm joint}(E_2 \times \mathcal{A}) = 0$, then
\[
\|\nu_k^{\rm joint} - \wt\pi_{\m L,\alpha}^{\rm joint}\|_{\rm TV}
\leq 1
= \nu_0^{\rm joint}\bigl((E_2 \times \mathcal{A})^c\bigr)
\leq M_0 p_2,
\]
which directly establishes the desired bound below.
From now on, assume that $\nu_0^{\rm joint}(E_2 \times \mathcal{A}) > 0$. Write $\left.\nu_0^{\rm joint}\right|_{E_2 \times \mathcal{A}}$ for the conditional initial distribution restricted to $E_2 \times \mathcal{A}$. Since $E_2 \subseteq E_3$, Assumption~C.2(b)(ii)B implies that $\left.\nu_0^{\rm joint}\right|_{E_2 \times \mathcal{A}}$ is an $\frac{M_0\exp(\varepsilon_1)}{\nu_0^{\rm joint}(E_2 \times \mathcal{A})}$-warm start with respect to $\wt\pi_{\m L,\alpha}^{\rm joint}$. Consider the decomposition
\[
\nu_k^{\rm joint}
=
\nu_0^{\rm joint}(E_2 \times \mathcal{A})
\left(\left.\nu_0^{\rm joint}\right|_{E_2 \times \mathcal{A}}\right)\m T^k
+
\nu_0^{\rm joint}\bigl((E_2 \times \mathcal{A})^c\bigr)
\left(\left.\nu_0^{\rm joint}\right|_{(E_2 \times \mathcal{A})^c}\right)\m T^k,
\]
where $\left(\left.\nu_0^{\rm joint}\right|_{E_2 \times \mathcal{A}}\right)\m T^k$ denotes the distribution after $k$ transitions of $\m T$ initialized at $\left.\nu_0^{\rm joint}\right|_{E_2 \times \mathcal{A}}$ (and analogously for the second term).

Applying the $s$-conductance bound (see, e.g., Corollary~1.6 in~\cite{lovasz1993random}), for any $k \geq 1$, the $k$-step distribution $\nu_k^{\rm joint}$ approximates the stationary distribution $\wt\pi_{\m L,\alpha}^{\rm joint}$ as follows:
\[
\begin{aligned}
&\|\nu_k^{\rm joint} - \wt\pi_{\m L,\alpha}^{\rm joint}\|_{\rm TV}\\
&\leq
\nu_0^{\rm joint}(E_2 \times \mathcal{A})
\left\|
\left(\left.\nu_0^{\rm joint}\right|_{E_2 \times \mathcal{A}}\right)\m T^k
- \wt\pi_{\m L,\alpha}^{\rm joint}
\right\|_{\rm TV}
+
\nu_0^{\rm joint}\bigl((E_2 \times \mathcal{A})^c\bigr)
\left\|
\left(\left.\nu_0^{\rm joint}\right|_{(E_2 \times \mathcal{A})^c}\right)\m T^k
- \wt\pi_{\m L,\alpha}^{\rm joint}
\right\|_{\rm TV}\\
&\leq
\nu_0^{\rm joint}(E_2 \times \mathcal{A}) \cdot
\frac{M_0\exp(\varepsilon_1)}{\nu_0^{\rm joint}(E_2 \times \mathcal{A})}
\left[
s + \left(1 - \frac{\Phi_s^2}{2}\right)^k
\right]
+ \nu_0^{\rm joint}\bigl((E_2 \times \mathcal{A})^c\bigr)\\
&\leq
M_0 p_2 +
M_0\exp(\varepsilon_1)
\left[
s + \left(1 - \frac{\Phi_s^2}{2}\right)^k
\right].
\end{aligned}
\]
Now set $s = \frac{8\exp(\varepsilon_1)}{\lambda\psi \wedge 1}(\exp(\varepsilon_1)p_1 + p_2)$. When $s < \frac{1}{2}$, substituting $\omega_1 = (1 - \zeta)\exp(-2\varepsilon_1)\frac{\omega}{8}$ with $\zeta = \frac{1}{2}$ yields
\[
\Phi_s
\geq \min\left\{
\frac{\omega}{64}\exp(-2\varepsilon_1),\,
\frac{\omega\lambda\psi}{512}\exp(-4\varepsilon_1)
\right\}
\geq
\frac{\omega}{512}\exp(-4\varepsilon_1)
\bigl((\lambda\psi) \wedge 1\bigr).
\]

Let $\nu_k$ denote the $\theta$-marginal of $\nu_k^{\rm joint}$. Combining the preceding conductance lower bound with the elementary inequality $(1 - u)^k \leq \exp(-ku)$ (for $u \in [0, 1]$) and the fact that
taking marginals cannot increase total variation distance, we obtain
\begin{equation*}
\begin{aligned}
\|\wt\pi_{\m L,\alpha} - \nu_k\|_{\rm TV}
&\leq
\|\wt\pi_{\m L,\alpha}^{\rm joint} - \nu_k^{\rm joint}\|_{\rm TV}\\
&\leq
\frac{24 M_0\exp(2\varepsilon_1)}{\lambda\psi \wedge 1}
\bigl(\exp(\varepsilon_1)p_1 + p_2\bigr)
+
M_0\exp(\varepsilon_1)
\exp\Bigl(
-C k \exp(-8\varepsilon_1) (\lambda\psi \wedge 1)^2 \omega^2
\Bigr)
\end{aligned}
\end{equation*}
for a universal constant $C > 0$, where the term $M_0 p_2$ has been absorbed into the first term.
Finally, if $s \geq \frac{1}{2}$, the first term on the right-hand side equals $3 M_0\exp(\varepsilon_1)s \geq 1$ (since $M_0 \geq 1$ and $\varepsilon_1 \geq 0$), so the same upper bound holds trivially as total variation distance is bounded above by $1$.
\end{proof}
\subsubsection{Proof of Lemma \ref{lemmaTVbound} }
\begin{proof}
  Consider any $(\theta,\wX^{(m)}),(\theta',\wX'^{(m)})\in (E_1\cap E_2)\times \m A$ satisfying $\|\theta-\theta'\|_2\leq\psi$. Then we have the following 
  \begin{equation*}
  {\small
\begin{aligned}
&\|\m T((\theta,\wX^{(m)}),\cdot)-\m T((\theta',\wX'^{(m)}),\cdot)\|_{\rm TV}\\
&\leq1-\frac{1-\zeta}{2}\int\sum_{\wX^*{}^{(m)}\in(\bX^{(n)})^{\otimes m}}A\big((\theta',\wX'^{(m)}),(\theta^{\triangle},\wX^*{}^{(m)})\big)\cdot\mu_n^{\otimes m}(\wX^*{}^{(m)})q\big(\theta^{\triangle}|\theta',\wX'^{(m)}\big)\,\dd\theta^{\triangle}\\
&-\frac{1-\zeta}{2}\int\sum_{\wX^*{}^{(m)}\in(\bX^{(n)})^{\otimes m}}A\big((\theta,\wX^{(m)}),(\theta^{\triangle},\wX^*{}^{(m)})\big)\cdot\mu_n^{\otimes m}(\wX^*{}^{(m)})q\big(\theta^{\triangle}|\theta,\wX^{(m)}\big)\,\dd\theta^{\triangle}\\
&+\frac{1-\zeta}{2}\int\sum_{\wX^*{}^{(m)}\in(\bX^{(n)})^{\otimes m}}\Big|A\big((\theta',\wX'^{(m)}),(\theta^{\triangle},\wX^*{}^{(m)})\big)\cdot\mu_n^{\otimes m}(\wX^*{}^{(m)})q\big(\theta^{\triangle}|\theta',\wX'^{(m)}\big)\\
&\qquad-A\big((\theta,\wX^{(m)}),(\theta^{\triangle},\wX^*{}^{(m)})\big)\cdot\mu_n^{\otimes m}(\wX^*{}^{(m)})q\big(\theta^{\triangle}|\theta,\wX^{(m)}\big)\Big|\,\dd\theta^{\triangle}\\
&\leq 1-(1-\zeta)\int_{E_3}\sum_{\wX^*{}^{(m)}\in\m A}\min\Big(A\big((\theta',\wX'^{(m)}),(\theta^{\triangle},\wX^*{}^{(m)})\big)\cdot\mu_n^{\otimes m}(\wX^*{}^{(m)})q\big(\theta^{\triangle}|\theta',\wX'^{(m)}\big),\\
&\qquad A\big((\theta,\wX^{(m)}),(\theta^{\triangle},\wX^*{}^{(m)})\big)\cdot\mu_n^{\otimes m}(\wX^*{}^{(m)})q\big(\theta^{\triangle}|\theta,\wX^{(m)}\big)\Big)\,\dd\theta^{\triangle}.
\end{aligned}}
\end{equation*}
By Assumption~$C.2$ , we have for $\theta\in E_1\cap E_2 \subset
E_3$ and $\wX ^{(m)}\in \m A$, 
{\small
\begin{equation*}
\big|\log\big( \wt\pi_{\m L,\alpha}^{\rm joint}\big((\theta,\wX^{(m)})\big)-\log \big(\pi_{\alpha}(\theta)\mu_n^{\otimes m}(\wX^{(m)})\big)\big|\leq \varepsilon_1,
\end{equation*}}
which implies
{\small
\begin{equation*}
    A((\theta',\wX'^{(m)}),(\theta^{\triangle},\wX^*{}^{(m)}))\geq \ov A((\theta',\wX'^{(m)}),(\theta^{\triangle},\wX^*{}^{(m)}))\exp(-
    2\varepsilon_1),
\end{equation*}}
hence
\begin{equation}\label{estimate 3.1}
{\footnotesize
\begin{aligned}
&\|\m T((\theta,\wX^{(m)}),\cdot)-\m T((\theta',\wX'^{(m)}),\cdot)\|_{\rm TV}\\
&\leq 1-(1-\zeta)\exp(-2\varepsilon_1)\int_{E_3}\sum_{\wX^*{}^{(m)}\in\m A}\min\Big(\ov A\big((\theta',\wX'^{(m)}),(\theta^{\triangle},\wX^*{}^{(m)})\big)\cdot\mu_n^{\otimes m}(\wX^*{}^{(m)})q\big(\theta^{\triangle}|\theta',\wX'^{(m)}\big),\\
&\qquad\ov A\big((\theta,\wX^{(m)}),(\theta^{\triangle},\wX^*{}^{(m)})\big)\cdot\mu_n^{\otimes m}(\wX^*{}^{(m)})q\big(\theta^{\triangle}|\theta,\wX^{(m)}\big)\Big)\,\dd\theta^{\triangle}\\
&\leq 1-(1-\zeta)\exp(-2\varepsilon_1)\int_{E_3}\sum_{\wX^*{}^{(m)}\in\m A}\min\Big(A^*\big(\theta',\theta^{\triangle}\big)\cdot\mu_n^{\otimes m}(\wX^*{}^{(m)})q^*\big(\theta^{\triangle}|\theta'\big), A^*\big(\theta,\theta^{\triangle}\big)\cdot\mu_n^{\otimes m}(\wX^*{}^{(m)})q^*\big(\theta^{\triangle}|\theta\big)\Big)\,\dd\theta^{\triangle}\\
&+(1-\zeta)\exp(-2\varepsilon_1)\int_{E_3}\sum_{\wX^*{}^{(m)}\in\m A}\Big|\ov A\big((\theta',\wX'^{(m)}),(\theta^{\triangle},\wX^*{}^{(m)})\big)\cdot\mu_n^{\otimes m}(\wX^*{}^{(m)})q\big(\theta^{\triangle}|\theta',\wX'^{(m)}\big)\\
&\qquad-A^*\big(\theta',\theta^{\triangle}\big)\cdot\mu_n^{\otimes m}(\wX^*{}^{(m)})q^*\big(\theta^{\triangle}|\theta'\big)\Big|\,\dd\theta^{\triangle}\\
&+(1-\zeta)\exp(-2\varepsilon_1)\int_{E_3}\sum_{\wX^*{}^{(m)}\in\m A}\Big|\ov A\big((\theta,\wX^{(m)}),(\theta^{\triangle},\wX^*{}^{(m)})\big)\cdot\mu_n^{\otimes m}(\wX^*{}^{(m)})q\big(\theta^{\triangle}|\theta,\wX^{(m)}\big)\\
&\qquad-A^*\big(\theta,\theta^{\triangle}\big)\cdot\mu_n^{\otimes m}(\wX^*{}^{(m)})q^*\big(\theta^{\triangle}|\theta\big)\Big|\,\dd\theta^{\triangle}.
\end{aligned}}
\end{equation}
To estimate the last two terms, we introduce the following lemma.
\begin{lemma}\label{strange}  
Under Assumption~$C.2$, for $\theta \in E_3$ and $\wX^{(m)}\in\m A$, we have 
\[
\int_{E_3} \big| q(\theta' \mid \theta, \wX^{(m)}) - q^{*}(\theta' \mid \theta) \big| \, \mathrm{d}\theta' \leq \frac{\omega}{8},
\]
and
\[
\int_{E_3}\sum_{\wX^*\in\m A}\Big|A^*(\theta,\theta')-\ov A\big((\theta,\wX^{(m)}),(\theta',\wX^*{}^{(m)})\big)\Big|q(\dd\theta'|\theta,\wX^{(m)})\mu_n^{\otimes m}(\wX^*)\leq\frac{\omega}{8}.
\]
\end{lemma}
Then for $\theta \neq \theta'$, we have
\begin{equation*}
\begin{aligned}
&\|\m T^*(\theta,\cdot)-\m T^*(\theta',\cdot)\|_{\rm TV}\\&= 1-\int\sum_{\wX^*{}^{(m)}\in(\bX^{(n)})^{\otimes m}
}\min\Big(A^*\big(\theta',\theta^{\triangle}\big)\mu_n^{\otimes m}(\wX^*{}^{(m)})q^*\big(\theta^{\triangle}|\theta'\big),A^*\big(\theta,\theta^{\triangle}\big)\mu_n^{\otimes m}(\wX^*{}^{(m)})q^*\big(\theta^{\triangle}|\theta\big)\Big)\,\dd\theta^{\triangle}\\
&=1-\int_{E_3}\sum_{\wX^*{}^{(m)}\in\m A}\min\Big(A^*\big(\theta',\theta^{\triangle}\big)\cdot\mu_n^{\otimes m}(\wX^*{}^{(m)})q^*\big(\theta^{\triangle}|\theta'\big),A^*\big(\theta,\theta^{\triangle}\big)\cdot\mu_n^{\otimes m}(\wX^*{}^{(m)})q^*\big(\theta^{\triangle}|\theta\big)\Big)\,\dd\theta^{\triangle}\\&-\int\sum_{\wX^*{}^{(m)}}\bold{1}((E_3\times\m A)^c)\min\Big(A^*\big(\theta',\theta^{\triangle}\big)\cdot\mu_n^{\otimes m}(\wX^*{}^{(m)})q^*\big(\theta^{\triangle}|\theta'\big), A^*\big(\theta,\theta^{\triangle}\big)\cdot\mu_n^{\otimes m}(\wX^*{}^{(m)})q^*\big(\theta^{\triangle}|\theta\big)\Big)\,\dd\theta^{\triangle}\\
&\geq 1-\int_{E_3}\sum_{\wX^*{}^{(m)}\in\m A}\min\Big(A^*\big(\theta',\theta^{\triangle}\big)\cdot\mu_n^{\otimes m}(\wX^*{}^{(m)})q^*\big(\theta^{\triangle}|\theta'\big),A^*\big(\theta,\theta^{\triangle}\big)\cdot\mu_n^{\otimes m}(\wX^*{}^{(m)})q^*\big(\theta^{\triangle}|\theta\big)\Big)\,\dd\theta^{\triangle}\\
&\quad-\int\sum_{\wX^*{}^{(m)}}\bold{1}((E_3\times\m A)^c)\mu_n^{\otimes m}(\wX^*{}^{(m)})q^*\big(\theta^{\triangle}|\theta\big)\,\dd\theta^{\triangle}\\
&=1-\int_{E_3}\sum_{\wX^*{}^{(m)}\in\m A}\min\Big(A^*\big(\theta',\theta^{\triangle}\big)\cdot\mu_n^{\otimes m}(\wX^*{}^{(m)})q^*\big(\theta^{\triangle}|\theta'\big),A^*\big(\theta,\theta^{\triangle}\big)\cdot\mu_n^{\otimes m}(\wX^*{}^{(m)})q^*\big(\theta^{\triangle}|\theta\big)\Big)\,\dd\theta^{\triangle}\\
&\quad-\Big(1-\int_{E_3}\sum_{\wX^*{}^{(m)}\in\m A}\mu_n^{\otimes m}(\wX^*{}^{(m)})q\big(\theta^{\triangle}|\theta,\wX^{(m)}\big)\,\dd\theta^{\triangle}+\int_{E_3}\sum_{\wX^*{}^{(m)}\in\m A}\mu_n^{\otimes m}(\wX^*{}^{(m)})q\big(\theta^{\triangle}|\theta,\wX^{(m)}\big)\,\dd\theta^{\triangle}\\
&\qquad-\int_{E_3}\sum_{\wX^*{}^{(m)}\in\m A}\mu_n^{\otimes m}(\wX^*{}^{(m)})q^*\big(\theta^{\triangle}|\theta\big)\,\dd\theta^{\triangle}\Big)\\
&\geq 1-\int_{E_3}\sum_{\wX^*{}^{(m)}\in\m A}\min\Big(A^*\big(\theta',\theta^{\triangle}\big)\cdot\mu_n^{\otimes m}(\wX^*{}^{(m)})q^*\big(\theta^{\triangle}|\theta'\big), A^*\big(\theta,\theta^{\triangle}\big)\cdot\mu_n^{\otimes m}(\wX^*{}^{(m)})q^*\big(\theta^{\triangle}|\theta\big)\Big)\,\dd\theta^{\triangle}\\
&\quad-\Big(\frac{\omega}{4}+\frac{\omega}{8}\Big),
\end{aligned}
\end{equation*}
where the last inequality is due to the fact that for $\theta\in E_1\cap E_2 \subset
E_3$ and $ \wX^{(m)} \in \mathcal{A}$, 
\[
\mu_n^{\otimes m}(\m A)\cdot\int_{E_3} q(\mathrm{d}\theta^{\triangle} \mid \theta, \wX^{(m)}) \geq 1 - \frac{\omega}{4} \quad 
\quad \text{and} \quad
\int_{E_3} \big| q(\theta^{\triangle} \mid \theta, \wX^{(m)}) - q^{*}(\theta^{\triangle} \mid \theta) \big| \, \mathrm{d}\theta^{\triangle} \leq \frac{\omega}{8}.
\]
Then, by the assumption that
$\|\m T^*(\theta,\cdot)-\m T^*(\theta',\cdot)\|_{\rm TV}\leq1-\omega,$
we get
\begin{equation}\label{estimate 3.2}
\begin{aligned}
    &\int_{E_3}\sum_{\wX^*{}^{(m)}\in\m A}\min\Big(A^*\big(\theta',\theta^{\triangle}\big)\cdot\mu_n^{\otimes m}(\wX^*{}^{(m)})q^*\big(\theta^{\triangle}|\theta'\big), A^*\big(\theta,\theta^{\triangle}\big)\cdot\mu_n^{\otimes m}(\wX^*{}^{(m)})q^*\big(\theta^{\triangle}|\theta\big)\Big)\,\dd\theta^{\triangle}\geq\omega-\frac{3\omega}{8}.\end{aligned}\end{equation}
When $\theta=\theta'$, choose $\bar\theta\in E_1$ such that
$0<\|\theta-\bar\theta\|_2\leq\psi$. Then we have
\[
\begin{aligned}
\int_\Theta A^*(\theta,\theta^{\triangle})q^*(\theta^{\triangle}|\theta)\,d\theta^{\triangle}
&\geq
\int_\Theta
\min\Big(
A^*(\theta,\theta^{\triangle})q^*(\theta^{\triangle}|\theta),
A^*(\bar\theta,\theta^{\triangle})q^*(\theta^{\triangle}|\bar\theta)
\Big)\,d\theta^{\triangle}\\
&=
1-\|\m T^*(\theta,\cdot)
-\m T^*(\bar\theta,\cdot)\|_{\rm TV}\\
&\geq\omega.
\end{aligned}
\]
Moreover, by Assumption C.2(b)(i) and Lemma~\ref{strange},
\[
\begin{aligned}
1-\mu_n^{\otimes m}(\mathcal A)
\int_{E_3}q^*(\theta^{\triangle}\mid\theta)\,d\theta^{\triangle}
&\leq
1-\mu_n^{\otimes m}(\mathcal A)
\int_{E_3}q(\theta^{\triangle}|\theta,\wt{\mathbf X}^{(m)})\,d\theta^{\triangle}\\&\quad+
\mu_n^{\otimes m}(\mathcal A)
\int_{E_3}
\left|
q(\theta^{\triangle}|\theta,\wt{\mathbf X}^{(m)})
-q^*(\theta^{\triangle}|\theta)
\right|\,d\theta^{\triangle}\\
&\leq
\frac{\omega}{4}+\frac{\omega}{8}
=\frac{3\omega}{8}.
\end{aligned}
\]
Using $0\leq A^*\leq1$, we therefore obtain \eqref{estimate 3.2} for $\theta=\theta'$
\[
\begin{aligned}
&\int_{E_3}\sum_{\wX^*{}^{(m)}\in\m A}\min\Big(A^*\big(\theta',\theta^{\triangle}\big)\cdot\mu_n^{\otimes m}(\wX^*{}^{(m)})q^*\big(\theta^{\triangle}|\theta'\big), A^*\big(\theta,\theta^{\triangle}\big)\cdot\mu_n^{\otimes m}(\wX^*{}^{(m)})q^*\big(\theta^{\triangle}|\theta\big)\Big)\,\dd\theta^{\triangle}\\
&=
\mu_n^{\otimes m}(\mathcal A)
\int_{E_3}A^*(\theta,\theta^{\triangle})q^*(\theta^{\triangle}|\theta)\,d\theta^{\triangle}\\
&\geq
\int_\Theta A^*(\theta,u)q^*(u\mid\theta)\,du
-
\left[
1-\mu_n^{\otimes m}(\mathcal A)
\int_{E_3}q^*(\theta^{\triangle}|\theta)\,d\theta^{\triangle}
\right]\\
&\geq
\omega-\frac{3\omega}{8}.
\end{aligned}
\]
Thus, the same restricted overlap lower bound holds when
$\theta=\theta'$.

And, by Lemma \ref{strange}, it holds
\begin{equation*}
\begin{aligned}
&\int_{E_3}\sum_{\wX^*{}^{(m)}\in\m A}\Big|\ov A\big((\theta,\wX^{(m)}),(\theta^{\triangle},\wX^*{}^{(m)})\big)\cdot\mu_n^{\otimes m}(\wX^*{}^{(m)})q\big(\theta^{\triangle}|\theta,\wX^{(m)}\big)\\& \quad-A^*\big(\theta,\theta^{\triangle}\big)\cdot\mu_n^{\otimes m}(\wX^*{}^{(m)})q^*\big(\theta^{\triangle}|\theta\big)\Big|\,\dd\theta^{\triangle}\\
&\leq\int_{E_3}\sum_{\wX^*{}^{(m)}\in\m A}\Big|\ov A\big((\theta,\wX^{(m)}),(\theta^{\triangle},\wX^*{}^{(m)})\big)-A^*\big(\theta,\theta^{\triangle}\big)\Big|\cdot\mu_n^{\otimes m}(\wX^*{}^{(m)})q\big(\theta^{\triangle}|\theta,\wX^{(m)}\big)\,\dd\theta^{\triangle}\\
&\quad+\int_{E_3}\sum_{\wX^*{}^{(m)}\in\m A}A^*(\theta,\theta^{\triangle})\Big|q\big(\theta^{\triangle}|\theta,\wX^{(m)}\big)-q^*\big(\theta^{\triangle}|\theta\big)\Big|\mu_n^{\otimes m}(\wX^*{}^{(m)})\,\dd\theta^{\triangle}\\
&\leq\frac{\omega}{4},
\end{aligned}
\end{equation*}
similarly,
\begin{equation*}
    \begin{aligned}
&\int_{E_3}\sum_{\wX^*{}^{(m)}\in\m A}\Big|\ov A\big((\theta',\wX'^{(m)}),(\theta^{\triangle},\wX^*{}^{(m)})\big)\cdot\mu_n^{\otimes m}(\wX^*{}^{(m)})q\big(\theta^{\triangle}|\theta',\wX'^{(m)}\big)\\&\quad-A^*\big(\theta',\theta^{\triangle}\big)\cdot\mu_n^{\otimes m}(\wX^*{}^{(m)})q^*\big(\theta^{\triangle}|\theta'\big)\Big|\,\dd\theta^{\triangle}\leq\frac{\omega}{4},
    \end{aligned}
\end{equation*}
 so we have 
 \begin{equation*}
\displaystyle
       \|\m T((\theta,\wX^{(m)}),\cdot)-\m T((\theta',\wX'^{(m)}),\cdot)\|_{\rm TV}\leq 1-(1-\zeta)\exp(-2\varepsilon_1)(\omega-\frac{7\omega}{8})=1-(1-\zeta)\exp(-2\varepsilon_1)\frac{\omega}{8}.
 \end{equation*}
 \end{proof}
\subsubsection{Proof of Lemma \ref{lemmaconduc}}
 \begin{proof}
     Let $s\geq \frac{8\exp(\varepsilon_1)}{\lambda\psi\wedge 1}(\exp(\varepsilon_1)p_1+p_2)$,  and let $S$ be any measurable set of $\Theta\times(\bX^{(n)})^{\otimes m}$ with $s\leq \wt\pi_{\m L,\alpha}^{\rm joint}(S)\leq \frac{1}{2}$. Define the following subsets:
\begin{equation*}
\begin{aligned}
&S_1:=\{x\in S\mid\m T(x,S^c)<\frac{\omega_1}{2}\},\\
&S_2:=\{x\in S^c\mid\m T(x,S)<\frac{\omega_1}{2}\},\\
&S_3:=(S_1\cup S_2)^c.
\end{aligned}
\end{equation*}
Then, as in~\cite{pmlr-v134-chewi21a}, if $\wt\pi_{\m L,\alpha}^{\rm joint}(S_1)\leq\frac{\wt\pi_{\m L,\alpha}^{\rm joint}(S)}{2}$ or $\wt\pi_{\m L,\alpha}^{\rm joint}(S_2)<\frac{\wt\pi_{\m L,\alpha}^{\rm joint}(S^c)}{2}$, then by stationarity of $\wt\pi_{\m L,\alpha}^{\rm joint}$ w.r.t $\m T$,
\begin{equation}\label{estimate 3.3}
\begin{aligned}
\int_S\m T(x,S^c)\wt\pi_{\m L,\alpha}^{\rm joint}(\dd x)&=\int\m T(x,S)\wt\pi_{\m L,\alpha}^{\rm joint}(\dd x)-\int_{S}\m T(x,S)\wt\pi_{\m L,\alpha}^{\rm joint}(\dd x)\\
&=\int_{S^c}\m T(x,S)\wt\pi_{\m L,\alpha}^{\rm joint}(\dd x)\\
&\geq\frac{\omega_1}{2}\cdot\max\{\wt\pi_{\m L,\alpha}^{\rm joint}(S\cap S_1^c),\wt\pi_{\m L,\alpha}^{\rm joint}(S^c\cap S_2^c)\}\\
&\geq\frac{\omega_1\cdot\wt\pi_{\m L,\alpha}^{\rm joint}(S)}{4}.
\end{aligned} 
\end{equation}
When $\wt\pi_{\m L,\alpha}^{\rm joint}(S_1)\wedge\wt\pi_{\m L,\alpha}^{\rm joint}(S_2)\geq\frac{
\wt\pi_{\m L,\alpha}^{\rm joint}(S)}{2}$, consider $x=(\theta,\wX^{(m)})\in \big((E_1\cap E_2)\times \m A\big)\cap S_1$ and $z=(\theta',\wX'^{(m)})\in \Big((E_1\cap E_2)\times \m A\Big)\cap S_2$. Then
\[
\|\m T_x-\m T_z\|_{\rm TV}\geq\m T(z,S^c)-\m T(x,S^c)> 1-\omega_1,
\]
thus $\|\theta-\theta'\|_2\geq\psi$. Define
\[
S_i^\Theta := \left\{
\theta \in \Theta :
\exists\,\widetilde X^{(m)} \in \mathcal A,\ 
(\theta,\widetilde X^{(m)}) \in S_i
\right\},
\qquad i=1,2.
\]
Then
\begin{equation*}
\inf_{\theta\in E_2\cap S_1^\Theta\cap E_1,\ \theta'\in E_2\cap S_2^\Theta\cap E_1}\|\theta-\theta'\|_2\geq\psi.
\end{equation*}
  Consider sets $\ E_2\cap K\cap E_1\cap S_1^\Theta$ and $E_2 \cap K\cap E_1\cap S_2^\Theta$. Using the isoperimetric inequality of $\pi_{\alpha}|_{K}$, we obtain
\begin{equation*}
\begin{aligned}
p_1+p_2+\pi_{\alpha}(K\cap S_1^{\Theta c}\cap S_2^{\Theta c})&\geq\pi_{\alpha}((E_2\cap E_1)^c\cup(K\cap S_1^{\Theta c}\cap S_2^{\Theta c}))\\&\geq
\pi_{\alpha}((K \cap (E_1 \cap E_2)^c\bigr)
\cup
\bigl(K \cap (S_1^\Theta)^c \cap (S_2^\Theta)^c))
\\
&\geq\lambda\psi\cdot\min\{\pi_{\alpha}(E_2\cap K\cap E_1\cap S_1^\Theta),\pi_{\alpha}(E_2\cap K\cap E_1\cap S_2^\Theta)\}.
\end{aligned}
\end{equation*}
Using Assumption~$C$,
\begin{equation*}
\begin{aligned}
\wt\pi_{\m L,\alpha}^{\rm joint}(S_1)&\leq\wt\pi_{\m L,\alpha}^{\rm joint}((S_1^{\Theta}\cap E_2)\times\m A)+p_2\\
&\leq\exp(\varepsilon_1)\pi_{\alpha}(S_1^{\Theta}\cap E_2)\mu_n^{\otimes m}(\m A)+p_2\\
&\leq\exp(\varepsilon_1)\pi_{\alpha}(E_1\cap K\cap S_1^{\Theta}\cap E_2)\mu_n^{\otimes m}(\m A)+p_2+2\exp(\varepsilon_1)p_1,
\end{aligned}
\end{equation*}
\begin{equation*}
\begin{aligned}\wt\pi^{\rm joint}_{\m L,\alpha}(S_2)&\leq\wt\pi^{\rm joint}_{\m L,\alpha}((S_2^{\Theta}\cap E_2)\times\m A)+p_2\\
&\leq\exp(\varepsilon_1)\pi_{\alpha}(S_2^{\Theta}\cap E_2)\mu_n^{\otimes m}(\m A)+p_2\\
&\leq\exp(\varepsilon_1)\pi_{\alpha}(E_1\cap K\cap S_2^{\Theta}\cap E_2)\mu_n^{\otimes m}(\m A)+p_2+2\exp(\varepsilon_1)p_1,
\end{aligned}
\end{equation*}
and
\begin{equation*}
\begin{aligned}
\wt\pi_{\m L,\alpha}^{\rm joint}((S_1\cup S_2)^c)&\geq\wt\pi^{\rm joint}_{\m L,\alpha}((S_1^{\Theta c}\cap S_2^{\Theta c})\times \m A)\\
&\geq\wt\pi_{\m L,\alpha}^{\rm joint}((K\cap E_2\cap S_1^{\Theta c}\cap S_2^{\Theta c})\times\m A)\\
&\geq\exp(-\varepsilon_1)\pi_{\alpha}(E_2\cap K\cap S_1^{\Theta c}\cap S_2^{\Theta c})\mu_n^{\otimes m}(\m A)\\
&\geq\exp(-\varepsilon_1)\cdot\big(\pi_{\alpha}(K\cap S_1^{\Theta c}\cap S_2^{\Theta c})-\pi_{\alpha}(E^c_2)\big)\mu_n^{\otimes m}(\m A)\\
&\geq\exp(-\varepsilon_1)\cdot\big(\pi_{\alpha}(K\cap S_1^{\Theta c}\cap S_2^{\Theta c})-p_2\big)\mu_n^{\otimes m}(\m A).
\end{aligned}
\end{equation*}

So,
\begin{equation}\label{estimate 3.4}
\begin{aligned}
\wt\pi_{\m L,\alpha}^{\rm joint}((S_1\cup S_2)^c)&\geq\exp(-\varepsilon_1)\pi_{\alpha}(K\cap S_1^{\Theta c}\cap S_2^{\Theta c})\mu_n^{\otimes m}(\m A)-\exp(-\varepsilon_1)p_2\\
&\geq\exp(-\varepsilon_1)\lambda\psi\min\{\pi_{\alpha}(E_2\cap K\cap E_1\cap S_1^\Theta),\pi_{\alpha}(E_2\cap K\cap E_1\cap S_2^\Theta)\}\mu_n^{\otimes m}(\m A)\\
&\quad-\exp(-\varepsilon_1)p_1-2\exp(-\varepsilon_1)p_2\\
&\geq\exp(-2\varepsilon_1)\lambda\psi\min\{\wt\pi^{\rm joint}_{\m L,\alpha}(S_1),\wt\pi_{\m L,\alpha}^{\rm joint}(S_2)\}\\
&\quad-(1+2\exp(-\varepsilon_1)\lambda\psi)p_1-(2\exp(-\varepsilon_1)+\lambda\psi\exp(-2\varepsilon_1))p_2\\
&\geq\frac{1}{4}\exp(-2\varepsilon_1)\lambda\psi\min\{\wt\pi_{\m L,\alpha}^{\rm joint}(S_1),\wt\pi_{\m L,\alpha}^{\rm joint}(S_2)\}.
\end{aligned}
\end{equation}

So,
\begin{equation*}
\begin{aligned}
\int_S\m T(x,S^c)\wt\pi_{\m L,\alpha}^{\rm joint}(\dd x)&\geq\frac{1}{2}\Big(\int_S\m T(x,S^c)\wt\pi_{\m L,\alpha}^{\rm joint}(\dd x)+\int_{S^c}\m T(x,S)\wt\pi_{\m L,\alpha}^{\rm joint}(\dd x)\Big)\\
&\geq\frac{\omega_1}{4}\wt\pi_{\m L,\alpha}^{\rm joint}(S_3)\\
&\geq\frac{\omega_1\lambda\psi}{32}\exp(-2\varepsilon_1)\wt\pi_{\m L,\alpha}^{\rm joint}(S).
\end{aligned}
\end{equation*}
Combining \eqref{estimate 3.3} and \eqref{estimate 3.4},
\begin{equation*}
\Phi_s\geq\min\{\frac{\omega_1}{4},\frac{\omega_1\lambda\psi}{32}\exp(-2\varepsilon_1)\}.
\end{equation*}
  \end{proof}

\subsection{Proof of Theorem \ref{th:4}}
 We  state a result that quantifies the deviation between $\wt\pi^{\rm joint}_{\m L,\alpha}$ and $\pi_{\alpha}\mu_n^{\otimes m}$, under the two constructions of $\m L$ (the zero-order and the first-order control variates).
\begin{lemma}[Deviation between joint distributions]\label{le:1} 
Let $C_{m,n}$ be a constant depending on $n$ and $m$. Then  there exists a subset $\m B\subset \m X^{\otimes n}$ with $\m P^*{}^{\otimes n}(\m B)\geq 1-\frac{1}{n^2}$ so that, for any $\bX^{(n)}\in \m B$ satisfying $\|\wh\theta(\bX^{(n)})-\theta^\dagger(\bX^{(n)})\|_2\leq C_{m,n}$, the following hold:
 \begin{enumerate}
    \item Consider $\m L=\m L_{\rm zcv}$, under the assumptions of Theorem~\ref{th:1}, fix any $c > 0$, there exist a subset $\m{A} \subset (\bX^{(n)})^{\otimes m}$ and a parameter subset $\Theta_1 = \{\theta \in \Theta : \|\theta - \wh\theta(\bX^{(n)})\|_2 \leq \delta_1 \sqrt{\frac{\log m}{\alpha}}\}$, where $\delta_1$ is a constant depending on $c$, such that, 
\begin{enumerate}
    \item $1 - \wt{\pi}_{\mathcal{L},\alpha}^{\rm joint}(\theta \in \Theta_1, \wX^{(m)} \in \m{A}) \lesssim \frac{1}{m^c}$ and $1 - \pi_\alpha(\theta \in \Theta_1) \mu_n^{\otimes m}(\wX^{(m)} \in \m{A}) \lesssim \frac{1}{m^c}$.
    \item For any $\theta \in \Theta_1$ and $\wX^{(m)} \in \m{A}$,
    \[
    \left|\log\Big(\wt{\pi}_{\mathcal{L},\alpha}^{\rm joint}(\theta, \wX^{(m)})\Big) - \log\Big(\pi_\alpha(\theta) \mu_n^{\otimes m}(\wX^{(m)})\Big)\right| \lesssim \sqrt{\frac{\alpha (\log m)^2}{m}} + \alpha C_{m,n} \sqrt{\frac{\log m}{m}}.
    \]
\end{enumerate}
\item Consider $\m L=\m L_{\rm fcv}$, under the assumptions of Theorem~\ref{th:2}, fix any  $c > 0$, there exist a subset $\m{A} \subset  (\bX^{(n)})^{\otimes m}$ and a parameter subset $\Theta_1 = \{\theta \in \Theta : \|\theta - \wh\theta(\bX^{(n)})\|_2 \leq \delta_1 \sqrt{\frac{\log m}{\alpha}}\}$  where $\delta_1$ is a constant depending on $c$, such that, 
\begin{enumerate}
    \item $1 - \widetilde{\pi}_{\m{L},\alpha}^{\rm joint}(\theta \in \Theta_1, \wX^{(m)} \in \mathcal{A}) \lesssim \frac{1}{m^c}$ and $1 - \pi_{\alpha}(\theta \in \Theta_1) \mu_n^{\otimes m}(\wX^{(m)} \in \mathcal{A}) \lesssim \frac{1}{m^c}$.
    \item For all $\theta \in \Theta_1$ and $\wX^{(m)} \in \mathcal{A}$,
    \begin{equation*}
    \begin{aligned}
               &\left| \log\left(\widetilde{\pi}_{\m{L},\alpha}^{\rm joint}(\theta, \wX^{(m)})\right) - \log\left(\pi_{\alpha}(\theta) \mu_n^{\otimes m}(\wX^{(m)})\right) \right| \\
               &\lesssim \alpha^{\frac{1-\beta_1}{2}}\frac{(\log m)^{\frac{\beta_1+2}{2}}}{m^{\frac{1}{2}}} + \alpha^{\frac{1-\beta_2}{2}}\frac{(\log m)^{\frac{\beta_2+3}{2}}}{m}. 
    \end{aligned}
        \end{equation*}
        \end{enumerate}
        \end{enumerate}
 \end{lemma}
\begin{proof}
Let
\[
\theta^{\diamond}(\bX^{(n)})
=
\theta^*-\m H(\theta^*)^{-1}
\frac{1}{n}\sum_{i=1}^n g(X_i,\theta^*).
\]
As in the proof of Theorem~\ref{th:1}, we use the same simplified notation.
To prove the theorem, it suffices to show that Assumption~$C.2$ holds
with parameters $\varepsilon_1=\mathcal O(1)$ and
$p_2=\mathcal O(p_1)$ with probability at least $1-\frac{1}{n^2}$
over $\mathbf X^{(n)}\in\mathcal X^{\otimes n}$.
The desired conclusion then follows directly from
Theorem~\ref{th:3}.

Define
\begin{equation*}
\begin{aligned}
&\mu(\bX^{(n)},\theta)
=\theta-\frac{h}{n}\sum_{i=1}^n g(X_i,\theta),\\
&\mu_1(\wX^{(m)},\theta)
=\theta-\frac{h}{m}\sum_{j=1}^m g(\wX_j,\theta),\\
&\mu_2(\wX^{(m)},\theta)
=\theta-\bigg[
\frac{h}{m}\sum_{j=1}^m g(\wt X_j,\theta)
-\frac{h}{m}\sum_{j=1}^m g(\wt X_j,\theta^\dagger)
+\frac{h}{n}\sum_{i=1}^n g(X_i,\theta^\dagger)
\bigg].
\end{aligned}
\end{equation*}
Also, define
\begin{equation*}
\begin{aligned}
&E_2=
\Big\{\theta\in\Theta:
\|\theta-\wh\theta\|_2
\leq\delta_1\sqrt{\frac{\log m}{\alpha}}\Big\},\\
&E_3=
\Big\{\theta\in\Theta:
\|\theta-\wh\theta\|_2
\leq C_{E}\delta_1\sqrt{\frac{\log m}{\alpha}}\Big\}.
\end{aligned}
\end{equation*}
Here $C_E > 1$ is a fixed constant to be determined.
\par\noindent (1):\par
Let $\mathcal B\subset\m X^{\otimes n}$ denote the event on which
all conclusions of Lemma~\ref{lemmaprobB1} hold, and let
$\m A$ denote the event on which all conclusions of
Lemma~\ref{lemmaprobA1} and part~1 of Lemma~\ref{le:1}
hold for $\widetilde{\bX}^{(m)}$. Taking the constants in the lower bound
$\alpha \gtrsim \log m$ sufficiently large ensures that  $E_2$ and $E_3$ lie inside $\Theta$, since $\theta^*$ is an interior point and $\widehat{\theta}$ concentrates near it.
By Lemmas~\ref{lemmaprobB1}, \ref{lemmaprobA1},
and~\ref{le:1}, for any given $\gamma_1$, we may choose
$c_0$ and $\delta_1$ sufficiently large so that
\[
\m P^*{}^{\otimes n}(\m B)\geq1-\frac{1}{n^2},
\qquad
\mb P_{\wX^{(m)}\sim\mu_n^{\otimes m}}(\m A)
\geq1-\frac{1}{m^2}.
\]
Moreover, the following bounds hold:
\begin{enumerate}
\item
$1-\wt{\pi}_{\mathcal L_{\rm zcv},\alpha}^{\rm joint}
(\theta\in E_2,\wX^{(m)}\in\m A)
\lesssim \frac{1}{m^{\gamma_1}}$ and
$1-\pi_\alpha(\theta\in E_2)
\mu_n^{\otimes m}(\wX^{(m)}\in\m A)
\lesssim \frac{1}{m^{\gamma_1}}$.

\item
For any $\theta\in E_3$ and $\wX^{(m)}\in\m A$,
\[
\left|
\log\Big(
\wt{\pi}_{\mathcal L_{\rm zcv},\alpha}^{\rm joint}
(\theta,\wX^{(m)})
\Big)
-\log\Big(
\pi_\alpha(\theta)\mu_n^{\otimes m}(\wX^{(m)})
\Big)
\right|
\lesssim
\sqrt{\frac{\alpha(\log m)^2}{m}}
+\alpha C_{m,n}\sqrt{\frac{\log m}{m}}.
\]
\end{enumerate}
Consequently, there exists a positive constant $\delta_1$
such that
\begin{equation}\label{estimate 4.1}
\begin{aligned}
1-\wt\pi_{\m L_{\rm zcv},\alpha}^{\rm joint}
(\theta\in E_2,\wX^{(m)}\in\m A)
&\lesssim\frac{1}{m^{\gamma_1}}\lesssim p_1,\\
1-\pi_\alpha(\theta\in E_2)
\mu_n^{\otimes m}(\wX^{(m)}\in\m A)
&\lesssim\frac{1}{m^{\gamma_1}}\lesssim p_1.
\end{aligned}
\end{equation}
Furthermore, for $\m\wX^{(m)}\in\m A$ and $\theta\in E_3$,
the choices of $\alpha$ and $m$ imply
\begin{equation}\label{estimate 4.2}
\begin{aligned}
&\Big|
\log(\wt\pi_{\m L_{\rm zcv},\alpha}^{\rm joint}
\big(\theta,\wX^{(m)}\big))
-\log(\pi_\alpha(\theta)
\mu_n^{\otimes m}(\wX^{(m)}))
\Big|\\
&\qquad\lesssim
\sqrt{\frac{\alpha(\log m)^2}{m}}
+\alpha C_{m,n}\sqrt{\frac{\log m}{m}}
\lesssim\varepsilon_1=\mathcal O(1).
\end{aligned}
\end{equation}
For $\|\theta-\wh\theta\|_2
\lesssim\sqrt{\frac{\log m}{\alpha}}$ and
$\wX^{(m)}\in\m A$, since $h$ is bounded above,
\eqref{estimate 0.3} and \eqref{estimate 0.4},
together with the third statement of
Lemma~\ref{lemmaprobB1} and first statement of
Lemma~\ref{lemmaprobA1}, yield
\begin{equation}\label{estimate 4.3}
\begin{aligned}
&\Big\|\mu_1(\wX^{(m)},\theta)-\theta\Big\|_2
=\Big\|\frac{h}{m}\sum_{j=1}^m
g(\wt X_j,\theta)\Big\|_2\\
&\leq
\Big\|
\frac{h}{m}\sum_{j=1}^m g(\wt X_j,\theta)
-\frac{h}{n}\sum_{i=1}^n g(X_i,\theta)
\Big\|_2\\
&\quad+
\Big\|
\frac{h}{n}\sum_{i=1}^n g(X_i,\theta)
-\frac{h}{n}\sum_{i=1}^n g(X_i,\theta^\diamond)
\Big\|_2
+\Big\|
\frac{h}{n}\sum_{i=1}^n g(X_i,\theta^\diamond)
\Big\|_2\\
&\leq
\Big\|
\frac{h}{m}\sum_{j=1}^m g(\wt X_j,\theta)
-\frac{h}{n}\sum_{i=1}^n g(X_i,\theta)
\Big\|_2\\
&\quad+
\Big\|
\frac{h}{n}\sum_{i=1}^n g(X_i,\theta)
-\frac{h}{n}\sum_{i=1}^n g(X_i,\theta^\diamond)
-h\mb E_{X\sim\m P^*}[g(X,\theta)]
+h\mb E_{X\sim\m P^*}[g(X,\theta^\diamond)]
\Big\|_2\\
&\quad+
\Big\|
h\mb E_{X\sim\m P^*}[g(X,\theta)]
-h\mb E_{X\sim\m P^*}[g(X,\theta^\diamond)]
\Big\|_2
+\Big\|
\frac{h}{n}\sum_{i=1}^n g(X_i,\theta^\diamond)
\Big\|_2\\
&\lesssim
\sqrt{\frac{\log m}{m}}
+\sqrt{\frac{\log n}{n}}
\|\theta-\theta^\diamond\|_2^{\beta_1}
+\frac{\log n}{n}
+\|\theta-\theta^\diamond\|_2
+\Big(\frac{\log n}{n}\Big)^{\frac{\beta_1+1}{2}}\\
&\lesssim
\sqrt{\frac{\log m}{m}}
+\sqrt{\frac{\log n}{n}}
\|\theta-\wh\theta\|_2^{\beta_1}
+\sqrt{\frac{\log n}{n}}
\|\wh\theta-\theta^\diamond\|_2^{\beta_1}
+\frac{\log n}{n}\\
&\quad+
\|\theta-\wh\theta\|_2
+\|\wh\theta-\theta^\diamond\|_2
+\Big(\frac{\log n}{n}\Big)^{\frac{\beta_1+1}{2}}\\
&\lesssim\sqrt{\frac{\log m}{\alpha}}.
\end{aligned}
\end{equation}
The preceding bound, restricted to $E_2$, gives
$\Big\|\mu_1(\widetilde{\mathbf{X}}^{(m)}, \theta) - \theta\Big\|_2
\leq C_{\delta_1}(\delta_1+1)\sqrt{\frac{\log m}{\alpha}}$
for a fixed $C_{\delta_1}$ depending on $\delta_1$. For $\delta_1$ large enough, we choose $C_E$ large enough so that $C_E\delta_1 \geq \delta_1 + C_{\delta_1}(\delta_1+1) + 3\sqrt{2}h$.
Then, for every $\theta \in E_2$,
\begin{equation*}
 B\Big(
\mu_1(\wX^{(m)},\theta),
3\sqrt{2h}\sqrt{\frac{\log m}{\alpha}}
\Big)
\subseteq E_3.
\end{equation*}
It follows that
\begin{equation*}
\begin{aligned}
\int_{E_3^c}
q_{\rm zcv}(\dd\theta'\mid\theta,\wX^{(m)})
&\leq
\int_{ B\left(
0,3\sqrt{2h}\sqrt{\frac{\log m}{\alpha}}
\right)^c}
N\Big(0,\frac{2h}{\alpha}I_d\Big)
\,\dd\theta'
\lesssim\frac{1}{m^2}.
\end{aligned}
\end{equation*}
Consequently, since $m\gtrsim\frac{1}{\sqrt{\omega}}$,
\begin{equation}\label{estimate 4.4}
\begin{aligned}
1-
\mb P_{\wX^{(m)}\sim\mu_n^{\otimes m}}(\m A)
\int_{E_3}
q_{\rm zcv}(\dd\theta'\mid\theta,\wX^{(m)})
\lesssim\frac{1}{m^2}\lesssim\omega.
\end{aligned}
\end{equation}
It remains to verify Assumption~$C.2.b.ii.A$.
For $\wX^{(m)}\in\m A$ and $\theta,\theta'\in E_3$,
we need to bound
\begin{equation*}
\begin{aligned}
\bigg|
\log\frac{
q^*(\theta\mid\theta')
}{
q_{\rm zcv}(\theta\mid\theta',\wX^{(m)})
}
\bigg|
=
\bigg|
\frac{\alpha}{4h}
\|\theta-\mu(\bX^{(n)},\theta')\|_2^2
-\frac{\alpha}{4h}
\|\theta-\mu_1(\wX^{(m)},\theta')\|_2^2
\bigg|.
\end{aligned}
\end{equation*}
For $\theta\in E_3$ and $\wX^{(m)}\in\m A$,
the first statement of Lemma~\ref{lemmaprobA1} gives
\begin{equation}\label{estimate 4.5}
\begin{aligned}
\Big\|
\mu_1(\wX^{(m)},\theta)
-\mu(\bX^{(n)},\theta)
\Big\|_2=
\Big\|
\theta-\frac{h}{m}\sum_{j=1}^m g(\wt X_j,\theta)
-\theta+\frac{h}{n}\sum_{i=1}^n g(X_i,\theta)
\Big\|_2
\lesssim\sqrt{\frac{\log m}{m}}.
\end{aligned}
\end{equation}
By \eqref{estimate 4.3} and \eqref{estimate 4.5},
and using the upper bound on $h$, we obtain,
for $\theta,\theta'\in E_3$,
\begin{equation*}
\begin{aligned}
&\bigg|
\frac{\alpha}{4h}
\bigl[
2\theta-\mu_1(\wX^{(m)},\theta')
-\mu(\bX^{(n)},\theta')
\bigr]
\cdot
\bigl[
\mu_1(\widetilde X^{(m)},\theta')
-\mu(\boldsymbol X^{(n)},\theta')
\bigr]
\bigg|\\
&\leq
\frac{\alpha}{4h}
\bigl\|
2\theta-\mu_1(\wX^{(m)},\theta')
-\mu(\bX^{(n)},\theta')
\bigr\|_2
\cdot\bigl\|
\mu_1(\wX^{(m)},\theta')
-\mu(\bX^{(n)},\theta')
\bigr\|_2\\
&=
\frac{\alpha}{4h}
\bigg\|
\frac{h}{m}\sum_{j=1}^m g(\widetilde X_j,\theta')
-\frac{h}{n}\sum_{i=1}^n g(X_i,\theta')
\bigg\|_2
\cdot\bigg\|
\theta'-\theta
-\frac{h}{n}\sum_{i=1}^n g(X_i,\theta')
+\theta'-\theta
-\frac{h}{m}\sum_{j=1}^m g(\widetilde X_j,\theta')
\bigg\|_2\\
&\lesssim
\alpha\sqrt{\frac{\log m}{m}}
\bigg\|
\theta'-\theta
-\frac{h}{n}\sum_{i=1}^n g(X_i,\theta')
+\theta'-\theta
-\frac{h}{m}\sum_{j=1}^m g(\widetilde X_j,\theta')
\bigg\|_2\\
&\lesssim
\alpha\sqrt{\frac{\log m}{m}}
\bigg(
\|\theta'-\theta\|_2
+\bigg\|
\frac{h}{n}\sum_{i=1}^n g(X_i,\theta')
\bigg\|_2
+\sqrt{\frac{\log m}{\alpha}}
\bigg)\\
&\lesssim
\alpha\sqrt{\frac{\log m}{m}}
\bigg(
\sqrt{\frac{\log m}{\alpha}}
+\bigg\|
\frac{h}{n}\sum_{i=1}^n g(X_i,\theta')
\bigg\|_2
\bigg).
\end{aligned}
\end{equation*}
Furthermore, the first statement of
Lemma~\ref{lemmaprobA1} yields
\begin{equation}\label{estimate 4.6}
\begin{aligned}
\Big\|
\frac{h}{n}\sum_{i=1}^n g(X_i,\theta')
\Big\|_2
&\leq
\Big\|
\frac{h}{n}\sum_{i=1}^n g(X_i,\theta')
-\frac{h}{m}\sum_{j=1}^m g(\wt X_j,\theta')
\Big\|_2
+\Big\|
\frac{h}{m}\sum_{j=1}^m g(\wt X_j,\theta')
\Big\|_2\lesssim\sqrt{\frac{\log m}{\alpha}}.
\end{aligned}
\end{equation}
Combining these bounds, we conclude that, for
$\wX^{(m)}\in\m A$ and $\theta,\theta'\in E_3$,
\begin{equation}\label{estimate 4.7}
\begin{aligned}
\bigg|
\log\frac{
q^*(\theta\mid\theta')
}{
q_{\rm zcv}(\theta\mid\theta',\wX^{(m)})
}
\bigg|
\lesssim
\sqrt{\frac{\alpha(\log m)^2}{m}}
\lesssim\omega.
\end{aligned}
\end{equation}
\par\noindent (2):\par
For this part, take $\mathcal B\subset\m X^{\otimes n}$
to be the event where the conclusions of both
Lemmas~\ref{lemmaprobB1} and~\ref{lemmaprobB2} are valid.
Similarly, define $\m A$ as the event where the conclusions
of Lemmas~\ref{lemmaprobA1} and~\ref{lemmaprobA2},
along with part~2 of Lemma~\ref{le:1}, are valid for
$\widetilde{\bX}^{(m)}$.
Invoking Lemmas~\ref{lemmaprobB1}, \ref{lemmaprobA1},
and~\ref{le:1}, we can, for each  $\gamma_2>0$,
take $c_0$ and $\delta_1$ large enough to ensure
\[
\m P^*{}^{\otimes n}(\m B)\geq1-\frac{1}{n^2},
\qquad
\mb P_{\wX^{(m)}\sim\mu_n^{\otimes m}}(\m A)
\geq1-\frac{1}{m^2}.
\]
Moreover, the following bounds hold:
\begin{enumerate}
\item
$1-\wt{\pi}_{\mathcal L_{\rm fcv},\alpha}^{\rm joint}
(\theta\in E_2,\wX^{(m)}\in\m A)
\lesssim \frac{1}{m^{\gamma_2}}$ and
$1-\pi_\alpha(\theta\in E_2)
\mu_n^{\otimes m}(\wX^{(m)}\in\m A)
\lesssim \frac{1}{m^{\gamma_2}}$.

\item
Whenever $\theta\in E_3$ and $\wX^{(m)}\in\m A$,
\[
\left|
\log\Big(
\wt{\pi}_{\mathcal L_{\rm fcv},\alpha}^{\rm joint}
(\theta,\wX^{(m)})
\Big)
-\log\Big(
\pi_\alpha(\theta)\mu_n^{\otimes m}(\wX^{(m)})
\Big)
\right|
\lesssim
\alpha^{\frac{1-\beta_1}{2}}
\frac{(\log m)^{\frac{\beta_1+2}{2}}}{m^{\frac12}}
+
\alpha^{\frac{1-\beta_2}{2}}
\frac{(\log m)^{\frac{\beta_2+3}{2}}}{m}.
\]
\end{enumerate}
 Hence, there exists a positive constant $\delta_1$ such that
\begin{equation}\label{estimate 4.9}
\begin{aligned}
1-\wt\pi_{\m L_{\rm fcv},\alpha}^{\rm joint}
(\theta\in E_2,\wX^{(m)}\in\m A)
&\lesssim\frac{1}{m^{\gamma_2}}\lesssim p_1,\\
1-\pi_\alpha(\theta\in E_2)
\mu_n^{\otimes m}(\wX^{(m)}\in\m A)
&\lesssim\frac{1}{m^{\gamma_2}}\lesssim p_1.
\end{aligned}
\end{equation}
And for $\theta\in E_3$ and $\wX^{(m)}\in\m A$, the prescribed parameter choices further yield
\begin{equation}\label{estimate 4.10}
\begin{aligned}
\Big|
\log(\wt\pi_{\m L_{\rm fcv},\alpha}^{\rm joint}
\big(\theta,\wX^{(m)}\big))
-\log(\pi_\alpha(\theta)
\mu_n^{\otimes m}(\wX^{(m)}))
\Big|
&\lesssim
\alpha^{\frac{1-\beta_1}{2}}
\frac{(\log m)^{\frac{\beta_1+2}{2}}}{m^{\frac12}}
+
\alpha^{\frac{1-\beta_2}{2}}
\frac{(\log m)^{\frac{\beta_2+3}{2}}}{m}\\
&\lesssim\varepsilon_1=\mathcal O(1).
\end{aligned}
\end{equation}
In the region
$\|\theta - \widehat{\theta}\|_2 \lesssim \sqrt{\frac{\log m}{\alpha}}$,
the first statement of Lemma~\ref{lemmaprobB1}, Assumption A.1, and
$\|\widehat{\theta} - \theta^*\|_2 \lesssim \sqrt{\frac{\log n}{n}}$
give the following bound directly.
Here $m \geq n^\gamma$ and $\alpha \leq n$, so
$\frac{\log n}{n} \lesssim \frac{\log m}{\alpha}$.
\begin{equation*}
\Big\|
\frac{1}{n}\sum_{i=1}^n g(X_i,\theta)
\Big\|_2\lesssim \|\theta-\theta^* \|_2+\sqrt{\frac{\log n}{n}}
\lesssim\sqrt{\frac{\log m}{\alpha}}.
\end{equation*}
Combining this bound with
$\|\wh\theta-\theta^\dagger\|_2
\leq C_{m,n}\lesssim1/\sqrt\alpha$
and the first statement of Lemma~\ref{lemmaprobA2}
leads to
\begin{equation}\label{estimate 4.11}
\begin{aligned}
&\Big\|
\mu_2(\wX^{(m)},\theta)-\theta
\Big\|_2\\
&=
\Big\|
\frac{h}{m}\sum_{j=1}^m g(\wt X_j,\theta)
-\frac{h}{m}\sum_{j=1}^m
g\big(\wt X_j,\theta^\dagger\big)
+\frac{h}{n}\sum_{i=1}^n
g\big(X_i,\theta^\dagger\big)
\Big\|_2\\
&\lesssim
\Big\|
\frac{1}{m}\sum_{j=1}^m g(\wt X_j,\theta)
-\frac{1}{m}\sum_{j=1}^m
g\big(\wt X_j,\theta^\dagger\big)
-\frac{1}{n}\sum_{i=1}^n g(X_i,\theta)
+\frac{1}{n}\sum_{i=1}^n
g\big(X_i,\theta^\dagger\big)
\Big\|_2\\
&\quad+
\Big\|
\frac{1}{n}\sum_{i=1}^n g(X_i,\theta)
\Big\|_2\\
&\lesssim
\sqrt{\frac{\log m}{m}}
\|\theta-\theta^\dagger\|_2^{\beta_1}
+
\sqrt{\frac{\log m}{m}}
\Big(\frac{\log n}{n}\Big)^{\frac14}
\|\theta-\theta^\dagger\|_2^{\frac{\beta_1+\beta_2}{2}}\\
&\quad+
\frac{\log m}{m}
\|\theta-\theta^\dagger\|_2^{\beta_2}
+
\sqrt{\frac{\log m}{m}}\frac{\log n}{n}
+
\Big\|
\frac{1}{n}\sum_{i=1}^n g\big(X_i,\theta\big)
\Big\|_2\\
&\lesssim
\sqrt{\frac{\log m}{m}}
\Big(\frac{\log m}{\alpha}\Big)^{\frac{\beta_1}{2}}
+
\sqrt{\frac{\log m}{m}}C_{m,n}^{\beta_1}\\
&\quad+
\sqrt{\frac{\log m}{m}}
\Big(\frac{\log n}{n}\Big)^{\frac14}
\Big(\frac{\log m}{\alpha}\Big)^{\frac{\beta_1+\beta_2}{4}}
+
\sqrt{\frac{\log m}{m}}
\Big(\frac{\log n}{n}\Big)^{\frac14}
C_{m,n}^{\frac{\beta_1+\beta_2}{2}}\\
&\quad+
\frac{\log m}{m}
\Big(\frac{\log m}{\alpha}\Big)^{\frac{\beta_2}{2}}
+
\frac{\log m}{m}C_{m,n}^{\beta_2}
+
\sqrt{\frac{\log m}{m}}\frac{\log n}{n}
+
\sqrt{\frac{\log m}{\alpha}}\\
&\lesssim\sqrt{\frac{\log m}{\alpha}}.
\end{aligned}
\end{equation}
As in the zero-order case, choosing $\delta_1$ and $C_E$ large enough ensures
\begin{equation*}
 B\Big(
\mu_2(\wX^{(m)},\theta),
3\sqrt{2h}\sqrt{\frac{\log m}{\alpha}}
\Big)
\subseteq E_3
\end{equation*}
whenever $\theta\in E_2$.
The probability of a proposal leaving $E_3$ can therefore
be bounded by
\begin{equation*}
\begin{aligned}
\int_{E_3^c}
q_{\rm fcv}(\dd\theta'\mid\theta,\wX^{(m)})
&\leq
\int_{ B\left(
0,3\sqrt{2h}\sqrt{\frac{\log m}{\alpha}}
\right)^c}
N\Big(0,\frac{2h}{\alpha}I_d\Big)
\,\dd\theta'
\lesssim\frac{1}{m^2}.
\end{aligned}
\end{equation*}
The condition $m\gtrsim1/\sqrt\omega$ now yields
\begin{equation}
\begin{aligned}
1-
\mb P_{\wX^{(m)}\sim\mu_n^{\otimes m}}(\m A)
\int_{E_3}
q_{\rm fcv}(\dd\theta'\mid\theta,\wX^{(m)})
\lesssim\frac{1}{m^2}\lesssim\omega.
\end{aligned}
\end{equation}

To compare the two proposal means, apply the first statement
of Lemma~\ref{lemmaprobA2} at $\theta'$.
For $\theta,\theta'\in E_3$ and $\wX^{(m)}\in\m A$,
the resulting estimate is
\begin{equation}\label{estimate 4.12}
\begin{aligned}
&\Big\|
\mu_2(\wX^{(m)},\theta')
-\mu(\bX^{(n)},\theta')
\Big\|_2\\
&=
\Big\|
\frac{h}{m}\sum_{j=1}^m g(\wt X_j,\theta')
-\frac{h}{m}\sum_{j=1}^m
g\big(\wt X_j,\theta^\dagger\big)
-\frac{h}{n}\sum_{i=1}^n g(X_i,\theta')
+\frac{h}{n}\sum_{i=1}^n
g\big(X_i,\theta^\dagger\big)
\Big\|_2\\
&\lesssim
\sqrt{\frac{\log m}{m}}
\|\theta'-\theta^\dagger\|_2^{\beta_1}
+
\sqrt{\frac{\log m}{m}}
\Big(\frac{\log n}{n}\Big)^{\frac14}
\|\theta'-\theta^\dagger\|_2^{\frac{\beta_1+\beta_2}{2}}\\
&\quad+
\frac{\log m}{m}
\|\theta'-\theta^\dagger\|_2^{\beta_2}
+
\sqrt{\frac{\log m}{m}}\frac{\log n}{n}\\
&\lesssim
\frac{(\log m)^{\frac{\beta_1+1}{2}}}
{m^{\frac12}\alpha^{\frac{\beta_1}{2}}}
+
\frac{(\log m)^{\frac{\beta_2+2}{2}}}
{m\alpha^{\frac{\beta_2}{2}}}.
\end{aligned}
\end{equation}
Meanwhile, estimates~\eqref{estimate 4.11}
and~\eqref{estimate 4.6} give:
\begin{equation*}
\begin{aligned}
&\|2\theta-\mu_2(\wX^{(m)},\theta')
-\mu(\bX^{(n)},\theta')\|_2\\
&\leq
2\|\theta-\theta'\|_2+
\Big\|
\frac{h}{m}\sum_{j=1}^m g(\wt X_j,\theta')
-\frac{h}{m}\sum_{j=1}^m
g\big(\wt X_j,\theta^\dagger\big)
+\frac{h}{n}\sum_{i=1}^n
g\big(X_i,\theta^\dagger\big)
\Big\|_2
+
\Big\|
\frac{h}{n}\sum_{i=1}^n g(X_i,\theta')
\Big\|_2\\
&\lesssim\sqrt{\frac{\log m}{\alpha}}.
\end{aligned}
\end{equation*}
Substituting this bound and~\eqref{estimate 4.12}
into the Gaussian log-density ratio gives
\begin{equation*}
\begin{aligned}
&\Bigg|
\log\frac{
q^*(\theta\mid\theta')
}{
q_{\rm fcv}(\theta\mid\theta',\wX^{(m)})
}
\Bigg|\\
&=
\Big|
\frac{\alpha}{4h}
\|\theta-\mu(\bX^{(n)},\theta')\|_2^2
-
\frac{\alpha}{4h}
\|\theta-\mu_2(\wX^{(m)},\theta')\|_2^2
\Big|\\
&=
\Big|
\frac{\alpha}{4h}
[2\theta-\mu_2(\wX^{(m)},\theta')
-\mu(\bX^{(n)},\theta')]
\cdot
[\mu_2(\wX^{(m)},\theta')
-\mu(\bX^{(n)},\theta')]
\Big|\\
&\lesssim
\alpha^{\frac{1-\beta_1}{2}}
\frac{(\log m)^{\frac{\beta_1+2}{2}}}{m^{\frac12}}
+
\alpha^{\frac{1-\beta_2}{2}}
\frac{(\log m)^{\frac{\beta_2+3}{2}}}{m}\\
&\lesssim\omega.
\end{aligned}
\end{equation*}
\end{proof}
\subsubsection{Proof of Lemma \ref{le:1}}
\begin{proof}
\par\noindent\textbf{(1)}
We first consider the zero-order control variate.

\par\noindent\textbf{(a)}
Let $\mathcal{B}\subset \m X^{\otimes n}$ denote the event on which
all conclusions of Lemma~\ref{lemmaprobB1} hold, and let $\m A$
denote the event on which all conclusions of
Lemma~\ref{lemmaprobA1} hold for $\widetilde{\bX}^{(m)}$.
By Lemmas~\ref{lemmaprobB1} and~\ref{lemmaprobA1}, for any given $c$,
we can choose $c_0$ sufficiently large so that
$\m P^*{}^{\otimes n}(\m B)\geq1-\frac{1}{n^2}$, and
$\mb P_{\wX^{(m)}\sim\mu_n^{\otimes m}}(\m A)
\geq 1-\frac{1}{m^{2c}\alpha^{d+1}}$.
It then follows from \eqref{estimate 1.9} that
\begin{equation}\label{estimate 1.16}
\begin{aligned}
&\wt{\pi}_{\m L_{zcv},\alpha}^{\rm joint}(\theta\in\Theta,\wX^{(m)} \in \m{A}^c)\\&\lesssim\alpha^{\frac{d}{2}}\int_{\Theta}\mb{E}_{\wX^{(m)}\sim\mu_n^{\otimes m}}\Big[\exp\Big(-\alpha\cdot\Big(\frac{1}{m}\sum_{j=1}^m\ell^*(\wt X_j,\theta)-\frac{1}{n}\sum_{i=1}^n\ell^*(X_i,\theta^{\diamond})\Big)\Big)\cdot\bold{1}(\wX^{(m)} \in \m{A}^c)\Big]\,\dd\theta\\
&\lesssim\alpha^{\frac{d}{2}}\int_{0}^{m^c\alpha^{\frac{d+1}{2}}}\mb{P}_{\wX^{(m)}\sim\mu_n^{\otimes m}}\big(\bold{1}(\wX^{(m)} \in \m{A}^c)>0\big)\,\dd t\\
&\quad+\alpha^{\frac{d}{2}}\sup_{\theta\in\Theta}\int_{m^c\alpha^{\frac{d+1}{2}}}^{+\infty}\mb{P}_{\wX^{(m)}\sim\mu_n^{\otimes m}}\Big(\exp\Big(-\alpha\cdot\Big(\frac{1}{m}\sum_{j=1}^m\ell^*(\wt X_j,\theta)-\frac{1}{n}\sum_{i=1}^n\ell^*(X_i,\theta^{\diamond})\Big)\Big)>t\Big)\,\dd t\\
&\lesssim\frac{1}{m^c\sqrt{\alpha}}+\alpha^{\frac{d}{2}}\int_{m^c\alpha^{\frac{d+1}{2}}}^{+\infty}\frac{1}{t^2}\,\dd t\\
&\lesssim\frac{1}{m^c}.
\end{aligned}
\end{equation}Next, define
\begin{equation*}
\Theta_1 = \Big\{\theta\in\Theta : \|\theta-\wh\theta\|_2 \leq \delta_1\sqrt{\frac{\log m}{\alpha}}\Big\},
\end{equation*}
Since
$\|\wh\theta-\theta^{\diamond}\|_2
\lesssim\Big(\frac{\log n}{n}\Big )^{\frac{1+\beta_1}{2}}$,
choosing $\delta_1$ sufficiently large yields
\begin{equation*}
\Big\{\theta\in\Theta : \|\theta-\theta^{\diamond}\|_2 \leq \frac{\delta_1}{2}\sqrt{\frac{\log m}{\alpha}}\Big\}\subseteq\Theta_1\subseteq\Big\{\theta\in\Theta : \|\theta-\theta^{\diamond}\|_2 \leq 2\delta_1\sqrt{\frac{\log m}{\alpha}}\Big\}.
\end{equation*}
The bounds in \eqref{step 1.2.1} and \eqref{step 1.3.1} then give
\begin{equation*}
\begin{aligned}
&\wt{\pi}_{\m L_{zcv},\alpha}^{\rm joint}(\theta\in\Theta^c_1,\wX^{(m)} \in \m{A})\\&\lesssim\alpha^{\frac{d}{2}}\mb{E}_{\wX^{(m)}\sim\mu_n^{\otimes m}}\Big[\int_{\Theta_1^c}\pi(\theta)\exp\Big(-\alpha\cdot\Big(\frac{1}{m}\sum_{j=1}^m\ell^*(\wt X_j,\theta)-\frac{1}{n}\sum_{i=1}^n\ell^*(X_i,\theta^{\diamond})\Big)\Big)\,\dd\theta\cdot\bold{1}(\wX^{(m)} \in \m{A})\Big]\\
&\lesssim\alpha^{\frac{d}{2}}\int_{\|\theta-\theta^{\diamond}\|_2\geq\frac{\delta_1}{2}\sqrt{\frac{\log m}{\alpha}}}\pi(\theta)\exp\Big(-\frac{\alpha M_1}{2}\|\theta-\theta^{\diamond}\|_2^2\Big)\,\dd\theta+\alpha^{\frac{d}{2}}\exp\Big(-\frac{\alpha K_1}{2}\Big)\\
&\lesssim\frac{1}{m^c}.
\end{aligned}
\end{equation*}
Combining this bound with \eqref{estimate 1.16}, we obtain
\begin{equation*}
    1-\wt{\pi}_{\m L_{zcv},\alpha}^{\rm joint}(\theta\in\Theta_1,\wX^{(m)} \in \m{A})\lesssim\frac{1}{m^c}.
\end{equation*}
An analogous argument gives
\begin{equation*}
    1-\pi_{\alpha}(\theta\in \Theta_1)\mu_n^{\otimes m}(\wX^{(m)} \in \m{A}) \lesssim\frac{1}{m^c}.
\end{equation*}

\par\noindent\textbf{(b)}
For $\theta\in \Theta_1$ and $\wX^{(m)} \in \m{A}$,
\eqref{estimate 1.4} gives
\[\Big|\alpha\cdot\Big(\frac{1}{m}\sum_{j=1}^m\ell^*(\wt X_j,\theta)-\frac{1}{n}\sum_{i=1}^n\ell^*(X_i,\theta)\Big)\Big|\lesssim\alpha\sqrt{\frac{\log m}{m}}\Big(\sqrt{\frac{\log m}{\alpha}}+C_{m,n}\Big)\]
Moreover, the proof of Theorem~\ref{th:1} yields
\begin{equation*}
    \begin{aligned}
        &\alpha^{\frac{d}{2}}\cdot\Bigg|\mb{E}_{\wX^{(m)}\sim\mu_n^{\otimes m}}\Big[\int\pi(\theta)\exp\Big(-\alpha\cdot\Big(\frac{1}{m}\sum_{j=1}^m\ell^*(\wt X_j,\theta)-\frac{1}{n}\sum_{i=1}^n\ell^*(X_i,\theta^{\diamond})\Big)\Big)\,\dd\theta\Big]\\&\qquad-\mb{E}_{\wX^{(m)}\sim\mu_n^{\otimes m}}\Big[\int\pi(\theta)\exp\Big(-\alpha\cdot\Big(\frac{1}{n}\sum_{i=1}^n\ell^*(X_i,\theta)-\frac{1}{n}\sum_{i=1}^n\ell^*
(X_i,\theta^{\diamond})\Big)\Big)\,\dd\theta\Big]\Bigg|\\&\lesssim
   \alpha\frac{(\log m)^2}{m}+\alpha^2\frac{\log m}{m}C^2_{m,n}
    \end{aligned}
\end{equation*}
Combining this estimate with the lower bound in
\eqref{lowerbound 1} and applying the mean value theorem
to the logarithm, we obtain
\begin{equation*}
    \begin{aligned}
        &\Bigg|\log\Bigg(\frac{\mb{E}_{\wX^{(m)}\sim\mu_n^{\otimes m}}\Big[\int\pi(\theta)\exp\Big(-\alpha\cdot\Big(\frac{1}{m}\sum_{j=1}^m\ell^*(\wt X_j,\theta)-\frac{1}{n}\sum_{i=1}^n\ell^*(X_i,\theta^{\diamond})\Big)\Big)\,\dd\theta\Big]}{\mb{E}_{\wX^{(m)}\sim\mu_n^{\otimes m}}\Big[\int\pi(\theta)\exp\Big(-\alpha\cdot\Big(\frac{1}{n}\sum_{i=1}^n\ell^*(X_i,\theta)-\frac{1}{n}\sum_{i=1}^n\ell^*(X_i,\theta^{\diamond})\Big)\Big)\,\dd\theta\Big]}\Bigg)\Bigg|\\&=\Bigg|\log\Bigg(\frac{\alpha^{\frac{d}{2}}\cdot\mb{E}_{\wX^{(m)}\sim\mu_n^{\otimes m}}\Big[\int\pi(\theta)\exp\Big(-\alpha\cdot\Big(\frac{1}{m}\sum_{j=1}^m\ell^*(\wt X_j,\theta)-\frac{1}{n}\sum_{i=1}^n\ell^*(X_i,\theta^{\diamond})\Big)\Big)\,\dd\theta\Big]}{\alpha^{\frac{d}{2}}\cdot\mb{E}_{\wX^{(m)}\sim\mu_n^{\otimes m}}\Big[\int\pi(\theta)\exp\Big(-\alpha\cdot\Big(\frac{1}{n}\sum_{i=1}^n\ell^*(X_i,\theta)-\frac{1}{n}\sum_{i=1}^n\ell^*(X_i,\theta^{\diamond})\Big)\Big)\,\dd\theta\Big]}\Bigg)\Bigg|
        \\&\lesssim \alpha\frac{(\log m)^2}{m}+\alpha^2\frac{\log m}{m}C^2_{m,n}
    \end{aligned}
\end{equation*}
Combining these estimates, we obtain
\begin{equation}\label{estimate 9.1}
\begin{aligned}
&\Big|\log\left(\wt{\pi}_{\m L_{zcv},\alpha}^{\rm joint}(\theta, \wX^{(m)})\right) - \log\left(\pi_\alpha(\theta) \mu_n^{\otimes m}(\wX^{(m)})\right)\Big|\\
&\leq\Big|\alpha\cdot\Big(\frac{1}{m}\sum_{j=1}^m\ell^*(\wt X_j,\theta)-\frac{1}{n}\sum_{i=1}^n\ell^*(X_i,\theta)\Big)\Big|\\
&+\Bigg|\log\Bigg(\frac{\mb{E}_{\wX^{(m)}\sim\mu_n^{\otimes m}}\Big[\int\pi(\theta)\exp\Big(-\alpha\cdot\Big(\frac{1}{m}\sum_{j=1}^m\ell^*(\wt X_j,\theta)-\frac{1}{n}\sum_{i=1}^n\ell^*(X_i,\theta^{\diamond})\Big)\Big)\,\dd\theta\Big]}{\mb{E}_{\wX^{(m)}\sim\mu_n^{\otimes m}}\Big[\int\pi(\theta)\exp\Big(-\alpha\cdot\Big(\frac{1}{n}\sum_{i=1}^n\ell^*(X_i,\theta)-\frac{1}{n}\sum_{i=1}^n\ell^*(X_i,\theta^{\diamond})\Big)\Big)\,\dd\theta\Big]}\Bigg)\Bigg|\\
&\lesssim\alpha\sqrt{\frac{\log m}{m}}\Big(\sqrt{\frac{\log m}{\alpha}}+C_{m,n}\Big)+\alpha\frac{(\log m)^2}{m}+\alpha^2\frac{\log m}{m}C^2_{m,n}\\
&\lesssim\sqrt{\frac{\alpha(\log m)^2}{m}}+\alpha C_{m,n}\sqrt{\frac{\log m}{m}}.
\end{aligned}
\end{equation}

\par\noindent\textbf{(2)}
We next consider the first-order control variate.
The proof of part \textbf{(a)} is analogous, with
$\wt{\pi}_{\m L_{zcv},\alpha}^{\rm joint}(\theta, \wX^{(m)})$
replaced by
$\wt{\pi}_{\m L_{fcv},\alpha}^{\rm joint}(\theta, \wX^{(m)})$.
For part \textbf{(b)}, the same reasoning, together with
\eqref{estimate 2.1} and \eqref{lowerbound 2}, yields
\begin{equation*}
    \begin{aligned}
        \Big|\log\left(\wt{\pi}_{\m L_{fcv},\alpha}^{\rm joint}(\theta, \wX^{(m)})\right) - \log\left(\pi_\alpha(\theta) \mu_n^{\otimes m}(\wX^{(m)})\right)\Big|\lesssim\alpha^{\frac{1-\beta_1}{2}}\frac{(\log m)^{\frac{\beta_1+2}{2}}}{m^{\frac{1}{2}}} + \alpha^{\frac{1-\beta_2}{2}}\frac{(\log m)^{\frac{\beta_2+3}{2}}}{m}
    \end{aligned}
\end{equation*}
\end{proof}
\section{Proof of technical lemmas}\label{Proof of technical lemmas}
 \subsection{Proof of Lemma \ref{strange}}
     \begin{proof}
To prove the first statement, we apply Assumption~$C.2.b.ii.A$ to obtain
\begin{equation}\label{omit 3.1}
    \begin{aligned}
\int_{E_3} \big| q(\theta' \mid \theta, \wX^{(m)}) - q^{*}(\theta' \mid \theta) \big| \, \mathrm{d}\theta'&=\int_{E_3} \Bigg|1 - \frac{q^{*}(\theta' \mid \theta)}{ q(\theta' \mid \theta, \wX^{(m)})
} \Bigg|q(\, \mathrm{d}\theta' \mid \theta, \wX^{(m)})
\\&\leq \Big(\exp\Big(\frac{\omega}{100}\Big)-1\Big)\cdot\int_{E_3} q(\, \mathrm{d}\theta' \mid \theta, \wX^{(m)})
\\&\leq \frac{\omega}{8}.
    \end{aligned}
\end{equation}

For the second statement, it suffices to show that for  $\theta \in E_3$, $\wX^{(m)},\wX^*{}^{(m)}\in\m A$, we have
\begin{equation*}
\int_{E_3}\Big|A^*(\theta,\theta')-\ov A\big((\theta,\wX^{(m)}),(\theta',\wX^*{}^{(m)})\big)\Big|q(\dd\theta'|\theta,\wX^{(m)})\leq\frac{\omega}{8}.
\end{equation*}

To simplify notation, we introduce the following conventions:
\begin{equation*}
\begin{aligned}
&\m U(\theta')=\pi_{\alpha}(\theta')\cdot q^*(\theta|\theta'),\\
&\m V(\theta')=\pi_{\alpha}(\theta)\cdot q^*(\theta'|\theta),\\
&\m W(\theta')=\pi_{\alpha}(\theta')\cdot q\big(\theta|\theta',\wX^*{}^{(m)}),\\
&\m X(\theta')=\pi_{\alpha}(\theta)\cdot q(\theta'|\theta,\wX^{(m)}),\\
&\m Y(\theta')=q(\theta'|\theta,\wX^{(m)}).
\end{aligned}
\end{equation*}
Proceeding as in \eqref{omit 3.1}, we have the following inequalities:
\begin{equation*}
\int_{E_3}\bigg|1-\frac{\mathcal{V}(\theta')}{\mathcal{X}(\theta')}\bigg|\mathcal{Y}(\theta')\,\mathrm{d}\theta'
=\int_{E_3}\bigg|1-\frac{q^*(\theta'|\theta)}{q(\theta'|\theta,\wX^{(m)})}\bigg|q(\theta'|\theta,\wX^{(m)})\,\mathrm{d}\theta'\leq\frac{\omega}{48},
\end{equation*}
\begin{equation*}
\int_{E_3}\bigg|1-\frac{\mathcal{X}(\theta')}{\mathcal{V}(\theta')}\bigg|\mathcal{Y}(\theta')\,\mathrm{d}\theta'
=\int_{E_3}\bigg|1-\frac{q(\theta'|\theta,\widetilde{X}^{(m)})}{q^*(\theta'|\theta)}\bigg|q(\theta'|\theta,\wX^{(m)})\,\mathrm{d}\theta'\leq\frac{\omega}{48},
\end{equation*}
\begin{equation*}
\int_{E_3}\bigg|1-\frac{\mathcal{U}(\theta')}{\mathcal{W}(\theta')}\bigg|\mathcal{Y}(\theta')\,\mathrm{d}\theta'
=\int_{E_3}\bigg|1-\frac{q^*(\theta|\theta')}{q\big(\theta|\theta',\wX^*{}^{(m)}\big)}\bigg|q(\theta'|\theta,\wX^{(m)})\,\mathrm{d}\theta'\leq\frac{\omega}{48},
\end{equation*}
and
\begin{equation*}
\int_{E_3}\bigg|1-\frac{\mathcal{W}(\theta')}{\mathcal{U}(\theta')}\bigg|\mathcal{Y}(\theta')\,\mathrm{d}\theta'
=\int_{E_3}\bigg|1-\frac{q\big(\theta|\theta',\wX^*{}^{(m)}\big)}{q^*(\theta|\theta')}\bigg|q(\theta'|\theta,\wX^{(m)})\,\mathrm{d}\theta'\leq\frac{\omega}{48}.
\end{equation*}
We have:
\begin{equation*}
\begin{aligned}
&\int_{E_3}\Big|A^*(\theta,\theta')-\ov A\big((\theta,\wX^{(m)}),(\theta',\wX^*{}^{(m)})\big)\Big|q(\dd\theta'|\theta,\wX^{(m)})\\
&=\int_{E_3}\Big|\frac{\m U(\theta')}{\m V(\theta')}-\frac{\m W(\theta')}{\m X(\theta')}\Big|\cdot\mathbf{1}(\m U(\theta')\leq\m V(\theta'))\cdot\mathbf{1}(\m W(\theta')\leq\m X(\theta'))\cdot\m Y(\theta')\,\dd\theta'\\
&+\int_{E_3}\Big|\frac{\m U(\theta')}{\m V(\theta')}-1\Big|\cdot\mathbf{1}(\m U(\theta')\leq\m V(\theta'))\cdot\mathbf{1}(\m W(\theta')>\m X(\theta'))\cdot\m Y(\theta')\,\dd\theta'\\
&+\int_{E_3}\Big|1-\frac{\m W(\theta')}{\m X(\theta')}\Big|\cdot\mathbf{1}(\m U(\theta')>\m V(\theta'))\cdot\mathbf{1}(\m W(\theta')\leq\m X(\theta'))\cdot\m Y(\theta')\,\dd\theta'.
\end{aligned}
\end{equation*}

For the first term, we have:
\begin{equation*}
\begin{aligned}
&\int_{E_3}\Big|\frac{\m U(\theta')}{\m V(\theta')}-\frac{\m W(\theta')}{\m X(\theta')}\Big|\cdot\mathbf{1}(\m U(\theta')\leq\m V(\theta'))\cdot\mathbf{1}(\m W(\theta')\leq\m X(\theta'))\cdot\m Y(\theta')\,\dd\theta'\\
&=\int_{E_3}\Big|\frac{\m U(\theta')\m X(\theta')-\m V(\theta')\m W(\theta')}{\m V(\theta')\m X(\theta')}\Big|\cdot\mathbf{1}(\m U(\theta')\leq\m V(\theta'))\cdot\mathbf{1}(\m W(\theta')\leq\m X(\theta'))\cdot\m Y(\theta')\,\dd\theta'\\
&\leq\int_{E_3}\Big(\Big|\frac{\m U(\theta')\m X(\theta')-\m U(\theta')\m V(\theta')}{\m V(\theta')\m X(\theta')}\Big|+\Big|\frac{\m U(\theta')\m V(\theta')-\m V(\theta')\m W(\theta')}{\m V(\theta')\m X(\theta')}\Big|\Big)\\
&\qquad\cdot\mathbf{1}(\m U(\theta')\leq\m V(\theta'))\cdot\mathbf{1}(\m W(\theta')\leq\m X(\theta'))\cdot\m Y(\theta')\,\dd\theta'\\
&\leq\int_{E_3}\Big(\Big|\frac{\m U(\theta')\m X(\theta')-\m U(\theta')\m V(\theta')}{\m U(\theta')\m X(\theta')}\Big|+\Big|\frac{\m U(\theta')\m V(\theta')-\m V(\theta')\m W(\theta')}{\m V(\theta')\m W(\theta')}\Big|\Big)\\
&\qquad\cdot\mathbf{1}(\m U(\theta')\leq\m V(\theta'))\cdot\mathbf{1}(\m W(\theta')\leq\m X(\theta'))\cdot\m Y(\theta')\,\dd\theta'\\
&\leq\int_{E_3}\Big(\Big|1-\frac{\m V(\theta')}{\m X(\theta')}\Big|+\Big|1-\frac{\m U(\theta')}{\m W(\theta')}\Big|\Big)\cdot\m Y(\theta')\,\dd\theta'\\
&\leq\frac{\omega}{24}.
\end{aligned}
\end{equation*}
For the second term:
\begin{equation*}
\begin{aligned}
&\int_{E_3}\Big|\frac{\m U(\theta')}{\m V(\theta')}-1\Big|\cdot\mathbf{1}(\m U(\theta')\leq\m V(\theta'))\cdot\mathbf{1}(\m W(\theta')>\m X(\theta'))\cdot\m Y(\theta')\,\dd\theta'\\
&=\int_{E_3}\Big(\frac{\m V(\theta')\m W(\theta')-\m U(\theta')\m W(\theta')}{\m V(\theta')\m W(\theta')}\Big)\cdot\mathbf{1}(\m U(\theta')\leq\m V(\theta'))\cdot\mathbf{1}(\m W(\theta')>\m X(\theta'))\cdot\m Y(\theta')\,\dd\theta'\\
&\leq\int_{E_3}\Big(\frac{\m V(\theta')\m W(\theta')-\m U(\theta')\m X(\theta')}{\m V(\theta')\m W(\theta')}\Big)\cdot\mathbf{1}(\m U(\theta')\leq\m V(\theta'))\cdot\mathbf{1}(\m W(\theta')>\m X(\theta'))\cdot\m Y(\theta')\,\dd\theta'\\
&\leq\int_{E_3}\Big(\Big|\frac{\m V(\theta')\m W(\theta')-\m U(\theta')\m V(\theta')}{\m V(\theta')\m W(\theta')}\Big|+\Big|\frac{\m U(\theta')\m V(\theta')-\m U(\theta')\m X(\theta')}{\m V(\theta')\m W(\theta')}\Big|\Big)\\
&\qquad\cdot\mathbf{1}(\m U(\theta')\leq\m V(\theta'))\cdot\mathbf{1}(\m W(\theta')>\m X(\theta'))\cdot\m Y(\theta')\,\dd\theta'\\
&\leq\int_{E_3}\Big(\Big|1-\frac{\m U(\theta')}{\m W(\theta')}\Big|+\Big|\frac{\m U(\theta')\m V(\theta')-\m U(\theta')\m X(\theta')}{\m U(\theta')\m X(\theta')}\Big|\Big)\cdot\m Y(\theta')\,\dd\theta'\\
&\leq\int_{E_3}\Big(\Big|1-\frac{\m U(\theta')}{\m W(\theta')}\Big|+\Big|\frac{\m V(\theta')}{\m X(\theta')}-1\Big|\Big)\cdot\m Y(\theta')\,\dd\theta'\\
&\leq\frac{\omega}{24}.
\end{aligned}
\end{equation*}

Similarly, the third term can be upper bounded by $\frac{\omega}{24}$, hence
\begin{equation*}
\int_{E_3}\Big|A^*(\theta,\theta')-\ov A\big((\theta,\wX^{(m)}),(\theta',\wX^*{}^{(m)})\big)\Big|q(\dd\theta'|\theta,\wX^{(m)})\leq\frac{\omega}{8}.
\end{equation*}  
\end{proof}
\subsection{Proof of the Concentration Inequalities}
\subsubsection{Proof of Lemma \ref{lemmaprobB1}}
 
\paragraph{Proof of the first statement.} 
Let $g_v(X_i, \theta)$ denote the $v$-th component of $g(X_i, \theta)$ for $1 \leq v \leq d$, and define the function class $\m F_v = \{g_v(\cdot, \theta) : \theta \in \Theta\}$. For $f_1,f_2:\mathcal{X}\to \mathbb R$, we define the pseudometric
\begin{equation}\label{mertic1}
d_{n}(f_1,f_2)=\sqrt{\frac{1}{n}\sum_{i=1}^{n}(f_1(X_{i})-f_2(X_{i}))^{2}}.
\end{equation} 
From this definition, it immediately follows that:
\begin{equation*}
    d_n\big(g_v(\cdot, \theta_1) , g_v(\cdot, \theta_2\big) = \sqrt{\frac{1}{n}\sum_{i=1}^n \left(g_v(X_i, \theta_1) - g_v(X_i, \theta_2)\right)^2} \leq d_n^g(\theta_1, \theta_2).
\end{equation*}
Consequently, the covering number is bounded by:
\begin{equation*}
    \log \mathcal{N}( \m F_v, d_n,\varepsilon) \leq \log \mathcal{N}( \Theta, d_n^g,\varepsilon) \lesssim \log\left( \frac{n}{\varepsilon} \right).
\end{equation*}
By Assumption~$B.2$, we have $\sup_{x\in\mathcal{X},\theta\in \Theta }|g_v(x,\theta)|\leq L$. To control the uniform deviation, we define $Z_{n,v} = \sup_{\theta \in \Theta} \left| \frac{1}{n}\sum_{i=1}^n g_v(X_i, \theta) - \mathbb{E}[g_v(X, \theta)] \right|$. Then, applying standard symmetrization techniques alongside the Dudley entropy integral (see for example Theorem 4.19 and Theorem 5.22 of~\cite{wainwright_2019}), the expectation is bounded as follows:
\begin{equation*}
    \begin{aligned}
         \mathbb{E}_{\bX^{n}\sim(\mathcal{P}^*)^{\otimes n}        
         }[Z_{n,v}] &\lesssim \frac{1}{\sqrt{n}} \int_0^{2L} \sqrt{\log \mathcal{N}(\m F_v, d_n,\varepsilon)} \, d\varepsilon \\
    &\lesssim \frac{1}{\sqrt{n}} \int_0^{2L} \sqrt{\log\left( \frac{n}{\varepsilon} \right)} \, d\varepsilon \\
    &\lesssim \sqrt{\frac{\log n}{n}}.    \end{aligned}
\end{equation*}
To establish the high probability bound, we apply the Talagrand concentration inequality (see for example Theorem 3.27 of~\cite{wainwright_2019}) to each component. Taking the union bound over all $d$ components, we conclude that for any constant $c > 0$, there exists a sufficiently large constant $c_0 > 0$ such that:
\begin{equation*}
    \mathbb{P}_{\bX^{n}\sim \m P^*{}^{\otimes n}}    \left( \sup_{\theta \in \Theta} \bigg\| \frac{1}{n}\sum_{i=1}^n g(X_i, \theta) - \mathbb{E}_{X\sim\mathcal{P}^*}[g(X, \theta)] \bigg\|_2 \geq c_0 \sqrt{\frac{\log n}{n}} \right) \leq \frac{1}{n^c}.
\end{equation*}
\paragraph{Proof of the second statement.} 

Because $d$ is a fixed constant independent of $(m,n)$ in our setting, this result follows as a direct corollary of the second assertion of Lemma 24 in \cite{MR4758128}, assuming the Lipschitz continuity of $\ell(x,\theta)$.

\paragraph{Proof of the third and fourth statements.} 
These two statements correspond to the special cases of the third and fourth assertions of Lemma \ref{lemmaprobB2} when $\beta_2=0$, and their proofs can be found in the proof of Lemma \ref{lemmaprobB2}.

\subsubsection{Proof of Lemma \ref{lemmaprobB2}}
 \par\noindent By compactness, choose $r > 0$ such that $\Theta \subseteq B(\theta^*, r)$.
\paragraph{Proof of the third statement.} We proceed via a peeling argument over the parameter space.
\par\noindent Define 
\[
\delta_n =
\begin{cases}
\left(\dfrac{\log n}{n}\right)^{\frac{3}{2\beta_2}}, & \beta_2 > 0, \\[6pt]
\left(\dfrac{\log n}{n}\right)^{\frac{1}{2\beta_1}}
, & \beta_2 = 0.
\end{cases}
\]
 For $k = 0, 1, \dots, \left\lfloor\log_2\frac{2r}{\delta_n}\right\rfloor+1$, we define
the set
\begin{equation}\label{partion A}
\mathcal{A}_k =
\begin{cases}
\{\theta, \theta' \in  \Theta: \|\theta - \theta'\|_2 \leq \delta_n\} & k=0, \\
\{\theta, \theta' \in  \Theta: 2^{k-1}\delta_n < \|\theta - \theta'\|_2 \leq 2^k\delta_n\} & k=1,2,\dots,\left\lfloor\log_2\frac{2r}{\delta_n}\right\rfloor, \\
\{\theta, \theta' \in  \Theta: 2^{k-1}\delta_n < \|\theta - \theta'\| _2\leq 2r\} & k=\left\lfloor\log_2\frac{2r}{\delta_n}\right\rfloor+1.
\end{cases}\end{equation}Then $\{\theta, \theta' \in  \Theta\} = \sum_{k=0}^{\left\lfloor\log_2\frac{2r}{\delta_n}\right\rfloor+1}\mathcal{A}_k$. Fix an integer $0 \leq k \leq \left\lfloor\log_2\frac{2r}{\delta_n}\right\rfloor+1$, for every component $v$, we consider the
function set
\[
\mathcal{L}_k = \left\{\frac{1}{(2^k\delta_n)^{\beta_2}+\frac{\log n}{n}}\bigl(g_v(x,\theta) - g_v(x,\theta^{\prime})\bigr): (\theta,\theta')\in\mathcal{A}_k\right\}.
\]
By Assumptions~$B.2-B.3$, it holds that for every $k$ and $f\in\mathcal{L}_k$, $\sup_{x\in\mathcal{X}}|f(x)|\leq L$ and \begin{equation}\label{estimate 10.1}
\mb E_{X\sim \m P^*}[f^2(X)]\leq \frac{L}{((2^{k}\delta_n)^{\beta_2}+\frac{\log n}{n})^2}(2^k\delta_n)^{2\beta_1
}.\end{equation} Then consider the star hull
$\overline{\mathcal{L}}_k$ of $\mathcal{L}_k$ and  the local Rademacher complexity associated with $\overline{\mathcal{L}}_k$,
\begin{equation*}
\m{\overline{R}}_{n}(\delta;\overline{\mathcal{L}}_k)=\mb{E}_
{X\sim \m P^*}\mb E_{\varepsilon}\Bigg[\sup_{\substack{f\in\overline{\mathcal{L}}_k\\ \sqrt{\mb E_{X\sim \m P^*}[f(X)^2]}\leq\delta}}\Big|\frac{1}{n}\sum_{i=1}^{n}\varepsilon_i f(X_i)\Big|\Bigg],
\end{equation*}
where \(\varepsilon_{i}\) are i.i.d. samples from the Rademacher distribution.
\par\noindent For $f,f^{\prime}:\mathcal{X}\to \mathbb R$, define the pseudometric
\begin{equation*}
d_{n}(f,f^{\prime})=\sqrt{\frac{1}{n}\sum_{i=1}^{n}(f(X_{i})-f^{\prime}(X_{i}))^{2}}.
\end{equation*}
Then by uniform boundedness of functions in class \(\overline{\mathcal{L}}_k\), we can obtain that the covering number satisfies 
\begin{equation*}
\begin{aligned}
\log \m N(\overline{\mathcal{L}}_k,d_{n},\varepsilon) 
&\lesssim \log\frac{1}{\varepsilon}+\log \m N (\mathcal{L}_k,d_{n},\varepsilon) .\\
\end{aligned}
\end{equation*}
By the assumption on the covering number of $d_n^g$, 
\begin{equation*}
    \log \m N (\mathcal{L}_k,d_{n},\varepsilon)\lesssim\log \Big(\m N ( B (\theta^*,r),nd_n^g,\varepsilon)\Big) \lesssim \log\frac{n}{\varepsilon}.   \end{equation*}
By (3.84) of \cite{wainwright_2019} and the Dudley entropy inequality, we can get
\begin{equation*}
    \overline{\m R}_{n}(\sqrt{\frac{\log n}{n}})\lesssim\frac{\log n}{n} .
    \end{equation*}
\begin{lemma}\label{localization}
(\cite{wainwright_2019}, Theorem 14.20) Given a uniformly $b$-bounded function class $\m F$ that is star shaped around $0$, let $(\delta^{*})^{2}\geq\frac{c}{n}$ be any solution to the inequality $\overline{\m{R}}_{n}(\delta;\mathcal{F})\leq\delta^{2}$, then we have
\begin{equation*}
\sup_{f\in\m{F}}{\frac{\bigg|\frac{1}{n}\sum_{i=1}^{n}f(X_{i})-\mathbb E[f(X)]\bigg|}{\sqrt{\mathbb E[f(X)^{2}]}+\delta^{*}}}\leq 10\delta^{*}
\end{equation*}
with probability greater than \(1-c_{1}\exp(-c_{2}n(\delta^{*})^{2})\).
\end{lemma}
Then by Lemma \ref{localization} and \eqref{estimate 10.1}, there exists a constant $c$ such that it holds with probability at least $1-\frac{1}{n^c}$ that for any $(\theta,\theta')\in\mathcal{A}_k$,
\[
\begin{aligned}
&\Big|\frac{1}{n}\sum_{i=1}^ng_v(X_i,\theta)-\frac{1}{n}\sum_{i=1}^ng_v(X_i,\theta')-\mb E_{X\sim \m P^*}[g_v(X,\theta)]+\mb E_{X\sim \m P^*}[g_v(X,\theta')]\Big|\\
&\lesssim\sqrt{\frac{\log n}{n}}\cdot(2^{k}\delta_n)^{\beta_1}+\frac{\log n}{n}\cdot\left((2^k\delta_n)^{\beta_2}+\frac{\log n}{n}\right)\\
&\lesssim\sqrt{\frac{\log n}{n}}\cdot(2^{k-1}\delta_n)^{\beta_1}+\frac{\log n}{n}(2^{k-1}\delta_n)^{\beta_2}+\Big(\frac{\log n}{n}\Big)^2\\
&\lesssim\sqrt{\frac{\log n}{n}}\cdot\big(\|\theta-\theta'\|_2+\delta_n\big)^{\beta_1}+\frac{\log n}{n}\cdot\big(\|\theta-\theta'\|_2+\delta_n\big)^{\beta_2}+\Big(\frac{\log n}{n}\Big)^2\\
&\lesssim\sqrt{\frac{\log n}{n}}\|\theta-\theta'\|_2^{\beta_1}+\frac{\log n}{n}\|\theta-\theta'\|_2^{\beta_2}+\Big(\frac{\log n}{n}\Big)^2.
\end{aligned}
\]
Then by $\log_2\frac{r}{\delta_n} \lesssim \log n$, consider the intersection of the above events for $k = 0, 1, \dots, \left\lfloor\log_2\frac{2r}{\delta_n}\right\rfloor+1$, and take the union bound for every component $v$, we can obtain the desired result.
\paragraph{Proof of the fourth statement.}\par\noindent Define $\delta_n = \left(\frac{\log n}{n}\right)^{\frac{5}{2(1+\beta_2)
}}$, for $\mathcal{A}_k$ defined in \eqref{partion A}, we consider the
function set
\[
\mathcal{L}_k = \left\{\frac{1}{(2^k\delta_n)^{1+\beta_2}+(\frac{\log n}{n})^2}\bigl(\ell(x,\theta) - \ell(x,\theta^{\prime}) - g(x,\theta^{\prime}) \cdot (\theta - \theta^{\prime})\bigr): (\theta,\theta')\in\mathcal{A}_k\right\}.
\]
By Assumptions~$B.2-B.3$, it holds that for every $k$ and $f\in\mathcal{L}_k$, $\sup_{x\in\mathcal{X}}|f(x)|\leq L$ and \begin{equation}\label{estimate 10.2}
\mb E_{X\sim \m P^*}[f^2(X)]\leq \frac{L}{((2^k\delta_n)^{1+\beta_2}+(\frac{\log n}{n})^2)^2}(2^k\delta_n)^{2(1+\beta_1)}.\end{equation}
By the Lipschitzness of $\ell(x,\theta)$ and assumption on the covering number of $d_n^g$, 
\begin{equation*}
    \log \m N (\mathcal{L}_k,d_{n},\varepsilon)\lesssim\log\frac{n^2}{\varepsilon}+\log \m N ( B (\theta^*,r),n^2d_n^g,\varepsilon) \lesssim \log\frac{n}{\varepsilon}.   \end{equation*}
As the second-order moment is bounded by \eqref{estimate 10.2}, Lemma \ref{localization} implies that there exists a constant $c$ such that it holds with probability at least $1-\frac{1}{n^c}$ that for any $(\theta,\theta')\in\mathcal{A}_k$,
\[
\begin{aligned}
&\Big|\frac{1}{n}\sum_{i=1}^n\ell(X_i,\theta)-\frac{1}{n}\sum_{i=1}^n\ell(X_i,\theta')-\frac{1}{n}\sum_{i=1}^n g(X_i,\theta')\cdot(\theta-\theta')\\
&\qquad-\mb E_{X\sim \m P^*}[\ell(X,\theta)]+\mb E_{X\sim \m P^*}[\ell(X,\theta')]+\mb E_{X\sim \m P^*}[g(X,\theta')\cdot(\theta-\theta')]\Big|\\
&\lesssim\sqrt{\frac{\log n}{n}}\cdot(2^{k}\delta_n)^{\beta_1+1}+\frac{\log n}{n}\Big((2^k\delta_n)^{1+\beta_2}+\Big(\frac{\log n}{n}\Big)^2\Big)\\
&\lesssim\sqrt{\frac{\log n}{n}}\cdot(2^{k-1}\delta_n)^{\beta_1+1}+\frac{\log n}{n}(2^{k-1}\delta_n)^{1+\beta_2}+\Big(\frac{\log n}{n}\Big)^3\\
&\lesssim\sqrt{\frac{\log n}{n}}\cdot\big(\|\theta-\theta'\|_2+\delta_n\big)^{\beta_1+1}+\frac{\log n}{n}\cdot\big(\|\theta-\theta'\|_2+\delta_n\big)^{\beta_2+1}+\Big(\frac{\log n}{n}\Big)^3\\
&\lesssim\sqrt{\frac{\log n}{n}}\|\theta-\theta'\|_2^{\beta_1+1}+\frac{\log n}{n}\|\theta-\theta'\|_2^{\beta_2+1}+\Big(\frac{\log n}{n}\Big)^3.
\end{aligned}
\]
Take the intersection of the above events for $k = 0, 1, \dots, \left\lfloor\log_2\frac{2r}{\delta_n}\right\rfloor+1$, the statement follows.
\paragraph{Proof of the first statement.}
\par\noindent Define 
\[
\delta_n =
\begin{cases}
\left(\dfrac{\log n}{n}\right)^{\frac{3}{4\beta_2}}, & \beta_2 > 0, \\[6pt]
\left(\dfrac{\log n}{n}\right)^{\frac{1}{2\beta_1}}, & \beta_2 = 0.
\end{cases}
\]
For component $v$ and $\mathcal{A}_k$ defined in \eqref{partion A}, we consider the
function set
\[
\mathcal{L}_k = \left\{\frac{1}{(2^k\delta_n)^{2\beta_2}+\frac{\log n}{n}}\Big|g_v(x,\theta)- g_v(x,\theta')\Big|^2: (\theta,\theta')\in\mathcal{A}_k\right\}.
\]
Then by Assumptions~$B.2-B.3$, for any $k$ and $f\in\mathcal{L}_k$, it holds that $\sup_{x\in\mathcal{X}}|f(x)|\leq L$ and \begin{equation}\label{estimate 10.3}
\mb E_{X\sim \m P^*}[f^2(X)]\leq \frac{L^2}{((2^k\delta_n)^{2\beta_2}+\frac{\log n}{n})^2}(2^k\delta_n)^{2(\beta_1+\beta_2)}.\end{equation}
\par\noindent Using the same argument to estimate the covering number as in the proof of the third statement of Lemma \ref{lemmaprobB2} and combining the second-order moment bound \eqref{estimate 10.3} with Lemma \ref{localization}, we obtain a constant $c$ such that it holds with probability at least $1-\frac{1}{n^c}$ that for any $(\theta,\theta')\in\mathcal{A}_k$,
\[
\begin{aligned}
&\left|\frac{1}{n}\sum_{i=1}^n\Big|g_v(X_i,\theta)- g_v(X_i,\theta')\Big|^2-
\mb E_{X\sim \m P^*}[(g_v(X,\theta)-g_v(X,\theta'))^2
]
\right| \\
&\lesssim\sqrt{\frac{\log n}{n}}\cdot(2^{k}\delta_n)^{\beta_1+\beta_2}+\frac{\log n}{n}\cdot\left((2^k\delta_n)^{2\beta_2}+\frac{\log n}{n}\right) \\
&\lesssim\sqrt{\frac{\log n}{n}}\cdot\bigl(\|\theta-\theta'\|_2+\delta_n\bigr)^{\beta_1+\beta_2}+\frac{\log n}{n}\cdot\bigl(\|\theta-\theta'\|_2+\delta_n\bigr)^{2\beta_2}+\left(\frac{\log n}{n}\right)^2\\
&\lesssim\sqrt{\frac{\log n}{n}}\|\theta-\theta'\|_2^{\beta_1+\beta_2}+\frac{\log n}{n}\|\theta-\theta'\|^{2\beta_2}_2+\left(\frac{\log n}{n}\right)^2.
\end{aligned}
\]
Moreover, we have
\begin{equation*}
    \mb E_{X\sim \m P^*}[(g_v(X,\theta)-g_v(X,\theta'))^2]\leq L\|\theta-\theta'\|^{2\beta_1}_2.
    \end{equation*}
We  consider the intersection of the above events for $k = 0, 1, \dots, \left\lfloor\log_2\frac{2r}{\delta_n}\right\rfloor+1$. Applying the union bound to every component $v$, we then derive the desired result.\paragraph{Proof of the second statement.}
\par\noindent Define $\delta_n = \left(\frac{\log n}{n}\right)^{\frac{7}{4(1+\beta_2)}}$. For $k = 0, 1, \dots, \left\lfloor\log_2\frac{2r}{\delta_n}\right\rfloor+1$, with $\mathcal{A}_k$ defined in \eqref{partion A}, we consider the function set 
\[
\mathcal{L}_k = \left\{\frac{1}{(2^k\delta_n)^{2(\beta_2+1)}+(\frac{\log n}{n})^3}\Big|\ell(x,\theta)- \ell(x,\theta')-g(x,\theta')\cdot(\theta-\theta')\Big|^2: (\theta,\theta')\in\mathcal{A}_k\right\}.
\]
Then by Assumptions~$B.2-B.3$, for any $k$ and $f\in\mathcal{L}_k$, it holds that $\sup_{x\in\mathcal{X}}|f(x)|\leq L$ and \begin{equation}\label{estimate 10.4}
\mb E_{X\sim \m P^*}[f^2(X)]\leq \frac{L^2}{((2^k\delta_n)^{2(\beta_2+1)}+(\frac{\log n}{n})^3)^2}(2^k\delta_n)^{2(\beta_1+\beta_2+2)}.\end{equation}
As in the proof of the third statement of Lemma \ref{lemmaprobB2}, we use the Lipschitzness of $\ell(x,\theta)$ and assumption on the covering number of $d_n^g$ to estimate the covering number. By \eqref{estimate 10.4} and Lemma \ref{localization}, there exists a constant $c$ such that it holds with probability at least $1-\frac{1}{n^c}$ that for any $(\theta,\theta')\in\mathcal{A}_k$,
\[
\begin{aligned}
&\left|\frac{1}{n}\sum_{i=1}^n\Big|\ell(X_i,\theta)- \ell(X_i,\theta')-g(X_i,\theta')\cdot(\theta-\theta')\Big|^2-\mb E_{X\sim \m P^*}[\Big|\ell(X,\theta)-\ell(X,\theta')- g(X,\theta')\cdot(\theta-\theta')\Big|^2]\right| \\
&\lesssim\sqrt{\frac{\log n}{n}}\cdot(2^{k-1}\delta_n)^{\beta_1+\beta_2+2}+\frac{\log n}{n}\cdot\left((2^k\delta_n)^{2(1+\beta_2)}+\Big(\frac{\log n}{n}\Big)^3\right) \\
&\lesssim\sqrt{\frac{\log n}{n}}\cdot\bigl(\|\theta-\theta'\|_2+\delta_n\bigr)^{\beta_1+\beta_2+2}+\frac{\log n}{n}\cdot\bigl(\|\theta-\theta'\|_2+\delta_n\bigr)^{2\beta_2+2}+\left(\frac{\log n}{n}\right)^4\\
&\lesssim\sqrt{\frac{\log n}{n}}\|\theta-\theta'\|_2^{\beta_1+\beta_2+2}+\frac{\log n}{n}\|\theta-\theta'\|^{2\beta_2+2}_2+\left(\frac{\log n}{n}\right)^4.
\end{aligned}
\]
Moreover, we have
\begin{equation*}
    \mathbb{E}_{X\sim P^*}[(\ell(X,\theta)-\ell(X,\theta')-g(X,\theta')\cdot(\theta-\theta'))^2]\leq L\|\theta-\theta'\|_2^{2(\beta_1+1)}.
\end{equation*}
Taking the intersection of the aforementioned events for $k = 0, 1, \dots, \left\lfloor\log_2\frac{2r}{\delta_n}\right\rfloor+1$, yields the desired result.

\subsubsection{Proof of Lemma \ref{lemmaprobA1}}
The first and second statements correspond to the special cases of the first and second statements in Lemma \ref{lemmaprobB1}, where we replace $\bX^{(n)} \sim \m P^{* \otimes n}$ with $\wX^{(m)} \sim \mu_n^{\otimes m}$.
\subsubsection{Proof of Lemma \ref{lemmaprobA2}}
\paragraph{Proof of the first statement.}
\par\noindent Define 
\[
\delta_n =
\begin{cases}
\left(\dfrac{\log n}{n}\right)^{\frac{1}{\beta_2}}, & \beta_2 > 0, \\[6pt]
\left(\dfrac{\log n}{n}\right)^{\frac{1}{2\beta_1}}, & \beta_2 = 0.
\end{cases}
\].\\For component $v$ and $\mathcal{A}_k$ defined in \eqref{partion A}, we consider the
function set
\[
\mathcal{L}_k = \left\{\frac{1}{(2^k\delta_n)^{\beta_2}+\frac{\log n}{n}}\Big(g_v(x,\theta)- g_v(x,\theta')\Big): (\theta,\theta')\in\mathcal{A}_k\right\}.
\]
Then by Assumption~$B.3$, for any $k$ and $f\in\mathcal{L}_k$, it holds that $\sup_{x\in \bX^{(n)}}|f(x)|\leq L$. And the first statement of Lemma \ref{lemmaprobB2} implies \begin{equation}
\label{estimate 10.5}
\begin{aligned}
\mb{E}_{\wX\sim \mu_n}[f^2(\wX)]&={\frac{\frac{1}{n}\sum_{i=1}^n\Big|g_v(X_i,\theta)- g_v(X_i,\theta')\Big|^2
}{\left((2^k\delta_n)^{\beta_2}+\frac{\log n}{n}\right)^2}}\\&
\lesssim{\frac{\|\theta-\theta'\|^{2\beta_1}_2+\sqrt{\frac{\log n}{n}} \,\|\theta-\theta'\|^{\beta_1+\beta_2}_2+\frac{\log n}{n}\,\|\theta-\theta'\|^{2\beta_2}_2+\left(\frac{\log n}{n}\right)^2
}{\left((2^k\delta_n)^{\beta_2}+\frac{\log n}{n}\right)^2}}
\\&\lesssim
{\frac{(2^k\delta_n)^{2\beta_1}+\sqrt{\frac{\log n}{n}} (2^k\delta_n)^{\beta_1+\beta_2}+\frac{\log n}{n}(2^k\delta_n)^{2\beta_2}+\left(\frac{\log n}{n}\right)^2
}{\left((2^k\delta_n)^{\beta_2}+\frac{\log n}{n}\right)^2}}.
\end{aligned}
\end{equation}
Then consider the star hull
$\overline{\mathcal{L}}_k$ of $\mathcal{L}_k$ and  the local Rademacher complexity associated with $\overline{\mathcal{L}}_k$,
\begin{equation*}
\m{\overline{R}}_{m}(\delta;\overline{\mathcal{L}}_k)=\mb{E}_
{\wX\sim \mu_n}\mb E_{\varepsilon}\Bigg[\sup_{\substack{f\in\overline{\mathcal{L}}_k\\ \sqrt{\mb E_{\wX\sim \mu_n}[f(\wX)^2]}\leq\delta}}\Big|\frac{1}{m}\sum_{j=1}^{m}\varepsilon_j f(\wt X_j)\Big|\Bigg],
\end{equation*}
where \(\varepsilon_{j}\) are i.i.d. samples from the Rademacher distribution.
\par\noindent For $f,f^{\prime}:\bX^{(n)}\to \mathbb R$, define the pseudometric
\begin{equation*}
d_{m}(f,f^{\prime})=\sqrt{\frac{1}{m}\sum_{j=1}^{m}(f(\wt X_{j})-f^{\prime}(\wt X_{j}))^{2}}.
\end{equation*}
Then by uniform boundedness of functions in class \(\overline{\mathcal{L}}_k\), we can obtain that the covering number satisfies 
\begin{equation*}
\begin{aligned}
\log \m N(\overline{\mathcal{L}}_k,d_{m},\varepsilon) 
&\lesssim \log\frac{1}{\varepsilon}+\log \m N (\mathcal{L}_k,d_{m},\varepsilon) .\\
\end{aligned}
\end{equation*}
By the assumption on the covering number of $d_m^g$ and $m\gtrsim n^{\gamma}$, 
\begin{equation*}
    \log \m N (\mathcal{L}_k,d_{m},\varepsilon)\lesssim\log \Big(\m N ( B (\theta^*,r),nd_m^g,\varepsilon)\Big) \lesssim \log\frac{m}{\varepsilon}.   \end{equation*}
By (3.84) of \cite{wainwright_2019} and the Dudley entropy inequality, we can get
\begin{equation*}
    \overline{\m R}_{m}(\sqrt{\frac{\log m}{m}})\lesssim\frac{\log m}{m} .
    \end{equation*}
As before, we can use \eqref{estimate 10.5} and Lemma \ref{localization} to obtain that there exists a constant $c$ such that it holds with probability at least $1-\frac{1}{m^c}$ that for any $(\theta,\theta')\in\mathcal{A}_k$,
\[
\begin{aligned}
&\Big|\frac{1}{m}\sum_{j=1}^m\big[g_v(\wt X_j,\theta)-g_v(\wt X_j,\theta')\big]-\frac{1}{n}\sum_{i=1}^n\big[g_v(X_i,\theta)-g_v(X_i,\theta')]\Big|\\
&\lesssim\sqrt{\frac{\log m}{m}}\|\theta-\theta'\|_2^{\beta_1}+\sqrt{\frac{\log m}{m}}\Big(\frac{\log n}{n}\Big)^{\frac{1}{4}}\|\theta-\theta'\|_2^{\frac{\beta_1+\beta_2}{2}}\\
&\quad+\frac{\log m}{m}\|\theta-\theta'\|_2^{\beta_2}+\sqrt{\frac{\log m}{m}}\frac{\log n}{n}.
\end{aligned}
\]
Finally, by considering the intersection of the above events for $k = 0, 1, \dots, \left\lfloor\log_2\frac{2r}{\delta_n}\right\rfloor+1$ and applying the union bound to each component $v$, and noting that $m\gtrsim n^{\gamma}$, the concentration inequality follows.
\paragraph{Proof of the second statement.}
\par\noindent Define $\delta_n = \left(\frac{\log n}{n}\right)^{\frac{2}{\beta_2+1}}$, for $\mathcal{A}_k$ defined in \eqref{partion A}, we consider the
function set
\[
\mathcal{L}_k = \left\{\frac{1}{(2^k\delta_n)^{1+\beta_2}+(\frac{\log n}{n})^2}\Big(\ell(x,\theta)- \ell(x,\theta')-g(x,\theta')\cdot(\theta-\theta')\Big): (\theta,\theta')\in\mathcal{A}_k\right\}.
\]
Then by Assumption~$B.3$, for any $k$ and $f\in\mathcal{L}_k$, it holds that $\sup_{x\in \bX^{(n)}}|f(x)|\leq L$. And  the second statement of Lemma \ref{lemmaprobB2} implies \begin{equation}\label{estimate 10.6}\begin{aligned}
\mb{E}_{\wX\sim \mu_n}[f^2(\wX)]&={\frac{\frac{1}{n}\sum_{i=1}^n\Big|\ell(X_i,\theta)- \ell(X_i,\theta')-g(X_i,\theta')\cdot(\theta-\theta')\Big|^2
}{\left((2^k\delta_n)^{1+\beta_2}+(\frac{\log n}{n})^2\right)^2}}\\&
\lesssim{\frac{\|\theta-\theta'\|^{2\beta_1+2}_2+\sqrt{\frac{\log n}{n}} \,\|\theta-\theta'\|^{\beta_1+\beta_2+2}_2+\frac{\log n}{n}\,\|\theta-\theta'\|^{2\beta_2+2}_2+\left(\frac{\log n}{n}\right)^4
}{\left((2^k\delta_n)^{1+\beta_2}+(\frac{\log n}{n})^2\right)^2}}
\\&\lesssim
{\frac{(2^k\delta_n)^{2\beta_1+2}+\sqrt{\frac{\log n}{n}} (2^k\delta_n)^{\beta_1+\beta_2+2}+\frac{\log n}{n}(2^k\delta_n)^{2\beta_2+2}+\left(\frac{\log n}{n}\right)^4
}{\left((2^k\delta_n)^{1+\beta_2}+(\frac{\log n}{n})^2\right)^2}}.
\end{aligned}
\end{equation}
By Lipschitzness of $\ell(x,\theta)$ and the assumption on the covering number of $d_m^g$, and using $m\gtrsim n^{\gamma}$, we can obtain that  the covering number satisfies 
\begin{equation*}
\begin{aligned}
\log \m N (\mathcal{L}_k,d_{m},\varepsilon) 
\lesssim\log\frac{m}{\varepsilon}.\\
\end{aligned}
\end{equation*}
As before, we can use \eqref{estimate 10.6} and Lemma \ref{localization} to obtain that there exists a constant $c$ such that it holds with probability at least $1-\frac{1}{m^c}$ that for any $(\theta,\theta')\in\mathcal{A}_k$,
\[
\begin{aligned}
&\Big|\frac{1}{m}\sum_{j=1}^m\big[\ell(\wt X_j,\theta)-\ell(\wt X_j,\theta')-g(\wt X_j,\theta')\cdot(\theta-\theta')\big]\\
&\qquad-\frac{1}{n}\sum_{i=1}^n\big[\ell(X_i,\theta)-\ell(X_i,\theta')-g(X_i,\theta')\cdot(\theta-\theta')\big]\Big|\\
&\lesssim\sqrt{\frac{\log m}{m}}\|\theta-\theta'\|_2^{\beta_1+1}+\sqrt{\frac{\log m}{m}}\Big(\frac{\log n}{n}\Big)^{\frac{1}{4}}\|\theta-\theta'\|_2^{\frac{\beta_1+\beta_2}{2}+1}\\
&\quad+\frac{\log m}{m}\|\theta-\theta'\|_2^{\beta_2+1}+\sqrt{\frac{\log m}{m}}\Big(\frac{\log n}{n}\Big)^2.
\end{aligned}
\]
Taking the intersection of the above events for $k = 0, 1, \dots, \left\lfloor\log_2\frac{2r}{\delta_n}\right\rfloor+1$ and using
$m\gtrsim n^{\gamma}$, the concentration inequality holds.

\end{document}